%% file: main.tex
\documentclass[onefignum,onetabnum]{siamart251216}

\usepackage{graphicx,psfrag,amsmath,amsfonts,verbatim}
\usepackage{indentfirst}
\usepackage{xcolor}
\usepackage{ragged2e}
\usepackage{booktabs,makecell, multirow, tabularx}
\usepackage{longtable}
\usepackage{enumerate}
\usepackage{amssymb}
\usepackage{setspace}
\usepackage{tikz}
\usetikzlibrary{arrows.meta, fit}
\usepackage[numbers]{natbib} 
\allowdisplaybreaks[4]

\renewcommand\arraystretch{1.2}

\input{shared}

\usepackage{hyperref}

\makeatletter
\AtBeginDocument{\long\def\cref@override@label@type#1\@nil#2{#1}}
\makeatother

\ifpdf
\hypersetup{
  pdftitle={Error Analysis of Deep Learning}
}
\fi

\begin{document}

\maketitle

\begin{abstract}
This review examines deep learning theory through the lens of error decomposition, systematically addressing approximation, statistical, and optimization errors. We trace the development of error analyses from bounds on individual error components to frameworks that account for all three. Using nonparametric regression and the Deep Ritz Method as dual running examples, we demonstrate how these analyses apply across fundamentally different loss structures. We then review representative analyses in scientific machine learning, generative learning, reinforcement learning, representation learning, and large language models, showing how learning errors affect solutions, distributions, and predictions through problem-specific arguments. Throughout, we clarify which error components are controlled, how they are combined, and what assumptions are needed to obtain guarantees for trained networks.
\end{abstract}




\makeatletter
\renewcommand{\tableofcontents}{%
  \begingroup
  \c@tocdepth=2\relax
  {\centering\large\bfseries Contents\par}%
  \vspace{0.5em}%
  \@starttoc{toc}%
  \endgroup
  \vspace{0.5em}%
}
\newcommand*\l@section{\@dottedtocline{1}{1.5em}{2.3em}}
\newcommand*\l@subsection{\@dottedtocline{2}{3.8em}{3.2em}}
\newcommand*\l@subsubsection{\@dottedtocline{3}{7.0em}{4.1em}}
\newcommand*\l@paragraph[2]{}    
\newcommand*\l@subparagraph[2]{} 
\makeatother
\tableofcontents

\input{sections/introduction}

\input{sections/foundations}

\input{sections/classical}

\input{sections/overparameterization}

\input{sections/complete}

\input{sections/applications}

\input{sections/horizon}

\appendix

\input{sections/appendix}

\bibliographystyle{siamplain}
\bibliography{references}

\end{document}

%% file: shared.tex
\usepackage{lipsum}
\usepackage{amsfonts}
\usepackage{graphicx}
\usepackage{epstopdf}
\usepackage{algorithmic}
\ifpdf
  \DeclareGraphicsExtensions{.eps,.pdf,.png,.jpg}
\else
  \DeclareGraphicsExtensions{.eps}
\fi

\newsiamremark{rem}{Remark}
\newsiamremark{hypothesis}{Hypothesis}
\crefname{hypothesis}{Hypothesis}{Hypotheses}
\newsiamthm{claim}{Claim}

\headers{Error Analysis of Deep Learning}{Y. Jiao, P. Wu, J.Z. Yang, and P. Zhang}

\title{Error Analysis of Deep Learning
}

\author{
Yuling Jiao\thanks{School of Artificial Intelligence, National Center for Applied Mathematics in Hubei, Hubei Key Laboratory of Computational Science, Wuhan University, Wuhan, China. Email: \texttt{yulingjiaomath@whu.edu.cn}}
\and Peiying Wu\thanks{School of Mathematics and Statistics, Wuhan University, Wuhan, China. Email: \texttt{peiyingwu@whu.edu.cn}}
\and Jerry Zhijian Yang\thanks{National Center for Applied Mathematics in Hubei, Wuhan Institute for Math \& AI, School of Mathematics and Statistics, Hubei Key Laboratory of Computational Science, Wuhan University, Wuhan, China. Email: \texttt{zjyang.math@whu.edu.cn}}
\and Pingwen Zhang\thanks{Wuhan Institute for Math \& AI, Wuhan University, Wuhan, School of Mathematical Sciences, Peking University, Beijing, China. Email: \texttt{pwzhang@whu.edu.cn}}
}

\usepackage{amsopn}

\usepackage{mathtools}

\newcommand{\argmin}{\mathop{\rm argmin}}
\newcommand{\argmax}{\mathop{\rm argmax}}

\newtheorem{asmp}{Assumption}[section]

\newtheorem{remark}{Remark}[section]

\def\approxcorrect{\checkmark\kern-1.1ex\raisebox{.89ex}{$\times$}}

\usepackage{amsmath,amsfonts,bm}

\def\1{\bm{1}}

\DeclareMathAlphabet{\mathsfit}{\encodingdefault}{\sfdefault}{m}{sl}
\SetMathAlphabet{\mathsfit}{bold}{\encodingdefault}{\sfdefault}{bx}{n}



%% file: sections/introduction.tex
\section{Introduction}\label{intro}

Deep learning is used across scientific and engineering disciplines for classification, regression, sequential decision making, generative modeling, and language processing. Neural networks are also used to approximate solutions of partial differential equations (PDEs), including high-dimensional problems \cite{Justin2018DGM,han2018solving}. The Deep Ritz Method (DRM) \cite{Weinan2017The} and Physics-Informed Neural Networks (PINNs) \cite{raissi2019physics} formulate such problems through variational energies and differential-equation residuals, respectively.

One approach to understanding the accuracy of trained neural networks separates three contributions: \textbf{approximation error}, which measures how well a network class represents the target \cite{devore2021neural}; \textbf{statistical error}, which accounts for replacing population quantities by finite-sample estimates \cite{bartlett2021deep}; and \textbf{optimization error}, which measures the empirical-risk gap left by the training algorithm \cite{bottou2018optimization,sun2020optimization}. Analyses organized around a function class often balance approximation and statistical errors while assuming exact empirical minimization. Analyses organized around training dynamics instead study a specified algorithm, using approaches such as the Neural Tangent Kernel (NTK), lazy training, or mean-field theory. These approaches describe different training regimes: kernel and linearized models use fixed features, whereas mean-field models can describe evolving features. Their statistical guarantees depend on the target class, data assumptions, and training regime. Combining a given network approximation theorem with a guarantee for the trained estimator therefore requires compatible statistical and optimization estimates.

Two developments are particularly relevant to this review. Norm-controlled analysis \cite{neyshabur2015norm,bartlett2017spectrally,golowich2018size,jiao2023approximation,ding2025convergence} bounds approximation and statistical errors using constraints on network weights, allowing meaningful guarantees when the parameter count exceeds the sample size. Separately, optimization analyses on structured architectures compare the training output with a constructed approximating network \cite{kohler2023opt,jiao2024drm}. Together with suitable complexity bounds, the latter approach yields population-risk guarantees that account for approximation, sampling, and optimization in the same argument.

This review presents the approximation, statistical, and optimization tools used in these analyses, with particular attention to how error decompositions combine them. We use two running problems to compare the roles of the tools and to explain representative complete three-error guarantees. We follow the line of work developed in \cite{jiao2021deep,jiao2022rate,jiao2023approximation,Duan2022,ding2025convergence,jiao2024drm}, with cross-references to neighboring approaches.

\subsection{Three approaches to error analysis}\label{intro:regimes}

The three analytical chapters distinguish how network complexity is controlled and whether the training algorithm is analyzed explicitly.

The classical error analysis (\S\ref{sec:classical_setting}) assumes the optimization error is zero and studies the approximation--statistical trade-off via size-based parameterizations of the network class. It draws on the analyses developed in \cite{schmidt2020nonparametric,Duan2022,jiao2021deep}. For bounded-noise regression, the fast-rate analysis presented below attains the classical H\"older exponent $2\beta/(2\beta+d)$ for the expected squared $L^2(P_X)$ error, up to logarithmic factors. The approximation tools include Yarotsky's $x^2$ sawtooth construction for ReLU networks and local-polynomial constructions with approximate partitions of unity for tanh networks in Sobolev norms. The statistical tools use Rademacher complexity and peeling.

The over-parameterized error analysis (\S\ref{sec:overparam_setting}) still assumes zero optimization error, but controls complexity through weight-norm constraints rather than parameter counts alone. The approximation results \cite{jiao2023approximation,ding2025convergence} control both the approximation accuracy and the weight norms of the constructed networks. Norm-dependent Rademacher bounds use Golowich's contraction-and-peeling argument \cite{golowich2018size}, with additional estimates for the derivative terms in DRM. These bounds allow arbitrarily large widths under the prescribed norm and depth conditions.

The complete error analysis (\S\ref{sec:complete_setting}) incorporates the optimization error via a complementary error decomposition that uses an explicitly constructed approximating network instead of the empirical risk minimizer as the analytical anchor, together with a structured architecture (parallel neural network, PNN) on which projected gradient descent (PGD) admits a convergence guarantee. The representative case is Theorem~\ref{thm:reg_complete} for nonparametric regression (following \citealp{kohler2023opt}); Theorem~\ref{thm:drm_complete} gives the corresponding result for the Deep Ritz Method (\citealp{jiao2024drm}).

This organization compares three analytical approaches rather than prescribing a chronological sequence. Size and norm constraints provide different ways to control statistical complexity, while a complete analysis additionally controls the training algorithm. We reuse the relevant tools across chapters and specify the additional estimates needed for each result.

\subsection{Two running problems}\label{intro:running}

We use two running problems throughout the review. The first is nonparametric regression: estimating an unknown function $f_0 : [0,1]^d \to \mathbb{R}$ from i.i.d.\ samples $\{(\boldsymbol{X}_i, Y_i)\}_{i=1}^n$, with the target quantity being the expected squared $L^2(P_X)$ error. The second is the DRM for elliptic PDEs, principally a Neumann boundary-value problem on $\Omega \subset [0,1]^d$ solved via Monte Carlo discretization of the energy functional, with the target quantity being the squared $H^1(\Omega)$ distance between the network estimator and the true weak solution. We also discuss related results with other boundary conditions and network classes. In both problems, statistical error is controlled through deviations between empirical and population objectives. Regression learns from noisy observations, whereas DRM approximates domain and boundary integrals by Monte Carlo sampling. The analyses also differ in their target norms and regularity assumptions: $L^2$ and H\"older regularity for regression, and $H^1$ and Sobolev regularity for DRM. Each analytical chapter (\S\ref{sec:classical_setting}--\S\ref{sec:complete_setting}) presents representative convergence results for these two types of problems.

\subsection{Paper organization}\label{intro:org}

The remainder of the review proceeds as follows.

\textbf{Chapter~\ref{sec:foundations}} (Foundations) establishes the technical preliminaries: neural-network notation (\S\ref{nnstructure}), symmetric and asymmetric error decompositions and their roles in the subsequent statistical analysis, and a preview of the complementary approximation-based decomposition (\S\ref{sec:decomp}); brief landscape surveys of the approximation, statistical, and optimization literatures (\S\ref{sec:approx_landscape}--\S\ref{sec:opt_landscape}); the formal setup of the two running problems (\S\ref{sec:running_examples}); and a notation summary that collects general mathematical conventions and function-space definitions (\S\ref{sec:notation_summary}).

\textbf{Chapter~\ref{sec:classical_setting}} (Classical Error Analysis) develops the approximation tools (\S\ref{sec:classical_approx}), the size-based statistical tools (\S\ref{sec:classical_stat}), and presents slow and fast convergence rates for the two running problems in the settings of the cited results (\S\ref{sec:classical_rates}). It also introduces minimax lower bounds under the corresponding statistical observation models (\S\ref{sec:minimax_lb}), and closes with a survey of mechanisms for escaping the curse of dimensionality (\S\ref{sec:low_intrinsic}).

\textbf{Chapter~\ref{sec:overparam_setting}} (Over-Parameterized Error Analysis) develops norm-controlled approximation and statistical bounds (\S\ref{sec:overparam_approx}--\S\ref{sec:overparam_stat}) and derives the resulting rates for the two running problems (\S\ref{sec:overparam_rates}).

\textbf{Chapter~\ref{sec:complete_setting}} (Complete Error Analysis) proves the approximation-based error decomposition (\S\ref{sec:new_decomp}), develops approximation and statistical bounds for parallel networks (\S\ref{sec:complete_approx}--\S\ref{sec:complete_stat}), and analyzes PGD with a detailed proof of the key optimization lemma (\S\ref{sec:complete_opt}). These estimates yield the complete regression and DRM guarantees in \S\ref{sec:complete_rates}.

\textbf{Chapter~\ref{sec:applications}} (Applications) introduces learning objectives and representative error analyses in scientific machine learning (\S\ref{section:app:sml}), generative learning (\S\ref{section:app:generative}), reinforcement learning (\S\ref{section:app:rl}), representation learning (\S\ref{section:app:repssl}), and large language models (\S\ref{section:app:llm}). It explains how approximation and finite-sample estimation affect the performance measures of each problem, distinguishing results that assume exact empirical minimization, retain a training residual, or directly analyze the training algorithm.

\textbf{Chapter~\ref{sec:open_horizons}} (Open Horizons) discusses improving the rates and resource requirements of complete analyses, reducing their dependence on particular architectures and training constraints, and analyzing specified training algorithms in the application models considered here.

\textbf{Appendix~\ref{app:notation}} collects notation tables organized by topic and a list of abbreviations. \textbf{Appendix~\ref{app:tool_index}} provides a tool dictionary and a convergence-rate index. We include selected proofs and proof sketches to explain the main techniques; further details are available in the cited sources.

%% file: sections/foundations.tex
\section{Foundations: Networks, Risks, and Error Decomposition}\label{sec:foundations}

The error analysis of deep learning decomposes the excess risk into three sources: approximation error, statistical error, and optimization error. This chapter establishes the foundations on which the rest of the review is built. We begin by fixing notation for neural networks (\S\ref{nnstructure}), present the error decomposition framework (\S\ref{sec:decomp}), survey the broad landscape of approximation, statistical, and optimization analyses (\S\ref{sec:approx_landscape}--\S\ref{sec:opt_landscape}), and introduce the two running examples that appear throughout the review (\S\ref{sec:running_examples}). Finally, we collect general mathematical conventions and function-space definitions in a notation summary (\S\ref{sec:notation_summary}); comprehensive symbol tables organized by topic appear in Appendix~\ref{app:notation}.

\subsection{Neural network notation}\label{nnstructure}
To ensure clarity and consistency, we first establish the standard notation used for neural networks throughout this review. Conceptually, a feedforward neural network realizes a function $\phi_{\boldsymbol{\theta}}(\boldsymbol{x})$ that maps an input vector $\boldsymbol{x} \in \mathbb{R}^d$ to an output. The vector $\boldsymbol{\theta}$ represents the collection of all learnable parameters (weights and biases) within the network. In our notation, we will consistently use bold lowercase letters for vectors (e.g., $\boldsymbol{x}, \boldsymbol{b}$), bold uppercase letters for matrices (e.g., $\boldsymbol{A}$), and italicized letters for scalar values (e.g., $L, W$).

More formally, we consider a function $\phi_{\boldsymbol{\theta}}: \mathbb{R}^{d} \rightarrow \mathbb{R}$ with depth $L \in \mathbb{N}$ and an element-wise activation function $\rho$. The network is the composition of affine maps and activations, with an input layer dimension $N_0 = d$ and output dimension $N_L=1$:
\begin{equation}
\label{eqnn1}
\phi_{\boldsymbol{\theta}}(\boldsymbol{x})
\;=\; \big(\boldsymbol{A}_{L-1}\,{\cdot}\, + \boldsymbol{b}_{L-1}\big) \circ \rho \circ \big(\boldsymbol{A}_{L-2}\,{\cdot}\, + \boldsymbol{b}_{L-2}\big) \circ \rho \circ \cdots \circ \rho \circ \big(\boldsymbol{A}_{0}\,{\cdot}\, + \boldsymbol{b}_{0}\big)\,(\boldsymbol{x}),
\end{equation}
where each factor $\big(\boldsymbol{A}_{\ell}\,{\cdot}\, + \boldsymbol{b}_{\ell}\big)$ denotes the affine map $\boldsymbol{z} \mapsto \boldsymbol{A}_{\ell}\boldsymbol{z} + \boldsymbol{b}_{\ell}$; for $L = 1$ the composition reduces to the single affine map $\boldsymbol{A}_0\boldsymbol{x} + \boldsymbol{b}_0$.
For each layer $\ell \in \{0, \ldots, L-1\}$, the terms $\boldsymbol{A}_{\ell} \in \mathbb{R}^{N_{\ell+1} \times N_{\ell}}$ and $\boldsymbol{b}_{\ell} \in \mathbb{R}^{N_{\ell+1}}$ represent the weight matrix and bias vector. For $L\ge2$, the network's width is defined as $W := \max \{N_{1}, \ldots, N_{L-1}\}$; the affine case $L=1$ has no hidden-width parameter. We denote the vector of all trainable parameters by $\boldsymbol{\theta} = ((\boldsymbol{A}_{0}, \boldsymbol{b}_{0}), \ldots, (\boldsymbol{A}_{L-1}, \boldsymbol{b}_{L-1}))$. The set of all possible parameter vectors constitutes the parameter space, denoted by $\Theta$. The family of all functions realized by this architecture is denoted by $\mathcal{NN}(W, L)$. For theoretical analysis it is common to constrain the magnitude of these parameters; we define the infinity norm of the parameter vector as
$$\|\boldsymbol{\theta}\|_{\infty} : = \max_{\ell = 0, \ldots, L-1} \max\{ \|\boldsymbol{A}_{\ell}\|_{\max}, \|\boldsymbol{b}_{\ell}\|_{\max} \}\, ,$$
where $\|\cdot\|_{\max}$ is the maximum absolute value of a matrix or vector's entries.

\begin{definition}[Bounded-parameter feedforward neural network classes]\label{def:nnclass}
Let $\rho$ be the activation function (specified in context) and let $\mathcal{NN}(W, L)$ be the family of feedforward networks of width $W$ and depth $L$ defined by \eqref{eqnn1}. This review uses two related parameterizations of the bounded-parameter sub-class of $\mathcal{NN}(W, L)$.

(i) Width-and-magnitude-indexed (used as the building block of the $\mathcal{PNN}$ class of Definition~\ref{def:pnn}):
\[
\mathcal{NN}_{\rho}(W, L, B_{\boldsymbol{\theta}}) := \big\{\phi_{\boldsymbol{\theta}} \in \mathcal{NN}(W, L) \,:\, \|\boldsymbol{\theta}\|_{\infty} \le B_{\boldsymbol{\theta}}\big\}.
\]

(ii) Depth-and-nonzero-count-indexed (used in DRM \S\ref{sec:drm_slow}--\S\ref{sec:drm_fast}, matching the codomain of the Gühring--Raslan tanh Sobolev approximation of Theorem~\ref{thm:A2}):
\[
\mathcal{NN}_{\rho}(L, \mathfrak{D}, B_{\boldsymbol{\theta}}) := \big\{\phi_{\boldsymbol{\theta}} \,:\, \mathrm{depth}\le L,\ \#\{\text{nonzero entries of }\boldsymbol{\theta}\}\le \mathfrak{D},\ \|\boldsymbol{\theta}\|_{\infty}\le B_{\boldsymbol{\theta}}\big\}.
\]
\end{definition}

For a fixed scalar-output architecture, the total parameter count is $\mathfrak D=\sum_{\ell=0}^{L-1}(N_\ell+1)N_{\ell+1}$, with $N_0=d$ and $N_L=1$. For $L\ge2$, padding every hidden layer to width $W$ gives
\begin{align} \label{eq:align_dim}
\mathfrak{D}(W, L, d) := (L-2)W^2 + L W + W d + 1\, 
\end{align}
parameters. For $L=1$, the affine map instead has $d+1$ parameters. Under the entrywise bound, the parameter space of a fixed architecture is the hypercube $\Theta = [-B_{\boldsymbol{\theta}}, B_{\boldsymbol{\theta}}]^{\mathfrak{D}}$; narrower architectures can be embedded in the padded parameter space by inserting zeros.

Throughout this review, $\mathcal{F}_{NN}$ denotes the neural network hypothesis class under consideration. The classical analyses use output-clipped, size-controlled ReLU networks for regression and bounded-parameter tanh networks for DRM. Chapter~\ref{sec:overparam_setting} uses norm-controlled ReLU networks and linear combinations of norm-controlled tanh sub-networks, with their respective constraints defined there. Chapter~\ref{sec:complete_setting} uses the parallel class $\mathcal{PNN}(\mathfrak{m},M,\{W,L,B_{\boldsymbol{\theta}}\})$ of Definition~\ref{def:pnn}.

\subsection{Error decomposition framework}\label{sec:decomp}

Given a pointwise loss function $\ell(y',y)$, define the population risk
\[
\mathcal{L}(\phi) := \mathbb{E}_{(\boldsymbol{X}, Y)}\bigl[\ell(\phi(\boldsymbol{X}), Y)\bigr].
\]
The Bayes risk is $\mathcal{L}^*=\inf_\phi\mathcal{L}(\phi)$, where the infimum is over admissible measurable predictors. Given a training sample $\{(\boldsymbol{X}_i,Y_i)\}_{i=1}^n$, the empirical risk is
\[
\widehat{\mathcal{L}}_n(\phi) := \frac{1}{n}\sum_{i=1}^n \ell(\phi(\boldsymbol{X}_i), Y_i).
\]

With the above objects fixed, we identify three primary sources of error.

\begin{itemize}
\item \textbf{Approximation error} measures the limitation imposed by the hypothesis class $\mathcal{F}_{NN}$:
\[
\mathcal{E}_{app} := \inf_{\phi\in\mathcal{F}_{NN}}\mathcal{L}(\phi)-\mathcal{L}^*.
\]

\item \textbf{Statistical error} arises because $\widehat{\mathcal{L}}_n$ is only an approximation of $\mathcal{L}$ on a finite sample. It is often captured by the uniform deviation
\[
\mathcal{E}_{sta} := \sup_{\phi \in \mathcal{F}_{NN}} \big|\mathcal{L}(\phi) - \widehat{\mathcal{L}}_n(\phi)\big|.
\]

\item \textbf{Optimization error} measures the empirical-risk gap of the algorithm output $\phi_{\boldsymbol{\theta}_{\mathcal{A}}}\in\mathcal{F}_{NN}$:
\[
\mathcal{E}_{opt} := \widehat{\mathcal{L}}_n(\phi_{\boldsymbol{\theta}_{\mathcal{A}}})-\inf_{\phi\in\mathcal{F}_{NN}}\widehat{\mathcal{L}}_n(\phi).
\]
\end{itemize}

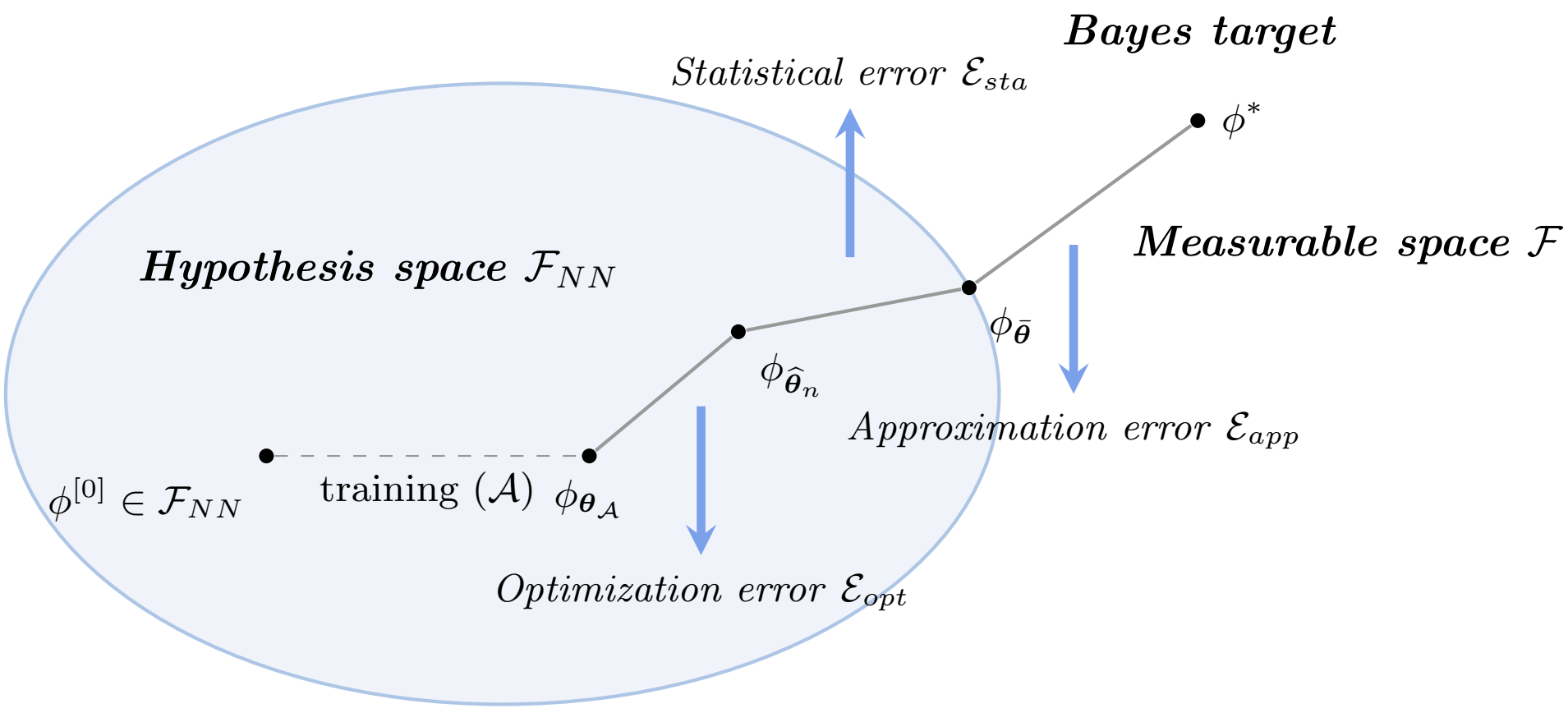
\begin{figure}[!t]
    \centering
    \begin{tikzpicture}[>=stealth, font=\small]

        \definecolor{webbg}{HTML}{F0F4FA}
        \definecolor{webborder}{HTML}{AEC6E6}
        \definecolor{webarrow}{HTML}{7CA1EB}

        \def\ellipseAngle{0}

        \draw[fill=webbg, draw=webborder, thick, rotate=\ellipseAngle] (0,0) ellipse (4cm and 2.5cm);
        \node[font=\bfseries\itshape, text=black] at (-1.0, 1.0) {Hypothesis space $\mathcal{F}_{NN}$};

        \coordinate (fbar) at ([rotate=\ellipseAngle] 20:4cm and 2.5cm);

        \node (n_fbar) at (fbar) [circle, fill=black, inner sep=1.2pt] {};
        \node[below right=1pt] at (n_fbar) {$\phi_{\bar{\boldsymbol{\theta}}}$};

        \node (n_fstar) at (5.6, 2.2) [circle, fill=black, inner sep=1.2pt, label=right:$\phi^{*}$] {};
        \node[font=\bfseries\itshape, text=black] at (5.6, 2.9) {Bayes target};

        \node (n_ftheta) at (1.9, 0.5) [circle, fill=black, inner sep=1.2pt, label=below right:$\phi_{\widehat{\boldsymbol{\theta}}_n}$] {};

        \node (n_fA) at (0.7, -0.5) [circle, fill=black, inner sep=1.2pt, label=below:$\phi_{\boldsymbol{\theta}_{\mathcal{A}}}$] {};
        \node (n_finit) at (-1.9, -0.5) [circle, fill=black, inner sep=1.2pt, label=below left:{$\phi^{[0]} \in \mathcal{F}_{NN}$}] {};

        \draw[dashed, gray!80] (n_finit) -- node[below, text=black] {training ($\mathcal{A}$)} (n_fA);
        \draw[thick, gray!80] (n_fA) -- (n_ftheta);
        \draw[thick, gray!80] (n_ftheta) -- (n_fbar);
        \draw[thick, gray!80] (n_fbar) -- (n_fstar);

        \node[font=\bfseries\itshape] at (6.8, 1.2) {Measurable space $\mathcal{F}$};

        \draw[->, line width=2pt, webarrow] (1.6, -0.1) -- (1.6, -1.3) node[below, text=black] {\textit{Optimization error} $\mathcal{E}_{opt}$};
        \draw[->, line width=2pt, webarrow] (2.8, 1.1) -- (2.8, 2.3) node[above, text=black] {\textit{Statistical error} $\mathcal{E}_{sta}$};
        \draw[->, line width=2pt, webarrow] (4.6, 1.2) -- (4.6, 0.0) node[below, text=black] {\textit{Approximation error} $\mathcal{E}_{app}$};

    \end{tikzpicture}
    \caption{Schematic diagram of error decomposition.}
    \label{fig:error_decom1}
\end{figure}

When the respective minima are attained, write $\phi^*$ for a Bayes predictor, $\phi_{\bar{\boldsymbol{\theta}}}$ for a population-risk minimizer in $\mathcal{F}_{NN}$, and $\phi_{\widehat{\boldsymbol{\theta}}_n}$ for an empirical risk minimizer (ERM), as illustrated in Figure~\ref{fig:error_decom1}. For any fixed comparison function $\phi\in\mathcal{F}_{NN}$, inserting empirical risks at $\phi$ and at the algorithm output bounds the excess risk by $\mathcal{L}(\phi)-\mathcal{L}^*+\mathcal{E}_{opt}+2\mathcal{E}_{sta}$. Taking the infimum over $\phi$ and then expectation gives
\begin{equation}\label{eq:sym_decomp}
\mathbb{E}\big[\mathcal{L}(\phi_{\boldsymbol{\theta}_{\mathcal{A}}})\big] - \mathcal{L}^* \;\le\; \underbrace{\inf_{\phi\in\mathcal{F}_{NN}}\mathcal{L}(\phi)-\mathcal{L}^*}_{\mathcal{E}_{app}} + \underbrace{2\,\mathbb{E}\sup_{\phi \in \mathcal{F}_{NN}} \big|\mathcal{L}(\phi) - \widehat{\mathcal{L}}_n(\phi)\big|}_{\mathcal{E}_{sta,\,slow}} + \mathbb{E}[\mathcal{E}_{opt}].
\end{equation}

We call \eqref{eq:sym_decomp} the symmetric decomposition because its statistical term is a two-sided uniform deviation. It is commonly used for slow-rate statistical bounds. Figure~\ref{fig:error_decom1} illustrates the roles of the comparison functions; the connecting segments are schematic, not distances equal to the risk errors.

An alternative asymmetric decomposition retains a negative excess-risk term in the statistical process. Using empirical near-optimality and $\mathbb{E}\widehat{\mathcal{L}}_n(\phi)=\mathcal{L}(\phi)$ for each fixed comparison function gives
\begin{equation}\label{eq:asym_decomp}
\mathbb{E}\big[\mathcal{L}(\phi_{\boldsymbol{\theta}_{\mathcal{A}}})\big] - \mathcal{L}^* \;\le\; 2\,\mathcal{E}_{app} + \underbrace{\mathbb{E}\sup_{\phi\in\mathcal{F}_{NN}}\big[\mathcal{L}(\phi) + \mathcal{L}^* - 2\widehat{\mathcal{L}}_n(\phi)\big]}_{\mathcal{E}_{sta,\,fast}} + 2\,\mathbb{E}[\mathcal{E}_{opt}].
\end{equation}

This inequality does not require strong convexity. Under suitable variance and complexity conditions, the negative excess-risk term supports fast-rate statistical bounds; these conditions are developed in \S\ref{sec:four_conditions}.

Both decompositions measure optimization error relative to the infimum of the empirical risk. The approximation-based decomposition in \S\ref{sec:new_decomp} instead compares the algorithm output with an explicitly constructed approximating network, which is used in the complete regression and DRM analyses.

\subsection{Approximation error: scope of this survey}\label{sec:approx_landscape}

This subsection sketches the approximation-theory literature that informs the deep-learning analyses developed downstream. The body of the review draws on only a small slice of this landscape; we record the rest for context.

Classical universal approximation results \cite{cybenko1989approximation,hornik1989multilayer,hornik1991approximation,pinkus1999approximation} establish that shallow neural networks can approximate continuous functions on compact sets to arbitrary precision under suitable activation assumptions. In particular, continuous non-polynomial activations have this property \cite{leshno1993multilayer}. These results are qualitative and do not specify approximation rates. Barron \cite{barron1993universal} obtained an $L^2$ approximation rate $O(n_{\mathrm{neurons}}^{-1/2})$ for functions satisfying a Fourier-moment condition, with the exponent independent of dimension. The constructions discussed below give quantitative rates on classical smoothness classes.

Beginning in the mid-2010s, attention shifted to deep ReLU networks and explicit rate bounds as a function of width $W$ and depth $L$. Yarotsky \cite{yarotsky2017error} showed that for $f$ in the unit ball of the integer-order Sobolev space $W^{s,\infty}([0,1]^d)$, $s\in\mathbb N^+$, ReLU networks achieve $L^\infty$ error $O(n_{\mathrm{params}}^{-s/d})$ up to logarithmic factors via the $x^2$ sawtooth construction (developed in detail in \S\ref{sec:classical_approx}). Subsequent work \cite{shen2020deep,yarotsky2020phase,yarotsky2018optimal,lu2021deep} refined the approximation bounds by separately tracking $W$ and $L$; for H\"older order $0<\beta\le1$, Shen et al.\ \cite{shen2022optimal} obtained a logarithmic improvement to the rate $(WL)^{-2\beta/d}$ and established optimality in width and depth separately. Jiao et al.\ \cite{jiao2021deep} give the higher-smoothness approximation bound with a polynomial-in-$d$ error prefactor that we adopt as Theorem~\ref{thm:A1}. Extensions to piecewise-smooth, shift-invariant, and band-limited targets are in \cite{petersen2018optimal,yang2022approximation,montanelli2019deep}; comprehensive surveys are \cite{petersen2020neural,devore2021neural}.

PDE methods such as DRM require approximation of derivatives as well as function values. ReLU approximation in first-order Sobolev norms is studied in \cite{guhring2020error}. For smooth activations, G\"uhring and Raslan \cite{guhring2021approximation} combine averaged Taylor polynomials with an approximate partition of unity; the tanh specialization used in \cite{jiao2021error} is stated in Theorem~\ref{thm:A2}. Related constructions appear in \cite{deryck2021approximation,jiao2024drm}. Piecewise-polynomial activations provide another route: $\mathrm{ReLU}^k(x):=\max\{0,x\}^k$ supports polynomial constructions and derivative approximation \cite{li2019better}. Duan et al.\ \cite{Duan2022} use $\mathrm{ReLU}^2$ for DRM, Jiao et al.\ \cite{jiao2022rate} use $\mathrm{ReLU}^3$ for PINNs, and Lu et al.\ \cite{lu2021machine} obtain localized statistical bounds with $\mathrm{ReLU}^3$. Piecewise-polynomial networks admit useful pseudo-dimension estimates \cite{bartlett2019nearly}; finite bounds also exist for sigmoidal networks \cite{karpinski1997polynomial}, while the tanh analyses presented here use parameter covering. Extensions to Besov and other smoothness classes are studied in \cite{suzuki2019adaptivity,schmidt2019deep,li2019better,montanelli2019new,jiao2023improved}.

Another approach controls approximation and statistical complexity through constraints on weight norms \cite{neyshabur2015norm,bartlett2017spectrally,golowich2018size}. When the parameter count is large relative to the sample size, such constraints can yield useful bounds even when parameter-count bounds are vacuous. Definition~\ref{sec:overparam_setup} defines the norm product $\kappa(\boldsymbol{\theta})$, including the treatment of biases. For $1$-Lipschitz activations, this product bounds the network's Lipschitz constant with respect to the $\ell^\infty$ input norm. The approximation results for ReLU networks \cite{jiao2023approximation} and linear combinations of tanh sub-networks \cite{ding2025convergence} control approximation accuracy together with the relevant norm budget; they are stated in \S\ref{sec:overparam_approx}.

Other constructions use activations such as sine, floor, and exponential functions to improve the dimension dependence of approximation bounds \cite{yarotsky2020phase,yarotsky2021elementary,shen2021neural,shen2021deepnn}. For example, Jiao et al.\ \cite{jiao2023deep} construct fixed-depth ReLU-sine-$2^x$ networks approximating functions in $\mathcal{H}^{\beta}([0,1]^d,B)$, $0<\beta\le1$, in $L^p$ for $1\le p<\infty$, with width growing logarithmically in the inverse accuracy at fixed $d,\beta,B$. Approximation results also cover convolutional \cite{zhou2020universality,6ZHOU2020319}, recurrent \cite{schafer2006recurrent,jiao2024approximationboundsrecurrentneural}, and self-attention \cite{yun2019transformers} architectures.

\medskip

The running regression example uses ReLU approximation in $L^\infty$, which also controls $L^2(P_X)$ error. The running DRM example uses tanh Sobolev approximation following \cite{jiao2021error,ding2025convergence,jiao2024drm}, with the $\mathrm{ReLU}^3$ result of \cite{lu2021machine} presented separately. Chapters~\ref{sec:classical_setting}--\ref{sec:complete_setting} track the size, weight constraints, and architecture required by each approximation result.

\subsection{Statistical error: scope of this survey}\label{sec:stat_landscape}

This subsection records the statistical tools developed in machine-learning theory that the rate analyses of \S\ref{sec:classical_setting}--\S\ref{sec:complete_setting} draw on, together with the broader empirical-process literature.

The classical theory of uniform convergence begins with Vapnik and Chervonenkis \cite{vapnik2015uniform}, who introduced the VC dimension as a combinatorial measure of capacity for binary classifiers; the pseudo-dimension extension to real-valued classes is in \cite{anthony2009neural}. Parallel developments in empirical process theory \cite{VanJhon,van2000empirical,gine2021mathematical,cucker2002mathematical} established symmetrization, contraction, and chaining as the standard toolkit; these techniques translate uniform-convergence questions into combinatorial complexity bounds (covering numbers, bracketing numbers, $\gamma_2$-functionals). Bartlett and Mendelson \cite{bartlett2002rademacher} subsequently developed the Rademacher complexity as a sample-dependent measure of class complexity, which can retain information about the sampled inputs; the framework reduces statistical analysis to bounding the Rademacher complexity, which can be controlled by Dudley's chaining inequality (Proposition~\ref{prop:dudley}) given a covering-number estimate. This is the path we follow in the classical-regime statistical analysis of \S\ref{sec:classical_stat}.

Fast statistical bounds use additional structure beyond global uniform convergence. A variance--mean condition, together with suitable complexity control, can yield a complexity-over-$n$ bound instead of its square root. Local Rademacher methods \cite{bartlett2005local,koltchinskii2006local} control classes near the risk minimizer through a fixed-point argument. Peeling divides the class into risk or variance shells and combines concentration bounds across them. Offset Rademacher methods \cite{liang2015learning} incorporate a negative quadratic term into the empirical process. Section~\ref{sec:four_conditions} describes the bounded-loss setting, and Section~\ref{sec:classical_rates} applies these tools to regression and DRM. The final convergence rate also depends on how the class complexity grows with the approximation accuracy.

Size-based bounds such as $\mathrm{Pdim}(\mathcal{NN}_{\mathrm{ReLU}}(W,L))\le C\mathfrak D L\log(e\mathfrak D)$ \cite{bartlett2019nearly} can become vacuous when the parameter count is large relative to the sample size. Norm-controlled bounds provide an alternative \cite{neyshabur2015norm,bartlett2017spectrally,golowich2018size}. The Rademacher bound of Golowich et al.\ used in \S\ref{sec:overparam_stat} is independent of width under prescribed norm and depth conditions, with depth dependence retained explicitly.

For DRM, the empirical loss contains both function values and gradients, so statistical bounds must control the corresponding derivative classes. Global bounds are developed for $\mathrm{ReLU}^2$ in \cite{Duan2022} and for tanh in \cite{jiao2021error}. The PINN analysis of \cite{jiao2022rate} instead treats a residual loss involving second derivatives. Lu et al.\ \cite{lu2021machine} use localized complexity bounds for DRM and PINNs in a Dirichlet setting, with the respective stability estimates stated in the norms appropriate to those objectives. Section~\ref{sec:drm_fast} discusses the tanh slow-rate analysis and presents the $\mathrm{ReLU}^3$ result separately; Section~\ref{sec:minimax_lb} distinguishes the observation models used in the lower-bound comparisons.

\medskip

The following chapters use size-based covering bounds (\S\ref{sec:classical_stat}), norm-controlled Rademacher bounds (\S\ref{sec:overparam_stat}), and complexity estimates for parallel networks (\S\ref{sec:complete_stat}).

\subsection{Optimization error: scope of this survey}\label{sec:opt_landscape}

Optimization analysis studies the empirical-risk gap attained by a specified training algorithm. We introduce NTK and mean-field approaches before describing the PGD-on-PNN analysis used in Chapter~\ref{sec:complete_setting}.

The empirical risk is generally non-convex as a function of the network parameters. Convergence to a stationary point alone does not control the gap to its global infimum. Optimization guarantees therefore depend on assumptions about the architecture, initialization, loss, and data. The rate analyses in Chapters~\ref{sec:classical_setting}--\ref{sec:overparam_setting} set $\mathcal{E}_{opt}=0$ and study empirical minimizers without specifying how they are computed.

In the NTK regime, suitable scaling and initialization keep the training dynamics close to those of the network linearized at initialization \cite{jacot2018neural,chizat2019lazy}. For squared-loss training, positive lower bounds on the appropriately normalized empirical NTK Gram matrix along the optimization trajectory yield a Polyak--{\L}ojasiewicz (PL) inequality of the form \cite{du2019gradient,allen2019convergence,zou2019improved,liu2022loss}
$$\big\|\nabla \widehat{\mathcal{L}}_n(\boldsymbol{\theta})\big\|^2 \;\ge\; 2c_{\mathrm{PL}}\big(\widehat{\mathcal{L}}_n(\boldsymbol{\theta})-\widehat{\mathcal{L}}_n^*\big).$$
Here $\widehat{\mathcal{L}}_n^*:=\inf_{\boldsymbol{\theta}}\widehat{\mathcal{L}}_n(\boldsymbol{\theta})$, and $c_{\mathrm{PL}}>0$ depends on the kernel lower bound and the normalization of the loss. For sufficiently wide networks, the cited analyses control kernel variation along training under their initialization and data assumptions. The PL inequality, together with an appropriate descent estimate and step size, gives geometric decay of the optimization error. Related analyses bound population risk \cite{arora2019fine,cao2019generalization,allen2019learning}; in particular, Arora et al.\ \cite{arora2019fine} obtain a data-dependent generalization bound independent of network size. The fixed-kernel limit describes training without feature evolution, whereas finite-width guarantees also control deviations from that limit.

The mean-field viewpoint represents the network parameters by a probability measure. Under suitable scaling, the dynamics of wide two-layer networks are described by an evolution equation for that measure, often a Wasserstein gradient flow \cite{mei2018mean,chizat2018global,rotskoff2018neural,sirignano2020mean}. This description allows features to evolve during training. Chizat et al.\ \cite{chizat2019lazy} show how model scaling controls proximity to linearized dynamics; width alone does not determine whether training is lazy. Convergence conclusions depend on the particular dynamics, regularization, and initialization assumptions.

NTK and mean-field analyses provide optimization and generalization guarantees under their respective assumptions. This review develops the PGD-on-PNN approach because the selected regression and DRM results combine constructive approximation, statistical estimates, and algorithmic guarantees within the same analysis \cite{kohler2023opt,jiao2024drm}.

A parallel neural network is a linear combination of independently parameterized sub-networks. For the squared and energy losses considered here, the empirical objective is convex in the outer coefficients when the inner parameters are fixed. The joint-training analysis combines this convexity with bounds on random initialization and inner-parameter movement, while projection controls the parameter constraints. Section~\ref{sec:complete_opt} develops these estimates.

Other optimization questions include convergence to stationary points and implicit regularization; see \cite{sun2020optimization,bottou2018optimization} for broader surveys.

\medskip

Chapter~\ref{sec:complete_setting} presents the PGD-on-PNN guarantees of \cite{kohler2023opt,jiao2024drm} for the two running problems.

\subsection{Running examples: regression and DRM}\label{sec:running_examples}

We now specify the regression and DRM objectives, sampling models, and target errors used in Chapters~\ref{sec:classical_setting}--\ref{sec:complete_setting}.

Throughout, $\mathcal{H}$ denotes the target-quantity space: $L^2(P_X)$ for Regression and $H^1(\Omega)$ for DRM.

\subsubsection{Regression}\label{sec:reg_setup}
The nonparametric regression model considered throughout this review is
\begin{equation}\label{eq:reg_model}
Y = f_0(\boldsymbol{X}) + \xi,
\end{equation}
where $\boldsymbol{X} \in [0,1]^d$ is a $d$-dimensional predictor, $Y \in \mathbb{R}$ is a response, $f_0 : [0,1]^d \to \mathbb{R}$ is the unknown regression function (the target), and $\xi$ is a noise term independent of $\boldsymbol{X}$ with $\mathbb{E}[\xi] = 0$ and finite second moment. The bounded-loss analyses of \S\ref{sec:reg_slow}--\S\ref{sec:reg_fast} additionally assume bounded noise and output-clipped networks, as specified there; extensions to sub-Gaussian noise require the corresponding unbounded-loss or truncation arguments. We observe an i.i.d.\ sample $S_n = \{(\boldsymbol{X}_i, Y_i)\}_{i=1}^n$, with marginal distribution $P_X$ for $\boldsymbol{X}$.

The population risk (squared error) and empirical risk are
\[
\mathcal{L}_{\mathrm{REG}}(f) := \mathbb{E}\big[(Y - f(\boldsymbol{X}))^2\big], \qquad \widehat{\mathcal{L}}_{\mathrm{REG}}(f) := \frac{1}{n}\sum_{i=1}^n (Y_i - f(\boldsymbol{X}_i))^2.
\]
The Bayes-optimal predictor is $f_0$ itself. When the empirical minimum on $\mathcal{F}_{NN}$ is attained, select an ERM
\begin{equation}\label{eq:reg_erm}
\widehat{f}_{\boldsymbol{\theta}} \in \arg\min_{f \in \mathcal{F}_{NN}} \widehat{\mathcal{L}}_{\mathrm{REG}}(f).
\end{equation}
Writing $f_{\boldsymbol{\theta}_{\mathcal{A}}}$ for the network actually produced by the training algorithm $\mathcal{A}$ (the role of $\phi_{\boldsymbol{\theta}_{\mathcal{A}}}$ in \S\ref{sec:decomp}), the target quantity to be bounded is its expected excess $L^2$ risk
\[
\mathbb{E}_{S_n}\big[\mathcal{L}_{\mathrm{REG}}(f_{\boldsymbol{\theta}_{\mathcal{A}}}) - \mathcal{L}_{\mathrm{REG}}(f_0)\big] = \mathbb{E}_{S_n}\|f_{\boldsymbol{\theta}_{\mathcal{A}}} - f_0\|^2_{L^2(P_X)}.
\]
The symmetric decomposition bounds this quantity through the three errors. In Chapters~\ref{sec:classical_setting}--\ref{sec:overparam_setting}, zero optimization error means that the algorithm output is an empirical minimizer, denoted by $\widehat{f}_{\boldsymbol{\theta}}$. Chapter~\ref{sec:complete_setting} analyzes the output of projected gradient descent.

For each fixed comparison function $\bar f\in\mathcal{F}_{NN}$, unbiasedness gives $\mathbb{E}_{S_n}[\widehat{\mathcal{L}}_{\mathrm{REG}}(\bar f)-\mathcal{L}_{\mathrm{REG}}(\bar f)]=0$. Using this identity before bounding the remaining fluctuation, and then taking the infimum over $\bar f$, sharpens the symmetric bound to
\begin{align}\label{eq:reg_three_decomp}
\mathbb{E}_{S_n}\|f_{\boldsymbol{\theta}_{\mathcal{A}}} - f_0\|^2_{L^2(P_X)} &\le \underbrace{\inf_{\bar{f} \in \mathcal{F}_{NN}} \|\bar{f} - f_0\|^2_{L^2(P_X)}}_{\mathcal{E}^{\mathrm{REG}}_{app}} + \underbrace{\mathbb{E}_{S_n}\!\!\sup_{f \in \mathcal{F}_{NN}}\big|\mathcal{L}_{\mathrm{REG}}(f) - \widehat{\mathcal{L}}_{\mathrm{REG}}(f)\big|}_{\mathcal{E}^{\mathrm{REG}}_{sta,\,slow}} \\
& ~~~~ + \underbrace{\mathbb{E}_{S_n}\!\big[\widehat{\mathcal{L}}_{\mathrm{REG}}(f_{\boldsymbol{\theta}_{\mathcal{A}}}) - \widehat{\mathcal{L}}_{\mathrm{REG}}(\widehat{f}_{\boldsymbol{\theta}})\big]}_{\mathcal{E}^{\mathrm{REG}}_{opt}}. \notag
\end{align}
For an empirical minimizer, $\mathcal{E}^{\mathrm{REG}}_{opt}=0$, and \eqref{eq:reg_three_decomp} reduces to
\begin{equation}\label{eq:reg_sym_decomp}
\mathbb{E}_{S_n}\|\widehat{f}_{\boldsymbol{\theta}} - f_0\|^2_{L^2(P_X)} \le \mathcal{E}^{\mathrm{REG}}_{app} + \mathcal{E}^{\mathrm{REG}}_{sta,\,slow},
\end{equation}
while the asymmetric decomposition \eqref{eq:asym_decomp} (likewise with $\mathcal{E}_{opt}=0$) gives the fast-rate form
\begin{equation}\label{eq:reg_asym_decomp}
\mathbb{E}_{S_n}\|\widehat{f}_{\boldsymbol{\theta}} - f_0\|^2_{L^2(P_X)} \le 2\,\mathcal{E}^{\mathrm{REG}}_{app} + \mathcal{E}^{\mathrm{REG}}_{sta,\,fast},
\end{equation}
with the offset statistical term
\begin{equation}\label{eq:reg_offset_def}
\mathcal{E}^{\mathrm{REG}}_{sta,\,fast} \;:=\; \mathbb{E}_{S_n}\sup_{f \in \mathcal{F}_{NN}}\big[\mathcal{L}_{\mathrm{REG}}(f_0) - 2\widehat{\mathcal{L}}_{\mathrm{REG}}(f) + \mathcal{L}_{\mathrm{REG}}(f)\big].
\end{equation}
For an empirical minimizer, \eqref{eq:reg_asym_decomp} follows from \eqref{eq:asym_decomp}. Proposition~\ref{prop:S4} controls the offset supremum under the conditions of Proposition~\ref{prop:fast_engine}. Chapter~\ref{sec:complete_setting} instead uses a constructed approximating network to analyze the optimization contribution.
Throughout, we assume the smoothness condition $f_0 \in \mathcal{H}^\beta([0,1]^d, B_{\mathrm{REG}})$ (Hölder class with parameter $\beta \ge 1$ and constant $B_{\mathrm{REG}}$). Convergence rates derived from these decompositions appear in \S\ref{sec:reg_slow} (slow rate) and \S\ref{sec:reg_fast} (fast rate).

\subsubsection{Deep Ritz method}\label{sec:drm_setup}

For PDE solving, we consider the elliptic equation on a bounded open set $\Omega \subset [0,1]^d$ with Neumann boundary conditions:
\begin{equation}\label{eq:drm_pde}
-\Delta u + \omega u = h \ \text{ on } \Omega, \qquad \frac{\partial u}{\partial \boldsymbol{n}} = g \ \text{ on } \partial\Omega,
\end{equation}
under the standard regularity assumptions $\partial\Omega \in C^{2+\alpha}$, $h \in L^\infty(\Omega)$, $g \in L^\infty(\partial\Omega) \cap H^{1/2}(\partial\Omega)$, $\omega \in C^{\alpha}(\bar\Omega)$, $\omega \ge c_0 > 0$ for some $0 < \alpha < 1$. Under these assumptions, \eqref{eq:drm_pde} admits a unique weak solution $u_0 \in H^2(\Omega)$ \cite{agmon1959estimates}, and the constant $B_{\mathrm{DRM}} := \max\{\|h\|_{L^\infty}, \|g\|_{L^\infty}, \|\omega\|_{L^\infty}\}$ is finite; the boundedness of $h, g, \omega$ is required for the Bernstein-type concentration steps invoked in the statistical analysis of \S\ref{sec:classical_stat}.

The Deep Ritz method minimizes the energy functional
\begin{equation}\label{eq:drm_energy}
\mathcal{L}_{\mathrm{DRM}}(u) := \int_\Omega \!\Big(\tfrac{1}{2}\|\nabla u\|^2 + \tfrac{1}{2}\omega u^2 - hu\Big) \mathrm{d}\boldsymbol{x} - \int_{\partial\Omega} g\, Tu \,\mathrm{d}s,
\end{equation}
where $T$ denotes the trace operator. Discretizing both integrals via Monte Carlo sampling with $\{\boldsymbol{X}_p\}_{p=1}^{N_{\mathrm{in}}} \sim U(\Omega)$ i.i.d., $\{\boldsymbol{Y}_p\}_{p=1}^{N_b} \sim U(\partial\Omega)$ i.i.d., we collect the combined MC sample set $S := \{\boldsymbol{X}_p\}_{p=1}^{N_{\mathrm{in}}} \cup \{\boldsymbol{Y}_p\}_{p=1}^{N_b}$ (and $\mathbb{E}_S$ denotes expectation over both draws); the empirical energy is
\begin{equation}\label{eq:drm_emp}
\widehat{\mathcal{L}}_{\mathrm{DRM}}(u) := \tfrac{|\Omega|}{N_{\mathrm{in}}}\!\sum_{p=1}^{N_{\mathrm{in}}}\!\Big[\tfrac{\|\nabla u(\boldsymbol{X}_p)\|^2}{2}\!+\!\tfrac{\omega(\boldsymbol{X}_p)u(\boldsymbol{X}_p)^2}{2}\!-\!u(\boldsymbol{X}_p)h(\boldsymbol{X}_p)\Big] - \tfrac{|\partial\Omega|}{N_b}\!\sum_{p=1}^{N_b}\! u(\boldsymbol{Y}_p)g(\boldsymbol{Y}_p).
\end{equation}
When the empirical minimum on $\mathcal{F}_{NN}$ is attained, select an ERM
\begin{equation}\label{eq:drm_erm}
\widehat{u}_{\boldsymbol{\theta}} \in \arg\min_{u \in \mathcal{F}_{NN}} \widehat{\mathcal{L}}_{\mathrm{DRM}}(u).
\end{equation}
Writing $u_{\boldsymbol{\theta}_{\mathcal{A}}}$ for the network produced by the training algorithm (the role of $\phi_{\boldsymbol{\theta}_{\mathcal{A}}}$ in \S\ref{sec:decomp}), the target quantity to be bounded is its expected $H^1$-distance squared,
\[
\mathbb{E}_S\big[\|u_{\boldsymbol{\theta}_{\mathcal{A}}} - u_0\|_{H^1(\Omega)}^2\big].
\]
The energy--norm equivalence below relates this error to the excess energy. In Chapters~\ref{sec:classical_setting}--\ref{sec:overparam_setting}, zero optimization error means that the algorithm output is an empirical minimizer, denoted by $\widehat{u}_{\boldsymbol{\theta}}$. Chapter~\ref{sec:complete_setting} analyzes the output of projected gradient descent.

For regression, the excess risk equals the squared $L^2(P_X)$ error. For DRM, the following energy--norm equivalence converts an excess-energy bound into an $H^1$ error bound.

\begin{proposition}[Energy--norm equivalence for DRM; \citealp{Duan2022,lu2021machine}]\label{prop:drm_eq}
Let $u_0$ be the weak solution of \eqref{eq:drm_pde}, and let $\mathcal{L}_{\mathrm{DRM}}$ be defined by \eqref{eq:drm_energy}. For any $u \in H^1(\Omega)$,
\[
\frac{c_0 \wedge 1}{2}\,\|u - u_0\|^2_{H^1(\Omega)} \;\le\; \mathcal{L}_{\mathrm{DRM}}(u) - \mathcal{L}_{\mathrm{DRM}}(u_0) \;\le\; \frac{B_{\mathrm{DRM}} \vee 1}{2}\,\|u - u_0\|^2_{H^1(\Omega)}.
\]
\end{proposition}

\begin{proof}[Proof sketch]
Direct calculation using the variational characterization of $u_0$: for any $u$,
\[
\mathcal{L}_{\mathrm{DRM}}(u) - \mathcal{L}_{\mathrm{DRM}}(u_0) = \tfrac{1}{2}\int_\Omega \big(\|\nabla(u - u_0)\|^2 + \omega(u - u_0)^2\big)\mathrm{d}\boldsymbol{x},
\]
which is bounded between $\frac{c_0 \wedge 1}{2}$ and $\frac{B_{\mathrm{DRM}} \vee 1}{2}$ times $\|u - u_0\|_{H^1}^2$ by definition of $\omega$.
\end{proof}

The left-hand inequality is the coercivity estimate that produces the factor $2/(c_0\wedge1)$ below; the right-hand inequality bounds the approximation contribution by $(B_{\mathrm{DRM}}\vee1)/2$ times the squared $H^1$ approximation error. Applying the decomposition with a fixed comparison function and using the unbiasedness of its empirical energy gives
\begin{align}\label{eq:drm_three_decomp}
& \mathbb{E}_S\big[\|u_{\boldsymbol{\theta}_{\mathcal{A}}} - u_0\|_{H^1}^2\big] \le \frac{2}{c_0 \wedge 1}\bigg\{\underbrace{\frac{B_{\mathrm{DRM}}\vee 1}{2}\inf_{u \in \mathcal{F}_{NN}}\|u - u_0\|_{H^1}^2}_{\mathcal{E}^{\mathrm{DRM}}_{app}} \\ \notag
& ~~~~~~ + \underbrace{\mathbb{E}_S\!\!\sup_{u \in \mathcal{F}_{NN}}\!\big|\mathcal{L}_{\mathrm{DRM}}(u) - \widehat{\mathcal{L}}_{\mathrm{DRM}}(u)\big|}_{\mathcal{E}^{\mathrm{DRM}}_{sta,\,slow}} + \underbrace{\mathbb{E}_S\!\big[\widehat{\mathcal{L}}_{\mathrm{DRM}}(u_{\boldsymbol{\theta}_{\mathcal{A}}}) - \widehat{\mathcal{L}}_{\mathrm{DRM}}(\widehat{u}_{\boldsymbol{\theta}})\big]}_{\mathcal{E}^{\mathrm{DRM}}_{opt}}\bigg\}.
\end{align}
For an empirical minimizer, $\mathcal{E}^{\mathrm{DRM}}_{opt}=0$, and this reduces to
\begin{equation}\label{eq:drm_sym_decomp}
\mathbb{E}_S\big[\|\widehat{u}_{\boldsymbol{\theta}} - u_0\|_{H^1}^2\big] \le \frac{2}{c_0 \wedge 1}\big\{\mathcal{E}^{\mathrm{DRM}}_{app} + \mathcal{E}^{\mathrm{DRM}}_{sta,\,slow}\big\}.
\end{equation}
The asymmetric decomposition also applies to the DRM energy, since Monte Carlo quadrature is unbiased at each fixed comparison function. Combining \eqref{eq:asym_decomp} with the left inequality of Proposition~\ref{prop:drm_eq} gives, under $\mathcal{E}^{\mathrm{DRM}}_{opt}=0$,
\begin{equation}\label{eq:drm_asym_decomp}
\mathbb{E}_S\big[\|\widehat{u}_{\boldsymbol{\theta}} - u_0\|_{H^1}^2\big] \;\le\; \frac{2}{c_0 \wedge 1}\big\{2\,\mathcal{E}^{\mathrm{DRM}}_{app} + \mathcal{E}^{\mathrm{DRM}}_{sta,\,fast}\big\},
\end{equation}
with the offset statistical term
\begin{equation}\label{eq:drm_offset_def}
\mathcal{E}^{\mathrm{DRM}}_{sta,\,fast} \;:=\; \mathbb{E}_S\sup_{u \in \mathcal{F}_{NN}}\big[\mathcal{L}_{\mathrm{DRM}}(u) + \mathcal{L}_{\mathrm{DRM}}(u_0) - 2\,\widehat{\mathcal{L}}_{\mathrm{DRM}}(u)\big].
\end{equation}
The supremand is twice the population--empirical fluctuation at $u$ minus its excess energy. Coercivity therefore supplies a negative quadratic term in the $H^1$ error. Obtaining a fast statistical rate additionally requires control of the fluctuations and local complexity of the energy-loss class. Section~\ref{sec:drm_fast} presents the $\mathrm{ReLU}^3$ fast-rate result of \cite{lu2021machine} under that paper's Dirichlet setting, separately from the Neumann setup used here.

Throughout, we assume the smoothness condition $u_0 \in H^s(\Omega)$ for some $s \ge 2$. Convergence rates derived from these decompositions appear in \S\ref{sec:drm_slow} (slow rate) and \S\ref{sec:drm_fast} (fast rate). The tanh slow-rate analysis uses $s=2$; higher regularity enters through the $\mathrm{ReLU}^3$ analysis of Theorem~\ref{thm:drm_fast} and the PNN analysis of \S\ref{sec:drm_complete}.

\subsection{Notation summary}\label{sec:notation_summary}

The review uses a moderately large set of symbols across approximation, statistical, and optimization analyses. General conventions and function-space definitions are collected here; detailed symbol tables organized by topic appear in Appendix~\ref{app:notation}.

\par\textbf{General mathematical conventions.}
Throughout the review, bold lowercase letters $\boldsymbol{x}, \boldsymbol{b}, \boldsymbol{\theta}$ denote vectors; bold uppercase letters $\boldsymbol{A}, \boldsymbol{X}$ denote matrices and vector-valued random variables; italic letters $L, W, n, d$ denote scalars and dimensions; calligraphic letters $\mathcal{F}, \mathcal{H}, \mathcal{L}, \mathcal{NN}, \mathcal{PNN}$ denote function classes and risks; sans-serif/Fraktur letters $\mathfrak{m}, \mathfrak{D}, \mathfrak{R}, \mathfrak{M}$ denote counts, complexities, and minimax quantities. The notation $a \asymp b$ means $a = \Theta(b)$, $a \lesssim b$ means $a = O(b)$, and $a \vee b = \max(a, b)$, $a \wedge b = \min(a, b)$. The notation $\tilde{\mathcal{O}}(\cdot)$ hides logarithmic factors; $\operatorname{polylog}(n)$ denotes an unspecified power of $\log n$. Multi-index notation $\boldsymbol{x}^{\boldsymbol{\alpha}} = x_1^{\alpha_1}\cdots x_d^{\alpha_d}$ with $\|\boldsymbol{\alpha}\|_1 = \sum_i \alpha_i$ is used throughout.

\par\textbf{Sobolev spaces.} For any $s\in\mathbb N_0$ and $1\le p\le\infty$, the \textit{Sobolev space} $W^{s,p}(\Omega)$ is defined using weak derivatives by
$$
W^{s,p}(\Omega):=\left\{f\in L^p(\Omega): D^{\boldsymbol{\alpha}}f\in L^p(\Omega), \ \forall \boldsymbol{\alpha}\in \mathbb{N}_0^{d} \text{ with } \|\boldsymbol{\alpha}\|_1 \leq s\right\}.
$$
When $p=2$, we write $H^s(\Omega):=W^{s,2}(\Omega)$ for any $s\in\mathbb{N}_0$. For $f\in W^{s,p}(\Omega)$ with $1\leq p <\infty$, the Sobolev norm is
$$\|f\|_{W^{s,p}(\Omega)}:=\left(\sum_{0\leq \|\boldsymbol{\alpha}\|_1\leq s} \|D^{\boldsymbol{\alpha}}f\|^p_{L^p(\Omega)}\right)^{\frac{1}{p}},$$
and for $p=\infty$,
$$\|f\|_{W^{s,\infty}(\Omega)} := \max_{0\leq \|\boldsymbol{\alpha}\|_1\leq s}\|D^{\boldsymbol{\alpha}}f\|_{L^{\infty}(\Omega)}.$$

\par\textbf{Hölder spaces.} Let $\beta=s+r>0$, with $s=\lceil\beta\rceil-1\in\mathbb N_0$ and $r=\beta-s\in(0,1]$. For a finite constant $B>0$, the bounded \textit{Hölder class} $\mathcal{H}^\beta(\Omega, B)$ is
\begin{equation}\label{eq0}
\big\{ f:\Omega \rightarrow \mathbb{R}, \max _{\|\boldsymbol{\alpha}\|_1 \leq s}\left\|D^{\boldsymbol{\alpha}} f\right\|_{L^{\infty}} \leq B, \max _{\|\boldsymbol{\alpha}\|_1=s} \sup _{\boldsymbol{x} \neq \boldsymbol{y}} \frac{\left|D^{\boldsymbol{\alpha}} f(\boldsymbol{x})-D^{\boldsymbol{\alpha}} f(\boldsymbol{y})\right|}{\|\boldsymbol{x}-\boldsymbol{y}\|_2^r} \leq B\big\},
\end{equation}
with the Hölder norm
$$\|f\|_{\mathcal{H}^{\beta}(\Omega)} := \sum_{\|\boldsymbol{\alpha}\|_1 \leq s}\|D^{\boldsymbol{\alpha}}f\|_{L^\infty(\Omega)} + \sum_{\|\boldsymbol{\alpha}\|_1=s} \sup_{\substack{\boldsymbol{x},\boldsymbol{y}\in\Omega\\\boldsymbol{x}\neq\boldsymbol{y}}} \frac{\left|D^{\boldsymbol{\alpha}} f(\boldsymbol{x})-D^{\boldsymbol{\alpha}} f(\boldsymbol{y})\right|}{\|\boldsymbol{x}-\boldsymbol{y}\|_2^r}.$$
For any $f \in \mathcal H^\beta(\Omega,B)$, all partial derivatives through order $s=\lceil\beta\rceil-1$ exist, and the derivatives of order $s$ are $r$-H\"older continuous. At an integer exponent $\beta=k$, this is the $C^{k-1,1}$ convention; classical derivatives of order $k$ need not exist everywhere.

%% file: sections/classical.tex
\section{Classical Error Analysis via Size Control}\label{sec:classical_setting}

This chapter establishes the theoretical baseline for deep learning by isolating the fundamental trade-off between approximation power and statistical generalization. The classical regime predates the modern focus on over-parameterization and training dynamics. It operates under two structural assumptions. First, we assume that the empirical risk minimizer is found exactly ($\mathcal{E}_{opt} = 0$). This decouples the analysis from gradient descent, so the network can be evaluated purely as a statistical hypothesis class. Second, the hypothesis class is controlled by network size (width $W$, depth $L$, and nonzero-parameter count $\mathfrak{D}$), with output clipping for regression and an entrywise weight bound for the tanh DRM analysis.

Under these two assumptions, the error decomposition reduces to a strict two-error balance. This chapter traces the complete analytical pipeline for this classical regime: we first develop the approximation tools (\S\ref{sec:classical_approx}) and the size-based statistical tools (\S\ref{sec:classical_stat}) in their general, example-independent form, then instantiate both on the two running examples and couple them into the slow and fast convergence rates of each example (\S\ref{sec:classical_rates}), with the sharper external $\mathrm{ReLU}^3$ benchmark for the Deep Ritz fast rate recorded alongside. Section~\ref{sec:minimax_lb} then presents minimax lower bounds under specified observation models and explains when they provide comparable statistical benchmarks. Claims of optimality require matching observation models, target classes, and error metrics; in particular, the lower bound for recovery from noisy right-hand-side observations is not directly a lower bound for the Monte Carlo DRM setting. The chapter closes by replacing the ambient dimension $d$ with an intrinsic dimension $d^*$ when the data concentrate on a low-dimensional manifold (\S\ref{sec:low_intrinsic}), a technique used in several of the applications in \S\ref{sec:applications}.

\subsection{Approximation tools for this regime}\label{sec:classical_approx}
The downstream rate analyses of \S\ref{sec:classical_rates} rest on two approximation theorems, drawn from existing work and adapted to the two running examples. For regression, where the error is measured in $L^2$, we use the quantitative ReLU $L^\infty$ approximation bound of Jiao et al.\ \cite{jiao2021deep} and review Yarotsky's $x^2$ sawtooth construction. For the Deep Ritz Method, where the error is measured in $H^1$, we use the Sobolev approximation framework of G\"uhring and Raslan, presented here in its tanh instantiation. The two constructions are juxtaposed to show how the choice of target norm ($L^2$ vs.\ $H^1$) leads to different approximation requirements. The subsection closes with a brief survey of other architectures (CNN, RNN, Transformer). The theorems are instantiated on the running examples at the point of use in \S\ref{sec:classical_rates}.

\subsubsection{ReLU approximation of Hölder functions}\label{sec:A1}

The fundamental ReLU approximation theorem in the form used throughout this review is due to \cite{jiao2021deep}. Its statement parameterizes the network by two integer parameters $J$ and $I$, intended to control the width and depth respectively, and gives an error prefactor polynomial in the input dimension $d$ at fixed smoothness.

\begin{theorem}[ReLU approximation; \citealp{jiao2021deep}]\label{thm:A1}
Let $\beta = s + r$ with $s \in \mathbb{N}_0$, $r \in (0,1]$, and $f \in \mathcal{H}^\beta([0,1]^d, B_{\mathrm{REG}})$. For any $J, I \in \mathbb{N}^+$, there exists a ReLU network $\phi \in \mathcal{NN}_{\mathrm{ReLU}}(W, L)$ with width
\[
W = 38(s + 1)^2 \cdot 3^d \cdot d^{s + 1} \cdot J \lceil \log_2(8J) \rceil
\]
and depth (including the output affine map, unlike the hidden-layer count in the source)
\[
L = 21(s + 1)^2 \cdot I \lceil \log_2(8I) \rceil + 2d + 1
\]
such that
\[
\|f - \phi\|_{L^\infty([0,1]^d)} \;\le\; 19 B_{\mathrm{REG}} (s + 1)^2 \cdot d^{s + (\beta \vee 1)/2} \cdot (JI)^{-2\beta/d}.
\]
\end{theorem}

The line of work culminating in Theorem~\ref{thm:A1} unfolds in four stages, each of which sharpens a distinct aspect of the rate-vs-architecture trade-off. \textbf{(i) Yarotsky 2017} \cite{yarotsky2017error} quantified the approximation benefit of depth: for $f$ in the unit ball of the integer-order Sobolev space $W^{s,\infty}([0,1]^d)$, $s\in\mathbb N^+$, $\epsilon$-approximation in $L^\infty$ is achievable by a ReLU network with $O(\epsilon^{-d/s}\log(1/\epsilon))$ parameters and depth $O(\log(1/\epsilon))$ (with a constant $c=c(d,s)$ in front), via the iterated $x^2$ sawtooth construction (Proposition~2 there) presented in detail below. Equivalently, expressed in terms of network size, a network with total parameter count $\mathfrak{D}$ achieves an approximation error of $\mathcal{O}(\mathfrak{D}^{-s/d})$ (up to logarithmic factors). \textbf{(ii) Shen-Yang-Zhang 2020} \cite{shen2020deep} treated Hölder $\mathcal{H}^\beta$ with $\beta\in(0,1]$ and gave a non-asymptotic rate that separately tracks width $W$ and depth $L$, of order $\sqrt{d}\,W^{-2\beta/d} L^{-2\beta/d}$ in $L^p$, combining a bit-extraction encoding with the sawtooth iteration. \textbf{(iii) Shen-Yang-Zhang 2022} \cite{shen2022optimal} sharpened this, still for $0<\beta\le1$, to $\sqrt{d}\,(W^2 L^2 \log_3(W+2))^{-\beta/d}$ and proved matching lower bounds in the width-fixed and depth-fixed regimes separately, certifying the rate as optimal in $W$ and $L$ individually. The $C^s$-smooth extension is given by Lu-Shen-Yang-Zhang \cite{lu2021deep} as $\|f\|_{C^s}\,J^{-2s/d} I^{-2s/d}$ at width $W = O(J\log(eJ))$ and depth $L = O(I\log(eI)) + 2d$, but with a smoothness prefactor $C(s,d) = 85(s+1)^d 8^s$ that is exponential in $d$. \textbf{(iv) Jiao-Shen-Lin-Huang 2023} \cite[Corollary~3.1]{jiao2021deep}, the source of Theorem~\ref{thm:A1} as stated, improves on the exponential-in-$d$ error prefactor of the higher-order smooth approximation result of Lu et al.: the factor $d^{s + (\beta\vee 1)/2}$ is polynomial in $d$ at fixed smoothness $\beta$. This concerns the error prefactor, not the entire architecture: the width in Theorem~\ref{thm:A1} still contains $3^d$. The displayed approximation bound is the form invoked downstream in \S\ref{sec:reg_slow}--\S\ref{sec:reg_fast}.

A parallel line of work studies approximation in Sobolev norms, including first-order norms for ReLU networks \cite{guhring2020error} and approximation with more general activations \cite{guhring2021approximation}. These results control derivatives of the approximant as well as function values, under higher-order smoothness assumptions on the target. The first-order norm $H^1=W^{1,2}$ is the relevant one for our DRM example.

The sawtooth construction of \cite{yarotsky2017error} provides an elementary building block for ReLU approximation. We first give its realization of $x^2$, then explain how approximate multiplication and local polynomials enter function approximation. The precise width--depth allocation in Theorem~\ref{thm:A1} is given in \cite{jiao2021deep}.

\medskip\noindent\textbf{Stage 1: Yarotsky's sawtooth construction of $x^2$.}

The construction begins with the hat function
\[
g(x) = \begin{cases} 2x, & 0 \le x \le 1/2, \\ 2(1-x), & 1/2 < x \le 1, \end{cases}
\]
which has the closed-form ReLU representation
\[
g(x) = 2\rho(x) - 4\rho(x - 1/2) + 2\rho(x - 1), \qquad \rho(z) := \max\{0, z\}.
\]
This expresses $g$ using one hidden ReLU layer of width $3$ and an affine output map, hence depth $L=2$.

Iterating $g$ produces sawtooth functions of doubled frequency: define 
\[
g_k := \underbrace{g \circ \cdots \circ g}_{k \text{ times}}\, .
\]
Each $g_k$ is a piecewise-linear function on $[0,1]$ with $2^{k-1}$ teeth, each of height $1$ and base width $2^{-k+1}$. As a ReLU network, $g_k$ has $k$ hidden layers of width $3$, hence depth $L=k+1$. The key identity is the following representation of $x^2$ as a convergent series.

\begin{lemma}[Yarotsky's sawtooth identity]\label{lem:saw}
For all $x \in [0,1]$,
\[
x^2 = x - \sum_{k=1}^{\infty} \frac{g_k(x)}{4^k}.
\]
The truncation
\[
\widetilde{f}_J(x) := x - \sum_{k=1}^J \frac{g_k(x)}{4^k}
\]
satisfies $\|x^2 - \widetilde{f}_J\|_{L^\infty([0,1])} \le 4^{-(J+1)}$.
\end{lemma}

\begin{proof}
On each dyadic sub-interval $\big[(j-1)\cdot 2^{-J},\, j\cdot 2^{-J}\big]$ with $j \in \{1, \ldots, 2^J\}$, the function $\widetilde{f}_J$ is linear and agrees with $x^2$ at the endpoints. Since $\widetilde{f}_J$ is therefore the piecewise-linear interpolant of $x^2$ on a uniform grid of mesh size $2^{-J}$, the standard interpolation-error formula with $f''(x) = 2$ yields $\|x^2 - \widetilde{f}_J\|_\infty \le (2^{-J})^2 \cdot 2 / 8 = 4^{-J}/4 = 4^{-(J+1)}$. The pointwise identity follows by letting $J \to \infty$.
\end{proof}

\begin{proposition}[ReLU network for $x^2$]\label{prop:relu_xsq}
For any $\epsilon \in (0, 1)$, there exists a ReLU network $\phi_{\mathrm{sq}}$ of depth at most $\lceil\log_2(1/\epsilon)/2\rceil + 1$ and width at most $5$ such that
\[
\|x^2 - \phi_{\mathrm{sq}}\|_{L^\infty([0,1])} \le \epsilon.
\]
\end{proposition}

\begin{proof}
By Lemma~\ref{lem:saw}, choose $J$ such that $4^{-(J+1)} \le \epsilon$, i.e., $J \ge \log_2(1/\epsilon)/2 - 1$. To realize $\widetilde{f}_J$ as a ReLU network, use the iterative structure $g_k = g \circ g_{k-1}$: a stack of $J$ hidden ReLU layers implementing the iterates of $g$ produces $g_1, \ldots, g_J$ at successive depths, while a parallel channel carries $x$ unchanged through the network. A final linear layer collects the combination $x - \sum_{k=1}^J g_k(x)/4^k$. Total depth: $J + 1$; total width: at most $5$ (three for the deepest $g$ layer, plus one channel each for $x$ and the running accumulator).
\end{proof}

This construction represents the dyadic piecewise-linear interpolant with only $O(J)$ parameters, although the interpolant has $2^J$ intervals. Its uniform error therefore decreases exponentially in the network depth.

\medskip\noindent\textbf{Stage 2: From $x^2$ to Hölder approximation.}

The following operations explain the local-polynomial approach to ReLU approximation. The sharper width--depth rate in Theorem~\ref{thm:A1} uses the specific network assembly of \cite{jiao2021deep}.

\begin{enumerate}
\item \textbf{Multiplication.} The polarization identity $xy = ((x+y)^2 - x^2 - y^2)/2$ realizes multiplication via three suitably rescaled copies of $\phi_{\mathrm{sq}}$ followed by a linear combination. Composing produces a ReLU network for $xy$ on $[0,1]^2$ with $L^\infty$ error $\epsilon$ at depth $O(\log_2(1/\epsilon))$ and constant width.
\item \textbf{Polynomials.} By repeatedly applying the multiplication operation, the network can realize any monomial $\boldsymbol{x}^{\boldsymbol{\alpha}} = x_1^{\alpha_1} \cdots x_d^{\alpha_d}$ satisfying $\|\boldsymbol{\alpha}\|_1 \le s$ within depth $O(\log d \cdot \log(1/\epsilon))$ and width $O(d)$.
\item \textbf{Local polynomial decomposition.} For $f \in \mathcal{H}^\beta([0,1]^d, B_{\mathrm{REG}})$, partition $[0,1]^d$ into $G^d$ overlapping cells, build an approximate partition of unity $\{\Psi_{\boldsymbol{m}}\}_{\boldsymbol{m}}$ compatible with ReLU \cite{guhring2021approximation}, and on each cell compute the order-$s$ Taylor polynomial of $f$. The global approximation is formed by weighting these local polynomials with the partition of unity, which expands into a linear combination of fundamental building blocks of the form $\Psi_{\boldsymbol{m}}(\boldsymbol{x})\cdot \boldsymbol{x}^{\boldsymbol{\alpha}}$. Each block is then realized by the neural network using the multiplication and polynomial constructions developed in moves~(1) and (2). Aggregating these localized outputs gives a constructive approximation; the precise architecture estimates for Theorem~\ref{thm:A1} follow \cite{jiao2021deep}.
\end{enumerate}

Theorem~\ref{thm:A1} tracks a polynomial-in-$d$ error prefactor. Its displayed width still contains $3^d$, so the error prefactor and the dimensional dependence of the architecture should be distinguished.

\subsubsection{Tanh Sobolev approximation (Gühring--Raslan)}\label{sec:A2}

For the Deep Ritz method, approximation is measured in $H^1$, so both values and gradients must be controlled. ReLU networks can provide such approximation \cite{guhring2020error}. For example, the dyadic interpolant $\widetilde f_J$ in Lemma~\ref{lem:saw} satisfies $\|\widetilde f_J'-2x\|_{L^\infty(0,1)}\le 2^{-J}$: its slope on each cell is the derivative of $x^2$ at the cell midpoint. Thus derivative jumps at the grid nodes do not prevent first-order Sobolev convergence.

We use tanh here to follow the DRM analysis of \cite{jiao2021error}. The approximation theorem below controls a prescribed $W^{k,p}$ norm under a higher-order regularity assumption on the target. Smoothness of the activation permits these derivative estimates; the estimates themselves follow from the construction.

\begin{theorem}[{Tanh Sobolev approximation; \citealp{guhring2021approximation}, Proposition~4.8; \citealp[arXiv version, Theorem~4.1]{jiao2021error}}]\label{thm:A2}
Let $s, k, d \in \mathbb{N}$ with $s \ge k+1$, $1 \le p \le \infty$, and $0<\mu<\min\{1,(s-k)/k\}$ sufficiently small. There exist constants $C = C(d, s, p, k, \mu) > 0$ and $\epsilon^* > 0$ such that for every $\epsilon \in (0, \epsilon^*]$ and every $f \in W^{s,p}([0,1]^d)$ with $\|f\|_{W^{s,p}} \le 1$, there is a $\tanh$-activated network $\phi \in \mathcal{NN}_{\tanh}(L, \mathfrak{D}, B_{\boldsymbol{\theta}})$ with $L \le C \log(d + s)$, $ \mathfrak{D} \le C \epsilon^{-d/(s - k - \mu k)}$ and
\[
B_{\boldsymbol{\theta}} \;\le\; C\,\epsilon^{-\vartheta}, \qquad \vartheta:=2+\frac{2(d/p+d+k+\mu k)+d/p+d}{s-k-\mu k},
\]
such that
\[
\|f - \phi\|_{W^{k,p}([0,1]^d)} \;\le\; \epsilon.
\]
\end{theorem}

Theorem~\ref{thm:A2} supplies the approximation term in the tanh DRM analysis of \cite{jiao2021error}. Its size and weight bounds also enter the statistical estimate, so both must be retained when balancing the errors.

The norm-constrained and parallel-network variants used in later chapters are stated separately in Theorems~\ref{thm:A4} and~\ref{thm:A5}, with their respective parameter bounds.

The construction combines an approximate partition of unity, averaged Taylor polynomials, and small tanh networks for multiplication. The following outline records these ingredients and the parameter dependence used below.

\begin{proof}[Construction outline]
\medskip\noindent\textbf{Step 1: Tanh-realized partition of unity.}
Define the one-dimensional $\tanh$ bump
\[
\psi^a(x) := \tfrac{1}{2}\big[\rho(a(x+3/2)) - \rho(a(x-3/2))\big], \qquad a \ge 1,\quad \rho=\tanh,
\]
where larger $a$ gives a tighter approximate partition at the cost of larger pre-activation weights. For $G \in \mathbb{N}^+$ and $\boldsymbol{m} \in \{0, \ldots, G\}^d$, the tensor-product bump
\[
\Psi^a_{\boldsymbol{m}}(\boldsymbol{x}) := \prod_{l=1}^d \psi^a\!\big(3G(x_l - m_l/G)\big)
\]
is concentrated near $\boldsymbol{m}/G$. The collection $\{\Psi^a_{\boldsymbol{m}}\}_{\boldsymbol{m}}$ forms an exponential partition of unity on $(0,1)^d$ \cite{guhring2021approximation}: there exist $C(d,k), D, S_0 > 0$ such that for all $a \ge S_0$ and $k \in \mathbb{N}_0$,
\begin{align*}
&\big\| \Psi^a_{\boldsymbol{m}} \big\|_{W^{k,\infty}(\mathbb{R}^d)} \le C(aG)^k, \qquad
\big\| \Psi^a_{\boldsymbol{m}} \big\|_{W^{k,\infty}(\Omega^c_{\boldsymbol{m}})} \le C(aG)^k e^{-Da}, \\
&\Big\|\mathbf{1}_{(0,1)^d} - \sum_{\boldsymbol{m}} \Psi^a_{\boldsymbol{m}} \Big\|_{W^{k,\infty}((0,1)^d)} \le C(aG)^k e^{-Da},
\end{align*}
where $\Omega^c_{\boldsymbol{m}} := \{\boldsymbol{x} : \|\boldsymbol{x} - \boldsymbol{m}/G\|_\infty \ge 1/G\}$. Each one-dimensional factor is a $\tanh$ network with one hidden layer and $\mathcal{O}(1)$ nonzero weights, with pre-activation weight $\mathcal{O}(aG)$. The $d$-dimensional product $\Psi^a_{\boldsymbol{m}}$ is assembled via the multiplication tree of Step~4.

\medskip\noindent\textbf{Step 2: Local polynomial approximation.}
Set $a := G^\mu$ with $\mu \in (0,1)$. The function $f$ is approximated by the partition-weighted local Taylor expansion
\[
f_G(\boldsymbol{x}) := \sum_{\boldsymbol{m} \in \{0, \ldots, G\}^d}\, \sum_{\|\boldsymbol{\alpha}\|_1 \le s-1} c_{f, \boldsymbol{m}, \boldsymbol{\alpha}}\, \Psi^a_{\boldsymbol{m}}(\boldsymbol{x})\, \boldsymbol{x}^{\boldsymbol{\alpha}},
\]
where $c_{f,\boldsymbol m,\boldsymbol\alpha}$ are coefficients of averaged Taylor polynomials on cells $\Omega_{\boldsymbol m,G}$ of diameter $O(G^{-1})$. For a Sobolev extension $\widetilde f$, they satisfy $|c_{f,\boldsymbol m,\boldsymbol\alpha}|\le C G^{d/p}\|\widetilde f\|_{W^{s,p}(\Omega_{\boldsymbol m,G})}$. The local approximation estimate of \cite{guhring2021approximation} gives
\[
\|f - f_G\|_{W^{k,p}((0,1)^d)} \;\le\; C(d, s, p, k)\, \|f\|_{W^{s,p}((0,1)^d)}\, G^{-(s - k - \mu k)}.
\]
The exponent $s - k - \mu k$ reflects three contributions: the local order-$(s-1)$ Taylor polynomial on cells of side $1/G$ gives $G^{-(s-k)}$ in $W^{k,p}$; the approximate PU incurs a residual bounded by $CG^{k(1+\mu)}e^{-DG^\mu}$, super-polynomially negligible; and the PU's $W^{k,\infty}$-norm scales as $G^{k(1+\mu)}$ under $a = G^\mu$, leaking $\mu k$ extra orders through the product rule. The parameter $\mu \in (0,1)$ thus trades a small loss of convergence order for super-polynomial control of the PU residual.

\medskip\noindent\textbf{Step 3: Atomic tanh approximation of $x$ and $x^2$.}
Fix $x_0$ with $\rho'(x_0)\ne0$ and $\rho''(x_0)\ne0$, where $\rho=\tanh$. For any fixed $B>0$ and derivative order $k$, Taylor expansion gives depth-$2$ networks approximating $x$ and $x^2$ on $[-B,B]$ in $W^{k,\infty}$:
\begin{align*}
\phi_1(x) &:= \tfrac{C}{\epsilon \rho^{(1)}(x_0)}\Big[\rho(x_0) - \rho\!\big(x_0 - \tfrac{\epsilon}{C}\, x\big)\Big], \\
\phi_2(x) &:= \tfrac{C^2}{\epsilon^2 \rho^{(2)}(x_0)}\Big[\rho(x_0) - 2\rho\!\big(x_0 - \tfrac{\epsilon}{C}\, x\big) + \rho\!\big(x_0 - \tfrac{2\epsilon}{C}\, x\big)\Big],
\end{align*}
with $C=C(B,k,x_0)$ sufficiently large, widths $1$ and $2$, weight magnitudes $\mathcal{O}(\epsilon^{-1})$ and $\mathcal{O}(\epsilon^{-2})$, and errors $\|x - \phi_1\|_{W^{k,\infty}}, \|x^2 - \phi_2\|_{W^{k,\infty}} \le \epsilon$. The mechanism is Taylor isolation: setting $\delta = -\epsilon/C$ (here $\delta$ denotes the Taylor-scaling step) in $\rho(x_0+\delta x) = \rho(x_0) + \delta x\,\rho^{(1)}(x_0) + \tfrac{(\delta x)^2}{2}\rho^{(2)}(x_0) + \cdots$, the first-order finite difference isolates $-x\,\rho^{(1)}(x_0)$ with remainder $\mathcal{O}(\epsilon)$, and the second-order difference isolates $x^2\,\rho^{(2)}(x_0)$ similarly; differentiation with respect to $x$ multiplies each term by $\delta$, preserving the isolation \cite{guhring2021approximation}.

Here accuracy is improved by increasing the weight magnitudes while keeping the atomic depth fixed. The derivative approximation gives a bounded function-level Lipschitz constant on $[-B,B]$, but this does not imply a uniformly bounded product of weight norms. The norm-controlled results in Chapter~\ref{sec:overparam_setting} track those norms separately.

\medskip\noindent\textbf{Step 4: Network construction and error decomposition.}
Let $\Phi_{\boldsymbol{m},\boldsymbol{\alpha}}$ denote the $\tanh$ subnetwork that approximates the product $\Psi^a_{\boldsymbol{m}}\boldsymbol{x}^{\boldsymbol{\alpha}}$: monomials $\boldsymbol{x}^{\boldsymbol{\alpha}}$ are built from $\phi_1, \phi_2$ via iterated polarization $xy = \tfrac{1}{4}((x+y)^2 - (x-y)^2)$, then multiplied by $\Psi^a_{\boldsymbol{m}}$. The assembled network mirrors the form of $f_G$:
\[
\phi(\boldsymbol{x}) := \sum_{\boldsymbol{m} \in \{0,\ldots,G\}^d}\, \sum_{\|\boldsymbol{\alpha}\|_1 \le s-1} c_{f, \boldsymbol{m}, \boldsymbol{\alpha}}\, \Phi_{\boldsymbol{m},\boldsymbol{\alpha}}(\boldsymbol{x}).
\]
The error decomposes as
\[
\|f - \phi\|_{W^{k,p}} \;\le\; \|f - f_G\|_{W^{k,p}} \;+\; \sum_{\boldsymbol{m},\boldsymbol{\alpha}} |c_{f,\boldsymbol{m},\boldsymbol{\alpha}}| \cdot \|\Psi^a_{\boldsymbol{m}}\boldsymbol{x}^{\boldsymbol{\alpha}} - \Phi_{\boldsymbol{m},\boldsymbol{\alpha}}\|_{W^{k,p}}.
\]
The first term is bounded by Step~2. For the second, each $\Phi_{\boldsymbol{m},\boldsymbol{\alpha}}$ is tuned to match $\Psi^a_{\boldsymbol{m}}\boldsymbol{x}^{\boldsymbol{\alpha}}$ to within the target accuracy; choosing the per-factor atomic accuracy appropriately and summing over the $\asymp (G+1)^d$ cells yields a total assembly error $\le \epsilon$.

\medskip\noindent\textbf{Step 5: Parameter count.}
Choose $G\asymp\epsilon^{-1/(s-k-\mu k)}$. For fixed $d,s,p,k,\mu$, there are $O((G+1)^d)$ cell--monomial pairs, each implemented by a constant-size subnetwork. The multiplication tree gives depth $O(\log(d+s))$, and the total nonzero-parameter count is $O(\epsilon^{-d/(s-k-\mu k)})$. Tracking the atomic accuracies and the coefficients gives the weight bound in Theorem~\ref{thm:A2}; see \cite{guhring2021approximation,jiao2021error}.
\end{proof}

\subsubsection{Parallel theories for alternative architectures}\label{sec:other_arch_classical}

Beyond fully connected ReLU and tanh networks, parallel approximation results exist for several other architectures. These form an important part of the classical size-controlled landscape.

\begin{itemize}
\item Convolutional networks (CNNs). Translation-invariant CNNs achieve approximation rates analogous to fully connected networks for translation-invariant targets \cite{zhou2020universality,6ZHOU2020319}. Moreover, the zero-padding, restricted-bias, and pooling constructions of \cite{11Han2023DeepCN} encode translation-equivariance and translation-invariance structurally, properties that unconstrained fully connected architectures do not guarantee automatically. For $(\beta, C)$-smooth targets with $\beta \ge 2$, \cite{Lin2023OptimalAA} sharpens earlier sub-optimal eDCNN bounds \cite{149714287} to the optimal rate $\inf_{g}\|f-g\|_{L^\infty([0,1]^d)} \lesssim (L/\log L)^{-2\beta/d}$ in terms of depth $L$. Notably, this is achieved with max-pooling controlling the overall network size and without requiring a fully connected head.

\item Recurrent networks (RNNs). Beyond classical universal-approximation results for shallow and deep-narrow RNNs \cite{schafer2006recurrent,hoon2023minimal}, \cite{jiao2024approximationboundsrecurrentneural} establishes rigorous H\"older approximation bounds for deep, general RNNs. Specifically, each step's output approximates a past-dependent target $f^{(t)}\in\mathcal{H}^\beta([0,1]^{d_x\times t})$ at a rate of $(JI)^{-2\beta/(d_x t)}$ via a network of width $W\lesssim J\log J$ and depth $L\lesssim I\log I$.

\item Self-attention and Transformers. Universal approximation for these architectures was initially established via the contextual-mapping argument of \cite{yun2019transformers}, with one-layer self-attention shown to be universal in \cite{kajitsuka2023transformers}. While these foundational results target continuous functions and incur the curse of dimensionality, recent work by \cite{jiao2025transformers} demonstrates that a single one-head self-attention layer with feedforward blocks can approximate H\"older targets while effectively overcoming the curse of dimensionality. This is achieved via a Kolmogorov--Arnold superposition argument that bypasses the need for contextual mapping, utilizing $O(\log\tfrac1\epsilon)$ ReLU--floor feedforward layers of width $O(\epsilon^{-2/\beta}\log\tfrac1\epsilon)$ (or constant width under richer activations).

\item Super-expressive networks. By employing carefully engineered activations (e.g., $\sin$ or mod-$1$), dimension-free approximation rates become achievable in principle \cite{yarotsky2020phase,shen2021neural}. For example, \cite{jiao2023deep} explicitly constructs ReLU-sine-$2^x$ networks of depth $O(d)$ and width polynomial in $\epsilon^{-1}$ that break the curse of dimensionality on $\mathcal{H}^\beta([0,1]^d)$. In the context of this review, these results delineate the theoretical limits of expressivity.
\end{itemize}

\subsection{Statistical tools for this regime}\label{sec:classical_stat}

The statistical bounds developed here combine empirical-process estimates with concentration inequalities. Symmetrization and contraction reduce the expected fluctuations to network-class complexity, which can be bounded through pseudo-dimension or parameter covering. McDiarmid's inequality controls deviations using boundedness, whereas Bousquet's inequality also uses variance information. Fast-rate bounds require a variance--mean relation together with suitable localized complexity control; changing the concentration inequality alone does not provide this improvement. Section~\ref{sec:four_conditions} states the assumptions used in these arguments. Their application to regression and to the function and gradient classes in the DRM energy is developed in \S\ref{sec:classical_rates}.

\subsubsection{Concentration inequalities, symmetrization, and Rademacher complexity}

We begin with the concentration layer: three inequalities that pin random quantities to their expectations, stated here once and used throughout \S\ref{sec:four_conditions}--\S\ref{sec:classical_rates}.

\begin{lemma}[McDiarmid's inequality; \citealp{mohri2018foundations}, Theorem D.8]\label{lem:mcdiarmid}
Let $Z_1, \ldots, Z_n$ be independent random variables on $\mathcal{Z}$. Suppose there exist $c_1, \ldots, c_n > 0$ such that $F: \mathcal{Z}^n \to \mathbb{R}$ satisfies
\[
|F(z_1, \ldots, z_i, \ldots, z_n) - F(z_1, \ldots, z_i', \ldots, z_n)| \;\le\; c_i
\]
for all $i$ and all points $z_1, \ldots, z_n, z_i' \in \mathcal{Z}$. Then, with probability at least $1 - \delta$,
\[
F(Z_1, \ldots, Z_n) \;\le\; \mathbb{E}\big[F(Z_1, \ldots, Z_n)\big] + \sqrt{\tfrac{1}{2}\sum_{i=1}^n c_i^2\, \log(1/\delta)}.
\]
\end{lemma}

\begin{lemma}[Bernstein's inequality; e.g.\ \citealp{gine2021mathematical}]\label{lem:bernstein}
Let $X_1, \ldots, X_n$ be i.i.d.\ with $|X_i - \mathbb{E}[X_i]| \le M$ almost surely and $\mathrm{Var}(X_i) = \sigma^2$. Then, with probability at least $1 - \delta$,
\[
\frac{1}{n}\sum_{i=1}^n \big(X_i - \mathbb{E}[X_i]\big) \;\le\; \sqrt{\frac{2\sigma^2 \log(1/\delta)}{n}} + \frac{2M \log(1/\delta)}{3n}.
\]
\end{lemma}

\begin{lemma}[Bousquet's inequality; \citealp{Bousquet2002Bennett}]\label{lem:bousquet}
Let $Z_1, \ldots, Z_n$ be independent random variables and $\mathcal{G}$ a countable function class with $|g(Z_i)| \le M$ for all $g \in \mathcal{G}$ and all $i$. Define
\[
S \;:=\; \sup_{g \in \mathcal{G}} \sum_{i=1}^n \big(g(Z_i) - \mathbb{E}[g(Z_i)]\big).
\]
Then, with probability at least $1 - \delta$,
\[
S \;\le\; \mathbb{E}[S] + \sqrt{2\Big(\sup_{g \in \mathcal{G}} \sum_{i=1}^n \mathrm{Var}\big(g(Z_i)\big) + 4M\, \mathbb{E}[S]\Big)\log(1/\delta)} + \frac{2M \log(1/\delta)}{3}.
\]
\end{lemma}

The three inequalities differ in what they see: McDiarmid's uses only boundedness, Bousquet's lets the variance of the summands enter explicitly, and Bernstein's is the single-function special case. In each case concentration reduces the statistical error to the expected supremum; the remainder of this subsubsection and the next assemble the tools that bound it.

\begin{definition}[Empirical Rademacher complexity]\label{def:rademacher}
Let $\mathcal{F}$ be a class of measurable functions $f:\mathcal{Z}\to\mathbb{R}$ and $\{Z_i\}_{i=1}^n$ an i.i.d.\ sample on $\mathcal{Z}$. Let $\boldsymbol{\sigma}=(\sigma_1,\ldots,\sigma_n)$ be i.i.d.\ Rademacher signs ($\Pr(\sigma_i=\pm 1)=1/2$), independent of $\{Z_i\}$. The empirical Rademacher complexity of $\mathcal{F}$ is
\begin{equation}\label{eq:rad_def_classical}
\widehat{\mathfrak{R}}_n(\mathcal{F}) := \mathbb{E}_{\boldsymbol{\sigma}}\sup_{f\in\mathcal{F}}\frac{1}{n}\sum_{i=1}^n \sigma_i f(Z_i),
\end{equation}
with fully expected version $\mathfrak{R}_n(\mathcal{F}) := \mathbb{E}_{\{Z_i\}}\widehat{\mathfrak{R}}_n(\mathcal{F})$.
\end{definition}

This is the standard form of \cite{bartlett2002rademacher}. A symmetric variant with $|\cdot|$ inside the supremum,
\begin{equation}\label{eq:rad_def_sym}
\tilde{\mathfrak{R}}_n(\mathcal{F}) := \mathbb{E}_{\{Z_i\}}\mathbb{E}_{\boldsymbol{\sigma}}\sup_{f\in\mathcal{F}}\frac{1}{n}\bigg|\sum_{i=1}^n \sigma_i f(Z_i)\bigg| \;=\; \mathfrak{R}_n\big(\mathcal{F}\cup(-\mathcal{F})\big),
\end{equation}
is used in the complete three-error analysis of \S\ref{sec:complete_stat}.

\begin{proposition}[Symmetrization]\label{prop:symmetrization}
Let $\mathcal{F}$ be a class of measurable functions $f: \mathcal{Z} \to \mathbb{R}$, and let $Z_1, \ldots, Z_n$ be i.i.d.\ samples from a distribution on $\mathcal{Z}$. Then
\[
\mathbb{E}\,\sup_{f \in \mathcal{F}}\Big[\mathbb{E}[f(Z)] - \tfrac{1}{n}\sum_{i=1}^n f(Z_i)\Big] \;\le\; 2\,\mathfrak{R}_n(\mathcal{F}).
\]
\end{proposition}

\begin{proof}
Introduce an independent ghost sample $Z'_1, \ldots, Z'_n$ with the same distribution as $Z_1, \ldots, Z_n$. By Jensen's inequality applied to the convex function $\sup_f(\cdot)$,
\begin{align*}
\mathbb{E}\,\sup_f \Big[\mathbb{E}[f] - \tfrac{1}{n}\sum_i f(Z_i)\Big] & = \mathbb{E}_Z\,\sup_f \mathbb{E}_{Z'}\Big[\tfrac{1}{n}\sum_i f(Z'_i) - \tfrac{1}{n}\sum_i f(Z_i)\Big] \\
& \le \mathbb{E}_{Z,Z'}\,\sup_f \tfrac{1}{n}\sum_i\big[f(Z'_i) - f(Z_i)\big].
\end{align*}
Because $(Z_i, Z'_i)$ are exchangeable, multiplying each pair $f(Z'_i) - f(Z_i)$ by an independent Rademacher sign $\sigma_i$ does not change the joint distribution; together with $\sup_f \sum_i \sigma_i[f(Z'_i)-f(Z_i)] \le \sup_f \sum_i \sigma_i f(Z'_i) + \sup_f \sum_i (-\sigma_i) f(Z_i)$, this gives
\[
\mathbb{E}_{Z,Z'}\sup_f \tfrac{1}{n}\sum_i [f(Z'_i) - f(Z_i)] = \mathbb{E}_{Z,Z',\boldsymbol{\sigma}}\sup_f \tfrac{1}{n}\sum_i \sigma_i [f(Z'_i) - f(Z_i)] \le 2\,\mathfrak{R}_n(\mathcal{F}).
\]
\end{proof}

The slow-rate two-sided supremum $\sup_f|\mathcal{L}(f) - \widehat{\mathcal{L}}_n(f)|$ that appears in \eqref{eq:reg_sym_decomp} and \eqref{eq:drm_sym_decomp} is bounded by applying Proposition~\ref{prop:symmetrization} to $\mathcal{F}\cup(-\mathcal{F})$:
\begin{equation}\label{eq:sym_twosided}
\mathbb{E}\sup_{f\in\mathcal{F}}\Big|\mathbb{E}[f(Z)] - \tfrac{1}{n}\sum_{i=1}^n f(Z_i)\Big| \;\le\; 2\,\mathfrak{R}_n\big(\mathcal{F}\cup(-\mathcal{F})\big).
\end{equation}

The symmetrization inequality is stated for the expected complexity $\mathfrak{R}_n$ because the left-hand side involves an expectation over the sample; symmetrization naturally integrates the Rademacher average over both signs and data. The tools that follow, namely Massart's lemma, Talagrand contraction, and Dudley's entropy integral, are stated for the empirical complexity $\widehat{\mathfrak{R}}_n$, because each bounds the Rademacher average pointwise for any fixed sample realization. The two are linked by $\mathfrak{R}_n = \mathbb{E}_{\{Z_i\}}\widehat{\mathfrak{R}}_n$ (Definition~\ref{def:rademacher}). Since the covering-number and norm-controlled estimates on $\widehat{\mathfrak{R}}_n$ are uniform in the sample (the uniform covering number bound \eqref{eq:dudley_uni} and the Golowich estimate of Lemma~\ref{lem:S6} do not depend on the particular realization of $\{Z_i\}$), the same rate carries to $\mathfrak{R}_n$ without loss. This two-level architecture is standard in empirical process theory \cite{bartlett2002rademacher,mohri2018foundations}; we make it explicit because Chapter~\ref{sec:complete_setting} returns to the expected level directly via Lemma~\ref{lem:S7}, and the transition between the two levels should be visible rather than implicit.

\begin{lemma}[Massart's lemma; \citealp{mohri2018foundations}, Theorem 3.7]\label{lem:massart}
Let $A \subseteq \mathbb{R}^n$ be a finite set and $\sigma_1, \ldots, \sigma_n$ independent Rademacher signs. Then
\[
\mathbb{E}\Bigg[\sup_{(x_1, \ldots, x_n) \in A}\frac{1}{n}\sum_{i=1}^n \sigma_i x_i\Bigg] \;\le\; \max_{(x_1, \ldots, x_n) \in A}\sqrt{\frac{1}{n}\sum_{i=1}^n x_i^2}\;\sqrt{\frac{2\log|A|}{n}}.
\]
In particular, for a finite function class $\mathcal{F}$, taking $A = \{(f(Z_1), \ldots, f(Z_n)) : f \in \mathcal{F}\}$ gives
\[
\widehat{\mathfrak{R}}_n(\mathcal{F}) \;\le\; \max_{f \in \mathcal{F}}\sqrt{\frac{1}{n}\sum_{i=1}^n f(Z_i)^2}\;\sqrt{\frac{2\log|\mathcal{F}|}{n}},
\]
and if moreover $\sup_{f \in \mathcal{F}}\|f\|_\infty \le \mathcal{B}$, then $\widehat{\mathfrak{R}}_n(\mathcal{F}) \le \mathcal{B}\sqrt{2\log|\mathcal{F}|/n}$.
\end{lemma}

Massart's lemma is the finite-class baseline; chaining (Proposition~\ref{prop:dudley} below) lifts it to infinite classes by replacing $\log|\mathcal{F}|$ with metric entropy $\log\mathcal{C}(\epsilon,\mathcal{F})$.

\begin{lemma}[Talagrand contraction principle; \citealp{ledoux2013probability}]\label{lem:talagrand_contraction}
Let $\psi:\mathbb{R}\to\mathbb{R}$ be $L$-Lipschitz with $\psi(0)=0$. Then for any class $\mathcal{F}$,
\[
\widehat{\mathfrak{R}}_n(\psi\circ\mathcal{F}) \;\le\; L\,\widehat{\mathfrak{R}}_n(\mathcal{F}).
\]
\end{lemma}

This is the standard device for transferring Rademacher complexity from a loss-composed class $\ell \circ \mathcal{F}$ to the underlying $\mathcal{F}$ when the loss is Lipschitz; the Lipschitz constants of the two running examples are computed at their points of use in \S\ref{sec:classical_rates}.

\subsubsection{Covering numbers, Dudley chaining, and the two covering routes}

Throughout this subsection, $\delta$ denotes the covering radius.

For an infinite class $\mathcal{F}\subset L^\infty(\mathcal{Z})$, Massart's lemma applies only after a discretization, replacing $\log|\mathcal{F}|$ by a metric entropy $\log\mathcal{C}(\epsilon, \mathcal{F}, d)$ at scale $\epsilon$ in some pseudo-metric $d$ on $\mathcal{F}$.

\begin{definition}[Covering numbers of a function class]\label{def:covering}
Let $\mathcal{F}\subset L^\infty(\mathcal{Z})$ be a class of measurable functions, $\epsilon>0$, and $d$ a pseudo-metric on $\mathcal{F}$ (most commonly induced by a (semi)norm $\|\cdot\|$ on an ambient function space containing $\mathcal{F}$ via $d(f,g)=\|f-g\|$, in which case we write $\mathcal{C}(\epsilon, \mathcal{F}, \|\cdot\|)$). The covering number (internal cover) of $\mathcal{F}$ at scale $\epsilon$ under $d$ is
\[
\mathcal{C}(\epsilon, \mathcal{F}, d) := \min\Big\{N : \exists f_1,\ldots,f_N\in\mathcal{F},\;\forall f\in\mathcal{F}\;\exists j,\; d(f, f_j) \le \epsilon\Big\}.
\]
Two metrics recur throughout the review. The $L^\infty$ metric is the domain sup-norm, giving covering numbers of the form $\mathcal{C}(\epsilon, \mathcal{F}, \|\cdot\|_{L^\infty(\mathcal{Z})})$. Given a sample $z^n=(z_1,\ldots,z_n)\in\mathcal{Z}^n$, the empirical metric $d_{z^n,\infty}(f,g):=\max_{i\le n}|f(z_i)-g(z_i)|$ yields the per-sample empirical covering number $\mathcal{C}(\epsilon, \mathcal{F}, d_{z^n,\infty})$. Taking the worst case over samples defines the uniform covering number
\[
\mathcal{C}_n(\epsilon,\mathcal{F}) := \sup_{z^n\in\mathcal{Z}^n}\mathcal{C}\big(\epsilon, \mathcal{F}, d_{z^n,\infty}\big).
\]
\end{definition}

\begin{proposition}[Dudley chaining; \citealp{anthony2009neural,jiao2021error}]\label{prop:dudley}
Let $\mathcal{F}$ be a class of real-valued functions with $\sup_f \|f\|_\infty \le \mathcal{B}$ and $0 \in \mathcal{F}$. The empirical Rademacher complexity is bounded by chaining in either the empirical metric or the $L^\infty(\mathcal{Z})$ metric:
\begin{align}
\widehat{\mathfrak{R}}_n(\mathcal{F}) &\;\le\; \inf_{\delta \in (0, \mathcal{B})}\Bigg(4\delta + \frac{12}{\sqrt{n}}\int_\delta^{\mathcal{B}} \sqrt{\log \mathcal{C}_n(\epsilon, \mathcal{F})}\, \mathrm{d}\epsilon\Bigg), \label{eq:dudley_emp}\\
\widehat{\mathfrak{R}}_n(\mathcal{F}) &\;\le\; \inf_{\delta \in (0, \mathcal{B})}\Bigg(4\delta + \frac{12}{\sqrt{n}}\int_\delta^{\mathcal{B}} \sqrt{\log \mathcal{C}(\epsilon, \mathcal{F}, \|\cdot\|_{L^\infty})}\, \mathrm{d}\epsilon\Bigg). \label{eq:dudley_uni}
\end{align}
\end{proposition}
The two forms share the same standard chaining argument (see, e.g., \citealp[Chapter~17]{anthony2009neural} or \citealp[Theorem~9.1]{gyorfi2002distribution}): at each scale $\epsilon_k$ a $\sqrt{\log\mathcal{C}}/\sqrt{n}$ sub-Gaussian increment is added via Massart's lemma (Lemma~\ref{lem:massart}), and integrating over scales yields the Dudley integral; they differ only in the metric used to build the chain.

Two covering routes are available, distinguished by the activation function; which one applies in a given analysis is indicated at the point of use in \S\ref{sec:four_conditions} and \S\ref{sec:classical_rates}.

\begin{itemize}

    \item \textbf{The pseudo-dimension route (piecewise-polynomial activations).} This route controls $\mathcal{C}_n(\epsilon,\mathcal{F})$ through the pseudo-dimension $\mathrm{Pdim}(\mathcal{F})$ (the real-valued generalization of the VC-dimension, \citealp{anthony2009neural}). The fundamental Anthony--Bartlett bound converts pseudo-dimension to uniform covering:

\begin{proposition}[Uniform covering bound from pseudo-dimension; \citealp{anthony2009neural}]\label{prop:pdim_cover}
Let $\mathcal{F}$ be a class with $\sup_f \|f\|_\infty \le \mathcal{B}$, and let $v\ge\max\{1,\mathrm{Pdim}(\mathcal F)\}$ be an integer. For $n\ge v$ and $0<\epsilon\le\mathcal B$,
\[
\mathcal{C}_n(\epsilon, \mathcal{F}) \;\le\; \Big(\frac{4e n \mathcal{B}}{\epsilon v}\Big)^v.
\]
\end{proposition}

For piecewise-linear activations, the pseudo-dimension is controlled by the following result of \cite{bartlett2019nearly}.

\begin{proposition}[Pseudo-dimension of piecewise-linear networks; \citealp{bartlett2019nearly}]\label{prop:relu_pdim}
Let $\rho$ be a fixed piecewise-linear activation function with finitely many pieces. For a fixed-architecture network class $\mathcal F$ of depth $L$ with $\mathfrak D$ parameters,
\[
\mathrm{Pdim}(\mathcal F) \;\leq\; C_\rho\,\mathfrak D L\log(e\mathfrak D).
\]
The constant $C_\rho$ depends only on the activation's number of pieces and is absolute for ReLU. In particular, for $L\ge2$ and $W\ge d$, the parameter count in \eqref{eq:align_dim} gives
\[
\mathrm{Pdim}\big(\mathcal{NN}_{\mathrm{ReLU}}(W,L)\big)
\le C W^2L^2\log(eWL).
\]
\end{proposition}

The pseudo-dimension bound follows by counting polynomial sign patterns in the network parameters \citep{bartlett2019nearly}. Once activation patterns are fixed, a ReLU network output is polynomial, not generally affine, in all layer parameters jointly. The argument extends to piecewise-polynomial activations such as $\mathrm{ReLU}^k$. For tanh, finite bounds are available through Pfaffian-function methods \cite{karpinski1997polynomial}; the analysis below instead uses parameter covering.

\item \textbf{The parameter-covering route (smooth activations).} For smooth activations such as $\tanh$, we bound $\mathcal{C}(\epsilon,\mathcal F,\|\cdot\|_{L^\infty})$ by transferring a finite-dimensional parameter-space cover to the function class through a parameter-to-function Lipschitz estimate.

\begin{lemma}[Covering transfer under Lipschitz maps]\label{lem:lip_transfer}
Let $(\Theta, d_\Theta)$ be a pseudo-metric space and $\Phi : \Theta \to \mathcal{F}$ a surjection onto a function class $\mathcal{F}$ with $d_{\mathcal{F}}\big(\Phi(\boldsymbol{\theta}), \Phi(\boldsymbol{\theta}')\big) \le L_{\mathrm{param}}\, d_\Theta(\boldsymbol{\theta}, \boldsymbol{\theta}')$. Then, for every $\epsilon > 0$,
\[
\mathcal{C}(\epsilon, \mathcal{F}, d_{\mathcal{F}}) \;\le\; \mathcal{C}\big(\epsilon / L_{\mathrm{param}},\, \Theta,\, d_\Theta\big).
\]
\end{lemma}

\begin{proof}
The image under $\Phi$ of an $(\epsilon/L_{\mathrm{param}})$-net of $\Theta$ is an $\epsilon$-net of $\mathcal{F}$, and its centers lie in $\mathcal{F}$.
\end{proof}

\begin{proposition}[Covering number of a finite-dimensional ball; standard, see e.g.\ \citealp{gine2021mathematical}]\label{prop:param_cover}
Let $\Theta = \{\boldsymbol{\theta} \in \mathbb{R}^{\mathfrak{D}} : \|\boldsymbol{\theta}\|_\infty \le B_{\boldsymbol{\theta}}\}$. For every $\delta > 0$,
\[
\log \mathcal{C}\big(\delta, \Theta, \ell^\infty\big) \;\le\; \mathfrak{D}\, \log\Big(\frac{2 B_{\boldsymbol{\theta}}}{\delta} + 1\Big).
\]
\end{proposition}

The $L^\infty$ covering bound for a network class follows by composing the two facts. Let $\mathcal{F} = \{\phi_{\boldsymbol{\theta}} : \boldsymbol{\theta} \in \Theta\}$ be parameterized by the ball $\Theta$ above, and suppose the realization map $\boldsymbol{\theta} \mapsto \phi_{\boldsymbol{\theta}}$ is $L_{\mathrm{param}}$-Lipschitz from $(\Theta, \ell^\infty)$ to $L^\infty(\mathcal{X})$. Lemma~\ref{lem:lip_transfer} transfers a covering of $\Theta$ at scale $\epsilon / L_{\mathrm{param}}$ to a covering of $\mathcal{F}$ at scale $\epsilon$, and Proposition~\ref{prop:param_cover} bounds the former, so that for every $\epsilon > 0$
\begin{equation}\label{eq:cinfty_param}
\log \mathcal{C}(\epsilon, \mathcal{F}, \|\cdot\|_{L^\infty}) \;\le\; \mathfrak{D}\, \log\Big(\frac{2 B_{\boldsymbol{\theta}} L_{\mathrm{param}}}{\epsilon} + 1\Big).
\end{equation}
For the $\tanh$ networks of Theorem~\ref{thm:A2}, chained Jacobian bounds yield $L_{\mathrm{param}} = \mathcal{O}\big((B_{\boldsymbol{\theta}} W)^L\big)$ \citep{jiao2021error}; in this form the entropy \eqref{eq:cinfty_param} enters the tanh analyses of \S\ref{sec:classical_rates}.

\end{itemize}

\subsubsection{From slow to fast rates: bounding the statistical error}\label{sec:four_conditions}

This subsubsection derives the two generic statistical bounds used throughout the review: a slow bound of order $1/\sqrt{n}$, which requires only boundedness and Lipschitz continuity of the loss (Assumptions~\ref{assum:bounded}--\ref{assum:class} below), and a fast bound of order $1/n$, which requires in addition a variance assumption (Assumption~\ref{assum:var}). Write $\ell(\phi, Z)$ for the pointwise loss at the sample point $Z$ and define the excess loss
\begin{equation}\label{eq:excess_loss}
g_\phi(Z) \;:=\; \ell(\phi, Z) - \ell(\phi^*, Z), \qquad \mathbb{E}[g_\phi] \;=\; \mathcal{L}(\phi) - \mathcal{L}^* \;\ge\; 0.
\end{equation}
If $\operatorname{Var}(\ell(\phi^*,Z))<\infty$, recentring changes the expected uniform deviation by at most $\sqrt{\operatorname{Var}(\ell(\phi^*,Z))/n}$. Thus it does not change the slow-rate order in the bounded-loss settings below. Everything that follows is organized by four assumptions on the pair $(\ell, \mathcal{F}_{NN})$.

\begin{asmp}[Bounded excess loss]\label{assum:bounded}
The excess loss is uniformly bounded over the class: there exists a constant $M_g > 0$ such that $\sup_Z |g_\phi(Z)| \le M_g$ for all $\phi \in \mathcal{F}_{NN}$.
\end{asmp}

\begin{asmp}[Lipschitz loss]\label{assum:lip}
The excess loss is Lipschitz in the evaluation of the network: there exists a constant $\Lambda_g > 0$ such that, for all $\phi_1, \phi_2 \in \mathcal{F}_{NN}$ and every sample point $Z$,
\[
|g_{\phi_1}(Z) - g_{\phi_2}(Z)| \;\le\; \Lambda_g\, |\phi_1(x) - \phi_2(x)|.
\]
\end{asmp}

\begin{asmp}[Bounded class]\label{assum:class}
The hypothesis class is uniformly bounded: there exists a constant $\mathcal{B} > 0$ such that $|\phi(x)| \le \mathcal{B}$ for all $\phi \in \mathcal{F}_{NN}$ and all $x$.
\end{asmp}

\begin{asmp}[Variance condition]\label{assum:var}
The second moment of the excess loss is controlled by its mean: there exists a constant $\alpha_g > 0$ such that, for every $\phi \in \mathcal{F}_{NN}$,
\[
\mathbb{E}\big[g_\phi^2\big] \;\le\; \alpha_g\, \mathbb{E}\big[g_\phi\big].
\]
\end{asmp}

Assumptions~\ref{assum:bounded}--\ref{assum:class} also enter the slow-rate analysis. The loss envelope controls the fluctuations on a finite net, while Assumption~\ref{assum:lip} transfers a cover of the network class to the excess-loss class. For squared loss, bounded predictions alone do not bound the loss unless the responses are also bounded.

\begin{proposition}[Slow-rate statistical bound]\label{prop:slow_engine}
Suppose Assumptions~\ref{assum:bounded}--\ref{assum:class} hold. For every $\epsilon \in (0,1)$,
\[
\begin{aligned}
\mathbb{E}\sup_{\phi\in\mathcal{F}_{NN}}\Big|\,\mathbb{E}[g_\phi] - \frac{1}{n}\sum_{i=1}^n g_\phi(Z_i)\Big| \leq 2M_g\sqrt{\frac{2\log\!\bigl(2\mathcal{C}_n(\mathcal{B}\epsilon, \mathcal{F}_{NN})\bigr)}{n}} + 2\Lambda_g\mathcal{B}\epsilon.
\end{aligned}
\]
For every $\delta\in(0,1)$, with probability at least $1-\delta$, the sample-dependent supremum is at most the same right-hand side plus $M_g\sqrt{2\log(1/\delta)/n}$.
\end{proposition}

\begin{proof}[Proof sketch]
By two-sided symmetrization \eqref{eq:sym_twosided},
\[
\begin{aligned}
\mathbb{E}\sup_{\phi\in\mathcal F_{NN}}\Big|\mathbb E[g_\phi]-\frac1n\sum_{i=1}^n g_\phi(Z_i)\Big| \leq 2\mathbb E_Z\mathbb E_{\boldsymbol\sigma}\sup_{\phi\in\mathcal F_{NN}}\Big|\frac1n\sum_{i=1}^n\sigma_i g_\phi(Z_i)\Big|.
\end{aligned}
\]
Conditional on the sample, choose an internal $\mathcal B\epsilon$-net $\{\phi_1,\ldots,\phi_{N_\epsilon}\}$ in the empirical sup-norm, with $N_\epsilon\le\mathcal C_n(\mathcal B\epsilon,\mathcal F_{NN})$. The net depends on the sample but not on the Rademacher signs. Assumption~\ref{assum:lip} gives
\[
\sup_{\phi\in\mathcal F_{NN}}\Big|\frac1n\sum_{i=1}^n\sigma_i g_\phi(Z_i)\Big|
\le\max_{j\le N_\epsilon}\Big|\frac1n\sum_{i=1}^n\sigma_i g_{\phi_j}(Z_i)\Big|+\Lambda_g\mathcal B\epsilon.
\]
Apply Massart's lemma (Lemma~\ref{lem:massart}) to the finite set
\[
\Big\{\pm\big(g_{\phi_j}(Z_1),\ldots,g_{\phi_j}(Z_n)\big):1\le j\le N_\epsilon\Big\}.
\]
It has at most $2N_\epsilon$ elements, each with empirical root mean square at most $M_g$ by Assumption~\ref{assum:bounded}. Thus
\[
\mathbb E_{\boldsymbol\sigma}\max_{j\le N_\epsilon}\Big|\frac1n\sum_{i=1}^n\sigma_i g_{\phi_j}(Z_i)\Big|
\le M_g\sqrt{\frac{2\log(2N_\epsilon)}{n}}.
\]
Combining these estimates proves the expectation bound. Changing one observation changes the supremum by at most $2M_g/n$, so McDiarmid's inequality (Lemma~\ref{lem:mcdiarmid}) yields the stated high-probability bound.
\end{proof}

The factor $2$ inside the logarithm accounts for both signs in the absolute supremum. The envelope $M_g$ controls the finite-net fluctuation, while $\Lambda_g$ controls the discretization remainder. The concrete slow-rate bounds below use pseudo-dimension covering for regression and parameter covering of the function and gradient classes for DRM, together with their respective loss bounds and Dudley chaining.

The global bounded-difference argument above yields a square-root-order fluctuation term. McDiarmid's inequality uses the envelope $M_g$ but does not exploit the smaller variance that may accompany a smaller excess risk.

Assumption~\ref{assum:var} relates the second moment of the excess loss to its mean, so that small excess risk implies small stochastic fluctuations. The relation is imposed on $g_\phi$ rather than the raw loss, whose variance need not vanish at $\phi^*$. One way to establish it combines a quadratic lower bound on excess risk with an upper bound on the second moment in the same norm. Concretely, if
\begin{equation}\label{eq:C4_source}
\mathcal{L}(\phi) - \mathcal{L}^* \;\ge\; c\, \|\phi - \phi^*\|^2_{\mathcal{H}}
\qquad\text{and}\qquad
\mathbb{E}\big[g_\phi^2\big] \;\le\; C\, \|\phi - \phi^*\|^2_{\mathcal{H}},
\end{equation}
then Assumption~\ref{assum:var} holds with $\alpha_g=C/c$. The quadratic growth condition relates the norm controlling the second moment to the excess risk. For bounded squared-loss regression this relation follows from the conditional-mean identity; for DRM it also requires suitable bounds on the solution and the candidate functions, as discussed below.

On a class localized by $\mathbb E[g_\phi]\le r$, Assumption~\ref{assum:var} bounds the variance by $\alpha_g r$. Bousquet's inequality (Lemma~\ref{lem:bousquet}) uses this risk-dependent bound. The following argument combines it with a covering estimate for the expected localized supremum.

\begin{proposition}[Fast-rate statistical bound]\label{prop:fast_engine}
Suppose Assumptions~\ref{assum:bounded}--\ref{assum:var} hold. There are absolute constants $c_1, c_2, c_3$ such that for all $\epsilon, \delta \in (0,1)$, with probability at least $1 - \delta$, simultaneously for every $\phi \in \mathcal{F}_{NN}$,
\begin{align*}
\mathcal{L}(\phi) - \mathcal{L}^* \;\le\; \frac{2}{n}\sum_{i=1}^n g_\phi(Z_i) \;&+\; c_1 (M_g + \alpha_g)\, \frac{\log\!\big(2\mathcal{C}_n(\mathcal{B}\epsilon, \mathcal{F}_{NN})\big)}{n} \\
&+\; c_2 \Lambda_g \mathcal{B}\epsilon \;+\; c_3 (M_g + \alpha_g)\, \frac{\log(1/\delta)}{n}.
\end{align*}
Integrating the uniform tail bound gives the corresponding expectation bound with an additional $O((M_g+\alpha_g)/n)$ term. This term is absorbed by the complexity term because its logarithm is at least $\log 2$.
\end{proposition}

\begin{proof}[Proof sketch]
Throughout this proof, $r$ denotes the localization radius. Localize by excess risk: $\mathcal{F}_r := \{\phi \in \mathcal{F}_{NN} : \mathbb{E}[g_\phi] \le r\}$, and write $F_r$ for the localized supremum
\[
F_r \;:=\; \sup_{\phi \in \mathcal{F}_r}\Big(\mathbb{E}[g_\phi] - \frac{1}{n}\sum_{i=1}^n g_\phi(Z_i)\Big).
\]
First construct an internal $2\mathcal B\epsilon$-net of $\mathcal F_r$, rather than of the full class. Such a net has at most $\mathcal C_n(\mathcal B\epsilon,\mathcal F_{NN})$ elements: retain the balls of a full-class cover that meet $\mathcal F_r$ and choose a representative from each intersection. Its centers therefore satisfy $\mathbb E[g_\phi^2]\le\alpha_g r$. Symmetrization, Massart's lemma, and contraction for the square map on $[-M_g,M_g]$ give
\[
\mathbb E[F_r]\le C\left\{\sqrt{\frac{\alpha_g r H_\epsilon}{n}}+
\frac{M_g H_\epsilon}{n}+\Lambda_g\mathcal B\epsilon\right\},
\qquad
H_\epsilon:=\log\!\big(2\mathcal C_n(\mathcal B\epsilon,\mathcal F_{NN})\big).
\]
Here the empirical squared radius is bounded by the supremum over the entire localized class before taking expectations. The empirical second-moment supremum is at most $\alpha_g r$ plus its fluctuation about the population second moment. A further symmetrization and contraction bound this fluctuation by the localized Rademacher complexity times a constant multiple of $M_g$. Solving the resulting quadratic inequality gives the displayed estimate.

Bousquet's inequality and Young's inequality then imply, for any fixed $\eta>0$, with probability at least $1-e^{-t}$,
\[
F_r\le \eta r+C_\eta\left\{\frac{(M_g+\alpha_g)(H_\epsilon+t)}{n}
+\Lambda_g\mathcal B\epsilon\right\}.
\]
Choose $r_0$ to be a sufficiently large constant times the expression in braces with $t=\log(2/\delta)$. Apply this bound at $r_j=2^j r_0$ using $t_j=\log(2/\delta)+j$, for $j=0,1,\ldots$. The failure probabilities sum to at most $\delta$, and $j/2^j$ is bounded. Taking $\eta$ sufficiently small and the constant in $r_0$ sufficiently large therefore ensures $F_{r_j}\le r_j/4$ simultaneously. If $\mathbb E[g_\phi]>r_0$, choose $j$ with $r_j/2<\mathbb E[g_\phi]\le r_j$; then $\mathbb E[g_\phi]-n^{-1}\sum_{i=1}^n g_\phi(Z_i)\le \mathbb E[g_\phi]/2$. If $\mathbb E[g_\phi]\le r_0$, the bound on $F_{r_0}$ gives $\mathbb E[g_\phi]\le2n^{-1}\sum_{i=1}^n g_\phi(Z_i)+r_0/2$. Together these prove the stated inequality. This is the bounded-loss localization and peeling argument underlying \cite{bartlett2005local,koltchinskii2006local}.
\end{proof}

\begin{remark}
Unlike the slow-rate bound (Proposition~\ref{prop:slow_engine}), which controls the supremum of the empirical process, the first display is a uniform oracle inequality: it bounds the excess risk of every $\phi$ simultaneously by its own empirical excess loss plus complexity, a stronger statement made possible by the variance condition. The derivation that the first display bounds $\mathcal{E}_{sta,\,fast}$ is as follows. Write $\widehat{\mathcal{L}}_n(\phi) := \frac{1}{n}\sum_{i=1}^n \ell(\phi, Z_i)$ and recall $g_\phi = \ell(\phi) - \ell(\phi^*)$, so $\frac{1}{n}\sum g_\phi = \widehat{\mathcal{L}}_n(\phi) - \widehat{\mathcal{L}}_n(\phi^*)$. Let $K$ collect the complexity terms of the first display; $K$ does not depend on $\phi$. Substituting gives $\mathcal{L}(\phi) - \mathcal{L}^* \le 2(\widehat{\mathcal{L}}_n(\phi) - \widehat{\mathcal{L}}_n(\phi^*)) + K$. Rearranging,
\[
\mathcal{L}(\phi) + \mathcal{L}^* - 2\widehat{\mathcal{L}}_n(\phi) \;\le\; K \;+\; 2\big(\mathcal{L}^* - \widehat{\mathcal{L}}_n(\phi^*)\big).
\]
The right-hand side is free of $\phi$, so the supremum over $\mathcal{F}_{NN}$ inherits the same bound. Integrating the uniform tail first, and then using $\mathbb{E}[\widehat{\mathcal{L}}_n(\phi^*)] = \mathcal{L}^*$ eliminates the second term, yielding $\mathcal{E}_{sta,\,fast} \lesssim (M_g+\alpha_g)H_\epsilon/n+\Lambda_g\mathcal B\epsilon$. The factor~$2$ on the empirical term comes from the self-normalizing form of the offset decomposition.
\end{remark}

\begin{remark}[Another localization mechanism: peeling]
Peeling partitions $\mathcal F_{NN}$ into risk shells $\mathcal F_k:=\{\phi:2^{k-1}r_0\le\mathbb E[g_\phi]\le2^kr_0\}$. Under Assumption~\ref{assum:var}, the variance on each shell is controlled by its risk scale. Applying concentration on each shell and summing the failure probabilities yields a localized bound, provided the corresponding complexity estimates close the argument. Local Rademacher methods express this balance through a fixed point \cite{bartlett2005local,koltchinskii2006local}. Neither a pseudo-dimension bound nor a parameter-covering bound alone guarantees that the required local complexity decreases with the shell radius.
\end{remark}

\subsection{Rates for the running examples}\label{sec:classical_rates}

This subsection applies the approximation and statistical bounds to the two running examples, with optimization error set to zero. For regression, we first bound the approximation error, then the statistical error, and finally balance them over the network width. For tanh-DRM, the balance also tracks the weight bound and the complexity of the gradient class. The fast-rate results use variance control and localization in their respective settings. Chapter~\ref{sec:complete_setting} additionally analyzes optimization error.

\subsubsection{Regression: slow and fast rates}\label{sec:reg_slow}\label{sec:reg_fast}

For Regression, the target is $\phi^* = f_0$, the risk $\mathcal{L} = \mathcal{L}_{\mathrm{REG}}$, and the target-quantity norm $\|\cdot\|_{L^2(P_X)}$. In this subsubsection, assume $|\xi|\le\sigma$ almost surely and fix $\mathcal B\ge B_{\mathrm{REG}}+\sigma$. Use the output-clipped class
\[
\mathcal F_{NN}=\{[\phi]_{\mathcal B}:\phi\in\mathcal{NN}_{\mathrm{ReLU}}(W,L)\},
\qquad [\phi]_{\mathcal B}:=\min\{\mathcal B,\max\{-\mathcal B,\phi\}\}.
\]
Then $|Y|,|f_0|,|f|\le\mathcal B$ for every $f\in\mathcal F_{NN}$. Clipping does not increase the approximation error, covering numbers, or pseudo-dimension. The four assumptions of \S\ref{sec:four_conditions} now hold: the excess loss factorizes as $g_f = (Y - f)^2 - (Y - f_0)^2 = (f_0 - f)(2Y - f - f_0)$, giving $M_g \le 8\mathcal{B}^2$ and $\Lambda_g \le 4\mathcal{B}$; and since $\mathbb{E}[Y \mid\boldsymbol X] = f_0(\boldsymbol X)$ makes the cross term vanish, $\mathbb{E}[g_f] = \|f - f_0\|^2_{L^2(P_X)}$ exactly, whence $\mathbb{E}[g_f^2] \le 16\mathcal{B}^2\, \mathbb{E}[g_f]$ and Assumption~\ref{assum:var} holds with $\alpha_g \le 16\mathcal{B}^2$. Phrased for the squared distance $h(f) := (f - f_0)^2$, the same computation gives the variance--mean relation
\begin{align}\label{eq:var_mean}
&h(f)^2 = (f - f_0)^4 \le 4\mathcal{B}^2 (f - f_0)^2 = 4\mathcal{B}^2\, h(f) \\\notag
& \Longrightarrow \qquad \mathrm{Var}(h(f)) \le \mathbb{E}[h(f)^2] \le 4\mathcal{B}^2\,\mathbb{E}[h(f)].
\end{align}
Since $\mathbb E[h(f)]=\mathbb E[g_f]$, the squared-distance and excess-risk localizations coincide. Both the slow and the fast analysis below are therefore available.

We begin with the slow rate, which combines the symmetric decomposition \eqref{eq:reg_sym_decomp} with the pseudo-dimension route of \S\ref{sec:classical_stat}.

\begin{theorem}[Regression slow rate; following \citealp{schmidt2020nonparametric,jiao2021deep}]\label{thm:reg_slow}
Suppose $f_0 \in \mathcal{H}^\beta([0,1]^d, B_{\mathrm{REG}})$ for some $\beta \ge 1$, $|\xi|\le\sigma$ almost surely, and $\mathcal{E}_{opt} = 0$. Fix $\mathcal B\ge B_{\mathrm{REG}}+\sigma$ and use the output-clipped class $\mathcal F_{NN}=\{[\phi]_{\mathcal B}:\phi\in\mathcal{NN}_{\mathrm{ReLU}}(W,L)\}$ with the fixed depth $L=L_0(d,\beta)$ prescribed by Theorem~\ref{thm:A1} with $I=1$, and $W \asymp n^{d/(2(4\beta+d))}$ (up to logarithmic factors). Then, for sufficiently large $n$, the ERM \eqref{eq:reg_erm} satisfies
\[
\mathbb{E}_{S_n}\big\|\widehat{f}_{\boldsymbol{\theta}} - f_0\big\|_{L^2(P_X)}^2 \;\le\; \widetilde{O}\!\left(n^{-\frac{2\beta}{4\beta + d}}\right).
\]
Here $\widetilde O$ suppresses logarithmic factors in $n$; its constants may depend on $d,\beta,B_{\mathrm{REG}}$ and $\mathcal B$.
\end{theorem}

\begin{proof}[Proof sketch]
The symmetric decomposition \eqref{eq:reg_sym_decomp} under $\mathcal{E}_{opt}=0$ reduces to $\mathbb{E}\|\widehat{f} - f_0\|_{L^2}^2 \le \mathcal{E}^{\mathrm{REG}}_{app} + \mathcal{E}^{\mathrm{REG}}_{sta,\,slow}$. The derivation follows the three-step trade-off pattern of \S\ref{sec:classical_rates}.

\textbf{Step 1: Approximation error.} By Theorem~\ref{thm:A1} \cite{jiao2021deep,yarotsky2017error}, for any $f_0 \in \mathcal{H}^\beta([0,1]^d, B_{\mathrm{REG}})$ and any $J, I \in \mathbb{N}^+$ there exists a ReLU network $\phi$ of width $W = O(J \log(eJ))$ and depth $L = O(I \log(eI))$ with $\|f_0 - \phi\|_{L^\infty} \le C_1(d, \beta, B_{\mathrm{REG}})(JI)^{-2\beta/d}$. Clipping, squaring, and using $\|\cdot\|_{L^2(P_X)} \le \|\cdot\|_{L^\infty}$ gives
\begin{equation}\label{eq:reg_slow_app}
\mathcal{E}^{\mathrm{REG}}_{app} \;=\; \inf_{\bar{f} \in \mathcal{F}_{NN}}\|\bar{f} - f_0\|_{L^2(P_X)}^2 \;\le\; C_2(d,\beta,B_{\mathrm{REG}})\,(JI)^{-4\beta/d}.
\end{equation}

\textbf{Step 2: Statistical error bound.} Symmetrization \eqref{eq:sym_twosided}, a finite loss net, and the pseudo-dimension covering bounds of Propositions~\ref{prop:pdim_cover}--\ref{prop:relu_pdim} yield the following.

\begin{proposition}[Regression slow-rate Rademacher bound]\label{prop:S3}
Let $\mathcal{F}=\{[\phi]_{\mathcal B}:\phi\in\mathcal{NN}_{\mathrm{ReLU}}(W,L)\}$ with $L\ge2$ and $W\ge d$. For squared loss with $|Y|\le\mathcal B$ almost surely,
\begin{align}
\mathcal{E}^{\mathrm{REG}}_{sta,\,slow}  &:=\; \mathbb{E}\sup_{f \in \mathcal{F}}\Big|\mathcal{L}_{\mathrm{REG}}(f) - \widehat{\mathcal{L}}_{\mathrm{REG}}(f)\Big| \notag\\
& \le\; C\mathcal B^2\,\sqrt{\frac{W^2 L^2 \log(eWL)\log(en)}{n}}. \label{eq:reg_slow_sta}
\end{align}
\end{proposition}

\begin{proof}[Proof sketch]
The raw squared loss has envelope $4\mathcal B^2$ and is $4\mathcal B$-Lipschitz in the prediction. Set $v=\max\{1,\mathrm{Pdim}(\mathcal F)\}$. If $n\ge v$, cover $\mathcal F$ at radius $\mathcal B/(2n)$ and apply two-sided symmetrization and Massart's lemma, as in the proof of Proposition~\ref{prop:slow_engine}. Proposition~\ref{prop:pdim_cover} then gives $\mathcal E^{\mathrm{REG}}_{sta,\,slow}\le C\mathcal B^2\sqrt{v\log(en)/n}$. Proposition~\ref{prop:relu_pdim} bounds $v$ by $CW^2L^2\log(eWL)$. If $n<v$, the uniform loss bound $\mathcal E^{\mathrm{REG}}_{sta,\,slow}\le4\mathcal B^2$ gives the displayed inequality after increasing $C$.
\end{proof}

\textbf{Step 3: Trade-off and balance.} Substituting \eqref{eq:reg_slow_app} and \eqref{eq:reg_slow_sta} into the symmetric decomposition and suppressing logarithmic factors yields $\mathbb{E}\|\widehat{f} - f_0\|_{L^2}^2 \lesssim (JI)^{-4\beta/d} + WL/\sqrt{n}$. Taking $I = O(1)$ so that $L = O(1)$ and $W \asymp J$ up to logarithmic factors, this becomes $J^{-4\beta/d} + J/\sqrt{n}$. Balancing $J^{-4\beta/d} \asymp J/\sqrt{n}$ gives $J \asymp n^{d/(2(4\beta+d))}$, which recovers the stated rate.
\end{proof}
The rate $n^{-2\beta/(4\beta+d)}$ has a slower exponent than the classical H\"older benchmark $n^{-2\beta/(2\beta+d)}$ because $\mathcal{E}_{sta,\,slow}$ scales as $1/\sqrt{n}$. The fast-rate analysis below recovers that exponent; the Gaussian-noise lower bound is reviewed separately in \S\ref{sec:minimax_lb}.

The fast rate improves the statistical term from $1/\sqrt{n}$ to $1/n$ by replacing the symmetric decomposition with the asymmetric one and invoking the variance condition (Assumption~\ref{assum:var}) through Proposition~\ref{prop:fast_engine}. The approximation step is unchanged.

\begin{theorem}[Regression fast rate; \cite{schmidt2020nonparametric,jiao2021deep,yang2024nonparametric,yang2024optimal}]\label{thm:reg_fast}
Suppose $f_0 \in \mathcal{H}^\beta([0,1]^d, B_{\mathrm{REG}})$ for some $\beta \ge 1$, $|\xi|\le\sigma$ almost surely, and $\mathcal{E}_{opt} = 0$. Fix $\mathcal B\ge B_{\mathrm{REG}}+\sigma$ and use the output-clipped class $\mathcal F_{NN}=\{[\phi]_{\mathcal B}:\phi\in\mathcal{NN}_{\mathrm{ReLU}}(W,L)\}$ with the fixed depth $L=L_0(d,\beta)$ prescribed by Theorem~\ref{thm:A1} with $I=1$, and $W \asymp n^{d/(2(2\beta+d))}$ (up to logarithmic factors). Then, for sufficiently large $n$, the ERM \eqref{eq:reg_erm} satisfies
\[
\mathbb{E}_{S_n}\big\|\widehat{f}_{\boldsymbol{\theta}} - f_0\big\|_{L^2(P_X)}^2 \;\le\; \widetilde{O}\!\left(n^{-\frac{2\beta}{2\beta + d}}\right).
\]
The exponent agrees with the classical H\"older rate. Theorem~\ref{thm:holder-regression-lower-bound} establishes the corresponding lower bound in a Gaussian-noise experiment, rather than the bounded-noise setting used here. Extensions to sub-Gaussian noise are discussed below.
Here $\widetilde O$ suppresses logarithmic factors in $n$; its constants may depend on $d,\beta,B_{\mathrm{REG}}$ and $\mathcal B$.
\end{theorem}

\begin{proof}[Proof sketch]
\textbf{Step 1: Decomposition.} The asymmetric (offset) decomposition \eqref{eq:reg_asym_decomp} under $\mathcal{E}_{opt}=0$ reads
\begin{equation}\label{eq:reg_fast_decomp}
\mathbb{E}_{S_n}\|\widehat{f}_{\boldsymbol{\theta}} - f_0\|^2_{L^2(P_X)} \le 2\,\mathcal{E}^{\mathrm{REG}}_{app} + \mathcal{E}^{\mathrm{REG}}_{sta,\,fast},
\end{equation}
where the fast-rate statistical error is defined by the offset supremum \eqref{eq:reg_offset_def}:
\[
\mathcal{E}^{\mathrm{REG}}_{sta,\,fast} \;:=\; \mathbb{E}_{S_n}\sup_{f \in \mathcal{F}_{NN}}\big[\mathcal{L}_{\mathrm{REG}}(f_0) - 2\widehat{\mathcal{L}}_{\mathrm{REG}}(f) + \mathcal{L}_{\mathrm{REG}}(f)\big].
\]

\textbf{Step 2: Statistical error bound.} The approximation step is unchanged from the slow rate; the four assumptions are verified above with $M_g \le 8\mathcal{B}^2$, $\Lambda_g \le 4\mathcal{B}$, $\mathcal{B}$, and $\alpha_g \le 16\mathcal{B}^2$. Instantiating Proposition~\ref{prop:fast_engine} through the pseudo-dimension route (Propositions~\ref{prop:pdim_cover}--\ref{prop:relu_pdim}, with $L\ge2$, $W\ge d$, and $\mathrm{Pdim} \le C W^2 L^2 \log(eWL)$) gives the following.

\begin{proposition}[Regression fast-rate Rademacher bound]\label{prop:S4}
Under the assumptions of Proposition~\ref{prop:S3} and the bounded-noise condition $|\xi| \le \mathcal{B}$ a.s.,
\begin{equation}
\mathcal{E}^{\mathrm{REG}}_{sta,\,fast} \;\le\; C(\mathcal{B})\,\frac{W^2 L^2 \log(eWL)}{n}\,\mathrm{polylog}(en). \label{eq:reg_fast_sta}
\end{equation}
\end{proposition}

\begin{proof}[Proof sketch]
Proposition~\ref{prop:fast_engine} bounds the offset supremum $\mathcal{E}^{\mathrm{REG}}_{sta,\,fast}$ (as shown in the remark following that proposition). For $n\ge\max\{1,\mathrm{Pdim}(\mathcal F)\}$, substitute the pseudo-dimension covering bound with $\mathrm{Pdim} \le C W^2 L^2 \log(eWL)$ and choose $\epsilon = 1/(2n)$. Absorbing the $\Lambda_g\mathcal{B}\epsilon$ term yields the stated bound; for smaller $n$, boundedness of the loss suffices. The localization step is supplied by Proposition~\ref{prop:fast_engine}.
\end{proof}

\textbf{Step 3: Trade-off and balance.} With the unchanged approximation upper bound $\mathcal{E}^{\mathrm{REG}}_{app} \lesssim J^{-4\beta/d}$ from Theorem~\ref{thm:A1} (taking $I=O(1)$, hence $W \asymp J$ up to logarithmic factors and $L=O(1)$), substitute into \eqref{eq:reg_fast_decomp}. Balancing $J^{-4\beta/d} \asymp J^2/n$ gives $J \asymp n^{d/(2(2\beta + d))}$, which recovers the stated rate.
\end{proof}

\begin{remark}[Peeling versus offset Rademacher complexity]
The bounded-loss argument above uses a variance--mean condition and uniform complexity control. With sub-Gaussian responses, truncation or a suitable unbounded-loss analysis is needed; offset Rademacher methods provide one such approach \cite{liang2015learning}. Norm-constrained regression results are discussed in \cite{yang2024nonparametric,yang2024optimal}. The DRM fast-rate result stated below uses the network class and localized analysis of \cite{lu2021machine}.
\end{remark}

\subsubsection{Deep Ritz method: slow and fast rates}\label{sec:drm_slow}\label{sec:drm_fast}

For DRM, the target is $\phi^*=u_0\in H^2(\Omega)$ and the error is measured in $H^1(\Omega)$. Proposition~\ref{prop:drm_eq} relates the excess energy to the squared $H^1$ error. Its coercive lower bound converts excess energy into solution error, producing the factor $2/(c_0\wedge1)$ in \eqref{eq:drm_sym_decomp}; its upper bound controls the approximation contribution by the squared $H^1$ approximation error. Both directions enter the slow-rate analysis.

The interior DRM loss depends on both $u$ and $\nabla u$, whereas the boundary term depends on the trace of $u$. The statistical analysis therefore covers the function and gradient classes, rather than only the scalar network outputs. For the tanh class of Theorem~\ref{thm:A2}, the chained Jacobian bound gives
\[
\sup_{u\in\mathcal F_{NN}}\max\{\|u\|_{L^\infty(\Omega)},\|\nabla u\|_{L^\infty(\Omega)}\}
\lesssim \mathfrak D^L B_{\boldsymbol\theta}^L,
\]
for parameter budgets at least one. These candidate-class bounds and the bounded PDE data suffice for the global statistical estimate below.

A variance--mean bound for DRM loss differences additionally requires control of the target solution and the candidate functions in the relevant norms. In particular, $u_0\in H^2(\Omega)$ alone does not give a uniform $W^{1,\infty}$ bound in general dimension. The tanh slow-rate argument does not use this additional condition; the localized result presented afterwards uses the assumptions of its own Dirichlet setting.

\begin{theorem}[{DRM slow rate; \citealp[arXiv version, proof of Theorem~6.1]{jiao2021error}}]\label{thm:drm_slow}
Suppose $u_0 \in H^2(\Omega)$ is the weak solution of \eqref{eq:drm_pde}, and $\mathcal{E}_{opt} = 0$. Choose the tanh network class $\mathcal{F}_{NN} = \mathcal{NN}_{\tanh}(L, \mathfrak{D}, B_{\boldsymbol{\theta}})$ with equal interior and boundary sample sizes, both denoted by $N$, and parameters scaled per the trade-off balance below. Then the ERM \eqref{eq:drm_erm} satisfies, for any $\mu\in(0,1)$,
\[
\mathbb{E}\big\|\widehat{u}_{\boldsymbol{\theta}} - u_0\big\|_{H^1(\Omega)}^2 \;\le\; \widetilde O\!\left(N^{-(1-\mu)/(Cd\log(d+1))}\right).
\]
This is a squared-error consequence of \cite[arXiv version, proof of Theorem~6.1]{jiao2021error}. The notation $\widetilde O$ suppresses logarithmic factors in $N$; the hidden constants may depend on $\mu,d,\Omega,c_0$, the PDE data bounds, and $\|u_0\|_{H^2(\Omega)}$.
\end{theorem}

\begin{proof}[Proof sketch]
The symmetric decomposition \eqref{eq:drm_sym_decomp} differs from Regression in that the energy--norm equivalence (Proposition~\ref{prop:drm_eq}) introduces a multiplicative factor $2/(c_0 \wedge 1)$, since $\mathcal{L}_{\mathrm{DRM}}(u) - \mathcal{L}_{\mathrm{DRM}}(u_0) \asymp \|u-u_0\|_{H^1}^2$ rather than equalling it exactly. The derivation otherwise follows the same three-step pattern.

\textbf{Step 1: Approximation error.} The error must be controlled in $H^1$, so we use Theorem~\ref{thm:A2} (tanh Sobolev approximation) with $k=1$, $p=2$, $s=2$, applied with target accuracy $\epsilon$ so that $\mathcal{E}^{\mathrm{DRM}}_{app} \le \epsilon^2$, with $\mathfrak{D}$ as the size-controlling parameter. The effective smoothness shifts from $\beta$ to $1$; the tanh framework follows \cite{jiao2021error}, and $\mathrm{ReLU}^k$ variants \cite{Duan2022,jiao2022rate} are compared below.

\textbf{Step 2: Statistical error bound.} We use parameter covering to control the function and gradient classes in the tanh DRM loss. A cover of the parameter space is transferred to these classes using their parameter-to-function Lipschitz bounds. Dudley's entropy integral and contraction for the bounded quadratic loss then yield the prefactor in Proposition~\ref{prop:S5}.

\begin{proposition}[{DRM slow-rate Rademacher bound; \citealp[arXiv version, Theorem~5.13]{jiao2021error}}]\label{prop:S5}
For the tanh-DRM setting with $u_0 \in H^2(\Omega)$ and $\mathcal{F} = \mathcal{NN}_{\tanh}(L, \mathfrak{D}, B_{\boldsymbol{\theta}})$:
\[
\mathcal{E}^{\mathrm{DRM}}_{sta,\,slow} \;\le\; C(d, |\Omega|, |\partial\Omega|, B_{\mathrm{DRM}}, c_0)\, \frac{d\sqrt{L}\,\mathfrak{D}^{2L}\,B_{\boldsymbol{\theta}}^{2L}}{\sqrt{N}}\,\sqrt{\log\!\big(d L \mathfrak{D} B_{\boldsymbol{\theta}} N\big)},
\]
where $N$ denotes the (equal) interior and boundary sample sizes.
\end{proposition}

\textbf{Step 3: Trade-off and balance.} Theorem~\ref{thm:A2} with $k=1$, $p=2$, $s=2$ allows the choices $\mathfrak{D}\asymp\epsilon^{-d/(1-\mu)}$, $B_{\boldsymbol{\theta}}\asymp\epsilon^{-\vartheta}$ with $\vartheta=(9d+8)/(2(1-\mu))$, and $L=O(\log(d+1))$. The prefactor $\mathfrak{D}^{2L}B_{\boldsymbol{\theta}}^{2L}$ is polynomial in $\epsilon^{-1}$ with depth-dependent exponent
\[
\mathfrak{a} \;:=\; 2L\Big(\tfrac{d}{1-\mu} + \vartheta\Big) \;=\; O_\mu\!\big(d\log(d+1)\big),
\]
because $L$ is logarithmic in $d$. Suppressing logarithmic factors in $N$, substituting $\mathcal{E}^{\mathrm{DRM}}_{app}\le\epsilon^2$ and the Proposition~\ref{prop:S5} bound $\lesssim \epsilon^{-\mathfrak{a}}/\sqrt{N}$ into the symmetric decomposition and balancing $\epsilon^2 \asymp \epsilon^{-\mathfrak{a}}/\sqrt{N}$ gives $\epsilon \asymp N^{-1/(2(2+\mathfrak{a}))}$, which gives the stated algebraic rate. Incorporating the logarithmic factor in Proposition~\ref{prop:S5} changes the logarithmic factors in the bound. For fixed $\mu$, the resulting exponent is of order $1/(d\log(d+1))$: the depth-exponential prefactor $\mathfrak{D}^{2L}B_{\boldsymbol{\theta}}^{2L}$, with $L = O(\log(d+1))$, inflates the statistical complexity.
\end{proof}

The tanh-DRM analysis cited above uses global complexity bounds. We next present a fast-rate result obtained through localization for a different network class.

For a zero-Dirichlet Schr\"odinger problem, Lu et al.\ \cite{lu2021machine} obtain a localized DRM bound using $\mathrm{ReLU}^3$ networks:
\begin{theorem}[DRM fast rate, $\mathrm{ReLU}^3$; \citealp{lu2021machine}, Theorem~B.10]\label{thm:drm_fast}
Consider the zero-Dirichlet Schr\"odinger problem and the sampling, boundedness, and network-class assumptions of \cite{lu2021machine}. Suppose its solution $u_0\in H^s(\Omega)$ for an integer $s\ge2$, and $\mathcal E_{opt}=0$. With the sparse $\mathrm{ReLU}^3$ architecture chosen as in that work, the empirical Deep Ritz minimizer based on $N$ interior observations satisfies, with high probability,
\[
\big\|\widehat{u}_{\boldsymbol{\theta}}-u_0\big\|_{H^1(\Omega)}^2 \;\lesssim\; N^{-\frac{2(s-1)}{2(s-1)+d}}\log N.
\]
\end{theorem}
The proof combines Sobolev approximation with a local Rademacher complexity bound for the Dirichlet energy. Coercivity converts excess energy into $H^1$ error, while the boundedness assumptions control the variance of the loss differences. This result concerns the boundary-constrained class and sampling model of \cite{lu2021machine}, rather than the Neumann ERM in \eqref{eq:drm_erm}.

\smallskip
In both examples, a variance--mean relation and localized complexity control improve the statistical bound. For PDE losses, the norm in which this yields a solution error depends on the stability estimate of the particular problem; DRM and PINN objectives need not control the same Sobolev norm.



\subsection{Minimax lower bounds for the two examples} \label{sec:minimax_lb}

The upper bounds in Section~\ref{sec:classical_rates} describe the accuracy attained by specified neural estimators. Minimax lower bounds provide complementary limits on worst-case accuracy over all estimators under a fixed observation model. We develop a common Fano argument for regression and elliptic recovery, making explicit both the target metric and the information supplied by the observations.

Let $\mathsf d$ denote the metric in which estimation accuracy is evaluated. For a statistical experiment
\[ \mathcal E_n=\{P_\theta^{(n)}:\theta\in\Theta_{\mathrm{stat}}\}, \]
in which $\theta$ determines a target function $\phi_\theta$, define the minimax squared-error risk by
\begin{equation} \mathfrak M_n(\Theta_{\mathrm{stat}},\mathsf d;\mathcal E_n) :=\inf_{\widehat\phi}\sup_{\theta\in\Theta_{\mathrm{stat}}} \mathbb E_\theta\!\left[ \mathsf d^2(\widehat\phi,\phi_\theta) \right]. \label{eq:minimax-risk-abstract} \end{equation}
The infimum is over all measurable estimators based on the $n$ observations. A lower bound
\[ \mathfrak M_n(\Theta_{\mathrm{stat}},\mathsf d;\mathcal E_n)\ge c r_n \]
means that every estimator has expected squared error at least $c r_n$ for some parameter in $\Theta_{\mathrm{stat}}$. It therefore rules out worst-case risk $o(r_n)$ uniformly over this class, independently of the algorithm or network architecture. We call $r_n$ the minimax rate only when a matching upper bound is available for the same target class, error metric, and observation model, including the design and noise law. Standard lower-bound methods include those of Le Cam, Fano, and Assouad~\cite{tsybakov2009lower}; the two examples below use Fano's method.

\subsubsection{A common Fano template}

Fano's method reduces function estimation to identifying one of many separated candidates whose data distributions are difficult to distinguish. Suppose that, for every sufficiently large resolution parameter $m$, one can select parameters $\theta_1,\ldots,\theta_{K_m}\in\Theta_{\mathrm{stat}}$ with target functions
\[ \Phi_m=\{\phi^{(1)},\ldots,\phi^{(K_m)}\}, \qquad \phi^{(j)}:=\phi_{\theta_j}. \]
Write $P_j$ for the law of one observation under $\theta_j$. We measure distinguishability by the Kullback--Leibler divergence $D_{\mathrm{KL}}(P\|Q):=\int\log(\mathrm dP/\mathrm dQ)\,\mathrm dP$ for $P\ll Q$. This divergence is generally neither symmetric nor a metric; small KL divergence limits the evidence available to distinguish the two data laws. The candidate family is chosen so that
\begin{align}
&\log K_m \asymp m^d, \label{eq:fano-cardinality}\\
&\mathsf d^2\!\left(\phi^{(i)},\phi^{(j)}\right) \gtrsim m^{-2\nu_{\mathrm{loss}}}, \qquad i\ne j, \label{eq:fano-separation}\\
&\overline D_{\mathrm{KL},m} :=\frac{1}{K_m^2}\sum_{i,j=1}^{K_m} D_{\mathrm{KL}}(P_i\|P_j) \lesssim m^{-2\nu_{\mathrm{KL}}}. \label{eq:fano-kl}
\end{align}
These conditions provide exponentially many candidates, separation in the target metric, and a bound on the information in one observation. Together, they make choosing the wrong candidate costly in the target metric, while limiting the evidence that each observation supplies for distinguishing among the many candidates. The exponents $\nu_{\mathrm{loss}}$ and $\nu_{\mathrm{KL}}$ distinguish the cost of estimating the wrong target from the difficulty of identifying it through the data.

For any fixed estimator, restricting the worst-case risk to the finite family and averaging gives
\begin{align} \sup_{\theta\in\Theta_{\mathrm{stat}}} \mathbb E_\theta \mathsf d^2(\widehat\phi,\phi_\theta) &\ge \max_{1\le j\le K_m} \mathbb E_{\theta_j} \mathsf d^2(\widehat\phi,\phi^{(j)}) \notag\\ &\ge \frac{1}{K_m}\sum_{j=1}^{K_m} \mathbb E_{\theta_j} \mathsf d^2(\widehat\phi,\phi^{(j)}). \label{eq:minimax-to-finite-average} \end{align}
To interpret this average as an expected risk, introduce an auxiliary experiment. First draw a hidden index $J\sim\operatorname{Unif}\{1,\ldots,K_m\}$, then generate the observations conditionally as $Z_{1:n}\mid\{J=j\}\sim P_j^{(n)}$, where $P_j^{(n)}:=P_{\theta_j}^{(n)}$. Each candidate is selected with probability $1/K_m$, so the average above equals $\mathbb E[\mathsf d^2(\widehat\phi,\phi^{(J)})]$ under this joint experiment. The random index is only a proof device; the original minimax problem does not assume that the unknown parameter is random.

To turn estimation into identification in this auxiliary experiment, convert $\widehat\phi$ into the nearest-candidate decoder
\[ \widehat J\in\arg\min_{1\le j\le K_m} \mathsf d(\widehat\phi,\phi^{(j)}). \]
On the event $\widehat J\ne J$, the triangle inequality and \eqref{eq:fano-separation} give $\mathsf d^2(\widehat\phi,\phi^{(J)})\gtrsim m^{-2\nu_{\mathrm{loss}}}$. Consequently, for a constant $c_1>0$ independent of $m$ and $n$,
\begin{equation} \mathbb E\!\left[ \mathsf d^2(\widehat\phi,\phi^{(J)}) \right] \ge c_1m^{-2\nu_{\mathrm{loss}}} \mathbb P(\widehat J\ne J). \label{eq:testing-to-estimation} \end{equation}
Thus a lower bound on identification error yields a lower bound on estimation risk for every estimator.

Fano's inequality controls this identification error:
\begin{equation} \mathbb P(\widehat J\ne J) \ge 1-\frac{I(J;Z_{1:n})+\log 2}{\log K_m}, \label{eq:fano-inequality} \end{equation}
where $I(J;Z_{1:n})$ is the mutual information between the hidden index and the data. Under the uniform auxiliary index, this mutual information is bounded by the average pairwise KL divergence between the joint data laws $P_j^{(n)}$; this step does not require independent observations. When the observations are independent and identically distributed under each candidate, KL additivity further expresses each joint divergence as $n$ times its one-observation counterpart. Combining these two facts gives
\[
\begin{aligned}
D_{\mathrm{KL}}(P_i^{(n)}\|P_j^{(n)})&=nD_{\mathrm{KL}}(P_i\|P_j),\\
I(J;Z_{1:n})&\le\frac{1}{K_m^2}\sum_{i,j=1}^{K_m}D_{\mathrm{KL}}(P_i^{(n)}\|P_j^{(n)})
=n\overline D_{\mathrm{KL},m}.
\end{aligned}
\]
For sufficiently large $m$, the testing error is bounded below by a positive constant whenever, for a sufficiently small constant $a\in(0,1)$,
\begin{equation} n\overline D_{\mathrm{KL},m} \le a\log K_m. \label{eq:fano-information-condition} \end{equation}
This condition compares the information supplied by the observations with the logarithmic number of candidates to be distinguished.

By \eqref{eq:fano-cardinality} and \eqref{eq:fano-kl}, these quantities have orders $n m^{-2\nu_{\mathrm{KL}}}$ (an upper bound on information) and $m^d$, respectively. Balancing the information bound against the logarithmic number of candidates gives $n m^{-2\nu_{\mathrm{KL}}}\asymp m^d$. Thus choosing
\[ m\asymp n^{1/(d+2\nu_{\mathrm{KL}})} \]
with a sufficiently large multiplicative constant ensures \eqref{eq:fano-information-condition}. Combining \eqref{eq:minimax-to-finite-average}, \eqref{eq:testing-to-estimation}, and Fano's inequality, and then taking the infimum over estimators, gives
\begin{equation} \mathfrak M_n(\Theta_{\mathrm{stat}},\mathsf d;\mathcal E_n) \gtrsim n^{-\frac{2\nu_{\mathrm{loss}}} {d+2\nu_{\mathrm{KL}}}}. \label{eq:generic-fano-rate} \end{equation}
The separation exponent determines the error incurred at this resolution, while the KL exponent determines which resolution the data can distinguish. For localized alternatives of smoothness $s$, a target Sobolev norm of order $t$, and observations through a differential operator of order $q$, the usual scaling is
\[ \nu_{\mathrm{loss}}=s-t, \qquad \nu_{\mathrm{KL}}=s-q. \]
The corresponding squared-error lower-bound scale is
\begin{equation} n^{-\frac{2(s-t)}{d+2(s-q)}}. \label{eq:sobolev-fano-rate} \end{equation}
Applying this template requires verifying the packing and KL bounds for the chosen function class, boundary conditions, observation operator, design distribution, and noise law. Table~\ref{tab:minimax-template} records the substitutions that the next two examples verify.

\begin{table}[htbp] \centering \caption{Substitution of the two examples into the common Fano template ($s,t,q$ are defined in the paragraph preceding the table).} \label{tab:minimax-template} \begin{tabular}{@{}lcccc@{}} \toprule Problem & $(s,t,q)$ & $\nu_{\mathrm{loss}}$ & $\nu_{\mathrm{KL}}$ & Squared-error lower bound \\ \midrule H\"older regression & $(\beta,0,0)$ & $\beta$ & $\beta$ & $n^{-2\beta/(d+2\beta)}$ \\ Noisy-RHS elliptic recovery & $(s,1,2)$ & $s-1$ & $s-2$ & $n^{-2(s-1)/(d+2s-4)}$ \\ \bottomrule \end{tabular} \end{table}

\subsubsection{Regression on a H\"older ball}

We first verify the template for regression with Gaussian noise. Consider
\[ Y_i=f_0(\boldsymbol X_i)+\xi_i, \qquad \boldsymbol X_i\stackrel{\mathrm{i.i.d.}}{\sim}P_X, \qquad \xi_i\stackrel{\mathrm{i.i.d.}}{\sim}\mathcal N(0,\sigma^2), \]
where $\sigma>0$, the noise is independent of the design, and $P_X$ has a density on $[0,1]^d$ bounded above and below by positive constants. The lower bound also applies when the admissible noise class contains $\mathcal N(0,\sigma^2)$, by restriction to this Gaussian submodel. For the bounded-noise results of \S\ref{sec:reg_slow}--\S\ref{sec:reg_fast}, it provides a comparison of exponents rather than a lower bound for the same observation model.

\begin{theorem}[Minimax lower bound for H\"older regression] \label{thm:holder-regression-lower-bound} Let $\beta>0$ and $B_{\mathrm{REG}}>0$. There is a constant $c>0$, independent of $n$, such that
\[ \inf_{\widehat f} \sup_{f_0\in\mathcal H^\beta([0,1]^d;B_{\mathrm{REG}})} \mathbb E_{f_0}\|\widehat f-f_0\|_{L^2(P_X)}^2 \ge c\,n^{-\frac{2\beta}{2\beta+d}}. \]
This is the classical nonparametric regression rate \citep{stone1982optimal}. \end{theorem} \begin{proof}[Proof sketch] Let $b\in C_c^\infty((0,1)^d)$ be a nonzero bump. Place rescaled, disjoint copies $b_1,\ldots,b_{m^d}$ in a regular partition of the cube and, for sufficiently small $c>0$, set
\[ f_{\boldsymbol\tau} =cB_{\mathrm{REG}}m^{-\beta}\sum_{j=1}^{m^d}\tau_j b_j, \qquad \boldsymbol\tau\in\{0,1\}^{m^d}. \]
The scaling ensures that the alternatives lie in the H\"older ball. The Varshamov--Gilbert lemma selects a subfamily with $\log K_m\asymp m^d$ and pairwise Hamming distance of order $m^d$. Hence
\[ \|f_{\boldsymbol\tau}-f_{\boldsymbol\tau'}\|_{L^2(P_X)}^2 \asymp m^{-2\beta}, \]
so $\nu_{\mathrm{loss}}=\beta$. Since the design law is common to all hypotheses, only the conditional Gaussian mean of $Y$ changes. The one-observation Gaussian KL divergence is therefore
\[ D_{\mathrm{KL}}(P_{\boldsymbol\tau}\|P_{\boldsymbol\tau'}) =\frac{1}{2\sigma^2} \|f_{\boldsymbol\tau}-f_{\boldsymbol\tau'}\|_{L^2(P_X)}^2 \lesssim m^{-2\beta}. \]
Thus $\nu_{\mathrm{KL}}=\beta$, and \eqref{eq:generic-fano-rate} gives the stated exponent. \end{proof}

\subsubsection{Noisy right-hand-side observations for an elliptic PDE}

The second example recovers a solution from observations of its right-hand side, so the observation operator enters the KL calculation. Following \cite{lu2021machine}, let $\Omega=(0,1)^d$ and consider the fixed operator
\[ Lu=-\Delta u+\omega u, \]
with homogeneous Dirichlet boundary conditions. Assume that $\omega$ is bounded and that $L$ is coercive and has the elliptic regularity required below. The boundary condition differs from the Neumann condition in the running DRM example (\S\ref{sec:drm_setup}); the resulting lower bound therefore supplies a benchmark rather than a theorem applying directly to the Neumann estimators in Chapters~\ref{sec:classical_setting}--\ref{sec:complete_setting}. The data are
\[ Y_i=(Lu_0)(\boldsymbol X_i)+\xi_i, \qquad \boldsymbol X_i\stackrel{\mathrm{i.i.d.}}{\sim}\operatorname{Unif}(\Omega), \qquad \xi_i\stackrel{\mathrm{i.i.d.}}{\sim}\mathcal N(0,\sigma^2), \]
with $\sigma>0$ and independent design and noise. Define
\begin{align*}
 \mathcal U_s(B) :=\left\{u\in H^s(\Omega)\cap H_0^1(\Omega): \|u\|_{H^s(\Omega)}\le B\right\}.     
\end{align*}
\begin{theorem}[Minimax lower bound for noisy-RHS elliptic recovery] \label{thm:noisy-rhs-lower-bound} Let $s\ge2$ and $B>0$. Under the preceding assumptions, there is a constant $c>0$, independent of $n$, such that
\[ \inf_{\widehat u} \sup_{u_0\in\mathcal U_s(B)} \mathbb E_{u_0}\|\widehat u-u_0\|_{H^1(\Omega)}^2 \ge c\,n^{-\frac{2(s-1)}{d+2s-4}}. \]
This is the noisy-RHS lower-bound exponent studied for the Deep Ritz method in \citealp{lu2021machine}. \end{theorem} \begin{proof}[Proof sketch] Using the same type of disjoint localized bumps, define
\[ u_{\boldsymbol\tau} =cBm^{-s}\sum_{j=1}^{m^d}\tau_j b_j, \qquad \boldsymbol\tau\in\{0,1\}^{m^d}, \]
where the supports lie strictly inside $\Omega$ and $c>0$ is small enough that $u_{\boldsymbol\tau}\in\mathcal U_s(B)$. A Varshamov--Gilbert subfamily has $\log K_m\asymp m^d$ and
\[ \|u_{\boldsymbol\tau}-u_{\boldsymbol\tau'}\|_{H^1(\Omega)}^2 \asymp m^{-2(s-1)}. \]
Therefore $\nu_{\mathrm{loss}}=s-1$.

The observations depend on $Lu_{\boldsymbol\tau}$ rather than directly on $u_{\boldsymbol\tau}$. The design distribution and Gaussian noise variance are common to all hypotheses, so for one observation,
\[ D_{\mathrm{KL}}(P_{\boldsymbol\tau}\|P_{\boldsymbol\tau'}) =\frac{1}{2\sigma^2|\Omega|} \|L(u_{\boldsymbol\tau}-u_{\boldsymbol\tau'})\|_{L^2(\Omega)}^2 \lesssim m^{-2(s-2)}. \]
(For the present unit cube, $|\Omega|=1$.) Thus $\nu_{\mathrm{KL}}=s-2$. Substitution into \eqref{eq:generic-fano-rate} yields
\[ \mathfrak M_n \gtrsim n^{-\frac{2(s-1)}{d+2(s-2)}} =n^{-\frac{2s-2}{d+2s-4}}. \]
\end{proof}

The two examples differ in how the target metric and observation law determine the exponents. Regression uses an $L^2$ distance between functions in both the target error and the Gaussian KL calculation. Elliptic recovery measures solution error in $H^1$, while the KL divergence is controlled by the $L^2$ distance between $Lu$ for a second-order operator. For the latter problem, Lu et al.\ \cite[Appendix~C]{lu2021machine} compare the lower bound with a modified DRM using a Fourier approximation space and separate quadrature samples. Their construction uses $n$ observations for the right-hand-side term and a larger independent sample for the gradient-energy term, with a specified relation between the sample sizes. This approximation space and sampling allocation differ from those of the neural-network DRM estimator above. The upper-bound calculation also retains a $1/n$ contribution. The displayed lower bound is therefore not sharp in every dimension: in dimension one its exponent exceeds one, whereas a one-parameter Gaussian submodel already gives an $n^{-1}$ lower bound.

For comparison with the preceding upper bounds, the observation model is as important as the smoothness class and target norm. The PDE result above concerns noisy random observations of the right-hand side. This differs from both a setting in which the PDE data are fully known and Monte Carlo quadrature is an algorithmic choice, and an information-based-complexity model in which an unknown black-box right-hand side is accessible only through finitely many noiseless point evaluations. These three information models impose different constraints on estimators and require separate lower-bound comparisons.

\subsection{Escaping the curse of dimensionality} \label{sec:low_intrinsic}

For classical nonparametric estimation over general isotropic H\"older or Sobolev classes under standard observation models, the ambient dimension $d$ enters the minimax statistical rates and cannot be removed without additional structure. At the approximation level the situation is more subtle: the dimension dependence can also be altered by changing the representation system itself.

We organize the literature into three broad mechanisms for weakening ambient-dimensional dependence. The first imposes structural restrictions on the target function or on the underlying problem, so that the relevant function class is substantially smaller than a generic $d$-variate smoothness class. The second exploits low-dimensional geometry of the data distribution, for example when the design is concentrated on or near a manifold or another set of small covering dimension. The third changes the representation system itself, using architectures or activation functions whose approximation power can exceed that of standard piecewise-polynomial networks.

These mechanisms capture different sources of complexity reduction, and the associated dimension parameters generally require separate justification. The data-side route is developed in the greatest detail below because it is the one used in the downstream GAN and reinforcement-learning analyses of \S\ref{sec:applications}.

\subsubsection{Structure of the target and the problem} \label{sec:target_problem_structure}

Target-side assumptions can weaken ambient-dimensional dependence in several different ways. It is useful to distinguish reductions in local input dimension from changes in function-space complexity and from algebraic or problem-specific tractability.

The most direct mechanism is a reduction in the number of variables that enter each constituent function. For example, a ridge-type target
\[ \phi^*(\boldsymbol{x}) = g(\boldsymbol A^\top\boldsymbol{x}), \qquad \operatorname{rank}(\boldsymbol A)=d_{\mathrm{act}}\ll d \]
depends only on $d_{\mathrm{act}}$ active directions even when the design distribution is fully $d$-dimensional. Additive and interaction models similarly take the form
\[ \phi^*(\boldsymbol{x}) = \sum_k g_k(\boldsymbol{x}_{I_k}), \qquad |I_k|\leq d_{\mathrm{loc}}\ll d, \]
while hierarchical composition models replace a single $d$-variate map by a composition of functions with small local input dimensions. Deep networks can exploit these local dependencies, so that the relevant rates are governed by componentwise smoothness and local arity rather than by the ambient dimension alone \citep{bauer2019deep,schmidt2020nonparametric,kohler2023opt}. In such models there need not be a single effective dimension: different layers or components may be governed by different local dimensions.

A second mechanism changes the complexity of the function class without making the input itself low-dimensional. Anisotropic and mixed-smoothness classes encode direction-dependent or mixed regularity, while Barron and related spectral classes control Fourier or variation norms rather than a generic isotropic H\"older modulus. Their approximation or estimation complexities can therefore have substantially weaker dependence on $d$ than those of an unrestricted $d$-variate H\"older ball \citep{suzuki2019adaptivity,barron1993universal}. Here the relevant complexity parameter is generally not a geometric dimension.

Finally, high-dimensional tractability may arise from algebraic or problem-specific structure. Low-rank tensor representations, symmetry or equivariance, and probabilistic representations of certain high-dimensional PDEs can reduce representational or computational complexity even when neither the target nor the data support is low-dimensional \citep{cohen2016expressive,elesedy2021provably,hutzenthaler2020proof}. Thus, reducing ambient-dimensional dependence does not always amount to identifying a single intrinsic dimension.

\subsubsection{Intrinsic dimension of the data} \label{sec:data_geometric_statistical_dimension}

Let $P_X$ be the distribution of the input $\boldsymbol{X}$ and let
\[ \mathcal{X} := \operatorname{supp}(P_X) \subseteq [0,1]^d. \]
The prediction error is measured only where the distribution places mass:
\[ \|\phi-\phi^*\|_{L^2(P_X)}^2 = \int_{\mathcal{X}} |\phi(\boldsymbol{x})-\phi^*(\boldsymbol{x})|^2 \,\mathrm{d}P_X(\boldsymbol{x}). \]
Consequently, the approximation and statistical complexity may depend on the geometry of $\mathcal{X}$ rather than on the full cube $[0,1]^d$.

For a compact set $\mathcal{X}$ with $\epsilon$-covering number $\mathcal{C}(\epsilon,\mathcal{X},\|\cdot\|_2)$, define its upper Minkowski dimension by
\[ \dim_{\mathrm{M}}(\mathcal{X}) := \limsup_{\epsilon\downarrow 0} \frac{\log\mathcal{C}(\epsilon,\mathcal{X},\|\cdot\|_2)} {\log(1/\epsilon)}. \]
Throughout this subsection, $\dim_{\mathrm{M}}$ refers to the upper Minkowski dimension, and we write $d^*=\dim_{\mathrm{M}}(\mathcal{X})$. By definition, $\mathcal{C}(\epsilon,\mathcal{X},\|\cdot\|_2)\lesssim\epsilon^{-d^\sharp}$ for every $d^\sharp>d^*$ and all sufficiently small $\epsilon$; obtaining the exact exponent $d^*$ uniformly requires an additional covering-regularity condition.

A particularly structured case is exact smooth-manifold support. Suppose that $P_X(\mathcal{M})=1$, where $\mathcal{M}\subseteq[0,1]^d$ is a compact $d_{\mathcal{M}}$-dimensional Riemannian submanifold with controlled reach and volume. Then $\dim_{\mathrm{M}}(\mathcal{M})=d_{\mathcal{M}}$. Local coordinate charts reduce approximation on $\mathcal{M}$ to approximation on subsets of $\mathbb{R}^{d_{\mathcal{M}}}$. Under the regularity assumptions of Jiao et al.\ \cite{jiao2021deep}, a ReLU regression estimator satisfies, up to logarithmic factors,
\[ \mathbb{E} \bigl\| \widehat f_n-f_0 \bigr\|_{L^2(P_X)}^2 \lesssim n^{-\frac{2\beta}{2\beta+d_{\mathcal{M}}}} \operatorname{polylog}(n). \]
The exponent depends on the manifold dimension, while the prefactor can still depend on the ambient dimension and on geometric quantities such as the reach and volume of $\mathcal{M}$. Under standard regularity conditions on the design measure, the corresponding nonparametric regression problem has the minimax rate $n^{-2\beta/(2\beta+d_{\mathcal{M}})}$ \citep{bickel2007local,yang2016bayesian}. Thus the ReLU estimator of \cite{jiao2021deep} is minimax-optimal up to logarithmic factors under exact manifold support.

Exact smooth-manifold assumptions impose substantial geometric regularity and exclude, for example, singular or fractal supports. A more general description of low-dimensional geometry is provided by Minkowski dimension: $d^*=\dim_{\mathrm{M}}(\mathcal{X})<d$. For such sets, Nakada and Imaizumi \cite{nakada2020adaptive} show that for every $d^\sharp>d^*$, their estimator attains $\|\widehat f_n-f_0\|_{L^2(P_X)}^2\lesssim n^{-2\beta/(2\beta+d^\sharp)}(1+\log n)^2$, recovering the intrinsic-dimensional exponent up to an arbitrarily small relaxation and logarithmic factors. Jiao et al.\ \cite{jiao2021deep} give an alternative Johnson--Lindenstrauss-type construction with projected dimension $d_0\gtrsim d^*/\delta^2$, where $\delta\in(0,1)$ denotes the projection distortion, that improves the ambient-dimensional prefactor.

Data may instead lie near, rather than exactly on, a manifold. One model takes $\operatorname{supp}(P_X)\subseteq\mathcal{M}_{\rho_{\mathrm{tube}}}$, the $\rho_{\mathrm{tube}}$-neighborhood of a $d_{\mathcal{M}}$-dimensional manifold $\mathcal{M}$ within $[0,1]^d$. Here $\rho_{\mathrm{tube}}\in(0,1)$ is a geometric radius, distinct from the projection distortion $\delta\in(0,1)$. Although this neighborhood has ambient Minkowski dimension $d$ for every fixed positive radius, its geometry can still be exploited at resolutions coarser than the tube width. The approximate-manifold analysis of Jiao et al.\ \cite[Theorem~6.1]{jiao2021deep} introduces an effective projection dimension $d_{\delta}=O(d_{\mathcal{M}}\log(d/\delta)/\delta^{2})$ and obtains, under its distributional assumptions and a condition relating $\rho_{\mathrm{tube}}$ to the approximation resolution,
\[ \mathbb{E} \bigl\| \widehat f_n-f_0 \bigr\|_{L^2(P_X)}^2 \lesssim n^{-\frac{2\beta}{2\beta+d_{\delta}}} \operatorname{polylog}(n). \]
For this approximate-manifold construction, removing the remaining ambient-dimensional dependence of the rate remains open \cite{jiao2021deep}. Across these cases, $d_{\mathcal{M}}$ and $d^*$ describe the support geometry, whereas $d_{\delta}$ and $d_0$ arise from specific projection constructions. These dimension parameters should therefore not be substituted for one another without a separate argument.

In settings where the geometric assumption yields both an intrinsic-dimensional approximation bound and the corresponding statistical-complexity bound, the generic balances of \S\ref{sec:reg_slow} take the schematic form (with $d_{\mathrm{eff}}$ understood as the appropriate dimension parameter from the cases above)
\[ \mathcal{E}^{\mathrm{REG}}_{\mathrm{slow}}(d_{\mathrm{eff}}) \lesssim n^{-\frac{2\beta}{4\beta+d_{\mathrm{eff}}}} \operatorname{polylog}(n), \qquad \mathcal{E}^{\mathrm{REG}}_{\mathrm{fast}}(d_{\mathrm{eff}}) \lesssim n^{-\frac{2\beta}{2\beta+d_{\mathrm{eff}}}} \operatorname{polylog}(n). \]
Matching these upper bounds to minimax lower bounds requires specifying the design class. For a fixed design $P_X$, support dimension alone does not guarantee the $L^2(P_X)$-separation needed in a Fano construction; a local mass-regularity condition such as $c_{P_X} r^{d^*}\leq P_X(B(\boldsymbol{x},r))\leq C_{P_X} r^{d^*}$ is natural. When the supremum is instead taken jointly over regression functions and design measures whose supports have Minkowski dimension at most $d^*$, Nakada and Imaizumi \cite{nakada2020adaptive} establish the lower bound $n^{-2\beta/(2\beta+d^*)}$.

The geometric refinement is used downstream in the generative learning and reinforcement learning analyses of \S\ref{section:app:generative} and \S\ref{section:app:rl}.

\subsubsection{Super-expressive activations} \label{sec:representation_dimension_escape}

A separate literature obtains unusually strong approximation rates by changing the activation family or the way information is encoded in the network. Periodic activations can implement efficient lookup operations and yield nearly exponential approximation rates in certain parameter regimes \citep{yarotsky2020phase}. Networks constructed from floor, exponential, and step activations can approximate H\"older functions with errors that decay exponentially in a width parameter \citep{shen2021neural}. Related constructions combine ReLU with oscillatory or exponential activations.

These results concern representational efficiency rather than a reduction of intrinsic dimension: they alter the relation between nominal network size and approximation power. Unusually fast approximation therefore need not imply comparable statistical or computational tractability. Parameter precision, hypothesis-class complexity, and the ability to recover the constructive network from data remain separate questions \citep{yarotsky2020phase}. We therefore distinguish approximation, statistical, and computational tractability, and do not invoke super-expressive constructions in the downstream rate analyses. Of the three routes surveyed in this section, only the data-geometry route of \S\ref{sec:data_geometric_statistical_dimension} is drawn on in the later applications (\S\ref{sec:applications}).

%% file: sections/overparameterization.tex
\section{Over-Parameterized Error Analysis via Norm Control}\label{sec:overparam_setting}

The previous chapter analyzed the two-error coupling under size-controlled parameterizations of the network class: width $W$ for Regression, nonzero-parameter count $\mathfrak{D}$ for DRM. In modern deep-learning practice, however, networks are typically heavily over-parameterized, with $\mathfrak{D}$ exceeding the sample size $n$ by orders of magnitude, and both size-based statistical bounds of \S\ref{sec:classical_stat} become vacuous. This chapter develops a parallel two-error analysis in which the controlling parameter is replaced by a norm on the network's weight matrices. For ReLU regression, norm control removes the width dependence from the statistical bound while retaining a depth factor. DRM requires a separate analysis of the derivative class.

Double descent \cite{belkin2019reconciling} and benign overfitting \cite{bartlett2020benign} motivate analyses beyond network size alone. Nagarajan and Kolter \cite{nagarajan2019uniform} also demonstrate limitations of uniform-convergence explanations of generalization in specific learning settings. As in the preceding chapter, the analysis here assumes exact empirical minimization, $\mathcal{E}_{opt}=0$, and focuses on approximation and statistical errors, now for norm-constrained network classes.

\begin{definition}[Norm-constrained network class]\label{sec:overparam_setup}
Define the augmented weight matrices
\[
\tilde{\boldsymbol{A}}_{L-1} := (\boldsymbol{A}_{L-1}, \boldsymbol{b}_{L-1}), \qquad
\tilde{\boldsymbol{A}}_{\ell} := \begin{pmatrix} \boldsymbol{A}_\ell & \boldsymbol{b}_\ell \\ \mathbf{0} & 1 \end{pmatrix}, \quad \ell = 0, \ldots, L-2,
\]
and the operator norm $\|\boldsymbol{A}\| := \max_i \sum_j |a_{ij}|$. The norm product $\kappa(\boldsymbol{\theta}) := \prod_{\ell=0}^{L-1}\|\tilde{\boldsymbol{A}}_\ell\|$ bounds the network's Lipschitz constant in the $\ell^\infty$ input norm when $\rho$ is $1$-Lipschitz, and the norm-constrained class is
\[
\mathcal{NN}_\rho(W, L, \mathcal{K}) := \big\{\phi_{\boldsymbol{\theta}} \in \mathcal{NN}_\rho(W, L) : \kappa(\boldsymbol{\theta}) \le \mathcal{K}\big\},
\]
where $\rho$ denotes the activation function.
\end{definition}
The chapter studies convergence rates as $\mathcal{K}$ varies, allowing large widths while tracking the remaining depth dependence.

\subsection{Norm-controlled approximation theorems}\label{sec:overparam_approx}

The norm-controlled refinements extend Theorems~\ref{thm:A1} and~\ref{thm:A2} to a setting where the weight-norm product $\kappa(\boldsymbol{\theta})$, rather than the network size, is the controlling parameter. For ReLU-activated networks, the fundamental result is a refinement of Theorem~\ref{thm:A1} in which the controlling parameter is the norm product $\mathcal{K}$ rather than the network size $W \cdot L$.

\begin{theorem}[Norm-controlled ReLU approximation; \citealp{jiao2023approximation}]\label{thm:A3}
Let $d \in \mathbb{N}$ and $\beta = s + r > 0$ with $s \in \mathbb{N}_0$, $r \in (0, 1]$. There exists a constant $C > 0$ such that for any $\mathcal{K} \ge 1$, taking $W \ge C \mathcal{K}^{(2d + \beta)/(2d + 2)}$ and $L \ge 2 \lceil \log_2(d + s)\rceil + 3$, the norm-constrained ReLU class $\mathcal{NN}_{\mathrm{ReLU}}(W, L, \mathcal{K})$ satisfies
\[
\sup_{f \in \mathcal{H}^\beta([0,1]^d, B_{\mathrm{REG}})}\, \inf_{\phi \in \mathcal{NN}_{\mathrm{ReLU}}(W, L, \mathcal{K})}\, \|f - \phi\|_{C([0,1]^d)} \;\lesssim\; \mathcal{K}^{-\frac{\beta}{d + 1}}.
\]
\end{theorem}

\begin{proof}[Proof sketch]
The construction uses local Taylor polynomials and norm-controlled multiplication networks. Their accuracy and norm budgets are chosen jointly: the width increases with the desired accuracy, whereas the required depth depends only on $d$ and $\beta$. This gives the stated relation between approximation error and $\mathcal K$; see \cite[Theorem~3.2]{jiao2023approximation}.
\end{proof}

For tanh-activated networks, the counterpart is due to \cite{ding2025convergence}. The Taylor-scaling construction of \S\ref{sec:A2} naturally produces a linear combination of $\mathfrak{m}$ norm-constrained tanh sub-networks, denoted $\mathcal{NN}^{\mathfrak{m}}_{\tanh}(W, L, \mathcal{K})$.

\begin{theorem}[Norm-controlled tanh approximation; \citealp{ding2025convergence}, Theorem~4.2]\label{thm:A4}
Let $s,k\in\mathbb N_0$, $d\in\mathbb N$, with $s \ge k + 1$, $1 \le p \le \infty$, and $0<\mu<1$ with $s-k-\mu k>0$. There exist constants $C = C(s,d,k,p,\mu,B) > 0$ and $\epsilon^* > 0$ such that for every $\epsilon \in (0, \epsilon^*]$ and every $f \in W^{s,p}([0,1]^d)$ with $\|f\|_{W^{s,p}} \le B$, there is a tanh-activated network $\phi \in \mathcal{NN}^{\mathfrak{m}}_{\tanh}(W, L, \mathcal{K})$ writable as $\phi = \sum_{j=1}^{\mathfrak{m}} c_j \phi^j$, with sub-networks $\phi^j \in \mathcal{NN}_{\tanh}(W_{\mathrm{sub}}, L, \kappa(\boldsymbol{\theta}^j))$ and aggregate norm budget $\sum_{j=1}^{\mathfrak{m}} |c_j|\, \kappa(\boldsymbol{\theta}^j)\le\mathcal{K}$. The parameters satisfy
\begin{align*}
\mathfrak{m} &= C\big(\lceil (\epsilon/2C)^{-1/(s - k - \mu k)}\rceil + 1\big)^d, \quad
W_{\mathrm{sub}} = 2^{\lceil \log_2(d + s - 1)\rceil + 1}, \quad
W \ge \mathfrak{m} \cdot W_{\mathrm{sub}}, \\
L &= \lceil \log_2(d + s - 1)\rceil + 2, \\
\mathcal{K} &= C\, \epsilon^{-\lceil \log_2(d + s - 1)\rceil - 1 - \big(d\lceil \log_2(d + s - 1)\rceil + 2d + 1 + \mu\big)/(s - k - \mu k)},
\end{align*}
and $\|f - \phi\|_{W^{k,p}([0,1]^d)} \le \epsilon$.
\end{theorem}

\begin{proof}[Proof sketch]
The proof reuses the Taylor-scaling construction of \S\ref{sec:A2}: an approximate partition of unity over a grid of $O(\epsilon^{-d/(s-k-\mu k)})$ cells yields a local Taylor expansion of $f$, and each cell--monomial pair is approximated by a small tanh network via the atomic $x, x^2$ construction. The construction produces one sub-network per cell--monomial pair, giving $\mathfrak{m} \asymp \epsilon^{-d/(s-k-\mu k)}$. Each sub-network $\phi^j$ is tracked individually for its norm product $\kappa(\boldsymbol{\theta}^j)$; the outer coefficients $\{c_j\}$ come from the Taylor coefficients of the local expansion. The aggregate norm budget $\mathcal{K}$ then inherits the polynomial scaling in $\epsilon^{-1}$. Full bookkeeping is in \cite{ding2025convergence}.
\end{proof}

The parallel neural network (PNN, Definition~\ref{def:pnn}) studied in Chapter~\ref{sec:complete_setting} also uses a linear combination of sub-networks. There, the outer coefficients and inner parameters satisfy separate constraints, $\sum_k |c_k|\le M$ and $\|\boldsymbol{\theta}^k\|_\infty\le B_{\boldsymbol{\theta}}$, which define the projection sets for the PGD analysis of \S\ref{sec:complete_opt}. These separate bounds do not follow from the aggregate norm budget alone; the corresponding approximation guarantee is given in Theorem~\ref{thm:A5}.

\subsection{Norm-controlled statistical bound}\label{sec:overparam_stat}

The norm-controlled statistical bound, due to \cite{golowich2018size}, replaces the network-size dependence of Proposition~\ref{prop:S3} with a norm-product dependence, via a contraction-and-peeling argument that bypasses covering numbers entirely. The central result is the following lemma \cite{golowich2018size}.

\begin{lemma}[Norm-controlled Rademacher bound; \citealp{jiao2023approximation}, Lemma~2.3]\label{lem:S6}
Let $\rho=\mathrm{ReLU}$, and let $\mathcal{F} = \mathcal{NN}_\rho(W, L, \mathcal{K})$. For any $\boldsymbol{x}_1, \ldots, \boldsymbol{x}_n \in [0, 1]^d$,
\[
\widehat{\mathfrak{R}}_n(\mathcal{F}) \;\le\; \frac{2 \mathcal{K} \sqrt{L + 2 + \log(d + 1)}}{\sqrt{n}}.
\]
\end{lemma}

\begin{proof}[Proof sketch]
Lemma~\ref{lem:S6} avoids covering numbers entirely by a direct layer-wise contraction. The key inequality, used here for the positively homogeneous ReLU activation, any matrix norm-bound $\|\boldsymbol{A}\| \le R$, and any convex monotone $g : \mathbb{R} \to [0, \infty)$, is
\begin{equation}\label{eq:golowich_contract}
\mathbb{E}_{\boldsymbol{\sigma}} \sup_{f \in \mathcal{F}',\, \|\boldsymbol{A}\| \le R} g\!\Big(\big\|\!\sum_i \sigma_i \rho(\boldsymbol{A} f(\boldsymbol{x}_i))\big\|_\infty\Big) \;\le\; 2\, \mathbb{E}_{\boldsymbol{\sigma}} \sup_{f \in \mathcal{F}'} g\!\Big(R \big\|\!\sum_i \sigma_i f(\boldsymbol{x}_i)\big\|_\infty\Big).
\end{equation}
The operator-norm bound $R$ appears outside the activation as a multiplicative factor, independently of the matrix width. Positive homogeneity first allows the ReLU network to be represented with normalized hidden-layer bounds and the norm budget absorbed into the output layer; see \cite[Proposition~2.4]{jiao2023approximation}. Applying \eqref{eq:golowich_contract} with $g(t)=\exp(\lambda t)$, $\lambda>0$, then introduces a factor of $2$ per layer in the exponential-moment bound. Taking logarithms converts these factors into a term of order $L$, while the maximum over the augmented input coordinates contributes $\log(d+1)$. Optimizing over $\lambda$ yields the factor $\sqrt{L+2+\log(d+1)}$, with $\mathcal K$ as a multiplicative norm bound. Full details of the exponential-moment argument are in \cite{golowich2018size}.
\end{proof}

The structural consequence is that the Rademacher complexity of an arbitrarily wide network is bounded by $O(\mathcal{K}\sqrt{L+\log(d+1)}/\sqrt{n})$, regardless of $W$: size-based bounds give $W L / \sqrt{n}$, which blows up; norm-based bounds remain independent of width for fixed depth and norm budget. Except for this single substitution, the empirical-process machinery of \S\ref{sec:classical_stat} (symmetrization, contraction, concentration) applies unchanged. The instantiation on the two running examples, including the derivative-aware DRM treatment, is carried out alongside the trade-off derivations in \S\ref{sec:overparam_rates}.

\subsection{Over-parameterized rates for the running examples}\label{sec:overparam_rates}

This subsection synthesizes the norm-controlled approximation and statistical components developed in \S\ref{sec:overparam_approx}--\S\ref{sec:overparam_stat}, integrating them through the trade-off framework established in \S\ref{sec:reg_slow}. \S\ref{sec:reg_overparam} treats Regression, following \cite{jiao2023approximation} with the minimax-optimal refinement of \cite{yang2024nonparametric,yang2024optimal}; \S\ref{sec:drm_overparam} treats DRM, following \cite{ding2025convergence}.

\subsubsection{Regression over-parameterized rate}\label{sec:reg_overparam}

For Regression, the three-step pattern of \S\ref{sec:reg_slow} carries over with two ingredients replaced by their norm-controlled counterparts: Theorem~\ref{thm:A3} \cite{jiao2023approximation} supplies the approximation bound, and Lemma~\ref{lem:S6} \cite{golowich2018size} replaces the size-dependent Rademacher analysis. The controlling parameter is the weight-norm product $\mathcal{K}$ rather than the network width.

\begin{theorem}[Bounded-response consequence of norm-controlled approximation]\label{thm:reg_overparam}
Suppose $f_0 \in \mathcal{H}^\beta([0,1]^d, B_{\mathrm{REG}})$ for some $\beta > 0$, the noise is bounded almost surely, and $\mathcal{E}_{opt} = 0$. Choose $\mathcal F_{NN}=\{\min(\mathcal B,\max(-\mathcal B,\phi)):\phi\in\mathcal{NN}_{\mathrm{ReLU}}(W,L,\mathcal K)\}$, where the fixed constant $\mathcal B$ bounds $|f_0|$ and $|Y|$, with $L$ fixed and satisfying Theorem~\ref{thm:A3} and $\mathcal{K} \asymp n^{(d+1)/(2(d+1+2\beta))}$. Then, for any $W\ge C\mathcal K^{(2d+\beta)/(2d+2)}$ with $C$ as in Theorem~\ref{thm:A3},
\[
\mathbb{E}_{S_n}\big\|\widehat f_{\boldsymbol{\theta}} - f_0\big\|_{L^2(P_X)}^2 \;\le\; C(d, \beta, B_{\mathrm{REG}},\mathcal B,L)\, n^{-\frac{\beta}{d+1+2\beta}}.
\]
\end{theorem}

\begin{proof}[Proof sketch]
Squaring Theorem~\ref{thm:A3} and using $\|\cdot\|_{L^2(P_X)}\le\|\cdot\|_\infty$ gives the approximation bound
\begin{equation}\label{eq:reg_op_app}
\mathcal{E}^{\mathrm{REG,} op}_{app} \;\le\; C(d, \beta, B_{\mathrm{REG}})\, \mathcal{K}^{-\frac{2\beta}{d+1}}.
\end{equation}
For the statistical error, two-sided symmetrization, squared-loss Talagrand contraction (constant $4\mathcal{B}$, where $\mathcal{B}$ bounds $|f|$ and $|Y|$), and Lemma~\ref{lem:S6} combine to give
\begin{equation}\label{eq:reg_op_sta}
\mathcal{E}^{\mathrm{REG,} op}_{sta} \;\le\; C(\mathcal{B})\, \frac{\mathcal{K}\sqrt{L + \log(d+1)}}{\sqrt{n}},
\end{equation}
which is linear in $\mathcal{K}$ and independent of $W$, the defining feature of the norm-controlled regime. Substituting \eqref{eq:reg_op_app} and \eqref{eq:reg_op_sta} into the symmetric decomposition \eqref{eq:reg_sym_decomp} and balancing $\mathcal{K}^{-2\beta/(d+1)} \asymp \mathcal{K}/\sqrt{n}$ gives $\mathcal{K}^* \asymp n^{(d+1)/(2(d+1+2\beta))}$, which recovers the stated rate.
\end{proof}

Theorem~\ref{thm:reg_overparam} is a bounded-response consequence of the norm-controlled approximation and Rademacher bounds, rather than a literal restatement of \cite[Theorem~4.1]{jiao2023approximation}. The cited theorem treats Gaussian observation noise and obtains the same polynomial exponent with an additional logarithmic factor (and an optimization residual). Under the bounded-response assumption used here, direct contraction of the clipped squared loss gives the displayed bound without that factor.

The exponent $\beta/(d+1+2\beta)$ is not minimax-optimal; it results from balancing the global $1/\sqrt n$ statistical bound with the $\mathcal K^{-\beta/(d+1)}$ approximation bound. Yang and Zhou \cite{yang2024nonparametric} study shallow ReLU networks with a variation-norm constraint and combine sharper approximation estimates with localized complexity bounds to obtain the minimax rate, up to logarithmic factors, when $0<\beta<(d+3)/2$. Related approximation results for $\mathrm{ReLU}^k$ networks appear in \cite{yang2024optimal}. These results concern empirical risk minimizers. Chapter~\ref{sec:complete_setting} studies algorithm outputs under the architectures and optimization assumptions specified there.

\subsubsection{DRM over-parameterized rate}\label{sec:drm_overparam}

For DRM, the over-parameterized analysis of \cite{ding2025convergence} uses norm-controlled Sobolev approximation and a statistical bound for the energy, including its derivative terms. The paper treats Robin boundary conditions; taking the diffusion matrix to be the identity and the Robin coefficient to be zero includes the Neumann problem of \S\ref{sec:drm_setup}. We give a version with explicit layerwise and outer-coefficient constraints and verify the uniform statistical bound for this candidate class below.

\begin{theorem}[Mean-square DRM bound with layerwise constraints; adapted from \citealp{ding2025convergence}]\label{thm:drm_overparam}
Let $\Omega\subset[0,1]^d$ be a bounded convex domain with smooth boundary, and let $u_0\in H^2(\Omega)$ solve the Neumann problem under the bounded coefficient and data assumptions of \S\ref{sec:drm_setup}. Assume $\mathcal E_{opt}=0$. For fixed $\mu\in(0,1/3]$, write $q_d=\lceil\log_2(d+1)\rceil$. For sufficiently large $\mathcal K$, choose
\begin{align*}
\mathfrak m &\asymp C(d)\big(\lceil\mathcal K^{1/(C_1d q_d)}\rceil\big)^d,
\qquad W\ge C(d)\mathfrak m,\\
L&=q_d+2,\qquad
N\asymp C(d,\Omega,B_{\mathrm{DRM}})\mathcal K^{C_2d/\log_2(d+1)},
\end{align*}
where $N$ is each of the interior and boundary sample sizes, and $C_1,C_2$ are sufficiently large absolute constants. Set
\[
S_{\mathcal K}=C(d)\mathcal K^{3/(2C_1q_d)},
\qquad
R_{\mathcal K}=C(d)\mathcal K^{(3d+4)/(C_1d q_d)}.
\]
Let the ERM candidate class be
\begin{equation}\label{eq:drm_op_candidate}
\begin{aligned}
\mathcal F_{NN}=\Big\{u={}&\sum_{\nu=1}^{\mathfrak m}c_\nu\phi^\nu:
 \ \phi^\nu\in\mathcal{NN}_{\tanh}(W_\nu,L),\
 \sum_{\nu=1}^{\mathfrak m}W_\nu\leq W,\\
&\sum_{\nu=1}^{\mathfrak m}|c_\nu|\kappa(\boldsymbol\theta^\nu)\leq\mathcal K,\
 \sum_{\nu=1}^{\mathfrak m}|c_\nu|\leq S_{\mathcal K},\
 \|\tilde{\boldsymbol A}_{\nu,\ell}\|\leq R_{\mathcal K}
 \quad\text{for all }\nu,\ell\Big\}.
\end{aligned}
\end{equation}
Here $\tilde{\boldsymbol A}_{\nu,\ell}$ is the augmented matrix of layer $\ell$ in sub-network $\nu$, with the norm of Definition~\ref{sec:overparam_setup}; $W$ bounds the total parallel width. Multiplicative constants may also depend on the fixed $\mu$, the domain, the data bounds, and $\|u_0\|_{H^2(\Omega)}$. Then the ERM over \eqref{eq:drm_op_candidate} satisfies
\[
\mathbb E\|\widehat u_{\boldsymbol\theta}-u_0\|_{H^1(\Omega)}^2
\lesssim \mathcal K^{-\frac{1-\mu}{C_1d\lceil\log_2(d+1)\rceil}}.
\]
Equivalently, the mean-square rate is
$N^{-\frac{(1-\mu)\log_2(d+1)}{C_1C_2d^2\lceil\log_2(d+1)\rceil}}$,
of order $N^{-c_\mu/d^2}$ for a positive constant $c_\mu$ independent of $N$. The bound is uniform over larger widths satisfying the same coefficient, layerwise, and aggregate norm budgets.
\end{theorem}

\begin{proof}[Proof sketch]
For the Neumann problem, Proposition~\ref{prop:drm_eq} and \eqref{eq:drm_sym_decomp} reduce the mean-square error to approximation and statistical terms. Set $a_{\mathrm{DRM}}=(1-\mu)/(C_1d q_d)$ and $\epsilon=\mathcal K^{-a_{\mathrm{DRM}}}$. Apply the local-polynomial tanh construction underlying \cite[Corollary~4.1]{ding2025convergence} to a bounded Sobolev extension of $u_0$, and restrict the resulting network to $\Omega$. To track its membership in \eqref{eq:drm_op_candidate}, choose grid resolution $G\asymp\epsilon^{-1/(1-\mu)}$ and local network precision
\[
\delta\asymp\epsilon G^{-(3d/2+1+\mu)}
\asymp\epsilon^{(3d/2+2)/(1-\mu)}.
\]
There are $O(G^d)$ sub-networks of width depending only on $d$, and each outer coefficient is bounded by $C(d)G^{d/2}$. Hence their absolute sum is at most $C(d)G^{3d/2}\leq S_{\mathcal K}$. Each sub-network has layer norms at most $C(d)\delta^{-2}\leq R_{\mathcal K}$. The aggregate norm budget is bounded by
\[
C(d)G^{3d/2+1+\mu}\delta^{-(q_d+1)}
\lesssim\epsilon^{-r_d/(1-\mu)},\qquad
r_d=(q_d+1)(1-\mu)+(3d/2+1+\mu)(q_d+2).
\]
Taking $C_1$ sufficiently large gives $r_d<C_1d q_d$, so this budget is at most $\mathcal K$ for sufficiently large $\mathcal K$. Thus the constructed approximant belongs to the ERM candidate class and has $H^1$ error at most $C\epsilon$. Squaring this bound gives
\begin{equation}\label{eq:drm_op_app}
\mathcal E^{\mathrm{DRM,}op}_{app}
\lesssim \mathcal K^{-\frac{2(1-\mu)}{C_1d\lceil\log_2(d+1)\rceil}}.
\end{equation}

For the statistical term, the layerwise bound in \eqref{eq:drm_op_candidate} holds for every candidate sub-network. The identities $\tanh'=1-\tanh^2$ and $xy=[(x+y)^2-(x-y)^2]/4$ represent its first derivatives by a network with activations in $\{\tanh,\mathrm{ReLU},\mathrm{ReLU}^2\}$. Carrying signed intermediate values through positive and negative ReLU channels gives an extension of depth at most $5L+3$, layer norms at most $C(d)R_{\mathcal K}^2$, and hidden-output magnitudes at most $C(d)R_{\mathcal K}^{2L+2}$. These bounds are uniform over the candidate class and do not depend on width.

At each extended layer, a preactivation is a linear combination of the previous coordinates and a constant, with bounded coefficient $\ell^1$ norm. Applying the contraction inequality on its bounded range, where $\mathrm{ReLU}^2$ has a controlled Lipschitz constant, and iterating over the layers bounds the symmetric Rademacher complexities of individual sub-network outputs and first derivatives by $C(d)R_{\mathcal K}^{Cq_d^2}/\sqrt N$. The outer-coefficient constraint multiplies these bounds by at most $S_{\mathcal K}$. Candidate outputs and first derivatives have envelope $C(d)S_{\mathcal K}R_{\mathcal K}^{L}$; square-function contraction and polarization therefore also control all quadratic energy terms, including cross terms between sub-networks. Symmetrization, with the bounded coefficients and data of \S\ref{sec:drm_setup}, yields
\begin{equation}\label{eq:drm_op_uniform_stat}
\begin{aligned}
\mathbb E\sup_{u\in\mathcal F_{NN}}
\big|\mathcal L_{\mathrm{DRM}}(u)-\widehat{\mathcal L}_{\mathrm{DRM}}(u)\big|
&\lesssim \frac{S_{\mathcal K}^2R_{\mathcal K}^{Cq_d^2}}{\sqrt N} \lesssim \frac{\mathcal K^{C_0q_d}}{\sqrt N},
\end{aligned}
\end{equation}
where $C_0$ is an absolute constant and the implicit constants are independent of width. Since $q_d^2\leq4d$, a sufficiently large absolute $C_2$ ensures
\[
\frac{C_2d}{2\log_2(d+1)}\geq C_0q_d+a_{\mathrm{DRM}}.
\]
The stated sample-size choice consequently gives
\begin{equation}\label{eq:drm_op_sta_raw}
\mathcal E^{\mathrm{DRM,}op}_{sta}
\lesssim \mathcal K^{-\frac{1-\mu}{C_1d\lceil\log_2(d+1)\rceil}}.
\end{equation}
The decomposition and $\mathcal E_{opt}=0$ therefore give
\[
\mathbb E\|\widehat u_{\boldsymbol\theta}-u_0\|_{H^1(\Omega)}^2
\lesssim \mathcal K^{-2a_{\mathrm{DRM}}}+\mathcal K^{-a_{\mathrm{DRM}}}
\lesssim \mathcal K^{-a_{\mathrm{DRM}}}.
\]
Substitution of the sample-size choice gives the stated rate.
\end{proof}

\begin{remark}[Why the over-parameterized DRM rate is slow]
The norm-controlled Sobolev approximation and the global statistical bound lead to a slower rate than the localized $\mathrm{ReLU}^3$ result in Theorem~\ref{thm:drm_fast}. These results use different activations and boundary-value settings; their rates should therefore be read with their respective assumptions.
\end{remark}

\begin{remark}[Scope of the over-parameterized release for DRM]
The regression bound \eqref{eq:reg_op_sta} removes width dependence but retains a depth factor. The DRM bound uses the additional layerwise and coefficient budgets in \eqref{eq:drm_op_candidate} to control network derivatives uniformly. In the PNN analysis of the next chapter, the outer $\ell^1$ constraint removes dependence on the number of parallel sub-networks; the width and depth of each sub-network still enter the constants.
\end{remark}

\subsection{Numerical illustration: network size and weight norms}\label{sec:overparam_numerical}

To compare network size with the norm quantities used above, consider the Robin problem on $\Omega=(0,1)^2$,
\[
-\Delta u+\pi^2u=2\pi^2u_0\quad\text{in }\Omega,
\qquad \partial_nu+u=g\quad\text{on }\partial\Omega,
\]
with $u_0(\boldsymbol{x})=\cos(\pi x_1)+\cos(\pi x_2)$ and $g=\partial_nu_0+u_0$. We use two tanh hidden layers of common width $W$, giving depth $L=3$ and $\mathfrak D=W^2+5W+1$ trainable parameters. The interior and boundary energies are discretized by midpoint quadrature with $N=128^2$ interior points and $N_b=512$ boundary points. Each network is trained in double precision for $10\,000$ full-batch Adam steps, without norm constraints, using a cosine learning-rate schedule from $10^{-3}$ to $10^{-5}$.

For $W\in\{16,32,64,96,128,192,256,512\}$, we repeat training with ten initialization seeds. We evaluate $\|u_{\boldsymbol\theta}-u_0\|_{H^1}/\|u_0\|_{H^1}$ by Gauss quadrature and compute $\kappa(\boldsymbol\theta)$ as in Definition~\ref{sec:overparam_setup}. We also estimate the $C^2$ norm by the maximum absolute partial derivative of order at most two on a $257\times257$ grid.

\begin{figure}[htbp]
\centering
\includegraphics[width=\linewidth]{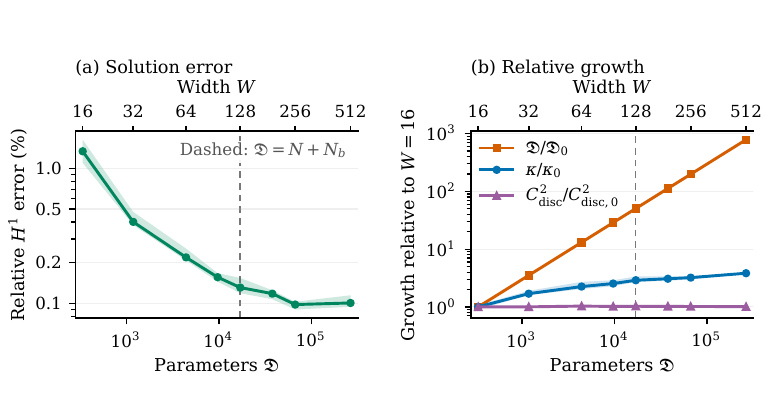}
\caption{Robin problem: relative $H^1$ error and growth of parameter count, weight-norm product, and the grid estimate of the $C^2$ norm. Growth factors are relative to the corresponding medians at $W=16$. Lines and shaded bands show medians and interquartile ranges over ten initializations. The dashed line marks $\mathfrak D=N+N_b$, a parameter-count reference.}
\label{fig:robin_width_norm_error}
\end{figure}

Figure~\ref{fig:robin_width_norm_error} shows that, as $W$ increases from 16 to 512, the parameter count grows by a factor of $785.5$, whereas the median $\kappa$ grows by a factor of $3.83$ and the median grid estimate of the $C^2$ norm increases by about $1.1\%$. The median relative $H^1$ error decreases overall from $1.3389\%$ to $0.1005\%$. Thus, in this example, improved accuracy accompanies very different growth in parameter count and network norms. These are measurements of trained networks; the statistical bounds above instead require prescribed budgets for an entire candidate class.

%% file: sections/complete.tex
\section{Complete Error Analysis via PNN and PGD}\label{sec:complete_setting}

Chapters~\ref{sec:classical_setting} and~\ref{sec:overparam_setting} analyzed two-error trade-offs under $\mathcal{E}_{opt}=0$, treating the ERM as a black box. This chapter surveys analyses that incorporate optimization into the error decomposition, bounding the population risk of the algorithm output. We focus on an approximation-based decomposition (\S\ref{sec:new_decomp}), parallel neural networks (PNNs, Definition~\ref{def:pnn}), and projected gradient descent (PGD, \S\ref{sec:complete_opt}). Work by Kohler and collaborators studies such complete error bounds for over-parameterized regression \cite{kohler2021over,kohler2022analysis,kohler2023opt,kohler2025rate,drews2024universal,drews2023analysis}; we present the ReLU PGD analysis of \cite{kohler2023opt} in detail. For smoothness $\beta\geq1/2$, \cite{kohler2025rate} also obtains a nearly minimax-optimal rate for logistic parallel networks trained by unprojected gradient descent under different initialization and parameter choices (see \S\ref{sec:reg_complete}). The DRM analysis of \cite{jiao2024drm} uses a related comparison-network argument with derivative estimates for the energy functional.

\subsection{Setup: PNN, PGD, and the approximation-based decomposition}\label{sec:new_decomp}

In the two-error analyses of \S\ref{sec:classical_setting} and \S\ref{sec:overparam_setting}, the ERM $\widehat\phi$ serves as the comparison point. For optimization analysis, a constructed approximating network $\bar\phi$ provides explicit coefficients and parameter bounds through Theorems~\ref{thm:A5r} and~\ref{thm:A5}. The analysis compares the algorithm's empirical risk with that of this network. The following subsections describe the architecture, algorithm, and error decomposition.

\subsubsection{The parallel neural network class} We consider the parallel architectures analyzed in \cite{kohler2023opt,jiao2024drm}. Following the notation of \cite{jiao2024drm}, we distinguish the outer coefficient budget $M$ from the sub-network parameter bound $B_{\boldsymbol{\theta}}$.

\begin{definition}[Parallel neural network class; \citealp{jiao2024drm}]\label{def:pnn}
Given parameters $\mathfrak{m}\in\mathbb{N}^+$, $M>0$ and $W, L \in \mathbb{N}^+$, $B_{\boldsymbol{\theta}} > 0$, the PNN class is
\begin{align}
\mathcal{PNN}& (\mathfrak{m}, M, \{W, L, B_{\boldsymbol{\theta}}\}) := \\\notag
& \bigg\{\Phi(\boldsymbol{x}) = \sum_{k=1}^{\mathfrak{m}} c_k\, \phi^k(\boldsymbol{x}) \;:\; \phi^k \in \mathcal{NN}_\rho(W, L, B_{\boldsymbol{\theta}}),\; \sum_{k=1}^{\mathfrak{m}} |c_k| \le M\bigg\},
\end{align}
where $\mathcal{NN}_\rho(\cdot,\cdot,\cdot)$ is the bounded-parameter class of Definition~\ref{def:nnclass}. All sub-networks have the same depth $L$; narrower networks are embedded in the common width-$W$ architecture by padding zeros. We split the total parameter vector as $\boldsymbol{\theta}_{\mathrm{total}} = (\boldsymbol{\theta}_{\mathrm{in}}, \boldsymbol{\theta}_{\mathrm{out}})$, where $\boldsymbol{\theta}_{\mathrm{in}} = (\boldsymbol{\theta}_1, \ldots, \boldsymbol{\theta}_{\mathfrak{m}})$ collects the sub-network weights and $\boldsymbol{\theta}_{\mathrm{out}} = (c_1, \ldots, c_{\mathfrak{m}})$ collects the linear-combination coefficients.
\end{definition}

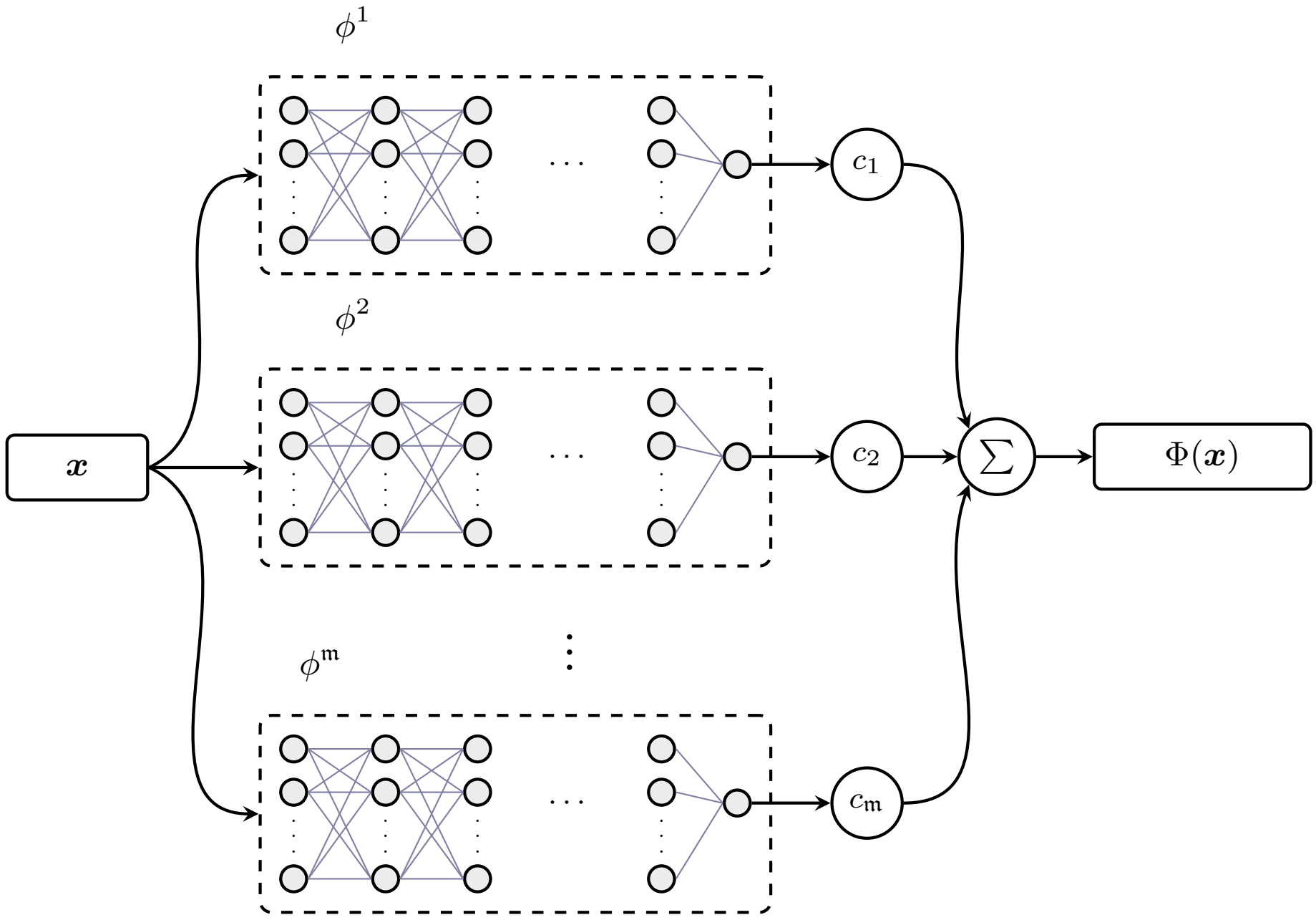
\begin{figure}[!t]
\centering
\begin{tikzpicture}[
  >=stealth,
  io/.style={draw=black, thick, rounded corners=2pt, minimum width=1.3cm, minimum height=0.6cm, align=center, text=black, font=\small},
  neuron/.style={draw=black, thick, circle, fill=gray!15, minimum size=0.24cm, inner sep=0pt, text=black, font=\tiny},
  subbox/.style={draw=black, thick, dashed, rounded corners=3pt, inner sep=5pt},
  coef/.style={draw=black, thick, circle, minimum size=0.65cm, inner sep=0pt, font=\footnotesize, text=black},
  sumop/.style={draw=black, thick, circle, minimum size=0.7cm, inner sep=0pt, font=\small, text=black},
  conn/.style={draw=blue!30!black!50, thin},
  arr/.style={->, thick, black},
  font=\small,
]
\node[io] (input) at (0, -0.1) {$\boldsymbol{x}$};

\foreach \i/\dy in {1/0.5, 2/0.1, 3/-0.7} {
  \node[neuron] (s1in\i)  at (2.0, 2.7+\dy) {};
  \node[neuron] (s1h1\i)  at (2.85, 2.7+\dy) {};
  \node[neuron] (s1h2\i)  at (3.7, 2.7+\dy) {};
  \node[neuron] (s1hL\i)  at (5.4, 2.7+\dy) {};
}
\node[text=black, font=\tiny] at (2.0, 2.50) {$\vdots$};
\node[text=black, font=\tiny] at (2.85, 2.50) {$\vdots$};
\node[text=black, font=\tiny] at (3.7, 2.50) {$\vdots$};
\node[text=black, font=\tiny] at (5.4, 2.50) {$\vdots$};
\node[text=black, font=\footnotesize] at (4.55, 2.7) {$\cdots$};
\node[neuron] (s1out) at (6.1, 2.7) {};
\foreach \i in {1,2,3} \foreach \j in {1,2,3}
  \draw[conn] (s1in\i.east) -- (s1h1\j.west);
\foreach \i in {1,2,3} \foreach \j in {1,2,3}
  \draw[conn] (s1h1\i.east) -- (s1h2\j.west);
\foreach \i in {1,2,3} \draw[conn] (s1hL\i.east) -- (s1out.west);
\node[subbox, fit=(s1in1)(s1in3)(s1out)] (sub1) {};
\node[text=black, font=\small, anchor=south, xshift=-15mm, yshift=5pt] at (sub1.north) {$\phi^1$};
\node[coef] (d1) at (7.3, 2.7) {$c_1$};

\foreach \i/\dy in {1/0.5, 2/0.1, 3/-0.7} {
  \node[neuron] (s2in\i)  at (2.0, 0+\dy) {};
  \node[neuron] (s2h1\i)  at (2.85, 0+\dy) {};
  \node[neuron] (s2h2\i)  at (3.7, 0+\dy) {};
  \node[neuron] (s2hL\i)  at (5.4, 0+\dy) {};
}
\node[text=black, font=\tiny] at (2.0, -0.2) {$\vdots$};
\node[text=black, font=\tiny] at (2.85, -0.2) {$\vdots$};
\node[text=black, font=\tiny] at (3.7, -0.2) {$\vdots$};
\node[text=black, font=\tiny] at (5.4, -0.2) {$\vdots$};
\node[text=black, font=\footnotesize] at (4.55, 0) {$\cdots$};
\node[neuron] (s2out) at (6.1, 0) {};
\foreach \i in {1,2,3} \foreach \j in {1,2,3}
  \draw[conn] (s2in\i.east) -- (s2h1\j.west);
\foreach \i in {1,2,3} \foreach \j in {1,2,3}
  \draw[conn] (s2h1\i.east) -- (s2h2\j.west);
\foreach \i in {1,2,3} \draw[conn] (s2hL\i.east) -- (s2out.west);
\node[subbox, fit=(s2in1)(s2in3)(s2out)] (sub2) {};
\node[text=black, font=\small, anchor=south, xshift=-15mm, yshift=5pt] at (sub2.north) {$\phi^2$};
\node[coef] (d2) at (7.3, 0) {$c_2$};

\node[text=black, font=\Large] at (4.55, -1.7) {$\vdots$};

\foreach \i/\dy in {1/0.5, 2/0.1, 3/-0.7} {
  \node[neuron] (sqin\i)  at (2.0, -3.2+\dy) {};
  \node[neuron] (sqh1\i)  at (2.85, -3.2+\dy) {};
  \node[neuron] (sqh2\i)  at (3.7, -3.2+\dy) {};
  \node[neuron] (sqhL\i)  at (5.4, -3.2+\dy) {};
}
\node[text=black, font=\tiny] at (2.0, -3.4) {$\vdots$};
\node[text=black, font=\tiny] at (2.85, -3.4) {$\vdots$};
\node[text=black, font=\tiny] at (3.7, -3.4) {$\vdots$};
\node[text=black, font=\tiny] at (5.4, -3.4) {$\vdots$};
\node[text=black, font=\footnotesize] at (4.55, -3.2) {$\cdots$};
\node[neuron] (sqout) at (6.1, -3.2) {};
\foreach \i in {1,2,3} \foreach \j in {1,2,3}
  \draw[conn] (sqin\i.east) -- (sqh1\j.west);
\foreach \i in {1,2,3} \foreach \j in {1,2,3}
  \draw[conn] (sqh1\i.east) -- (sqh2\j.west);
\foreach \i in {1,2,3} \draw[conn] (sqhL\i.east) -- (sqout.west);
\node[subbox, fit=(sqin1)(sqin3)(sqout)] (subq) {};
\node[text=black, font=\small, anchor=south, xshift=-18mm, yshift=5pt] at (subq.north) {$\phi^{\mathfrak{m}}$};
\node[coef] (dq) at (7.3, -3.2) {$c_{\mathfrak{m}}$};

\draw[arr] (input.east) to[out=25, in=180] (sub1.west);
\draw[arr] (input.east) -- (sub2.west);
\draw[arr] (input.east) to[out=-25, in=180] (subq.west);

\draw[arr] (s1out.east) -- (d1.west);
\draw[arr] (s2out.east) -- (d2.west);
\draw[arr] (sqout.east) -- (dq.west);

\node[sumop] (sum) at (8.5, 0) {$\sum$};
\draw[arr] (d1.east) to[out=0, in=110] (sum.north west);
\draw[arr] (d2.east) -- (sum.west);
\draw[arr] (dq.east) to[out=0, in=-110] (sum.south west);

\node[io, minimum width=2.0cm] (out) at (10.4, 0) {$\Phi(\boldsymbol{x})$};
\draw[arr] (sum.east) -- (out.west);
\end{tikzpicture}
\caption{Structure of the parallel neural network class $\mathcal{PNN}(\mathfrak{m}, M, \{W, L, B_{\boldsymbol{\theta}}\})$ of Definition~\ref{def:pnn}. The $\mathfrak{m}$ sub-networks share a width-$W$, depth-$L$ architecture and receive the same input $\boldsymbol{x}$. Their weighted sum $\Phi(\boldsymbol{x}) = \sum_{k=1}^{\mathfrak{m}} c_k\phi^k(\boldsymbol{x})$ forms the output, with $\sum_k|c_k|\le M$ and each sub-network's parameters bounded by $B_{\boldsymbol{\theta}}$.}
\label{fig:pnn_structure}
\end{figure}

The outer $\ell^1$ constraint yields a Rademacher complexity bound proportional to $M$ and independent of the sub-network count $\mathfrak m$ (\S\ref{sec:complete_stat}).

\subsubsection{The PGD algorithm}\label{sec:pgd_algorithm} Given a $\mathcal{PNN}$ class with parameters $(\mathfrak{m}, M, \{W, L, B_{\boldsymbol{\theta}}\})$ and target accuracy parameters $\eta, \zeta > 0$, $T \in \mathbb{N}^+$, $\lambda > 0$, the PGD algorithm of \cite{kohler2023opt,jiao2024drm} proceeds in three stages.

(i) Initialization. The outer coefficients are initialized at zero and the inner sub-network weights are initialized i.i.d.\ uniform:
\[
(\boldsymbol{\theta}_{\mathrm{out}})^{[0]} := \boldsymbol{0}, \qquad
\big(\boldsymbol{\theta}_{\mathrm{in}}\big)^{[0]} \sim \text{i.i.d.\ } \mathcal{U}[-B_{\boldsymbol{\theta}}, B_{\boldsymbol{\theta}}].
\]

(ii) Constraint set. Two projection radii $\eta, \zeta > 0$ define the constraint set $A_\eta \times B_\zeta$ where
\[
A_\eta := \big\{\boldsymbol{\theta}_{\mathrm{in}} : \|\boldsymbol{\theta}_{\mathrm{in}} - (\boldsymbol{\theta}_{\mathrm{in}})^{[0]}\|_2 \le \eta\big\},
\qquad
B_\zeta := \Big\{\boldsymbol{\theta}_{\mathrm{out}} : \sum_{k=1}^{\mathfrak{m}} |c_k| \le \zeta\Big\}.
\]
The set $A_\eta$ constrains the inner weights to a Euclidean ball around initialization, while $B_\zeta$ is the $\ell^1$-ball on outer coefficients. Thus $B_{\boldsymbol\theta}$ controls initialization, $\eta$ controls subsequent inner-weight displacement, and $\zeta$ bounds the outer coefficients; below we set $\zeta=\bar M$. The projection onto $A_\eta$ has a closed form, and the projection onto $B_\zeta$ is computable in linear time \cite{duchi2008efficient}.

(iii) Iterative update. For $t = 0, 1, \ldots, T-1$,
\begin{equation}\label{eq:pgd_update}
\big(\boldsymbol{\theta}_{\mathrm{in}}, \boldsymbol{\theta}_{\mathrm{out}}\big)^{[t+1]} = \mathrm{Proj}_{A_\eta \times B_\zeta} \big\{(\boldsymbol{\theta}_{\mathrm{in}}, \boldsymbol{\theta}_{\mathrm{out}})^{[t]} - \lambda \cdot \nabla_{(\boldsymbol{\theta}_{\mathrm{in}}, \boldsymbol{\theta}_{\mathrm{out}})}\widehat{\mathcal{L}}\big\}.
\end{equation}
Here $\mathrm{Proj}$ is the Euclidean projection onto $A_\eta\times B_\zeta$. The generic iteration bound below uses the iterate of smallest empirical risk among $t=0,\ldots,T-1$. The DRM result takes the last iterate, with the additional step-size condition ensuring descent. In the regression result, each ReLU sub-network is clipped before the outer summation, and initialization uses layer-dependent ranges; these choices are specified in \S\ref{sec:reg_complete}. At ReLU breakpoints, the inner update uses the chain-rule convention $\rho'(0)=0$ of \cite[Section~2.3]{kohler2023opt}, also for the ReLU representation of clipping. All iterates have inner parameter bound $B_{\boldsymbol\theta}+\eta$.

\begin{remark}
For fixed inner weights, squared regression loss and the coercive DRM energy are convex in the outer coefficients. The outer update can therefore be compared with fixed reference coefficients by a convex PGD argument. The inner projection controls the change of the corresponding reference network, while smoothness of the tanh empirical energy supplies the descent condition used for the last-iterate DRM bound.
\end{remark}

\subsubsection{The approximation-based decomposition} Choose a network $\bar\phi$ supplied by an approximation theorem such that
\begin{equation}\label{eq:bar_phi}
\bar\phi\in\mathcal{PNN}^{\prime},\qquad \|\bar\phi-\phi^*\|_{\mathcal H}\le\epsilon.
\end{equation}
Here $\mathcal{PNN}^{\prime}$ is the comparison class from Theorem~\ref{thm:A5r} or~\ref{thm:A5}, with parameters $(\bar{\mathfrak m},\bar M,\{\bar W,\bar L,B_{\bar{\boldsymbol\theta}}\})$ chosen for accuracy $\epsilon$. The norm $\mathcal H$ is $L^2(P_X)$ for regression and $H^1(\Omega)$ for DRM. The network need not minimize the approximation error. The larger class used for the uniform statistical bound contains both $\bar\phi$ and the PGD iterates; the initialization argument below constructs a separate feasible comparison point for the random projection set.

\begin{lemma}[Approximation-based decomposition; \citealp{jiao2024drm}]\label{lem:D1}
Let $\phi_{\boldsymbol{\theta}_{\mathcal{A}}}$ be the output of an algorithm $\mathcal{A}$, and suppose that both $\phi_{\boldsymbol{\theta}_{\mathcal{A}}}$ and the comparison network $\bar\phi$ in \eqref{eq:bar_phi} belong to the $\mathcal{PNN}$ class over which the uniform deviation is taken. Then
\begin{equation}\label{eq:D1}
\mathcal{L}(\phi_{\boldsymbol{\theta}_{\mathcal{A}}}) - \mathcal{L}(\phi^*) \;\le\; \underbrace{\big[\mathcal{L}(\bar{\phi}) - \mathcal{L}(\phi^*)\big]}_{\mathcal{E}^{\mathrm c}_{app}}
+ \underbrace{\big[\widehat{\mathcal{L}}(\phi_{\boldsymbol{\theta}_{\mathcal{A}}}) - \widehat{\mathcal{L}}(\bar{\phi})\big]}_{\mathcal{E}_{opt,-}}
+ \underbrace{2 \sup_{\Phi \in \mathcal{PNN}}\big|\mathcal{L}(\Phi) - \widehat{\mathcal{L}}(\Phi)\big|}_{\mathcal{E}^{\mathrm c}_{sta}}.
\end{equation}
Here the superscript $\mathrm c$ distinguishes these terms from those in \S\ref{sec:decomp}: the approximation gap uses the constructed network $\bar\phi$, and the statistical term includes the factor of two. Their DRM instances are denoted by $\mathcal E^{\mathrm{DRM,c}}_{app}$ and $\mathcal E^{\mathrm{DRM,c}}_{sta}$.
Moreover, $\mathcal{E}_{opt,-} \le \mathcal{E}_{opt}$, where $\mathcal{E}_{opt} := \widehat{\mathcal{L}}(\phi_{\boldsymbol{\theta}_{\mathcal{A}}}) - \widehat{\mathcal{L}}(\widehat{\phi})$ is the standard ERM-anchored optimization error of \S\ref{sec:decomp}.
\end{lemma}

\begin{proof}
The decomposition is obtained by inserting $\bar{\phi}$ into the chain
\begin{align*}
&\mathcal{L}(\phi_{\boldsymbol{\theta}_{\mathcal{A}}}) - \mathcal{L}(\phi^*) \\
& ~~~~ = \big[\mathcal{L}(\phi_{\boldsymbol{\theta}_{\mathcal{A}}}) - \widehat{\mathcal{L}}(\phi_{\boldsymbol{\theta}_{\mathcal{A}}})\big] + \big[\widehat{\mathcal{L}}(\phi_{\boldsymbol{\theta}_{\mathcal{A}}}) - \widehat{\mathcal{L}}(\bar{\phi})\big] + \big[\widehat{\mathcal{L}}(\bar{\phi}) - \mathcal{L}(\bar{\phi})\big] + \big[\mathcal{L}(\bar{\phi}) - \mathcal{L}(\phi^*)\big].
\end{align*}
Bounding the first and third bracketed terms by the uniform supremum yields the statistical contribution $\mathcal{E}^{\mathrm c}_{sta}$; the remaining two are $\mathcal{E}_{opt,-}$ and $\mathcal{E}^{\mathrm c}_{app}$ respectively. For the inequality $\mathcal{E}_{opt,-} \le \mathcal{E}_{opt}$, write
\[
\mathcal{E}_{opt,-} = \big[\widehat{\mathcal{L}}(\phi_{\boldsymbol{\theta}_{\mathcal{A}}}) - \widehat{\mathcal{L}}(\widehat{\phi})\big] + \big[\widehat{\mathcal{L}}(\widehat{\phi}) - \widehat{\mathcal{L}}(\bar{\phi})\big] = \mathcal{E}_{opt} + \big[\widehat{\mathcal{L}}(\widehat{\phi}) - \widehat{\mathcal{L}}(\bar{\phi})\big].
\]
The second bracket is non-positive because $\widehat{\phi}$ is the empirical risk minimizer over the $\mathcal{PNN}$ class containing $\bar{\phi}$. Hence $\mathcal{E}_{opt,-} \le \mathcal{E}_{opt}$.
\end{proof}

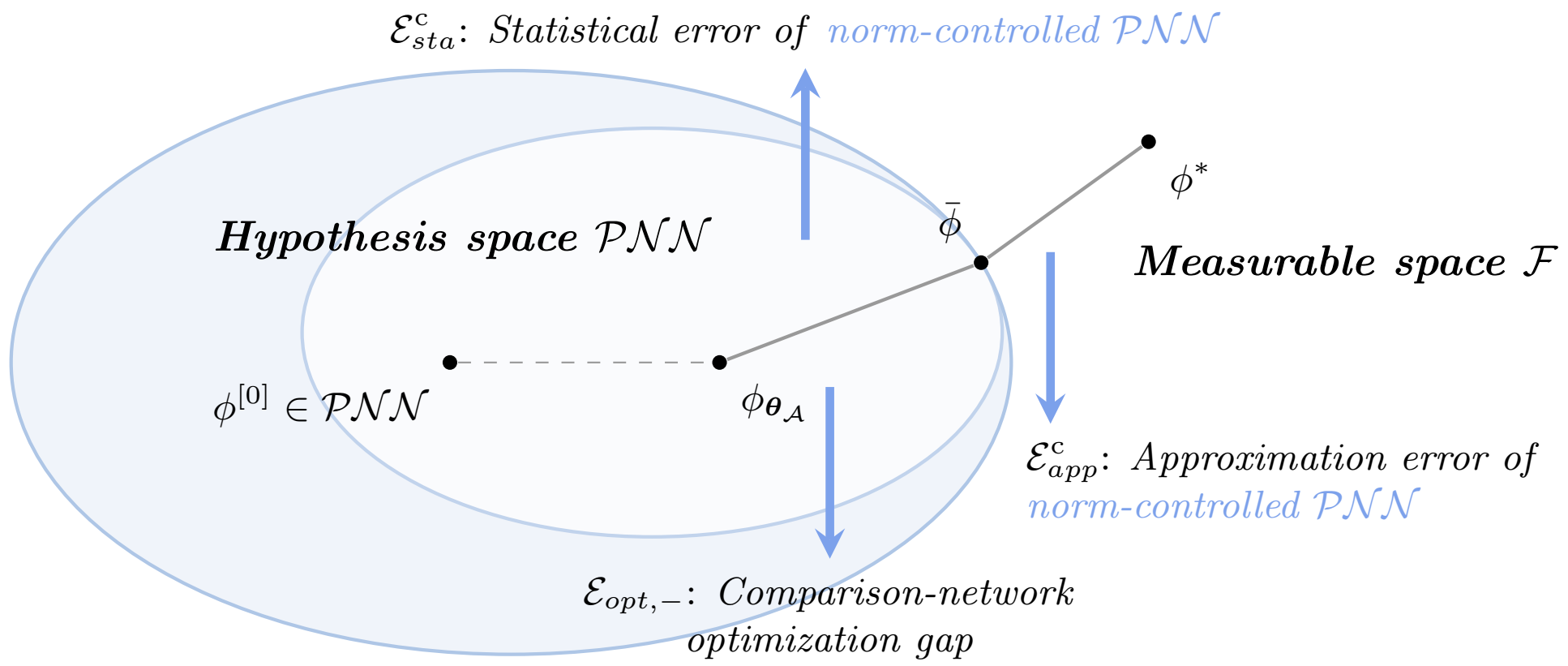
\begin{figure}[!t]
    \centering
    \begin{tikzpicture}[>=stealth, font=\small]

        \definecolor{webbg}{HTML}{F0F4FA}
        \definecolor{webborder}{HTML}{AEC6E6}
        \definecolor{webarrow}{HTML}{7CA1EB}

        \draw[fill=webbg, draw=webborder, thick] (0,0) ellipse (4.08cm and 2.38cm);

        \draw[fill=white, draw=webborder, thick, opacity=0.7] (1.150, 0.244) ellipse (2.856cm and 1.666cm);

        \node[font=\bfseries\itshape, text=black] at (-0.4, 1.0) {Hypothesis space $\mathcal{PNN}$};

        \coordinate (phibar) at (3.834, 0.814);

        \node (n_phibar) at (phibar) [circle, fill=black, inner sep=1.2pt] {};
        \node[above left=1pt] at (n_phibar) {$\bar{\phi}$};

        \node (n_phistar) at (5.2, 1.8) [circle, fill=black, inner sep=1.2pt, label=below right:$\phi^{*}$] {};

        \node (n_phiA) at (1.7, 0.0) [circle, fill=black, inner sep=1.2pt, label=below right:$\phi_{\boldsymbol{\theta}_{\mathcal{A}}}$] {};

        \node (n_phiinit) at (-0.5, 0.0) [circle, fill=black, inner sep=1.2pt, label=below left:{$\phi^{[0]} \in \mathcal{PNN}$}] {};

        \draw[dashed, gray!80] (n_phiinit) -- (n_phiA);
        \draw[thick, gray!80] (n_phiA) -- (n_phibar);
        \draw[thick, gray!80] (n_phibar) -- (n_phistar);

        \node[font=\bfseries\itshape] at (6.8, 0.8) {Measurable space $\mathcal{F}$};

        \draw[->, line width=2pt, webarrow] (2.6, -0.2) -- (2.6, -1.6)
            node[below, text=black, align=center] {$\mathcal{E}_{opt,-}$: \textit{Comparison-network}\\\textit{optimization gap}};

        \draw[->, line width=2pt, webarrow] (2.4, 1.0) -- (2.4, 2.4)
            node[above, text=black] {$\mathcal{E}^{\mathrm c}_{sta}$: \textit{Statistical error of} \textcolor{webarrow}{\textit{norm-controlled $\mathcal{PNN}$}}};

        \draw[->, line width=2pt, webarrow] (4.4, 0.9) -- (4.4, -0.5)
            node[below right, xshift=-10pt, yshift=0pt, text=black, align=left] {$\mathcal{E}^{\mathrm c}_{app}$: \textit{Approximation error of} \\ \textcolor{webarrow}{\textit{norm-controlled $\mathcal{PNN}$}}};

    \end{tikzpicture}
    \caption{Schematic diagram of the approximation-based decomposition of \S\ref{sec:new_decomp}. Compared with Figure~\ref{fig:error_decom1}, the analytical anchor on the right of $\phi_{\boldsymbol{\theta}_{\mathcal{A}}}$ is the constructed approximating network $\bar{\phi}$ (constructible from the target $\phi^*$ via the PNN approximation theorem of \S\ref{sec:complete_approx}) rather than the unobservable ERM. The optimization error becomes $\mathcal{E}_{opt,-} = \widehat{\mathcal{L}}(\phi_{\boldsymbol{\theta}_{\mathcal{A}}}) - \widehat{\mathcal{L}}(\bar{\phi})$, replacing the ERM-anchored $\mathcal{E}_{opt}$ of \S\ref{sec:decomp}. The statistical and approximation contributions are bounded by the PNN-specific tools developed in \S\ref{sec:complete_stat} and \S\ref{sec:complete_approx}; the inner ellipse depicts the sub-class $\mathcal{PNN}^\prime$ of \eqref{eq:bar_phi} that hosts $\bar{\phi}$.}
    \label{fig:error_decom2}
\end{figure}

Figure~\ref{fig:error_decom2} gives the geometric counterpart of \eqref{eq:D1} and should be read alongside Figure~\ref{fig:error_decom1} of \S\ref{sec:decomp}. The remaining sections develop the three ingredients (\S\ref{sec:complete_approx}--\S\ref{sec:complete_opt}) needed to instantiate \eqref{eq:D1}, and \S\ref{sec:complete_rates} performs the final trade-off algebra.

\subsection{PNN approximation}\label{sec:complete_approx}

We present the ReLU approximation result used in the regression analysis of \cite{kohler2023opt} and the tanh Sobolev approximation result used for DRM in \cite{jiao2024drm}. Both specify an outer coefficient budget and an inner parameter bound separately, as required by the PNN analysis below.

\begin{theorem}[PNN ReLU approximation for H\"older targets; \citealp{kohler2023opt}]\label{thm:A5r}
Let $\beta=s+r>0$, with $s=\lceil\beta\rceil-1$ and $r\in(0,1]$, and let $f_0\in\mathcal H^\beta([0,1]^d,B_{\mathrm{REG}})$. For any partition resolution $G \in \mathbb{N}$ sufficiently large, there exists a ReLU-activated parallel network $\bar{f} \in \mathcal{PNN}(\bar{\mathfrak{m}}, \bar{M}, \{\bar{W}, \bar{L}, B_{\bar{\boldsymbol{\theta}}}\})$ with
\begin{align*}
\bar{\mathfrak{m}} &= \binom{d+s}{d}\,(G+1)^d, &\qquad
\bar{M} &\asymp (G+1)^d, \\
\bar{W} &= 18(s+d)^{d+1}, &\qquad
\bar{L} &= O(\log G), \\
B_{\bar{\boldsymbol{\theta}}} &= O(G),
\end{align*}
such that
\[
\|\bar{f} - f_0\|_{L^\infty([0,1]^d)} \;\le\; C_{\mathrm{ReLU}}\, G^{-\beta},
\]
where $C_{\mathrm{ReLU}}$ and the constants in the parameter bounds depend on $d,\beta,B_{\mathrm{REG}}$. In particular, the bound on the outer coefficients uses the bounds on the Taylor coefficients supplied by the H\"older ball.
\end{theorem}

The construction \cite[Lemma~6]{kohler2023opt} proceeds via a partition-of-unity Taylor expansion on a grid of $(G+1)^d$ cells. On each cell, a Taylor polynomial of degree $s$ is multiplied by a ReLU-realized hat function; the resulting local approximants are assembled into the PNN using one sub-network for each cell--monomial pair, giving $\bar{\mathfrak{m}} \asymp \bar{M} \asymp (G+1)^d$ for fixed $d$ and $s$. The sub-network width is the fixed polynomial in $d$ and $s$ required to implement multiplication via the $x^2$ sawtooth iteration (\S\ref{sec:classical_approx}); the depth grows only logarithmically with $G$; the weight bound $B_{\bar{\boldsymbol{\theta}}} = O(G)$ is dominated by the input-layer weights encoding the cell centres, all deeper weights being bounded by a fixed constant.

\begin{theorem}[PNN tanh approximation; \citealp{jiao2024drm}, Theorem~7]\label{thm:A5}
Let $s,d\in\mathbb N$, $s\ge2$, $1\le p\le\infty$, and $\|f\|_{W^{s,p}([0,1]^d)}\le B$. Fix $0<\mu<\min\{1,s-1\}$. For every sufficiently small $\epsilon>0$, there exists a tanh parallel network $\Phi\in\mathcal{PNN}(\mathfrak m,M,\{W,L,B_{\boldsymbol\theta}\})$ with
\begin{align*}
\mathfrak{m} &= C_1\, \epsilon^{-d/(s - 1 - \mu)}, &\qquad
M &= C_2\, \epsilon^{-d(p+1)/((s-1-\mu)p)}, \\
W &= 2^{\lceil \log_2(d + s - 1)\rceil + 1}, &\qquad
L &= \lceil \log_2(d + s - 1)\rceil + 2, \\
B_{\boldsymbol{\theta}} &= C_3\, \epsilon^{-(2d+2s)/(s-1-\mu)}, &&
\end{align*}
such that
\[
\|f - \Phi\|_{W^{1,p}([0,1]^d)} \;\le\; \epsilon.
\] The constants depend on $s,d,p,\mu,B$; for $p=\infty$, $(p+1)/p$ is interpreted as $1$, and integer architecture sizes are rounded up.
\end{theorem}

The two theorems share the same $\mathcal{PNN}$ architecture template but differ in activation, target function space, and parameter scaling. The ReLU construction controls the $L^\infty$ norm, sufficient for the squared $L^2$ regression risk; the tanh construction controls the $W^{1,p}$ Sobolev norm required by the DRM energy. Both supply the approximation building blocks for the complete analyses of \S\ref{sec:complete_rates}.

\subsection{PNN statistical tools}\label{sec:complete_stat}

The statistical analysis of the PNN class exploits the outer $\ell^1$-constraint $\sum_k|c_k|\le M$. Expanding the PNN representation $\sum_k c_k \phi^k$ inside the Rademacher average and applying H\"older's inequality decouples the outer coefficients from the inner sub-networks; the $\ell^1$-constraint then bounds the coefficient sum by $M$, and the $\sup$ over sub-network indices is absorbed into the symmetric Rademacher complexity of a single sub-network. The sub-network count $\mathfrak{m}$ therefore drops out of the bound. The remaining single-sub-network complexity still depends on its width, depth, and parameter bound, and is controlled by the source-specific estimates described below. The instantiation on the two running examples is carried out in the proof sketches of \S\ref{sec:reg_complete} and \S\ref{sec:drm_complete}.

\begin{lemma}[$\mathcal{PNN}$ Rademacher reduction]\label{lem:S7}
Let $\mathcal{PNN}(\mathfrak{m}, M, \{W, L, B_{\boldsymbol{\theta}}\})$ be the parallel network class of Definition~\ref{def:pnn}, for a fixed activation $\rho$ with finite single-sub-network Rademacher complexity. For any sample size $n$,
\[
\mathfrak{R}_n\big(\mathcal{PNN}(\mathfrak{m}, M, \{W, L, B_{\boldsymbol{\theta}}\})\big) \;\le\; M \cdot \tilde{\mathfrak{R}}_n\big(\mathcal{NN}_\rho(W, L, B_{\boldsymbol{\theta}})\big).
\]
\end{lemma}

\begin{proof}
For a fixed sample and signs, write $\mathcal F_{\mathrm{sub}}=\mathcal{NN}_\rho(W,L,B_{\boldsymbol\theta})$. Then
\begin{align*}
\sup_{\Phi\in\mathcal{PNN}}\frac1n\sum_{i=1}^n\sigma_i\Phi(\boldsymbol x_i)
&=\sup_{\substack{\sum_k|c_k|\le M\\ \phi^k\in\mathcal F_{\mathrm{sub}}}}
\sum_k c_k\left(\frac1n\sum_i\sigma_i\phi^k(\boldsymbol x_i)\right)\\
&\le \sup_{\substack{\sum_k|c_k|\le M\\ \phi^k\in\mathcal F_{\mathrm{sub}}}}
\sum_k|c_k|\sup_{\phi\in\mathcal F_{\mathrm{sub}}}
\left|\frac1n\sum_i\sigma_i\phi(\boldsymbol x_i)\right|\\
&\le M\sup_{\phi\in\mathcal F_{\mathrm{sub}}}
\left|\frac1n\sum_i\sigma_i\phi(\boldsymbol x_i)\right|.
\end{align*}
Taking expectations over the signs and sample gives the stated symmetric-complexity bound.
\end{proof}

The single-sub-network complexity is estimated separately for each example. For regression, \cite[Lemma~5]{kohler2023opt} uses the Bartlett--Harvey--Liaw--Mehrabian pseudo-dimension bound \cite{bartlett2019nearly} to obtain a covering estimate for clipped ReLU sub-networks, followed by Hoeffding's inequality and integration of the tail bound. For DRM, \cite{jiao2024drm} separates the energy into gradient-squared, value-squared, source-linear, and boundary-linear terms. The linear-combination argument of Lemma~\ref{lem:S7} applies to both the value and derivative classes; contraction for the squared terms contributes an additional factor of $M$. These estimates remove dependence on $\mathfrak m$ while retaining dependence on the sub-network architecture and parameter bounds.

\subsection{PNN optimization tools}\label{sec:complete_opt}

The PGD algorithm was described in \S\ref{sec:new_decomp}; this subsection provides its convergence analysis. The optimization error targeted by the analysis is
\begin{equation}\label{eq:opt_def}
\mathcal{E}_{opt,-} := \widehat{\mathcal{L}}(\phi_{\boldsymbol{\theta}_{\mathcal{A}}}) - \widehat{\mathcal{L}}(\bar{\phi}).
\end{equation}
Here $\bar\phi$ is the constructed comparison network with parameters $(\bar{\mathfrak m},\bar M,\{\bar W,\bar L,B_{\bar{\boldsymbol\theta}}\})$. Training uses $\mathfrak m\gg\bar{\mathfrak m}$ sub-networks of the same width and depth, initialization bound $B_{\bar{\boldsymbol\theta}}$, and outer radius $\zeta=\bar M$. The projection radius $\eta$ enlarges the inner parameter bound at the iterates to $B_{\bar{\boldsymbol\theta}}+\eta$, which is used in the statistical analysis.

The central structural observation is that the empirical risk $\widehat{\mathcal{L}}$ is convex in the outer coefficients $\boldsymbol{\theta}_{\mathrm{out}}$ when the inner weights $\boldsymbol{\theta}_{\mathrm{in}}$ are held fixed, because the PNN output depends linearly on $\boldsymbol{\theta}_{\mathrm{out}}$ and the loss is convex in its first argument. This partial convexity is what makes the convergence analysis tractable: the inner weights are confined to a ball $A_\eta$ around initialization where their perturbation is bounded, while the outer coefficients enjoy a classical convex PGD guarantee.

\subsubsection{Factorisation and decomposition of the optimization error}\label{sec:factorisation}

Throughout this subsection, $\delta$ denotes the cover radius.

The proof strategy splits $\mathcal{E}_{opt,-}$ into two pieces, each controlled by a separate argument. The transition network $\phi_q^*$, constructed below, bridges the random initialization to the anchor $\bar{\phi}$. Its gap to $\bar{\phi}$, the initialization error, is controlled by a covering argument: with enough random draws, some sub-networks land near each target parameter, and the replica averaging keeps the outer coefficients bounded. Its gap to the PGD output, the iteration error, is controlled by a convex optimization argument on the outer coefficients, with the inner-weight drift controlled by the projection radius. The factorisation $\mathfrak{m} = \bar{\mathfrak{m}}\cdot R\cdot Q_{\mathrm{trial}}$ supplies the parameters for both mechanisms simultaneously.

The total sub-network count is factorised as
\begin{equation}\label{eq:factorisation}
\mathfrak{m} = \bar{\mathfrak{m}} \cdot R \cdot Q_{\mathrm{trial}},
\end{equation}
where $\bar{\mathfrak{m}}$ is the number of sub-networks in the anchor $\bar{\phi}$, fixed by the approximation theorem; $R$ is the number of replicas per target, whose role is to shrink the $\ell_2$-norm of the outer parameter vector by a factor $1/\sqrt{R}$ when each target coefficient $\bar{c}_k$ is evenly distributed across its $R$ replicas as $\bar{c}_k/R$; and $Q_{\mathrm{trial}}$ is the number of independent random initializations drawn per replica, chosen large enough that with high probability every target is covered by at least one initial draw within distance $\delta$. The total sub-network count $\mathfrak{m} = \bar{\mathfrak{m}} R Q_{\mathrm{trial}}$ is the product of these three quantities; over-parameterization enters through $Q_{\mathrm{trial}}$, whose sufficient size is polynomial in $1/\delta$ for fixed sub-network parameter dimension, but can become super-polynomial when that dimension grows.

The analysis proceeds by bounding the two parts of the decomposition
\begin{equation}\label{eq:opt_decomp}
\mathcal{E}_{opt,-}
= \underbrace{\big[\widehat{\mathcal{L}}(\phi_{\boldsymbol{\theta}_{\mathcal{A}}}) - \widehat{\mathcal{L}}(\phi_q^*)\big]}_{\text{iteration error}}
  + \underbrace{\big[\widehat{\mathcal{L}}(\phi_q^*) - \widehat{\mathcal{L}}(\bar{\phi})\big]}_{\text{initialization error}}.
\end{equation}
The initialization and iteration terms are bounded in \S\ref{sec:init_error} and \S\ref{sec:iter_error}, respectively. PGD itself starts with zero outer coefficients; $\phi_q^*$ is a comparison network, not the initial iterate.

\subsubsection{Bounding the initialization error}\label{sec:init_error}

For the uniform initialization used in DRM, let $J_{\max}=\bar W(\bar W+1)\bar L$ bound the number of trainable parameters in each comparison sub-network (with $\bar W\ge d$). Draw $\mathfrak m=\bar{\mathfrak m}RQ_{\mathrm{trial}}$ independent initializations. Assign a disjoint block of $Q_{\mathrm{trial}}$ draws to each target-replica pair $(k,v)$.

\begin{lemma}[Cover initialization; \citealp{jiao2024drm}]\label{lem:O2a}
With the notation above, for $0<\delta\le2B_{\bar{\boldsymbol\theta}}$,
\[
\Pr\big(\mathcal{E}_{\mathrm{cover}}(\delta)\big) \;\ge\; 1 - \bar{\mathfrak{m}} R\Big[1 - \Big(\frac{\delta}{2B_{\bar{\boldsymbol{\theta}}}}\Big)^{J_{\max}}\Big]^{Q_{\mathrm{trial}}},
\]
where $\mathcal{E}_{\mathrm{cover}}(\delta)$ is the event that for every $k \in [\bar{\mathfrak{m}}]$ and $v \in [R]$ there exists an index $s_{k,v}$ in its assigned block such that $\|(\boldsymbol{\theta}_{s_{k,v}})^{[0]} - \bar{\boldsymbol{\theta}}_k\|_\infty \le \delta$.
\end{lemma}

\begin{proof}
A single draw lands within $\ell^\infty$ distance $\delta$ of the target with probability at least $p_\delta=(\delta/(2B_{\bar{\boldsymbol\theta}}))^{J_{\max}}$, including targets on the boundary of the initialization cube. A block misses its target with probability at most $(1-p_\delta)^{Q_{\mathrm{trial}}}$. The union bound over the $\bar{\mathfrak m}R$ disjoint blocks proves the claim and ensures distinct selected indices.
\end{proof}

On the event $\mathcal{E}_{\mathrm{cover}}(\delta)$, fix indices $s_{k,v}$ as in Lemma~\ref{lem:O2a}, with all $s_{k,v}$ distinct for different $(k,v)$ pairs. Define the transition parameter $\boldsymbol{\theta}_{\mathrm{total}}^{\mathfrak{m},*} = (\boldsymbol{\theta}_{\mathrm{in}}^{\mathfrak{m},*}, \boldsymbol{\theta}_{\mathrm{out}}^{\mathfrak{m},*})$ by
\[
\boldsymbol{\theta}_{\mathrm{in}}^{\mathfrak{m},*} := (\boldsymbol{\theta}_{\mathrm{in}}^{\mathfrak{m}})^{[0]}, \qquad
c_i^* :=
\begin{cases}
\bar{c}_k / R & \text{if } i = s_{k,v} \text{ for some } k \in [\bar{\mathfrak{m}}],\, v \in [R], \\
0 & \text{otherwise}.
\end{cases}
\]
Denote by $\phi_q^*$ the PNN parameterized by $\boldsymbol{\theta}_{\mathrm{total}}^{\mathfrak{m},*}$. This construction freezes the inner weights at their random initialization and distributes each target coefficient $\bar{c}_k$ uniformly across its $R$ replicas, giving
\begin{align*}
\|\boldsymbol{\theta}_{\mathrm{out}}^{\mathfrak{m},*}\|_1 &= \sum_{i=1}^{\mathfrak{m}} |c_i^*| = \sum_{k=1}^{\bar{\mathfrak{m}}} |\bar{c}_k| \le \bar{M},\\
\|\boldsymbol{\theta}_{\mathrm{out}}^{\mathfrak{m},*}\|_2 &= \bigg(\sum_{i=1}^{\mathfrak{m}} |c_i^*|^2\bigg)^{\!1/2} = \bigg(\frac{1}{R}\sum_{k=1}^{\bar{\mathfrak{m}}} |\bar{c}_k|^2\bigg)^{\!1/2} \le \frac{\bar{M}}{\sqrt{R}}.
\end{align*}
The first bound ensures $\boldsymbol{\theta}_{\mathrm{total}}^{\mathfrak{m},*} \in A_\eta \times B_\zeta$ (since $\zeta = \bar{M}$); the second reduces the reference coefficient vector's Euclidean norm and hence the corresponding terms in the iteration bound below.

The next lemma bounds the initialization error $\widehat{\mathcal{L}}(\phi_q^*) - \widehat{\mathcal{L}}(\bar{\phi})$. It relies on the following structural assumption on the empirical risk.
\begin{asmp}[Lipschitz dependence on network evaluations]\label{assum:L_lip}
Use $\|v\|_{\mathrm{eval}}=\|v\|_{L^\infty}$ for regression and $\|v\|_{\mathrm{eval}}=\max\{\|v\|_{L^\infty(\bar\Omega)},\|\nabla v\|_{L^\infty(\bar\Omega;\mathbb R^d)}\}$ for DRM. On a bounded class containing the comparison networks and the PGD iterates, assume
$|\widehat{\mathcal{L}}(\phi) - \widehat{\mathcal{L}}(\psi)| \le L_{\widehat{\mathcal{L}}}\,\|\phi - \psi\|_{\mathrm{eval}}$
for any two PNNs $\phi, \psi$. The value of $L_{\widehat{\mathcal{L}}}$ depends on the pointwise loss $\ell$ and must be estimated for each concrete instance.
\end{asmp}
The required sub-network bounds depend on the activation and the evaluation norm. Write $\mathcal B_{\mathrm{sub}}(B)$ and $\mathfrak L(B)$ for bounds satisfying
\begin{equation}\label{eq:subnet_props}
\|\phi_{\boldsymbol\theta}\|_{\mathrm{eval}}\le\mathcal B_{\mathrm{sub}}(B),\qquad \|\phi_{\boldsymbol\theta}-\phi_{\tilde{\boldsymbol\theta}}\|_{\mathrm{eval}}\le\mathfrak L(B)\|\boldsymbol\theta-\tilde{\boldsymbol\theta}\|_\infty,
\end{equation}
for inner parameter vectors bounded by $B$. In particular, set
\begin{equation}\label{eq:frakturL}
\mathfrak L:=\mathfrak L(B_{\bar{\boldsymbol\theta}}),\qquad B\in\{B_{\bar{\boldsymbol\theta}},B_{\bar{\boldsymbol\theta}}+\eta\}.
\end{equation}
For tanh DRM subnetworks, these bounds include derivatives; for the regression architecture, the subnetworks are clipped. The tanh value bound $(\bar W+1)B$ alone does not control either an unclipped deep ReLU network or the DRM gradient term.

\begin{lemma}[Initialization error]\label{lem:O2init}
On the event $\mathcal{E}_{\mathrm{cover}}(\delta)$,
\[
\widehat{\mathcal{L}}(\phi_q^*) - \widehat{\mathcal{L}}(\bar{\phi}) \;\le\; L_{\widehat{\mathcal{L}}}\,\bar{M}\,\mathfrak{L}\,\delta,
\]
where $L_{\widehat{\mathcal{L}}}$ is the Lipschitz constant from Assumption~\ref{assum:L_lip} and $\mathfrak{L}$ is given by \eqref{eq:frakturL}.
\end{lemma}

\begin{proof}
By Assumption~\ref{assum:L_lip}, $|\widehat{\mathcal{L}}(\phi_q^*) - \widehat{\mathcal{L}}(\bar{\phi})| \le L_{\widehat{\mathcal{L}}}\,\|\phi_q^* - \bar{\phi}\|_{\mathrm{eval}}$. Writing both networks as sums over sub-networks,
\begin{align*}
\|\phi_q^* - \bar{\phi}\|_{\mathrm{eval}}
&= \Big\|\sum_{k=1}^{\bar{\mathfrak m}}\sum_{v=1}^{R}\frac{\bar c_k}{R}\big(\phi_{s_{k,v}}^{[0]}-\bar\phi^k\big)\Big\|_{\mathrm{eval}} \\
&\le \sum_{k,v} \frac{|\bar{c}_k|}{R} \cdot \mathfrak{L}\,\|(\boldsymbol{\theta}_{s_{k,v}})^{[0]} - \bar{\boldsymbol{\theta}}_k\|_\infty
\le \bar{M}\,\mathfrak{L}\,\delta,
\end{align*}
using $\sum_k |\bar{c}_k| \le \bar{M}$ and the cover guarantee.
\end{proof}

\subsubsection{Bounding the iteration error}\label{sec:iter_error}

The iteration error $\widehat{\mathcal{L}}(\phi_{\boldsymbol{\theta}_{\mathcal{A}}}) - \widehat{\mathcal{L}}(\phi_q^*)$ is controlled by a convex PGD argument on the outer coefficients. For fixed inner weights, the map $\boldsymbol{\theta}_{\mathrm{out}} \mapsto \widehat{\mathcal{L}}(\boldsymbol{\theta}_{\mathrm{in}}, \boldsymbol{\theta}_{\mathrm{out}})$ is convex because the PNN output is linear in $\boldsymbol{\theta}_{\mathrm{out}}$ and the loss is convex.

\begin{lemma}[Convex outer PGD; \citealp{kohler2023opt}, Lemma~1]\label{lem:O2b}
Let $\mathrm X$ and $\mathrm Y$ be nonempty closed convex sets, and suppose that $\boldsymbol y\mapsto F(\boldsymbol x,\boldsymbol y)$ is differentiable and convex for each fixed $\boldsymbol x\in\mathrm X$. For any sequence $\boldsymbol x_0,\ldots,\boldsymbol x_{T-1}\in\mathrm X$ and $\boldsymbol y_0\in\mathrm Y$, define, for $t=0,\ldots,T-1$,
\[
\boldsymbol y_{t+1}=\operatorname{Proj}_{\mathrm Y}\big\{\boldsymbol y_t-\lambda\nabla_{\boldsymbol y}F(\boldsymbol x_t,\boldsymbol y_t)\big\},
\]
with $\lambda = T^{-1}$. Suppose (a)~$\|\nabla_{\boldsymbol{y}} F(\boldsymbol{x}, \boldsymbol{y})\|_2 \le V$ on $\mathrm{X} \times \mathrm{Y}$, and (b)~there exists $\boldsymbol{y}^* \in \mathrm{Y}$ such that $|F(\boldsymbol{x}_t, \boldsymbol{y}^*) - F(\boldsymbol{x}_0, \boldsymbol{y}^*)| \le U_0 \|\boldsymbol{y}^*\|_2 \|\boldsymbol{x}_t - \boldsymbol{x}_0\|_2$ for all $t$. Then
\[
\min_{0\le t<T} F(\boldsymbol{x}_t, \boldsymbol{y}_t) \;\le\; F(\boldsymbol{x}_0, \boldsymbol{y}^*) + U_0 \|\boldsymbol{y}^*\|_2 \operatorname{diam}(\mathrm{X}) + \frac{\|\boldsymbol{y}^* - \boldsymbol{y}_0\|_2^2}{2} + \frac{V^2}{2T}.
\]
\end{lemma}

\begin{proof}[Proof sketch]
By convexity of $\boldsymbol{y} \mapsto F(\boldsymbol{x}_t, \boldsymbol{y})$, the projection inequality, and the polarisation identity,
\[
F(\boldsymbol{x}_t, \boldsymbol{y}_t) - F(\boldsymbol{x}_t, \boldsymbol{y}^*)
\le \frac{1}{2\lambda}\big[\|\boldsymbol{y}_t - \boldsymbol{y}^*\|_2^2 - \|\boldsymbol{y}_{t+1} - \boldsymbol{y}^*\|_2^2 + \lambda^2 V^2\big].
\]
Telescoping over $t = 0, \ldots, T-1$ with $\lambda = T^{-1}$, adding $\frac{1}{T}\sum_t F(\boldsymbol{x}_t, \boldsymbol{y}^*)$, and using condition~(b) together with $\min_t \le \text{average}$ yields the claim.
\end{proof}

To bound the iteration error $\widehat{\mathcal{L}}(\phi_{\boldsymbol{\theta}_{\mathcal{A}}}) - \widehat{\mathcal{L}}(\phi_q^*)$, we instantiate Lemma~\ref{lem:O2b} with $\boldsymbol{x} = \boldsymbol{\theta}_{\mathrm{in}}$, $\boldsymbol{y} = \boldsymbol{\theta}_{\mathrm{out}}$, reference point $\boldsymbol{y}^* = \boldsymbol{\theta}_{\mathrm{out}}^{\mathfrak{m},*}$, and $(\boldsymbol{x}_0, \boldsymbol{y}_0) = ((\boldsymbol{\theta}_{\mathrm{in}}^{\mathfrak{m}})^{[0]}, \boldsymbol{0})$. Conditions (a) and (b) are verified as follows.

Condition (a): gradient bound.
For regression, differentiation with respect to $c_k$ gives a sample average of the loss derivative times the clipped sub-network. For DRM, it also contains the corresponding gradient evaluations. Let $B_{\ell'}$ bound the appropriate dual norm of the derivative of the sampled loss, including the quadrature weights and PDE coefficients. The evaluation bound then yields
\begin{equation}\label{eq:grad_bound}
\|\nabla_{\boldsymbol\theta_{\mathrm{out}}}\widehat{\mathcal L}\|_2\le B_{\ell'}\,\mathcal B_{\mathrm{sub}}(B_{\bar{\boldsymbol\theta}}+\eta)\sqrt{\mathfrak m}=:V.
\end{equation}

Condition (b): inner-weight drift.
For fixed outer coefficients $\boldsymbol{\theta}_{\mathrm{out}}^{\mathfrak{m},*}$, the change of $\widehat{\mathcal{L}}$ when the inner weights evolve from the initialization $\boldsymbol{\theta}_{\mathrm{in}}^{[0]}$ to the $t$-th PGD iterate $\boldsymbol{\theta}_{\mathrm{in}}^{[t]}$ is bounded via Assumption~\ref{assum:L_lip} and the sub-network Lipschitz property \eqref{eq:subnet_props}. Because the weights at step~$t$ are bounded by $B_{\bar{\boldsymbol{\theta}}}+\eta$ rather than the initial bound $B_{\bar{\boldsymbol{\theta}}}$, the effective sub-network Lipschitz constant is $\mathfrak{L}(B_{\bar{\boldsymbol{\theta}}}+\eta)$, evaluated on the enlarged parameter set. Writing $\delta \Phi_t(\boldsymbol{x}) := \Phi_{\boldsymbol{\theta}_{\mathrm{in}}^{[t]}, \boldsymbol{\theta}_{\mathrm{out}}^{\mathfrak{m},*}}(\boldsymbol{x}) - \Phi_{\boldsymbol{\theta}_{\mathrm{in}}^{[0]}, \boldsymbol{\theta}_{\mathrm{out}}^{\mathfrak{m},*}}(\boldsymbol{x})$ and applying Cauchy--Schwarz,
\begin{align*}
\|\delta\Phi_t\|_{\mathrm{eval}}
&\le \sum_{i=1}^{\mathfrak m}|c_i^*|\,\|\phi_i^{[t]}-\phi_i^{[0]}\|_{\mathrm{eval}}\le \|\boldsymbol\theta_{\mathrm{out}}^{\mathfrak m,*}\|_2 \Big(\sum_{i=1}^{\mathfrak m}\|\phi_i^{[t]}-\phi_i^{[0]}\|_{\mathrm{eval}}^2\Big)^{1/2}\\
&\le \|\boldsymbol{\theta}_{\mathrm{out}}^{\mathfrak{m},*}\|_2\,
\mathfrak{L}(B_{\bar{\boldsymbol{\theta}}}+\eta)\,
\|\boldsymbol{\theta}_{\mathrm{in}}^{[t]} - \boldsymbol{\theta}_{\mathrm{in}}^{[0]}\|_2.
\end{align*}
Combined with Assumption~\ref{assum:L_lip}, this gives
\begin{equation}\label{eq:drift_bound}
\big|\widehat{\mathcal{L}}(\boldsymbol{\theta}_{\mathrm{in}}^{[t]}, \boldsymbol{\theta}_{\mathrm{out}}^{\mathfrak{m},*}) - \widehat{\mathcal{L}}(\boldsymbol{\theta}_{\mathrm{in}}^{[0]}, \boldsymbol{\theta}_{\mathrm{out}}^{\mathfrak{m},*})\big|
\;\le\; L_{\widehat{\mathcal{L}}}\,\mathfrak{L}(B_{\bar{\boldsymbol{\theta}}}+\eta)\;\|\boldsymbol{\theta}_{\mathrm{out}}^{\mathfrak{m},*}\|_2\;\|\boldsymbol{\theta}_{\mathrm{in}}^{[t]} - \boldsymbol{\theta}_{\mathrm{in}}^{[0]}\|_2.
\end{equation}
Since $\|\boldsymbol{\theta}_{\mathrm{in}}^{[t]} - \boldsymbol{\theta}_{\mathrm{in}}^{[0]}\|_2 \le \eta$ on $A_\eta$, condition~(b) of Lemma~\ref{lem:O2b} holds with $U_0 = L_{\widehat{\mathcal{L}}}\,\mathfrak{L}(B_{\bar{\boldsymbol{\theta}}}+\eta)$ and $\boldsymbol{y}^* = \boldsymbol{\theta}_{\mathrm{out}}^{\mathfrak{m},*}$. Substituting $\|\boldsymbol{\theta}_{\mathrm{in}}^{[t]} - \boldsymbol{\theta}_{\mathrm{in}}^{[0]}\|_2 \le \eta$ for $\|\boldsymbol{x}_t - \boldsymbol{x}_0\|_2$ in the drift bound of Lemma~\ref{lem:O2b} yields the stated factor $\eta$.

With zero initialization, $\|\boldsymbol y^*-\boldsymbol y_0\|_2^2\le\bar M^2/R$. Choose the output with smallest empirical risk among $0\le t<T$. Then its risk equals the minimum controlled by Lemma~\ref{lem:O2b}, giving the following bound. A last-iterate conclusion additionally requires empirical-risk descent.

\begin{lemma}[Iteration error]\label{lem:O2iter}
Under the hypotheses above, with $\zeta = \bar{M}$, zero initialization of $\boldsymbol{\theta}_{\mathrm{out}}$, $\lambda=T^{-1}$, and the best-iterate output,
\[
\widehat{\mathcal{L}}(\phi_{\boldsymbol{\theta}_{\mathcal{A}}}) - \widehat{\mathcal{L}}(\phi_q^*) \;\le\; \frac{L_{\widehat{\mathcal{L}}}\,\bar{M}\,\mathfrak{L}(B_{\bar{\boldsymbol{\theta}}}+\eta)\,\eta}{\sqrt{R}} + \frac{\bar{M}^2}{2R} + \frac{V^2}{2T},
\]
where $V$ is given by \eqref{eq:grad_bound}, $\mathfrak{L}(B_{\bar{\boldsymbol{\theta}}}+\eta)$ denotes the sub-network Lipschitz constant \eqref{eq:frakturL} evaluated with the weight bound $B_{\bar{\boldsymbol{\theta}}}+\eta$, and $L_{\widehat{\mathcal{L}}}$ is the Lipschitz constant from Assumption~\ref{assum:L_lip}.
\end{lemma}

The three terms reflect the three sources of suboptimality in the convex PGD argument: the first is the inner-weight drift (the $\boldsymbol{x}$-variable in Lemma~\ref{lem:O2b}); the second is the initial distance to the reference point (zero initialization vs.\ $\boldsymbol{\theta}_{\mathrm{out}}^{\mathfrak{m},*}$); the third is the gradient-descent convergence rate on the outer convex problem.

\subsubsection{Assembling the complete optimization bound}\label{sec:master_lemma}

Inserting the bounds of Lemmas~\ref{lem:O2iter} and~\ref{lem:O2init} into the decomposition \eqref{eq:opt_decomp} yields the complete optimization error bound.

\begin{lemma}[Combined initialization and iteration bound]\label{lem:O2}
Under the hypotheses of Lemmas~\ref{lem:O2init} and~\ref{lem:O2iter}, on the event $\mathcal{E}_{\mathrm{cover}}(\delta)$,
\begin{align}\label{eq:O2_bound}
\mathcal{E}_{opt,-} \;\le\; L_{\widehat{\mathcal{L}}}\,\bar{M}\,\mathfrak{L}\,\delta
+ \frac{L_{\widehat{\mathcal{L}}}\,\bar{M}\,\mathfrak{L}(B_{\bar{\boldsymbol{\theta}}}+\eta)\,\eta}{\sqrt{R}}
+ \frac{\bar{M}^2}{2R}
+ \frac{V^2}{2T},
\end{align}
where $L_{\widehat{\mathcal{L}}}$ is the Lipschitz constant from Assumption~\ref{assum:L_lip}, $\mathfrak{L}$ is given by \eqref{eq:frakturL}, $\mathfrak{L}(B_{\bar{\boldsymbol{\theta}}}+\eta)$ denotes the same formula evaluated with the weight bound $B_{\bar{\boldsymbol{\theta}}}+\eta$, and $V$ is given by \eqref{eq:grad_bound}.
\end{lemma}

\begin{remark}
The constants $\mathcal B_{\mathrm{sub}}$, $\mathfrak L$, $L_{\widehat{\mathcal L}}$, and $B_{\ell'}$ must be estimated for the chosen architecture and loss. In DRM, this includes first derivatives of the network; the source's quantitative optimization bound is given in Remark~\ref{rem:drm_opt}.

The first term in \eqref{eq:O2_bound} is controlled by the cover accuracy $\delta$, the next two by the replica count $R$, and the last by the iteration budget $T$. For fixed sub-network architecture, the covering lemma permits $Q_{\mathrm{trial}}$ of order $(2B_{\bar{\boldsymbol\theta}}/\delta)^{J_{\max}}\log(\bar{\mathfrak m}R/\xi)$ for failure probability $\xi$. Thus covering is polynomial in $1/\delta$ at fixed $J_{\max}$; a growing sub-network depth can instead produce super-polynomial over-parameterization.
\end{remark}

\subsection{Complete rates for the running examples}\label{sec:complete_rates}

The two examples use the approximation-based decomposition \eqref{eq:D1}, the PNN complexity reduction, and an outer-convexity optimization argument, with the architecture and algorithm choices required by their respective source theorems.

\subsubsection{Nonparametric regression}\label{sec:reg_complete}

For nonparametric regression, the comparison network is supplied by Theorem~\ref{thm:A5r}. The analysis of \cite{kohler2023opt} combines approximation, the statistical complexity of clipped parallel networks, and a best-iterate PGD bound. The sub-network initialization is layer-dependent, so the covering argument uses the corresponding rectangular parameter ranges.

\begin{theorem}[Regression complete error analysis; \citealp{kohler2023opt}, Thm.~1]\label{thm:reg_complete}
Let $\beta = s + r$ with $s = \lceil\beta\rceil-1 \in \mathbb{N}_0$, $r \in (0,1]$, $B_{\mathrm{REG}}>0$, and let $\{(\boldsymbol X_i, Y_i)\}_{i=1}^n$ be i.i.d.\ $\mathbb{R}^d\times\mathbb{R}$-valued with $\mathrm{supp}(\boldsymbol X)\subseteq[0,1]^d$, a sub-Gaussian response satisfying $\mathbb{E}\{e^{c_2 Y^2}\}<\infty$, and $(\beta, B_{\mathrm{REG}})$-smooth regression function $f_0(\boldsymbol x) = \mathbb{E}\{Y\mid\boldsymbol X=\boldsymbol x\}$. Run PGD on the clipped parallel architecture $f_{\boldsymbol\theta}=\sum_k c_k[\phi^k_{\boldsymbol\theta}]_{b_n}$, where $[z]_b=\min\{b,\max\{-b,z\}\}$ and $b_n=c_1\log n$. Initialize the outer coefficients at zero, the first-layer parameters uniformly on $[-n,n]$, and the remaining trainable inner parameters uniformly on $[-1,1]$, as in \cite[Section~2]{kohler2023opt}. Use the parameter choices
\begin{alignat*}{2}
\frac{\mathfrak{m}}{n^{c_3(\log n+1)}} \;&\longrightarrow\; \infty, &\qquad
W \;&=\; 18(s + d)^{d+1}, \\
L \;&=\; c_5\log n, &\qquad
B_{\boldsymbol{\theta}} \;&\sim\; n, \\
\eta \;&=\; 1, &\qquad
\zeta \;&=\; c_4\, n^{d/(2(2\beta + d))}, \\
T \;&=\; \mathfrak{m}^2 n^2, &\qquad
\lambda \;&=\; 1/T,
\end{alignat*}
together with output truncation at level $b_n = c_1\log n$: $\widehat{f}_{\boldsymbol{\theta}_{\mathcal{A}}} = [f_{\boldsymbol{\theta}^{[\widehat t]}}]_{b_n}$, where $[x]_b := \min(\max(x,-b),\,b)$ and $\widehat t = \arg\min_{0\le t\le T}\frac{1}{n}\sum_i|f_{\boldsymbol{\theta}^{[t]}}(\boldsymbol X_i) - Y_i|^2$. Then for $n$ sufficiently large,
\[
\mathbb{E}_{S_n,\boldsymbol\theta^{[0]}}\int |\widehat{f}_{\boldsymbol{\theta}_{\mathcal{A}}}(\boldsymbol x) - f_0(\boldsymbol x)|^2\, P_X(\mathrm d\boldsymbol x) \;\le\; c_8\, (\log n)^5\, n^{-\beta/(2\beta + d)}.
\]
The initialization is independent of the training sample, and the expectation is taken jointly over both sources of randomness.
\end{theorem}

\begin{proof}[Proof sketch]
The argument follows the comparison-network decomposition in \eqref{eq:D1}, with the truncation steps of \cite[Section~5.4]{kohler2023opt}. Barred quantities describe the comparison network; the training network has many more parallel sub-networks. Sub-network clipping bounds the functions entering the complexity estimate, whereas final-output truncation bounds the estimator. On the event $\max_i|Y_i|\le b_n$,
\[
\widehat{\mathcal L}([f]_{b_n})\le\widehat{\mathcal L}(f),
\]
so the empirical comparison for the untruncated PGD output also controls its final truncation. The statistical supremum is taken over the final-output-truncated class, using truncated responses $[Y]_{b_n}$; the exponential-moment assumption controls the response-tail remainders. For sufficiently large $n$, truncation also cannot increase the population squared distance to the bounded regression function $f_0$.

Theorem~\ref{thm:A5r} gives approximation error $O(G^{-2\beta})$ and outer radius $\bar M\asymp G^d$, with constants depending on the target's H\"older norm. For fixed sub-network width, clipping and the complexity bound in \cite[Lemma~5]{kohler2023opt} give statistical error $\widetilde O(\bar M/\sqrt n)$; the logarithmically growing depth contributes logarithmic factors. On the joint response-truncation and initialization-cover event, the comparison argument controls the optimization contribution by $\widetilde O(1/n)$ under the stated parameter choices. The event's complement contributes $\widetilde O(1/n)$ because the final estimator is bounded by $b_n$ and the failure probability is $O(1/n)$. Together with the response-tail bounds, this yields total error at most $C G^{-2\beta}+\widetilde O(G^d/\sqrt n)+\widetilde O(1/n)$. Balancing the first two terms gives $G\asymp n^{1/(2(2\beta+d))}$ and the stated rate; see \cite[Section~5.4]{kohler2023opt}.
\end{proof}

\begin{remark}
The rate is slower than the H\"older minimax rate $n^{-2\beta/(2\beta+d)}$ \cite{stone1982optimal} by a factor of two in the exponent. It results from balancing the approximation bound with the global statistical estimate $\widetilde O(\bar M/\sqrt n)$; Lemma~\ref{lem:S7} itself only reduces the PNN complexity to the single-sub-network complexity. Improving the generalization estimate is identified as a possible route to sharper rates in \cite[Remark~1]{kohler2023opt}. For $\beta\geq1/2$, \cite{kohler2025rate} obtains the nearly minimax-optimal rate $n^{-2\beta/(2\beta+d)+\epsilon}$ for any $\epsilon>0$ using logistic parallel networks trained by unprojected full-batch gradient descent, with the initialization and parameter choices prescribed in that work.
\end{remark}

\subsubsection{Deep Ritz method}\label{sec:drm_complete}

The DRM complete analysis, due to \cite{jiao2024drm}, applies the same framework to the DRM example. Four differences from \S\ref{sec:reg_complete} determine the final rate. First, the target norm is $H^1(\Omega)$ rather than $L^2(P_X)$, introducing an energy-norm equivalence prefactor $\frac{2}{c_0\wedge 1}$ (Proposition~\ref{prop:drm_eq}) on the right-hand side of \eqref{eq:D1}. Second, approximation uses the tanh PNN form of Theorem~\ref{thm:A5} (with $p=2$) rather than the ReLU PNN of Theorem~\ref{thm:A5r}, shifting the effective smoothness from $\beta$ to $s-1$. Third, the statistical bound uses the derivative-aware decomposition of \S\ref{sec:complete_stat} in place of the simple Talagrand contraction, introducing a $(B_{\bar{\boldsymbol{\theta}}}+\eta)^{2\bar L}$ Lipschitz-blow-up prefactor. Fourth, the optimization bound uses derivative-sensitive constants and a step size ensuring descent of the empirical energy, and the parameterization is expressed through a target accuracy $\epsilon$ rather than directly in $n$.

\begin{theorem}[DRM complete error analysis; \citealp{jiao2024drm}, Theorem~1]\label{thm:drm_complete}
Suppose $u_0 \in H^s(\Omega)$ for some integer $s \ge 2$ is the weak solution of \eqref{eq:drm_pde}. Fix $0<\mu<1$. For sufficiently small $\epsilon>0$, run PGD on a $\mathcal{PNN}(\mathfrak{m}, M, \{W, L, B_{\boldsymbol{\theta}}\})$ class with equal interior and boundary sample sizes, both denoted by $N$, and parameter choices
\begin{alignat*}{2}
\mathfrak{m} &= \lceil C\, \epsilon^{-\tilde{C}_1(\mu, d, \iota_0, s)}\rceil, &\qquad
W &= 2^{\lceil \log_2(d + s - 1)\rceil + 1}, \\
L &= \lceil \log_2(d + s - 1)\rceil + 2, &\qquad
B_{\boldsymbol{\theta}} &= C\, \epsilon^{-(2d+2s)/(s-1-\mu)}, \\
\eta &= \epsilon^{-\iota}, &\qquad
\zeta &= C\, \epsilon^{-3d/(2(s-1-\mu))}, \\
T &= \lceil C\, \epsilon^{-\tilde{C}_2(\mu, d, \iota_0, s)}\rceil, &\qquad
\lambda &= T^{-1}\asymp \epsilon^{\tilde{C}_2(\mu, d, \iota_0, s)}, \\
N &= \lceil C\, \epsilon^{-\tilde{C}_3(\mu, d, \iota_0, s)}\rceil,
\end{alignat*}
where $M=\zeta$, $\iota > 0$ is a free parameter controlling the inner radius $\eta = \epsilon^{-\iota}$, $\iota_0 = \max\{\iota, (2d+2s)/(s-1-\mu)\}$, and $\tilde{C}_1, \tilde{C}_2, \tilde{C}_3$ are the explicit positive exponents given in \cite[Theorem~1]{jiao2024drm}. Then with probability at least $1 - 2\,\epsilon^{\tilde{C}_3}$,
\[
\big\|u_{\boldsymbol{\theta}_{\mathcal{A}}} - u_0\big\|^2_{H^1(\Omega)} \;\le\; C\, \epsilon\, \log^{1/2}(C\epsilon^{-1}) \;=\; \tilde{\mathcal{O}}(\epsilon),
\]
where $u_{\boldsymbol{\theta}_{\mathcal{A}}}$ is the last PGD iterate. The constants in the parameter choices are chosen as in the cited theorem; in particular the step size satisfies its descent condition, and the output class has inner bound $B_{\boldsymbol\theta}+\eta$.
\end{theorem}

\begin{proof}[Proof sketch]
The argument follows the same pattern as Theorem~\ref{thm:reg_complete}, with the four differences recorded above.

Apply Theorem~\ref{thm:A5} with $p=2$ to a bounded Sobolev extension of $u_0$ to $[0,1]^d$, and restrict the resulting tanh PNN to $\Omega$. This supplies the comparison network $\bar u$, with constants depending on the extension operator and $\|u_0\|_{H^s(\Omega)}$. The energy-norm equivalence (Proposition~\ref{prop:drm_eq}) converts its $H^1$ approximation guarantee into $\mathcal{E}^{\mathrm{DRM,c}}_{app} \le C\epsilon^2$.

The statistical analysis separates the energy into gradient-squared, value-squared, source-linear, and boundary-linear terms. The linear-combination argument of Lemma~\ref{lem:S7} reduces the value and derivative complexities to their single-sub-network counterparts. Contraction handles the squared terms, and parameter covering controls the remaining complexities. By \cite[Theorem~5]{jiao2024drm}, with probability at least $1-\xi$,
\[
\mathcal E^{\mathrm{DRM,c}}_{sta}
\le C\bar M^2(B_{\bar{\boldsymbol\theta}}+\eta)^{2\bar L}N^{-1/2}
\left[\sqrt{\log((B_{\bar{\boldsymbol\theta}}+\eta)\bar W\bar L N)}+\sqrt{\log(1/\xi)}\right].
\]
Taking $\xi=\epsilon^{\tilde C_3}$ and the stated sample size gives $\mathcal E^{\mathrm{DRM,c}}_{sta}\le C\epsilon\sqrt{\log(C/\epsilon)}$ on this event.

The optimization contribution is controlled by \cite[Theorem~4]{jiao2024drm}, stated in Remark~\ref{rem:drm_opt}. Its initialization, replica, and iteration terms can all be made $O(\epsilon)$ by polynomial choices of $(\delta^{-1},R,T,Q_{\mathrm{trial}})$ because the sub-network width and depth are fixed for fixed $d,s$.

Choose the initialization parameters so that the covering event also has failure probability at most $\epsilon^{\tilde C_3}$. The union bound gives probability at least $1-2\epsilon^{\tilde C_3}$ for the intersection of the sampling and covering events. On this intersection, substituting the three bounds into Lemma~\ref{lem:D1} and applying the energy-norm equivalence yields the stated bound.
\end{proof}

The cited DRM theorem uses a target-accuracy, high-probability format, whereas Theorem~\ref{thm:reg_complete} bounds expected error in terms of sample size. Inverting $N = \lceil C\epsilon^{-\tilde{C}_3}\rceil$ gives the DRM rate $\tilde{\mathcal{O}}(N^{-\gamma})$, where $\gamma = \tilde{C}_3^{-1}$.

\begin{remark}[DRM instantiation of the optimization bound]\label{rem:drm_opt}
The optimization argument is applied to the full empirical energy, not separately to its four terms. For fixed sub-network width and depth, tanh parameter estimates give $\mathcal B_{\mathrm{sub}}(B)=O(B^{\bar L})$ and $\mathfrak L(B)=O(B^{2\bar L})$ in the derivative-sensitive evaluation norm. The resulting energy constants have orders
\[
L_{\widehat{\mathcal L}},\ B_{\ell'}=O\!\left(\bar M(B_{\bar{\boldsymbol\theta}}+\eta)^{\bar L}\right),
\]
and $V=O(\bar M(B_{\bar{\boldsymbol\theta}}+\eta)^{2\bar L}\sqrt{\mathfrak m})$. Choose
\[
T\ge\frac{C_{\mathrm{Lip}}}{2}\,\mathfrak m\bar M^2(B_{\bar{\boldsymbol\theta}}+\eta)^{4\bar L},
\qquad \lambda=T^{-1}.
\]
Then the step size satisfies the empirical-energy descent condition and $\lambda T=1$. With this choice, \cite[Theorem~4]{jiao2024drm} gives the last-iterate bound
\[
\mathcal{E}_{opt,-} \;\le\; C_0 \bar{M}^2 B_{\bar{\boldsymbol{\theta}}}^{3\bar{L}} \delta
+ \frac{C_1 \bar{M}^2 (B_{\bar{\boldsymbol{\theta}}}+\eta)^{3\bar{L}} \eta}{\sqrt{R}}
+ \frac{\bar{M}^2}{2R}
+ \frac{C_V \bar{M}^2 (B_{\bar{\boldsymbol{\theta}}}+\eta)^{4\bar{L}} \mathfrak{m}}{2T},
\]
where $C_0, C_1, C_V$ absorb $\bar{W}, \bar{L}, d$ and the constants from the remaining three sub-functionals. The first term involves only the initial weight bound $B_{\bar{\boldsymbol{\theta}}}$, while the second and fourth involve $B_{\bar{\boldsymbol{\theta}}}+\eta$ because the sub-networks at PGD iterates have weights bounded by the initial radius plus the inner drift. Compared with the squared-loss regression case (where the corresponding exponents are $O(\bar{L})$ for the initialization and drift terms and $O(1)$ for the iteration term), the exponents on $B_{\bar{\boldsymbol{\theta}}}$ here are inflated to $O(3\bar{L})$ and $O(4\bar{L})$ respectively. This inflation is the source of the large polynomial exponents $\tilde{C}_1, \tilde{C}_2$ appearing in the parameter choices of Theorem~\ref{thm:drm_complete}.
\end{remark}

\begin{remark}[Sub-optimality and limitations of the three-error framework]
The complete regression result is slower than the H\"older minimax rate, and the DRM result guarantees convergence with large polynomial resource requirements. The two analyses share a comparison-network decomposition and outer-layer convexity, but require different initialization, clipping, and derivative estimates. Quantitative comparisons with minimax bounds must also use the same observation model and target class.

The PGD analyses above use the parallel architecture, prescribed initialization, and explicit projection steps. A PNN can be represented as a feedforward network with a particular sparse connectivity pattern, but these proofs do not cover arbitrary dense or convolutional architectures with their usual training procedures. The logistic-network result cited in \S\ref{sec:reg_complete} dispenses with projection for full-batch gradient descent; it does not establish the same guarantee for unprojected stochastic gradient descent. These distinctions inform the open directions in \S\ref{sec:open_horizons}.
\end{remark}

\subsection{Numerical illustration: parallel networks and parameter movement}\label{sec:complete_numerical}

We next examine the inner-parameter displacement and outer-coefficient norm that enter the PNN analysis. On $\Omega=(0,1)$, consider
\[
-u''+u=h,\qquad u'(0)=u'(1)=0,\qquad h=-u_0''+u_0,
\]
where
\[
u_0(x)=\frac{5\cos(2\pi x)+\cos(4\pi x)}{\sqrt{13+58\pi^2}},
\qquad \|u_0\|_{H^1(0,1)}=1.
\]
Each of the $\mathfrak m$ tanh sub-networks has two hidden layers of width $W=4$ and depth $L=3$, giving $\mathfrak D=34\mathfrak m$ trainable parameters including the outer coefficients. All sub-network weights and biases are initialized independently and uniformly on $[-4,4]$, and the outer coefficients are initialized to zero. We discretize the energy by midpoint quadrature with $N=512$ points and train in double precision by full-batch gradient descent without projection for $T=40\,000$ steps, with step size $\lambda=10^{-3}\min\{1,64/\mathfrak m\}$.

\begin{figure}[htbp]
\centering
\includegraphics[width=\linewidth]{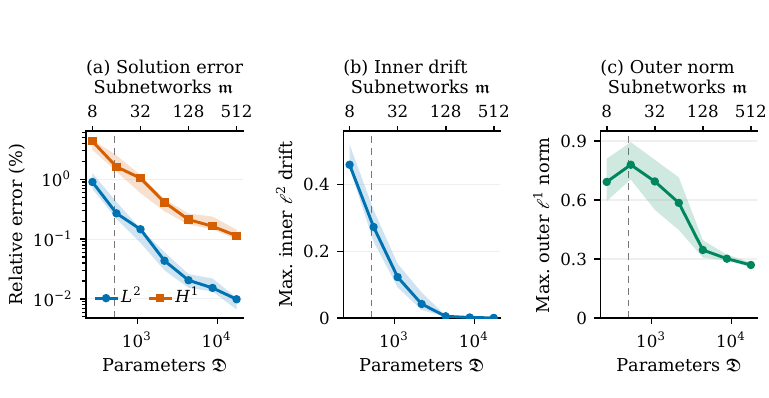}
\caption{One-dimensional Neumann problem: relative solution errors, maximum inner-parameter displacement, and maximum outer-coefficient $\ell^1$ norm as the number of sub-networks varies. Lines and shaded bands show medians and interquartile ranges over ten initializations. The dashed line marks $\mathfrak D=N$.}
\label{fig:pnn_size_parameter_movement}
\end{figure}

For $\mathfrak m\in\{8,16,32,64,128,256,512\}$, we use the same ten initialization seeds and evaluate relative $L^2$ and $H^1$ errors on $2^{14}$ test points. Figure~\ref{fig:pnn_size_parameter_movement} also records the two trajectory quantities
\[
\max_{0\le t\le T}\|\boldsymbol\theta_{\mathrm{in}}^{[t]}-\boldsymbol\theta_{\mathrm{in}}^{[0]}\|_2,
\qquad
\max_{0\le t\le T}\|\boldsymbol\theta_{\mathrm{out}}^{[t]}\|_1.
\]

The median relative $H^1$ error decreases from $4.339\%$ at $\mathfrak m=8$ to $0.114\%$ at $\mathfrak m=512$. Under the stated initialization and step-size choices, the median maximum inner displacement decreases, while the outer-coefficient norm initially rises and then decreases. Accurate solutions can therefore coexist with moderate values of the two quantities controlled in the analysis. Since this experiment uses a fixed grid and unprojected training, it illustrates their behavior rather than testing the Monte Carlo and PGD rate in Theorem~\ref{thm:drm_complete}.

%% file: sections/applications.tex
\section{Applications}\label{sec:applications}

\par This section reviews error analysis in scientific machine learning (Section~\ref{section:app:sml}), generative learning (Section~\ref{section:app:generative}), reinforcement learning (Section~\ref{section:app:rl}), representation learning (Section~\ref{section:app:repssl}), and large language models (Section~\ref{section:app:llm}). For each application, we introduce the learning problem and training criterion, then use representative studies to explain how function approximation and finite-sample estimation affect the computed solution, generated distribution, or prediction performance. The selected convergence results quantify these effects under the assumptions of the respective models.

\subsection{Scientific machine learning}
\label{section:app:sml}

Scientific machine learning uses mathematical models together with data to solve scientific problems. Here we focus on neural network methods for partial differential equations: forward problems seek a solution for given equations and boundary data, whereas inverse problems recover unknown sources or coefficients from observations. The Deep Ritz Method (DRM)~\cite{Weinan2017The}, analyzed in Sections~\ref{sec:drm_fast}, \ref{sec:drm_overparam}, and \ref{sec:drm_complete}, approximates the solution by minimizing a variational energy. Physics-informed neural networks (PINNs)~\cite{raissi2019physics} instead penalize violations of the governing equation and boundary conditions at sampled collocation points, together with discrepancies from observations when available. These residual losses involve different approximation quantities and empirical-process classes from those of DRM. We first review a forward elliptic problem and then an inverse source problem, examining how stability estimates convert training-loss bounds into solution or recovery guarantees.

\subsubsection{Physics-informed neural networks}

Following~\cite{jiao2022rate}, consider the Dirichlet problem
\begin{equation}
\left\{
\begin{aligned}
-\Delta u + u &= h, && x\in\Omega,\\
u &= \psi, && x\in\partial\Omega,
\end{aligned}
\right.
\label{eq:elliptic}
\end{equation}
on a bounded domain $\Omega\subset\mathbb{R}^{d}$ with sufficiently regular boundary. Let $u_0$ denote its exact solution. The population physics-informed loss is
\begin{equation}
\mathcal L_{\mathrm{PINN}}(u) \coloneqq \|\Delta u-u+h\|_{L^2(\Omega)}^2 + \|u-\psi\|_{L^2(\partial\Omega)}^2.
\label{eq:pinn_population}
\end{equation}

Given independent samples $\{\boldsymbol X_i\}_{i=1}^{N_{\mathrm{in}}}$ from $\mathrm{unif}(\Omega)$ and $\{\boldsymbol Y_j\}_{j=1}^{N_b}$ from $\mathrm{unif}(\partial\Omega)$, its empirical counterpart is
\begin{equation}
\begin{aligned}
\widehat{\mathcal L}_{\mathrm{PINN}}(u) = \frac{|\Omega|}{N_{\mathrm{in}}} \sum_{i=1}^{N_{\mathrm{in}}} \bigl( \Delta u(\boldsymbol X_i) -u(\boldsymbol X_i) +h(\boldsymbol X_i) \bigr)^2 + \frac{|\partial\Omega|}{N_b} \sum_{j=1}^{N_b} \bigl( u(\boldsymbol Y_j)-\psi(\boldsymbol Y_j) \bigr)^2.
\end{aligned}
\label{eq:pinn_erm}
\end{equation}

Let $\widehat u = \argmin_{u\in\mathcal F_{NN}}\widehat{\mathcal L}_{\mathrm{PINN}}(u)$ be the empirical risk minimizer over the network class $\mathcal F_{NN}$. The first bridge from residual minimization to solution recovery is a two-sided comparison between the loss excess and the solution error, due to \citealp[Proposition 2.2]{jiao2022rate} (see also its proof). Since $\mathcal L_{\mathrm{PINN}}(u_0)=0$, for any $u\in H^2(\Omega)\cap C^2(\bar\Omega)$,
\begin{equation}
\|u-u_0\|_{H^{1/2}(\Omega)}^2 \;\lesssim\; \mathcal L_{\mathrm{PINN}}(u)-\mathcal L_{\mathrm{PINN}}(u_0) \;\lesssim\; \|u-u_0\|_{C^2(\bar\Omega)}^2,
\label{eq:pinn_equiv}
\end{equation}
with constants depending only on $\Omega$ and $d$. The lower bound is the a-priori stability estimate of \citealp{jiao2022rate}; the upper bound holds since both residuals are $C^2$-controlled. The two-sided estimate plays the role of the energy--norm equivalence of DRM: it converts risk control into the $H^{1/2}$-error, while the loss excess is controlled by the $C^2$ approximation error.

For any candidate class $\mathcal F$, define the sample-dependent uniform deviation
\[
\Delta_{\mathrm{PINN}}(\mathcal F)
:=\sup_{u\in\mathcal F}
\left|\mathcal L_{\mathrm{PINN}}(u)-\widehat{\mathcal L}_{\mathrm{PINN}}(u)\right|.
\]
Inserting the ERM comparison of \S\ref{sec:decomp} into the two-sided estimate and taking the best comparison function gives the two-error decomposition of \citealp[Proposition 2.2]{jiao2022rate}:
\begin{align}
\mathbb E \|\widehat u-u_0\|_{H^{1/2}(\Omega)}^2 \lesssim \underbrace{ \inf_{u\in\mathcal F_{NN}} \|u-u_0\|_{C^2(\bar\Omega)}^2 }_{\mathcal E_{\mathrm{app}}^{\mathrm{PINN}}} + \underbrace{\mathbb E\Delta_{\mathrm{PINN}}(\mathcal F_{NN})}_{\mathcal E_{\mathrm{sta}}^{\mathrm{PINN}}}.
\label{eq:pinn_decomp}
\end{align}
As in the classical regime of Section~\ref{sec:classical_setting}, the optimization error is set to zero by assuming access to the empirical risk minimizer.

\begin{theorem}[PINN convergence rate; \citealp{jiao2022rate}, Theorem~5.1]
\label{thm:pinn_forward}
Suppose that $u_0\in C^3(\bar\Omega)$ and that the data in \eqref{eq:elliptic} satisfy the regularity assumptions of \cite{jiao2022rate}. For any $\epsilon>0$, choose the norm-restricted $\mathrm{ReLU}^3$ class
\[
\mathcal F_{NN}=\{u\in\mathcal{NN}_{\mathrm{ReLU}^3}(W,L):\|u\|_{C^2(\bar\Omega)}\le 2\|u_0\|_{C^2(\bar\Omega)}\}
\]
with depth $L=\lceil\log_2d\rceil+2$ and width $W\asymp C(d,\Omega,\|u_0\|_{C^3(\bar\Omega)})\,\epsilon^{-d}$, and take $N_{\mathrm{in}},\,N_b\asymp C(d,\Omega,\|u_0\|_{C^3(\bar\Omega)})\,\epsilon^{-(2d+4+\nu)}$ interior and boundary samples, where $\nu>0$ is arbitrarily small. Then
\[
\mathbb E \|\widehat u-u_0\|_{H^{1/2}(\Omega)} \le \epsilon.
\]
\end{theorem}

In \eqref{eq:pinn_decomp}, the approximation term is controlled by the $C^2$ approximation theorem for $\mathrm{ReLU}^3$ networks (\citealp[Theorem~3.1]{jiao2022rate}), while the statistical term is controlled through pseudo-dimension and Rademacher-complexity bounds for the residual class (\citealp[Theorem~4.13]{jiao2022rate}). The use of $\mathrm{ReLU}^3$ is technically important because the strong-form PINN loss contains second derivatives, and hence the approximation analysis must control the network in $C^2$. The resulting rate nevertheless retains the curse of dimensionality through the width $W=O(\epsilon^{-d})$. Converted to the sample size, with $N_{\mathrm{in}}$ and $N_b$ of the same order, the guarantee reads 
\[
\mathbb E\|\widehat u-u_0\|_{H^{1/2}(\Omega)} \lesssim N_{\mathrm{in}}^{-1/(2d+4+\nu)}.
\]

Several subsequent results extend this baseline. The over-parameterized setting is studied in~\cite{chen2025Convergence}, where size restrictions are replaced by norm control in the spirit of Section~\ref{sec:overparam_setting}. Under compositional sparsity, \cite{jiao2023improved} obtains dimension-mitigated approximation and generalization rates, while stabilized PINNs for the wave equation are analyzed in~\cite{jiao2024Stabilized}. More recently, \cite{jiao2026Complete} established a complete three-error analysis for elliptic PINNs using projected stochastic gradient descent (PSGD). Write $\mathcal F_{NN}^{\mathrm{PSGD}}$ for the admissible parallel tanh network class in that analysis, containing the algorithm output $u_{\boldsymbol\theta_{\mathcal A}}$ and the comparison network $\bar u$. Replacing the ERM by this algorithm output adds the optimization term:
\begin{align*}
&\mathbb{E}_{\mathrm{SGD}}\big\|u_{\boldsymbol\theta_{\mathcal A}}-u_0\big\|_{H^{1/2}(\Omega)}^2
\\
&\leq C_{\mathrm{stab}}\Big[
2\Delta_{\mathrm{PINN}}(\mathcal F_{NN}^{\mathrm{PSGD}})+\underbrace{\mathbb{E}_{\mathrm{SGD}}\big[\widehat{\mathcal L}_{\mathrm{PINN}}(u_{\boldsymbol\theta_{\mathcal A}})\big]-\widehat{\mathcal L}_{\mathrm{PINN}}(\bar u)}_{\mathcal E_{\mathrm{opt},-}^{\mathrm{PINN}}} +C_{\mathrm{app}}\underbrace{\|\bar u-u_0\|_{C^2(\bar\Omega)}^2}_{\mathcal E_{\mathrm{app}}^{\mathrm{PINN,c}}}\Big],
\end{align*}
where $C_{\mathrm{stab}}$ and $C_{\mathrm{app}}$ are the stability and approximation constants in \eqref{eq:pinn_equiv}. The network $\bar u$ serves as the approximation benchmark and the reference for the optimization gap. Here $\mathbb E_{\mathrm{SGD}}$ averages over the SGD sampling conditional on the training set and initialization.

Suppose that the Dirichlet problem~\eqref{eq:elliptic} on a bounded convex domain $\Omega\subseteq[0,1]^d$ has a unique solution $u_0\in C^m(\bar\Omega)$ for an integer $m\geq3$. For each sufficiently small $\epsilon>0$, Jiao et al.~\cite[Theorem~3.5]{jiao2026Complete} choose a parallel tanh network, projection radii, and a PSGD step size for which the numbers of sub-networks, training points, and iterations grow polynomially in $\epsilon^{-1}$. With high probability over the training samples and initialization, the squared approximation error is bounded by $C\epsilon^2$, the optimization gap by $C\epsilon$, and the statistical deviation by $C\epsilon\sqrt{\log(C/\epsilon)}$. Substitution into the preceding decomposition gives
\[
\mathbb E_{\mathrm{SGD}}
\bigl\|u_{\boldsymbol\theta_{\mathcal A}}-u_0\bigr\|_{H^{1/2}(\Omega)}^2
\lesssim \epsilon\sqrt{\log(C/\epsilon)}.
\]
Methodologically, this result follows the same structured-architecture approach as the complete DRM analysis in Chapter~\ref{sec:complete_setting}, but its optimization argument must additionally accommodate the second derivatives appearing in the PINN residual.

\subsubsection{Physics-informed neural networks for inverse problems}

PINNs also extend to inverse problems, but the decomposition acquires several new ingredients. We illustrate them through elliptic source identification, following~\cite{Duan2024Recovering}. Relative to the forward problem, three structural changes occur: both the state and the source are unknown; noisy interior measurements introduce an additional data-fidelity term; and the analysis uses Tikhonov regularization to stabilize the ill-posed source recovery problem. These changes alter both the population stability estimate and the oracle inequality.

Consider
\begin{equation}
\left\{
\begin{aligned}
-\Delta u+u &= h, && x\in\Omega,\\
\partial_{\mathbf n}u &=0, && x\in\partial\Omega,
\end{aligned}
\right.
\label{eq:pinn_inverse_pde}
\end{equation}
where both $u$ and $h$ are unknown. The original analysis of \cite{Duan2024Recovering} treats the more general operator $-\Delta u+V(x)u$ with nonzero Neumann data; we state the specialization $V\equiv1$ and homogeneous boundary data to keep the notation aligned with the model problems used elsewhere in the review.

Let $(u_0,h_0)\in H^2(\Omega)\times H^1(\Omega)$ denote the ground-truth pair and suppose that noisy interior measurements are available:
\begin{equation}
Y=u_0(\boldsymbol X)+e_y, \qquad \boldsymbol Z = \nabla u_0(\boldsymbol X)+\boldsymbol e_z, \qquad \boldsymbol X\sim\mathrm{unif}(\Omega),
\label{eq:pinn_inverse_data}
\end{equation}
with the componentwise noise bounds
\[
|e_y|\leq\delta,\qquad |(\boldsymbol e_z)_j|\leq\delta,\quad j=1,\ldots,d,
\qquad\text{almost surely}.
\]
Here $\delta$ is the noise level, as in \citealp[Eq.~(1.3)]{Duan2024Recovering}. When gradients are not observed directly, they may themselves be estimated from function values using Sobolev-regularized regression; see~\cite{ding2025Semi}.

The population measurement-fidelity functional is
\begin{equation}
\mathcal J(u) \coloneqq |\Omega| \mathbb E \left[ (Y-u(\boldsymbol X))^2 + \|\boldsymbol Z-\nabla u(\boldsymbol X)\|_2^2 \right].
\label{eq:pinn_fid}
\end{equation}

The physics residuals are
\[
\mathcal R_{\mathrm{int}}(u,h) \coloneqq \|\Delta u-u+h\|_{L^2(\Omega)}^2, \qquad \mathcal R_{\mathrm{bdy}}(u) \coloneqq \|\partial_{\mathbf n}u\|_{L^2(\partial\Omega)}^2,
\]
and the Tikhonov regularizer is
\[
\mathcal R(h) \coloneqq \|h\|_{H^1(\Omega)}^2.
\]

The regularized population objective is therefore
\begin{equation}
\mathcal L_{\mathrm{PINN}}^{\lambda_{\mathrm{reg}}}(u,h) \coloneqq \underbrace{\mathcal J(u)}_{\text{data fidelity}} + \underbrace{ \mathcal R_{\mathrm{int}}(u,h) + \mathcal R_{\mathrm{bdy}}(u) }_{\text{physics consistency}} + \underbrace{ \lambda_{\mathrm{reg}}\mathcal R(h) }_{\text{regularization}}.
\label{eq:pinn_inv_risk}
\end{equation}

The following stability estimate converts control of the population risk into recovery guarantees, with an explicit dependence on the noise level.

\begin{proposition}[{Stability estimate; adapted from \citealp[Theorem~2.1]{Duan2024Recovering}}]
\label{proposition:pinns:inverse:stability}
Suppose that $\Omega$ is a bounded, simply connected domain with sufficiently smooth boundary and that $(u_0,h_0)\in H^2(\Omega)\times H^1(\Omega)$ satisfies~\eqref{eq:pinn_inverse_pde}. Under the bounded-noise model above, for $(u,h)\in H^2(\Omega)\times H^1(\Omega)$ and $\lambda_{\mathrm{reg}}>0$,
\begin{equation}
\|u-u_0\|_{H^1(\Omega)} \lesssim \bigl[\mathcal L_{\mathrm{PINN}}^{\lambda_{\mathrm{reg}}}(u,h)\bigr]^{1/2}+\delta,
\label{eq:pinn_inverse_state_stability}
\end{equation}
and
\begin{equation}
\begin{aligned}
\|h-h_0\|_{L^2(\Omega)} \lesssim{}\left(1+\lambda_{\mathrm{reg}}^{-1/4}\bigl[\mathcal L_{\mathrm{PINN}}^{\lambda_{\mathrm{reg}}}(u,h)\bigr]^{1/4}\right)\cdot\left(\bigl[\mathcal L_{\mathrm{PINN}}^{\lambda_{\mathrm{reg}}}(u,h)\bigr]^{1/4}+\delta^{1/2}\right).
\end{aligned}
\label{eq:pinn_inverse_source_stability}
\end{equation}
The implicit constants depend only on $\Omega$, $d$, and $\|h_0\|_{H^1(\Omega)}$.
\end{proposition}

The source estimate follows by interpolating between weak-norm control from the state and physics residuals and the $H^1(\Omega)$ bound supplied by Tikhonov regularization.

We next pass to empirical risk minimization. Let $\{(\boldsymbol X_i,Y_i,\boldsymbol Z_i)\}_{i=1}^{N_{\mathrm{in}}}$ be i.i.d.\ copies of \eqref{eq:pinn_inverse_data}, and let $\{\bar{\boldsymbol X}_j\}_{j=1}^{N_b}$ be i.i.d.\ uniform samples on $\partial\Omega$, independent of the interior observations. The empirical version of \eqref{eq:pinn_inv_risk} replaces the population integrals by sample averages:
\begin{align}
\widehat{\mathcal L}_{\mathrm{PINN}}^{\lambda_{\mathrm{reg}}}(u,h)
\;=\;& \frac{|\Omega|}{N_{\mathrm{in}}}\sum_{i=1}^{N_{\mathrm{in}}}\Big[(Y_i-u(\boldsymbol X_i))^2+\big\|\boldsymbol Z_i-\nabla u(\boldsymbol X_i)\big\|_2^2\Big] \nonumber\\
&+\frac{|\Omega|}{N_{\mathrm{in}}}\sum_{i=1}^{N_{\mathrm{in}}}\big(\Delta u(\boldsymbol X_i)-u(\boldsymbol X_i)+h(\boldsymbol X_i)\big)^2
+\frac{|\partial\Omega|}{N_b}\sum_{j=1}^{N_b}\partial_{\mathbf n}u(\bar{\boldsymbol X}_j)^2 \nonumber\\
&+\lambda_{\mathrm{reg}}\frac{|\Omega|}{N_{\mathrm{in}}}\sum_{i=1}^{N_{\mathrm{in}}}\Big(h(\boldsymbol X_i)^2+\big\|\nabla h(\boldsymbol X_i)\big\|_2^2\Big).
\label{eq:pinn_inverse_empirical}
\end{align}

For this inverse problem, let $\mathcal F_{NN}$ and $\mathcal G_{NN}$ denote candidate neural network classes for the state $u$ and source $h$, respectively. Define the joint empirical risk minimizer by
\[
(\widehat u,\widehat h)\in\argmin_{(u,h)\in\mathcal F_{NN}\times\mathcal G_{NN}}
\widehat{\mathcal L}_{\mathrm{PINN}}^{\lambda_{\mathrm{reg}}}(u,h).
\]
For the following oracle inequality, take $N_{\mathrm{in}}=N_b$, as in \citealp[Theorem~3.1]{Duan2024Recovering}. It bounds the population risk by approximation and finite-sample errors, regularization bias, and observation noise.

\begin{proposition}[{Oracle inequality; adapted from \citealp[Theorem~3.1]{Duan2024Recovering}}]
\label{proposition:pinns:inverse:oracle}
Suppose that the functions in $\mathcal F_{NN}\cup\{u_0\}$ and $\mathcal G_{NN}\cup\{h_0\}$, together with the derivatives appearing in the loss, are uniformly bounded, with envelope bounds $\mathcal B_{\mathcal F},\mathcal B_{\mathcal G}$ as in \citealp[Section~3.1]{Duan2024Recovering}. For $\lambda_{\mathrm{reg}}\in(0,1)$ and $\tau\in(0,1)$, with probability at least $1-(13+7d)\tau$,
\begin{align}
\mathcal L_{\mathrm{PINN}}^{\lambda_{\mathrm{reg}}} (\widehat u,\widehat h) &\lesssim{} \underbrace{ \inf_{u\in\mathcal F_{NN}} \|u-u_0\|_{H^2(\Omega)}^2 + \inf_{h\in\mathcal G_{NN}} \|h-h_0\|_{H^1(\Omega)}^2 }_{\mathcal E_{\mathrm{app}}^{\mathrm{PINN},\mathrm{inv}}} \nonumber\\
&+ \underbrace{ \mathcal E_{\mathrm{gen}}^{\mathrm{PINN}} ( \mathcal F_{NN},\mathcal G_{NN}, N_{\mathrm{in}},N_b ) }_{\text{finite-sample error}} + \underbrace{ C\lambda_{\mathrm{reg}} \|h_0\|_{H^1(\Omega)}^2 }_{\text{regularization bias}} + \underbrace{ C\delta^2 }_{\text{observation-noise term}}.
\label{eq:pinn_oracle}
\end{align}
\end{proposition}

The data term involves $\nabla u$, the interior and boundary residuals involve $\Delta u$ and $\partial_{\mathbf n}u$, respectively, and the source penalty involves $\nabla h$. The finite-sample term in \eqref{eq:pinn_oracle} therefore includes empirical Rademacher complexities for the state and source classes and their derivative classes, together with a concentration term proportional to $\sqrt{\log(1/\tau)/N_{\mathrm{in}}}$. Bounds for these complexities are given in~\cite[Lemma~3.2 and Corollary~3.1]{Duan2024Recovering} using parameter covering arguments of the type discussed in Section~\ref{sec:classical_stat}.

To quantify the four terms in \eqref{eq:pinn_oracle}, suppose in addition that
\[
(u_0,h_0)\in H^{\alpha+2}(\Omega)\times H^{\alpha+1}(\Omega),\qquad \alpha\geq1.
\]
For noise level $0<\delta<1$, Duan et al.~\cite[Theorem~3.2]{Duan2024Recovering} take the approximation accuracy $\varepsilon=\delta$ and use bounded-depth tanh classes for the state and source, with $O(\delta^{-d/(\alpha-\mu)})$ nonzero parameters per network, where $0<\mu<\alpha$ is a sufficiently small approximation slack. Their construction approximates $u_0$ in $H^2(\Omega)$ and $h_0$ in $H^1(\Omega)$ to accuracy $O(\delta)$ using Sobolev approximation of the type in Theorem~\ref{thm:A2}. The squared approximation errors therefore contribute $O(\delta^2)$, as does the observation noise.

To control the finite-sample term, take
\[
N_{\mathrm{in}}=N_b\gtrsim
\delta^{-4-c\frac{d+1}{\alpha-\mu}\log(d+2)},
\]
where $c>0$ is an absolute constant. With this sample size, the derivative-class complexity and concentration estimates give the terms containing $\sqrt{\log(1/\delta)}$ and $\sqrt{\log(1/\tau)}$, respectively. For fixed $h_0$, the regularization term is $O(\lambda_{\mathrm{reg}})$. The oracle inequality thus gives, for each fixed $\lambda_{\mathrm{reg}}\in(0,1)$ and $0<\tau<(13+7d)^{-1}$, with probability at least $1-(13+7d)\tau$,
\begin{equation}\label{eq:pinn_inverse_risk_scale}
\mathcal L_{\mathrm{PINN}}^{\lambda_{\mathrm{reg}}}(\widehat u,\widehat h)
\lesssim
\underbrace{\delta^2\sqrt{\log(1/\delta)}+\delta^2
+\delta^{\frac{d+4(\alpha-\mu)}{2(\alpha-\mu)}}\sqrt{\log(1/\tau)}}_{B_{\delta,\tau}}
+\lambda_{\mathrm{reg}}.
\end{equation}
Here $B_{\delta,\tau}$ collects the approximation, finite-sample, and noise contributions.
Fix $c_\lambda>0$ and choose $\lambda_{\mathrm{reg}}=c_\lambda B_{\delta,\tau}$, taking $\delta$ sufficiently small that $\lambda_{\mathrm{reg}}<1$. On the same event,
\[
\mathcal L_{\mathrm{PINN}}^{\lambda_{\mathrm{reg}}}(\widehat u,\widehat h)\lesssim B_{\delta,\tau},
\qquad
\lambda_{\mathrm{reg}}^{-1/4}
\bigl[\mathcal L_{\mathrm{PINN}}^{\lambda_{\mathrm{reg}}}(\widehat u,\widehat h)\bigr]^{1/4}\lesssim1.
\]
The second bound controls the amplification factor in the source stability estimate. Since $B_{\delta,\tau}\geq\delta^2$, we also have $\delta\leq B_{\delta,\tau}^{1/2}$ and $\delta^{1/2}\leq B_{\delta,\tau}^{1/4}$. Substitution into \eqref{eq:pinn_inverse_state_stability}--\eqref{eq:pinn_inverse_source_stability} then gives
\[
\|\widehat u-u_0\|_{H^1(\Omega)}\lesssim B_{\delta,\tau}^{1/2},
\qquad
\|\widehat h-h_0\|_{L^2(\Omega)}\lesssim B_{\delta,\tau}^{1/4}.
\]

Expanding these powers and absorbing $\delta$ and $\delta^{1/2}$ into the respective leading terms gives
\begin{align*}
\big\|\widehat u-u_0\big\|_{H^1(\Omega)}
&\lesssim \delta\log^{1/4}(1/\delta)+\delta^{(d+4(\alpha-\mu))/(4(\alpha-\mu))}\log^{1/4}(1/\tau), \\
\big\|\widehat h-h_0\big\|_{L^2(\Omega)}
&\lesssim \delta^{1/2}\log^{1/8}(1/\delta)+\delta^{(d+4(\alpha-\mu))/(8(\alpha-\mu))}\log^{1/8}(1/\tau).
\end{align*}
For fixed confidence, $\lambda_{\mathrm{reg}}\asymp B_{\delta,\tau}\asymp\delta^2\sqrt{\log(1/\delta)}$. The slower recovery of the source reflects the weaker inverse stability estimate.

For the forward elliptic problem, the results reviewed above include approximation--statistical bounds, norm-controlled analysis, and an analysis incorporating PSGD optimization. For the inverse source problem, the cited result assumes empirical risk minimization and combines approximation and statistical errors with regularization bias and observation noise. The inverse stability estimate converts this risk bound into recovery guarantees for the state and source.

\subsection{Generative learning}\label{section:app:generative}

Generative learning aims to construct a sampling mechanism whose output distribution approximates an unknown data distribution $p_0$. Generative adversarial networks (GANs) learn a map from a simple source distribution by comparing generated and observed samples through a discriminator~\cite{Goodfellow2014Generative}. Diffusion models learn to reverse a gradual noising process~\cite{song2021score}, while flow matching learns a velocity field that transports a source distribution toward the data distribution~\cite{lipman2023flow}.

\subsubsection{Generative adversarial networks}\label{section:app:generative:gan}

\par The original GAN trains a discriminator using a logistic classification objective~\cite{Goodfellow2014Generative}. Here we study GANs that compare the target and generated distributions using an integral probability metric (IPM)~\cite{Huang2022Error}. The IPM compares the expectations of each test function under the two distributions and takes the supremum of the resulting absolute differences over a prescribed class of functions.

\par Let $p_0$ be an unknown target distribution on $\Omega \subseteq \mathbb{R}^d$, and let $p_z$ be an easy-to-sample source distribution on $\mathbb{R}^{d_Z}$. A generator $g : \mathbb{R}^{d_Z} \to \Omega$ transforms a latent draw $\boldsymbol Z\sim p_z$ into a generated sample $g(\boldsymbol Z)$, whose distribution is the pushforward $g_\sharp p_z$. The population objective seeks a generator for which $g_\sharp p_z$ is close to $p_0$ according to the most discriminating test function in the training discriminator class $\mathcal F$:
\begin{equation*}
\min_{g \in \mathcal{G}}\, \max_{f \in \mathcal{F}}\, \mathbb{E}_{\boldsymbol{X} \sim p_0}[f(\boldsymbol{X})] - \mathbb{E}_{\boldsymbol{Z} \sim p_z}[f \circ g(\boldsymbol{Z})],
\end{equation*}
where the generator class $\mathcal{G}$ and training discriminator class $\mathcal{F}$ are both neural network classes. If $\mathcal F$ is symmetric, it contains $-f$ whenever it contains $f$. Replacing $f$ by $-f$ reverses the sign of the expectation difference, so taking the supremum of the signed difference gives the same IPM
\begin{equation}\label{app:gan:ipm}
d_{\mathcal{F}}(p_0, g_\sharp p_z) \coloneqq \sup_{f \in \mathcal{F}} \bigl(\mathbb{E}_{p_0}[f] - \mathbb{E}_{g_\sharp p_z}[f]\bigr)
\end{equation}
between the target distribution $p_0$ and the generated distribution $g_\sharp p_z$. Given i.i.d.\ samples $\{\boldsymbol{X}_i\}_{i=1}^n \sim p_0$ and $\{\boldsymbol{Z}_j\}_{j=1}^m \sim p_z$, replace both population expectations by sample averages. The resulting empirical minimax estimator is
\begin{equation}\label{app:gan:erm}
\widehat{g}_{n, m} \in \arg\min_{g \in \mathcal{G}} \sup_{f \in \mathcal{F}}\, \frac{1}{n}\sum_{i=1}^n f(\boldsymbol{X}_i) - \frac{1}{m}\sum_{j=1}^m f \circ g(\boldsymbol{Z}_j),
\end{equation}
which minimizes the empirical IPM $d_{\mathcal{F}}(\widehat{p}_{0,n}, g_\sharp \widehat{p}_{z,m})$ between the empirical data distribution $\widehat p_{0,n}$ and the empirical generated distribution $g_\sharp\widehat p_{z,m}$. The training discriminator class $\mathcal F$ need not coincide with the class used to evaluate the learned distribution. We therefore measure the final distributional error with a benchmark test-function class $\mathcal H$, typically the H\"older class $\mathcal{H}^\beta([0,1]^d)$ for some smoothness $\beta$:
\begin{equation*}
d_{\mathcal{H}}(p_0, (\widehat{g}_{n,m})_\sharp p_z) \coloneqq \sup_{h \in \mathcal{H}} \mathbb{E}_{p_0}[h] - \mathbb{E}_{(\widehat{g}_{n,m})_\sharp p_z}[h].
\end{equation*}
Throughout the GAN analysis, $\mathcal H^\beta(\Omega)$ abbreviates the unit H\"older ball $\mathcal H^\beta(\Omega,1)$.

\par For the GAN empirical minimizer, the error in the benchmark IPM admits a four-term decomposition: two approximation errors and two statistical errors.

\begin{lemma}[GAN oracle inequality~{\cite[Lemma 9]{Huang2022Error}}]\label{lem:gan_oracle}
Suppose $\mathcal{F}$ is symmetric and that $p_0, g_\sharp p_z$ are supported on $\Omega \subseteq \mathbb{R}^d$. Then for any function class $\mathcal{H}$ on $\Omega$,
\begin{equation}\label{eq:gan_decomp}
\begin{aligned}
d_{\mathcal{H}}(p_0, (\widehat{g}_{n,m})_\sharp p_z)
\;\le\;& 2\underbrace{\mathcal{E}(\mathcal{H}, \mathcal{F}, \Omega)}_{\text{disc.\ approx.}}
\;+\; \underbrace{\inf_{g \in \mathcal{G}} d_{\mathcal{F}}(\widehat{p}_{0,n}, g_\sharp p_z)}_{\text{gen.\ approx.}} \\
&\;+\; \underbrace{d_{\mathcal{F}}(p_0, \widehat{p}_{0,n}) \wedge d_{\mathcal{H}}(p_0, \widehat{p}_{0,n})}_{\text{statistical (i)}}
\;+\; 2\underbrace{d_{\mathcal{F} \circ \mathcal{G}}(p_z, \widehat{p}_{z,m})}_{\text{statistical (ii)}},
\end{aligned}
\end{equation}
where the discriminator approximation error is
\[
\mathcal{E}(\mathcal{H}, \mathcal{F}, \Omega) \coloneqq \sup_{h \in \mathcal{H}}\, \inf_{f \in \mathcal{F}}\, \|h - f\|_{L^\infty(\Omega)}.
\]
\end{lemma}

The four terms have distinct roles: \emph{generator approximation} controls how well $\mathcal{G}$ can push $p_z$ to match the empirical $\widehat{p}_{0,n}$; \emph{discriminator approximation} controls how well $\mathcal{F}$ can imitate the benchmark class $\mathcal{H}$ in $L^\infty$; \emph{statistical (i)} measures the empirical-distribution fluctuation on the $p_0$-side (crucially, the minimum $\wedge$ allows the benchmark class $\mathcal{H}$ to control this data fluctuation even when the discriminator class $\mathcal{F}$ is large); \emph{statistical (ii)} measures it on the $p_z$-side, jointly over the composition class $\mathcal{F} \circ \mathcal{G}$. The minimum follows from two valid comparisons, inserting the empirical target before or after passing from the benchmark class to the discriminator class. The source coefficient two arises from replacing the empirical source by $p_z$ once for the fitted generator and once for the comparison generator. The lemma assumes exact empirical minimization. If the empirical IPM attained by the generator exceeds its infimum over $\mathcal G$ by at most $\epsilon_{\mathrm{opt}}$, then \eqref{eq:gan_decomp} gains the additive term $\epsilon_{\mathrm{opt}}$. This tolerance bounds the empirical optimization error $\mathcal E_{\mathrm{opt}}$ of Section~\ref{sec:decomp}.

For the generator side, a sufficiently large generator class can approximate the empirical distribution $\widehat{p}_{0,n}$ arbitrarily well in Wasserstein distance.

\begin{lemma}[Generator capacity; \citealp{Yang2022capacity}, Lemma 3.2]\label{lem:gan_gen}
Let $p_z$ be an absolutely continuous probability measure on $\mathbb{R}$ and let $\mathrm{supp}(\widehat{p}_{0,n}) \subseteq [0,1]^d$. Set $\mathcal{G} = \{g \in \mathcal{NN}_{\rm ReLU}(W_{\mathcal{G}}, L_{\mathcal{G}}) : g(\mathbb{R}) \subseteq [0,1]^d\}$ with $W_{\mathcal{G}} \ge 7d + 1$, $L_{\mathcal{G}} \ge 2$, and $W_{\mathcal{G}}^2 L_{\mathcal{G}} \gtrsim n$. Then for any $\mathcal{F} \subseteq \mathrm{Lip}([0,1]^d, \mathcal{K}_{\mathcal{F}})$,
\[
\inf_{g \in \mathcal{G}} d_{\mathcal{F}}(\widehat{p}_{0,n}, g_\sharp p_z) \;\le\; \mathcal{K}_{\mathcal{F}}\, \inf_{g \in \mathcal{G}} W_1(\widehat{p}_{0,n}, g_\sharp p_z) \;=\; 0.
\]
\end{lemma}

The lemma's role is to remove generator approximation as a source of error provided the generator class is large enough; the price is the requirement that $\mathcal{F}$ be uniformly Lipschitz with constant $\mathcal{K}_{\mathcal{F}}$, which becomes a parameter to track in the discriminator analysis. The zero infimum means that the empirical atomic law can be approximated arbitrarily closely; it does not require a generator that realizes this law exactly.

Discriminator approximation requires control of both the $L^\infty$ error for each $h\in\mathcal H^\beta$ and the Lipschitz constant used in Lemma~\ref{lem:gan_gen}.

\begin{lemma}[Discriminator capacity under Lipschitz constraint; \citealp{Huang2022Error}, Lemma 11]\label{lem:gan_disc}
For every $h\in\mathcal H^\beta([0,1]^d)$ and sufficiently large width and depth parameters $W_{\mathcal F},L_{\mathcal F}$, the cited construction gives a bounded ReLU network $f_h$ such that
\[
\|h-f_h\|_{L^\infty([0,1]^d)}
\lesssim
\left(\frac{W_{\mathcal F}L_{\mathcal F}}
{\log W_{\mathcal F}\log L_{\mathcal F}}\right)^{-2\beta/d},
\qquad
\operatorname{Lip}(f_h)\leq \mathcal K_{\mathcal F}^{\mathrm{adm}},
\]
where $\mathcal K_{\mathcal F}^{\mathrm{adm}}=\mathcal K_{\mathcal F}^{\mathrm{adm}}(W_{\mathcal F},L_{\mathcal F},d,\beta)$ is the explicit upper bound supplied by Lemma~11 of the reference. Consequently, any discriminator class containing these constructed networks obeys the displayed approximation bound.
\end{lemma}

The construction controls both the approximation error and the Lipschitz constant, extending the size-based approximation analysis of Section~\ref{sec:classical_approx}.

The two IPM-statistical terms in \eqref{eq:gan_decomp} are bounded using the symmetrization and entropy tools of Section~\ref{sec:classical_stat}, applied to two function classes: $\mathcal{H}^\beta([0,1]^d)$ and $\mathcal{F} \circ \mathcal{G}$.

\begin{lemma}[Statistical IPM bounds~{\cite[Theorem 5 and Corollary~35]{Huang2022Error}}]\label{lem:gan_stat}
Let $\mathcal{H} = \mathcal{H}^\beta([0,1]^d)$, and suppose that $|f|\leq B_{\mathcal F}$ uniformly on the target support and on all generator outputs. If the relevant pseudo-dimensions are finite, then
\begin{align*}
\mathbb{E}\big[d_{\mathcal{H}}(p_0, \widehat{p}_{0,n})\big]
&\lesssim_{d,\beta} n^{-\beta/d}\vee
n^{-1/2}(\log n)^{\mathbf 1_{\{2\beta=d\}}}, \\[1ex]
\mathbb{E}\big[d_{\mathcal{F}}(p_0, \widehat{p}_{0,n})\big]
&\lesssim B_{\mathcal F}
\sqrt{\frac{\operatorname{Pdim}(\mathcal F)\log n}{n}}, \\[1ex]
\mathbb{E}\big[d_{\mathcal{F} \circ \mathcal{G}}(p_z, \widehat{p}_{z,m})\big]
&\lesssim B_{\mathcal F}
\sqrt{\frac{\operatorname{Pdim}(\mathcal F\circ\mathcal G)\log m}{m}}.
\end{align*}
\end{lemma}

The first bound uses the Hölder-class metric entropy $\log \mathcal{C}(\epsilon, \mathcal{H}^\beta) \lesssim \epsilon^{-d/\beta}$ \cite{Kolmogorov1961entropy}; the second and third use the pseudo-dimension bound for ReLU networks (Proposition~\ref{prop:relu_pdim}, due to \citealp{bartlett2019nearly}). The latter bounds depend on the output envelope and on pseudo-dimension through parameter count, depth, and logarithmic factors.

Combining the preceding estimates in \eqref{eq:gan_decomp} yields the following non-asymptotic rate.

\begin{theorem}[GAN convergence; \citealp{Huang2022Error}, Theorem 5]\label{thm:gan_main}
Suppose $p_z$ is absolutely continuous on $\mathbb{R}$ and $p_0$ is supported on $[0,1]^d$. Let $\mathcal{H} = \mathcal{H}^\beta([0,1]^d)$. Take $\mathcal G$ to be the ReLU generator class with outputs in $[0,1]^d$, and take
\[
\mathcal F=\left\{f\in\mathcal{NN}_{\mathrm{ReLU}}(W_{\mathcal F},L_{\mathcal F}):
\|f\|_{L^\infty(\mathbb R^d)}\leq1,\quad
\operatorname{Lip}(f)\leq\mathcal K_{\mathcal F}\right\}.
\]
Choose the widths, depths, and Lipschitz radius according to the construction in the reference, satisfying
\begin{align*}
W_{\mathcal{G}}^2 L_{\mathcal{G}} &\asymp n, &
W_{\mathcal{F}} L_{\mathcal{F}} &\asymp n^{1/2} \log^2 n, &
\mathcal{K}_{\mathcal{F}} &\asymp \mathcal K_{\mathcal F}^{\mathrm{adm}}.
\end{align*}
If $m \gtrsim n^{2 + 2\beta/d} \log^6 n$, then the GAN ERM \eqref{app:gan:erm} satisfies
\[
\mathbb{E}\big[d_{\mathcal{H}}(p_0, (\widehat{g}_{n,m})_\sharp p_z)\big]
\;\lesssim\; n^{-\beta/d}\vee n^{-1/2}(\log n)^{\mathbf 1_{\{2\beta=d\}}}.
\]
\end{theorem}

The displayed rate recovers the minimax-optimal rate, up to the critical logarithm, for IPM estimation of compactly supported distributions in $d$ dimensions \cite{Liang2021How}, certifying that the GAN ERM is statistically efficient in this regime. Here $\beta$ measures test-function regularity, not density smoothness. An absolutely continuous higher-dimensional source can also be accommodated by projecting onto one coordinate \cite[Remark~7]{Huang2022Error}.

When the training objective uses the exact source expectation, namely $d_{\mathcal F}(\widehat p_{0,n},g_\sharp p_z)$, the source-sampling term in \eqref{eq:gan_decomp} is absent. The same generator and discriminator choices then give the displayed rate in $n$, plus $\epsilon_{\mathrm{opt}}$ for approximate empirical minimization, without a condition on $m$.

Two extensions concern low-dimensional target distributions and norm-controlled discriminators.

\emph{(i) Low-intrinsic-dimensional target.} Following Assumption~18 of \cite{Huang2022Error}, let
\[
\boldsymbol X=\widetilde{\boldsymbol X}+\boldsymbol\xi\sim p_0,
\qquad \mathbb E\boldsymbol\xi=\boldsymbol0,
\qquad V_{\mathrm{noise}}:=\mathbb E\|\boldsymbol\xi\|_2^2<\infty,
\]
where $\widetilde{\boldsymbol X}$ and $\boldsymbol\xi$ are independent, and $\widetilde{\boldsymbol X}$ is supported on a compact set $\mathcal X\subseteq[0,1]^d$ with $\dim_{\mathrm M}(\mathcal X)=d^*>0$. To make the exact covering exponent explicit, we additionally assume $\mathcal C(r,\mathcal X,\|\cdot\|_2)\lesssim r^{-d^*}$ for $0<r\le1$, as discussed in Section~\ref{sec:low_intrinsic}. Evaluate the noisy target using the bounded class $\mathcal H=\mathcal H^\beta(\mathbb R^d,1)$, $\beta>0$.

For an absolutely continuous source $p_z$ on $\mathbb R$, use the generator and bounded Lipschitz discriminator classes constructed in \cite[Theorem~19]{Huang2022Error}. If
\[
m\gtrsim
\begin{cases}
n^{(3d+4\beta)/(2d^*)}(\log n)^6,& d^*\le d/2,\\
n^{1+(d+2\beta)/d^*}(\log n)^4,& d^*>d/2,
\end{cases}
\]
the approximate empirical GAN estimator, with objective tolerance $\epsilon_{\mathrm{opt}}$, satisfies for $n\ge2$
\begin{equation}\label{eq:gan_intrinsic_noise}
\mathbb E\,d_{\mathcal H}\bigl(p_0,(\widehat g_{n,m})_\sharp p_z\bigr)
\le\epsilon_{\mathrm{opt}}
+2\sqrt d\,V_{\mathrm{noise}}^{(\beta\wedge1)/2}
+C\bigl(n^{-\beta/d^*}\vee n^{-1/2}\bigr)\log n.
\end{equation}
The low-dimensional entropy bound controls the statistical term. Two applications of H\"older continuity, to the population and empirical noise perturbations, produce the noise term. When $V_{\mathrm{noise}}=0$, only the intrinsic-dimensional statistical rate and optimization residual remain. Constants and architecture sizes can still depend on $d$. Related manifold results appear in \cite{Dahal2022Deep}.

\emph{(ii) Over-parameterized discriminator} \cite{jiao2023approximation,golowich2018size}. Rather than parameterizing the discriminator by size, one can parameterize it by the weight-norm budget $\mathcal K$ of Section~\ref{sec:overparam_setting}. This budget constrains the network weights and supplies a Lipschitz bound, whereas the radius $\mathcal K_{\mathcal F}$ used above constrains the discriminator functions directly. Take $\mathcal F=\mathcal{NN}_{\mathrm{ReLU}}(W,L,\mathcal K)$, with width and depth satisfying the conditions of Theorem~\ref{thm:A3}. Its norm-controlled approximation bound gives
\[
\mathcal E(\mathcal H^\beta,\mathcal F,[0,1]^d)
\lesssim_{d,\beta}\mathcal K^{-\beta/(d+1)}.
\]
Choosing $\mathcal K\asymp n^{(d+1)/d}$ makes this approximation term no larger than $n^{-\beta/d}$. The resulting GAN rate is obtained through the benchmark branch $d_{\mathcal H}(p_0,\widehat p_{0,n})$ in the minimum of \eqref{eq:gan_decomp}, together with vanishing generator approximation. This yields the same benchmark rate as Theorem~\ref{thm:gan_main} for population-source GANs while allowing large discriminator widths. Finite latent sampling still requires separate control of $\mathcal F\circ\mathcal G$ and a sufficiently large $m$; see \citealp[Remark~4.9]{jiao2023approximation}.

\subsubsection{Diffusion models}\label{section:app:generative:diffusion}

Diffusion models perturb data with progressively stronger noise and generate new samples by approximately reversing this process. In the score-based formulation, the reverse dynamics depend on the score, the gradient of the log density of the perturbed data, which is estimated by a neural network~\cite{song2021score}. For a target distribution $p_0$ on $\mathbb R^d$, we consider the Ornstein--Uhlenbeck forward process:
\begin{equation}\label{eq:forward:base}
\mathrm{d}\boldsymbol X_t=-\boldsymbol X_t\,\mathrm{d}t+\sqrt{2}\,\mathrm{d}\boldsymbol B_t, \quad t\in(0,T_{\mathrm{end}}),~\boldsymbol X_0\sim p_0,
\end{equation}
where $\boldsymbol B_t$ is a $d$-dimensional standard Brownian motion, and $T_{\mathrm{end}}>0$ is the terminal time. The transition distribution of the forward process can be expressed as:
\begin{equation}\label{eq:forward:solution}
(\boldsymbol X_t\mid \boldsymbol X_0=\boldsymbol x_0)\sim \mathcal{N}(\mu_t\boldsymbol x_0,\sigma_t^2\boldsymbol I_d),
\end{equation}
where the mean and variance coefficients are given, respectively, as $\mu_t\coloneqq\exp(-t)$ and $\sigma_t^2\coloneqq 1-\exp(-2t)$. Here $\boldsymbol I_d$ denotes the $d\times d$ identity matrix. The forward process~\eqref{eq:forward:base} is commonly referred to as the variance-preserving (VP) SDE as $\mu_t^2+\sigma_t^2=1$ for each $t\in(0,T_{\mathrm{end}})$. Denote by $p_t$ the marginal density of $\boldsymbol X_t$ for each $t\in(0,T_{\mathrm{end}})$, which satisfies 
\begin{equation}\label{eq:marginal:density}
p_t(\boldsymbol x_t)=\int\varphi_d(\boldsymbol x_t;\mu_t\boldsymbol x_0,\sigma_t^2\boldsymbol I_d)p_0(\mathrm d\boldsymbol x_0),
\end{equation}
where $\varphi_d(\cdot;\mu_t\boldsymbol x_0,\sigma_t^2\boldsymbol I_d)$ denotes the density function of the Gaussian distribution $\mathcal{N}(\mu_t\boldsymbol x_0,\sigma_t^2\boldsymbol I_d)$. Let $\tau=T_{\mathrm{end}}-t$ denote reverse time and define $\boldsymbol Y_\tau=\boldsymbol X_{T_{\mathrm{end}}-\tau}$. By the time-reversal formula~\cite{anderson1982reverse}, this process satisfies 
\begin{equation}\label{eq:base:reversal}
\begin{aligned}
\mathrm{d}\boldsymbol Y_\tau&=\bigl(\boldsymbol Y_\tau+2\nabla\log p_{T_{\mathrm{end}}-\tau}(\boldsymbol Y_\tau)\bigr)\mathrm{d}\tau+\sqrt{2}\mathrm{d}\overline{\boldsymbol B}_\tau, \quad \tau\in(0,T_{\mathrm{end}}), \\
\boldsymbol Y_0&\sim p_{T_{\mathrm{end}}}.
\end{aligned}
\end{equation}
Here $\overline{\boldsymbol B}_\tau$ is a standard Brownian motion in the reverse-time filtration. At reverse time $\tau$, the marginal law of $\boldsymbol Y_\tau$ is $p_{T_{\mathrm{end}}-\tau}$.

\par Consequently, sampling from $p_0$ reduces to two ingredients: an initialization that is close to $p_{T_{\mathrm{end}}}$, and access to the score functions $\{\nabla\log p_t\}_{t\in(0,T_{\mathrm{end}})}$. Gaussian initialization is motivated by the fact that the Ornstein--Uhlenbeck process converges exponentially fast to its stationary distribution $\gamma\coloneqq\mathcal{N}(\boldsymbol 0,\boldsymbol I_d)$, so that $p_{T_{\mathrm{end}}}\approx\gamma$ for moderately large $T_{\mathrm{end}}$. The second ingredient is where learning enters: the score is unknown and must be estimated from the samples of $p_0$. The following classical identity, usually attributed to Tweedie, expresses the score as a denoising conditional expectation and underlies denoising score matching~\cite{Vincent2011DenoisingScoreMatching}.

\begin{proposition}[Tweedie's formula]\label{proposition:tweedie}
For each $t\in(0,T_{\mathrm{end}})$,
\begin{equation}\label{eq:tweedie}
\nabla\log p_t(\boldsymbol x)=-\frac{1}{\sigma_t^2}\boldsymbol x+\frac{\mu_t}{\sigma_t^2}\mathbb{E}[\boldsymbol X_0 \mid \boldsymbol X_t=\boldsymbol x], \quad \boldsymbol x\in\mathbb{R}^d.
\end{equation}
\end{proposition}

\par At each fixed $t$, we can generate a noisy training input $\boldsymbol X_t=\mu_t\boldsymbol X_0+\sigma_t\boldsymbol\varepsilon$ using an independent draw $\boldsymbol\varepsilon\sim\gamma$. The added noise is known during training. By~\eqref{eq:tweedie}, its least-squares predictor given $\boldsymbol X_t=\boldsymbol x$ is $-\sigma_t\nabla\log p_t(\boldsymbol x)$. A candidate score $s$ thus predicts the noise by $-\sigma_t s(t,\boldsymbol x)$, giving the denoising loss below.

\begin{corollary}\label{corollary:dsm}
Let $s^*(t,\boldsymbol x)\coloneqq\nabla\log p_t(\boldsymbol x)$ be the score. The score function minimizes the following population risk, uniquely up to equality almost everywhere with respect to $\mathrm dt\,p_t(\mathrm d\boldsymbol x)$:
\begin{equation}\label{eq:population:sm}
\mathcal{L}_{\mathrm{DSM}}(s) \coloneqq \int_0^{T_{\mathrm{end}}}\mathbb{E}_{(\boldsymbol X_0,\boldsymbol\varepsilon)\sim p_0\otimes\gamma}\bigl[\|\sigma_ts(t,\mu_t\boldsymbol X_0+\sigma_t\boldsymbol\varepsilon)+\boldsymbol\varepsilon\|_2^2\bigr]\,\mathrm{d}t.
\end{equation}
Further, for any measurable function $s\colon (0,T_{\mathrm{end}})\times\mathbb{R}^d\to\mathbb{R}^d$ with finite denoising risk, we have
\begin{equation}\label{eq:excess:sm}
\mathcal{L}_{\mathrm{DSM}}(s)-\mathcal{L}_{\mathrm{DSM}}(s^*) = \int_0^{T_{\mathrm{end}}}\sigma_t^2\mathbb{E}_{\boldsymbol X_t\sim p_t}\bigl[\|s(t,\boldsymbol X_t)-s^*(t,\boldsymbol X_t)\|_2^2\bigr]\,\mathrm{d}t.
\end{equation}
\end{corollary}

\par In~\eqref{eq:excess:sm}, the factor $\sigma_t^2$ comes from measuring squared error in the predicted noise. Thus the same squared score error contributes less to the training loss at times when less noise has been added. Sampling an independent time and noise draw for each data point gives the empirical risk
\begin{equation}\label{eq:empirical:sm}
\widehat{\mathcal{L}}_{\mathrm{DSM}}(s) \coloneqq \frac{T_{\mathrm{end}}}{n}\sum_{i=1}^n\|\sigma_{t^i}s(t^i,\mu_{t^i}\boldsymbol X_0^i+\sigma_{t^i}\boldsymbol\varepsilon^i)+\boldsymbol\varepsilon^i\|_2^2, 
\end{equation}
where $\{(t^i,\boldsymbol X_0^i,\boldsymbol\varepsilon^i)\}_{i=1}^n$ are i.i.d.\ from $\mathrm{unif}(0,T_{\mathrm{end}})\otimes p_0\otimes\gamma$. Let $\mathcal S$ be a score network class and let $\widehat s\in\mathcal S$ minimize the empirical risk over $\mathcal S$ up to a tolerance $\epsilon_{\mathrm{opt}}$. For finite risks and finite expected uniform deviation, the empirical-risk comparison of Section~\ref{sec:decomp} and identity~\eqref{eq:excess:sm} give
\begin{equation}\label{eq:structural:score:oracle}
\begin{aligned}
\mathcal E_{\mathrm{score}}^2 &
\coloneqq \mathbb E\int_0^{T_{\mathrm{end}}}\sigma_t^2
\|\widehat s(t,\cdot)-s^*(t,\cdot)\|_{L^2(p_t)}^2\,\mathrm dt\\
&\leq
\underbrace{\inf_{s\in\mathcal S}\int_0^{T_{\mathrm{end}}}\sigma_t^2
\|s(t,\cdot)-s^*(t,\cdot)\|_{L^2(p_t)}^2\,\mathrm dt}_{\mathcal E_{\mathrm{app}}^2}\\
&+2\underbrace{\mathbb E\sup_{s\in\mathcal S}
\big|\mathcal L_{\mathrm{DSM}}(s)-\widehat{\mathcal L}_{\mathrm{DSM}}(s)\big|}_{\mathcal E_{\mathrm{sta}}^2}
+\underbrace{\epsilon_{\mathrm{opt}}}_{\text{optimization tolerance}}.
\end{aligned}
\end{equation}
The approximation term measures how closely $\mathcal S$ can fit the true score in the weighted norm. The statistical term measures the fluctuation of the denoising loss over $\mathcal S$ as the data, times, and noise draws vary. Both estimates must account for the tails of the noisy inputs, and the statistical estimate also involves the Gaussian training response.

\par It remains to convert the score estimation error into a guarantee for the law of the generated samples. In practice one replaces the unknown score in~\eqref{eq:base:reversal} by $\widehat{s}$, initializes at the stationary distribution $\gamma$ rather than at $p_{T_{\mathrm{end}}}$, and stops at reverse time $\tau=T_{\mathrm{end}}-\delta$, corresponding to the forward time $t=\delta>0$. This avoids the singularity of the score as $t\downarrow0$. The resulting generative process satisfies
\begin{equation}\label{eq:approx:reversal}
\begin{aligned}
\mathrm{d}\widehat{\boldsymbol Y}_\tau&=\bigl(\widehat{\boldsymbol Y}_\tau+2\widehat{s}(T_{\mathrm{end}}-\tau,\widehat{\boldsymbol Y}_\tau)\bigr)\mathrm{d}\tau+\sqrt{2}\,\mathrm{d}\overline{\boldsymbol B}_\tau, \quad \tau\in(0,T_{\mathrm{end}}-\delta), \\
\widehat{\boldsymbol Y}_0&\sim\gamma.
\end{aligned}
\end{equation}
Fix the trained score $\widehat s$ and use a fresh Gaussian initial draw and Brownian motion in \eqref{eq:approx:reversal}, independent of the training sample. Write $\widehat p_\delta$ for the resulting conditional law of $\widehat{\boldsymbol Y}_{T_{\mathrm{end}}-\delta}$. The chain rule for Kullback--Leibler divergence and Girsanov's theorem on the reverse path space then give the following bound.

\begin{proposition}\label{proposition:kl:score}
Assume that $D_{\mathrm{KL}}(p_0\,\|\,\gamma)<\infty$ and that $\widehat{s}$ is regular enough for Girsanov's theorem to apply along~\eqref{eq:approx:reversal}. Then, for every $\delta\in(0,T_{\mathrm{end}})$,
\begin{equation}\label{eq:kl:score}
D_{\mathrm{KL}}(p_\delta\,\|\,\widehat{p}_\delta)\leq\underbrace{\exp(-2T_{\mathrm{end}})\,D_{\mathrm{KL}}(p_0\,\|\,\gamma)}_{\text{initialization}}+\underbrace{\int_\delta^{T_{\mathrm{end}}}\mathbb{E}_{\boldsymbol X_t\sim p_t}\bigl[\|\widehat{s}(t,\boldsymbol X_t)-s^*(t,\boldsymbol X_t)\|_2^2\bigr]\,\mathrm{d}t}_{\text{score estimation}}.
\end{equation}
\end{proposition}

\par The first term reflects the exponential convergence of the Ornstein--Uhlenbeck semigroup: since $\gamma$ satisfies a logarithmic Sobolev inequality with constant one, $D_{\mathrm{KL}}(p_{T_{\mathrm{end}}}\,\|\,\gamma)\leq\exp(-2T_{\mathrm{end}})D_{\mathrm{KL}}(p_0\,\|\,\gamma)$, so that a terminal time $T_{\mathrm{end}}\asymp\log(1/\epsilon)$ already makes the initialization error negligible. The second term is an unweighted score error. The lower bound on $\sigma_t^2$ relates it to the weighted training error: for $t\geq\delta$ we have $\sigma_t^2\geq\sigma_\delta^2=1-\exp(-2\delta)\geq\delta$ whenever $\delta\leq1/2$, whence
\begin{equation}\label{eq:weight:conversion}
\int_\delta^{T_{\mathrm{end}}}\mathbb{E}_{\boldsymbol X_t\sim p_t}\bigl[\|\widehat{s}(t,\boldsymbol X_t)-s^*(t,\boldsymbol X_t)\|_2^2\bigr]\mathrm{d}t\leq\frac{1}{\delta}\int_0^{T_{\mathrm{end}}}\sigma_t^2\|\widehat{s}(t,\cdot)-s^*(t,\cdot)\|_{L^2(p_t)}^2\,\mathrm{d}t.
\end{equation}
Writing $m_2=(\mathbb E_{\boldsymbol X_0\sim p_0}\|\boldsymbol X_0\|_2^2)^{1/2}<\infty$, the early-stopping bias can be bounded in Wasserstein distance: since $\boldsymbol X_\delta=\mu_\delta\boldsymbol X_0+\sigma_\delta\boldsymbol\varepsilon$, the coupling $(\boldsymbol X_0,\boldsymbol X_\delta)$ gives $W_2(p_\delta,p_0)\lesssim\delta m_2+\sqrt{\delta d}$.

\par While the Wasserstein distance is controlled by finite second moments, establishing total variation (TV) distance requires density regularity. A singular target can remain at TV distance one from its smoothed law. Conditional on the training sample, let $\bar p_\delta$ be the output law of the same learned reverse dynamics as \eqref{eq:approx:reversal}, but initialized from $p_{T_{\mathrm{end}}}$ rather than $\gamma$. For a numerical output law $\widehat p^{\mathrm{num}}$, define
\[
\mathcal E_{\mathrm{early}}^{\mathrm{TV}}:=d_{\mathrm{TV}}(p_0,p_\delta),\quad
\mathcal E_{\mathrm{init}}^{\mathrm{TV}}:=d_{\mathrm{TV}}(p_{T_{\mathrm{end}}},\gamma),\quad
\mathcal E_{\mathrm{disc}}^{\mathrm{TV}}:=d_{\mathrm{TV}}(\widehat p_\delta,\widehat p^{\mathrm{num}}).
\]
For well-posed reverse processes satisfying the Girsanov hypotheses, using the common initial law $p_{T_{\mathrm{end}}}$ and then Pinsker's inequality bounds $d_{\mathrm{TV}}(p_\delta,\bar p_\delta)$ by the square root of half the integrated score error. The laws $\bar p_\delta$ and $\widehat p_\delta$ use the same learned transition kernel, so TV contraction bounds their distance by $d_{\mathrm{TV}}(p_{T_{\mathrm{end}}},\gamma)$. Inserting these intermediate laws gives
\begin{equation}\label{eq:diffusion_tv_decomp}
\begin{aligned}
&d_{\mathrm{TV}}(p_0,\widehat p^{\mathrm{num}})\\
&\leq d_{\mathrm{TV}}(p_0,p_\delta)+d_{\mathrm{TV}}(p_\delta,\bar p_\delta) +d_{\mathrm{TV}}(\bar p_\delta,\widehat p_\delta)
+d_{\mathrm{TV}}(\widehat p_\delta,\widehat p^{\mathrm{num}})\\
&\leq\underbrace{\mathcal E_{\mathrm{early}}^{\mathrm{TV}}}_{\text{early stopping}} + \underbrace{\mathcal E_{\mathrm{init}}^{\mathrm{TV}}}_{\text{initialization}} + \underbrace{\left(\frac{1}{2}\int_\delta^{T_{\mathrm{end}}}\|\widehat s_t-s_t^*\|_{L^2(p_t)}^2\,\mathrm dt\right)^{1/2}}_{\text{score estimation}} + \underbrace{\mathcal E_{\mathrm{disc}}^{\mathrm{TV}}}_{\text{discretization}}.
\end{aligned}
\end{equation}
For $0<\delta<T_{\mathrm{end}}$ with $\delta\leq 1/2$, taking expectations over the training sample in~\eqref{eq:diffusion_tv_decomp} gives
\[
\begin{aligned}
\mathbb E\,d_{\mathrm{TV}}(p_0,\widehat p^{\mathrm{num}})
&\leq \mathcal E_{\mathrm{early}}^{\mathrm{TV}}
+\mathcal E_{\mathrm{init}}^{\mathrm{TV}}
+\mathbb E\,\mathcal E_{\mathrm{disc}}^{\mathrm{TV}}+\left(\frac{\mathcal E_{\mathrm{app}}^2+2\mathcal E_{\mathrm{sta}}^2
+\epsilon_{\mathrm{opt}}}{2\delta}\right)^{1/2}.
\end{aligned}
\]
The last term follows by applying Jensen's inequality to the score contribution, using~\eqref{eq:weight:conversion} to bound the unweighted score integral by the weighted one, and then applying the error decomposition~\eqref{eq:structural:score:oracle}.

Oko et al.~\cite[Theorems~3.1 and~4.3]{Oko2023Diffusion} give an explicit rate for the score term in the same TV decomposition~\eqref{eq:diffusion_tv_decomp}. They consider densities supported on $[-1,1]^d$, bounded above and away from zero there, and lying in a fixed-radius ball of $B^\beta_{p,q}([-1,1]^d)$, where $0<p,q\leq\infty$ and $\beta>d(1/p-1/2)_+$. Here $p$ and $q$ are the integrability and summability indices of the Besov space. On a strip inside the boundary of the cube, the density is additionally $C^\infty$ with uniformly bounded derivatives. Their denoising loss averages the squared norm of the score residual $s(t,\boldsymbol X_t)+\boldsymbol\varepsilon/\sigma_t$ over the added noise and integrates over $[\delta,T_{\mathrm{end}}]$ without multiplying the integrand by $\sigma_t^2$, so its population excess risk equals the unweighted integrated squared score error in the sampling bound. Combining sparse ReLU approximation with a statistical estimate for $n$ independent data points gives, for their exact empirical minimizer $\widehat s_n$,
\begin{equation}\label{eq:diffusion_unweighted_score_rate}
\mathbb E\int_\delta^{T_{\mathrm{end}}}
\|\widehat s_n(t,\cdot)-s^*(t,\cdot)\|_{L^2(p_t)}^2\,\mathrm dt
\lesssim n^{-2\beta/(2\beta+d)}(\log n)^{16}.
\end{equation}
The score contribution to TV is bounded in expectation by the square root of half the bound in~\eqref{eq:diffusion_unweighted_score_rate}. Their choice of sampling times bounds the early-stopping and initialization errors by $O(n^{-\beta/(2\beta+d)})$, and the continuous sampler has zero discretization error. For their truncated output law $\widehat p^{\mathrm{tr}}$, Oko et al.~\cite[Theorem~5.1]{Oko2023Diffusion} thus obtain
\[
\mathbb E\,d_{\mathrm{TV}}(p_0,\widehat p^{\mathrm{tr}})
\lesssim n^{-\beta/(2\beta+d)}(\log n)^8.
\]

The rate above describes the effect of sample size for a continuous reverse process. For a finite-step sampler, Li and Yan~\cite{liyan2025diffusion} measure TV error against the target after one forward noising step. Their bound retains a discretization term even when the score is exact. This raises the question of how fine a grid suffices to preserve the statistical accuracy. In conditional diffusion, Chang et al.~\cite[Section~3.1]{chang2026deep} choose the mesh to keep the discretization term within the error bound due to score estimation. Adding the early-stopping bias then gives a TV rate for the conditional distribution, averaged over covariates.

In latent diffusion, accuracy in sampling the latent distribution controls only part of the error in the decoded output law. For their latent Schr\"odinger bridge model, Jiao et al.~\cite{jiao2024latent} bound the $W_2$ error of the decoded output in terms of the latent sampling error, the autoencoder reconstruction error, and the discrepancy between the pretraining and target distributions. Even exact sampling in latent space would leave the latter two contributions in the bound.

Diffusion models also admit deterministic sampling through a probability-flow ODE driven by the same score as the reverse SDE~\cite{song2021score}. In the reverse-time convention of \eqref{eq:base:reversal}, this ODE is
\begin{equation}\label{eq:diffusion:probability:flow}
\frac{\mathrm d}{\mathrm d\tau}\boldsymbol Y_\tau^{\mathrm{pf}}
=\boldsymbol Y_\tau^{\mathrm{pf}}
+s^*(T_{\mathrm{end}}-\tau,\boldsymbol Y_\tau^{\mathrm{pf}}),
\qquad
\boldsymbol Y_0^{\mathrm{pf}}\sim p_{T_{\mathrm{end}}}.
\end{equation}
Assume that the ODE generates a unique flow and that the corresponding continuity equation has a unique solution on each interval $[0,T_{\mathrm{end}}-\delta]$, $0<\delta<T_{\mathrm{end}}$. With the exact score and the displayed initialization, its marginal law is $p_{T_{\mathrm{end}}-\tau}$, as for the reverse SDE. Replacing $s^*$ by the same learned score $\widehat s$ and initializing from $\gamma$ gives an approximate ODE sampler. The score-error decomposition can then be combined with ODE stability: under suitable spatial Lipschitz and moment bounds, the coupling and Gr\"onwall argument developed below for flow matching (Proposition~\ref{proposition:gronwall}) provides a route to Wasserstein error control.

\bigskip

\subsubsection{Stochastic interpolants and flow matching}\label{section:app:generative:flow_matching}

Flow matching trains a time-dependent velocity field by regression along prescribed probability paths; samples are then generated by solving the learned ordinary differential equation~\cite{lipman2023flow}. Stochastic interpolants construct such paths from random source and target samples and characterize the corresponding velocity through conditional expectations~\cite{Albergo2023Stochastic}. We use an independent Gaussian source and a linear interpolant.

Let $\gamma\coloneqq\mathcal{N}(\boldsymbol 0,\boldsymbol I_d)$ be an initial distribution, and let $p_0$ be a target distribution with samples $\boldsymbol X_{\mathrm{data}}^1,\ldots,\boldsymbol X_{\mathrm{data}}^n\sim^{\mathrm{i.i.d.}}p_0$. Define a linear interpolant for $(\boldsymbol Z,\boldsymbol X_{\mathrm{data}}) \sim \gamma \otimes p_0$ by 
\begin{equation*}
\boldsymbol X_t=t\boldsymbol X_{\mathrm{data}}+(1-t)\boldsymbol Z, \quad t\in(0,1).
\end{equation*}
The interpolant moves from the Gaussian law $\gamma$ at time zero to the target law $p_0$ at time one. For $0<t<1$, its marginal density $p_t$ satisfies the continuity equation
\begin{equation*}
\partial_tp_t(\boldsymbol x_t)+\nabla\cdot(u(t,\boldsymbol x_t)p_t(\boldsymbol x_t))=0, \quad t\in(0,1), ~ \boldsymbol x_t\in\mathbb{R}^d, 
\end{equation*}
where the velocity field $u:(0,1)\times\mathbb{R}^d\to\mathbb{R}^d$ is defined as
\begin{equation*}
u(t,\boldsymbol x_t)\coloneqq\mathbb{E}[\boldsymbol X_{\mathrm{data}}-\boldsymbol Z \mid \boldsymbol X_t=\boldsymbol x_t]=-\frac{1}{1-t} \boldsymbol x_t+\frac{1}{1-t}\mathbb{E}[\boldsymbol X_{\mathrm{data}} \mid \boldsymbol X_t=\boldsymbol x_t],
\end{equation*}
see~\cite{liu2023flow}. For $a>0$, let $\gamma_a$ denote the density of $\mathcal N(\boldsymbol 0,a\boldsymbol I_d)$. Conditioned on $\boldsymbol X_{\mathrm{data}}=\boldsymbol x_{\mathrm{data}}$, the density of the interpolant $\boldsymbol X_t$ is given as 
\begin{equation*}
p_{\boldsymbol X_t\mid \boldsymbol X_{\mathrm{data}} = \boldsymbol x_{\mathrm{data}}} = \gamma_{(1-t)^2}(\cdot - t\boldsymbol x_{\mathrm{data}}),
\end{equation*}
and by the law of total probability we have
\begin{equation*}
p_t
\coloneqq\int p_{\boldsymbol X_t\mid \boldsymbol X_{\mathrm{data}}= \boldsymbol x_{\mathrm{data}}}\, p_0(\mathrm d\boldsymbol x_{\mathrm{data}})=
\int\gamma_{(1-t)^2}(\cdot - t\boldsymbol x_{\mathrm{data}}) p_0(\mathrm d\boldsymbol x_{\mathrm{data}}).
\end{equation*}
Moreover, we have the important relation between the score of $p_t$ and the velocity field
\begin{equation*}
\nabla\log p_t(\boldsymbol x_t)= \frac{t}{1-t} u(t,\boldsymbol x_t)-\frac{1}{1-t}\boldsymbol x_t, \quad t \in (0,1), ~ \boldsymbol x_t\in\mathbb{R}^d,
\end{equation*}
see, e.g.,~\cite[Proposition 3.1]{ding2024characteristic}. 

\par The characteristic curves of the continuity equation define a flow map $\psi:[0,1)\times\mathbb{R}^d\rightarrow\mathbb{R}^d$, specified by the following probability flow ODE:
\begin{equation*}
\begin{aligned}
\frac{\mathrm{d}}{\mathrm{d}t}\psi(t,\boldsymbol z)&=u(t,\psi(t,\boldsymbol z)), \quad t\in (0,1), ~ \boldsymbol z\in\mathbb{R}^d, \\
\psi(0,\boldsymbol z)&=\boldsymbol z. 
\end{aligned}
\end{equation*}
Assume that this ODE generates a unique flow and that the continuity equation has a unique solution with initial law $\gamma$ on each interval $[0,t]$, $t<1$. The flow then transports the interpolant marginals: $\psi(t,\cdot)_{\#}\gamma=p_t$ for $0<t<1$. Sampling therefore amounts to drawing $\boldsymbol Z\sim \gamma$ and integrating the probability flow ODE, which requires an estimate of the velocity field $u$. For a fixed velocity field and initial point, the trajectory is deterministic. ODE stability then controls sampling error, playing the role of the Girsanov argument used for the reverse SDE.

\par Exactly as in Corollary~\ref{corollary:dsm}, the velocity field is characterized as an $L^2$-projection, since it is a conditional expectation. This yields the flow matching objective
\begin{equation*}
\mathcal{L}_{\mathrm{FM}}(v)\coloneqq\int_0^1\mathbb{E}_{(\boldsymbol Z,\boldsymbol X_{\mathrm{data}})\sim\gamma\otimes p_0}\bigl[\|v(t,\boldsymbol X_t)-(\boldsymbol X_{\mathrm{data}}-\boldsymbol Z)\|_2^2\bigr]\,\mathrm{d}t,
\end{equation*}
where $\boldsymbol X_t$ is the interpolant. The crucial point, and the reason the objective is implementable, is that the regression target $\boldsymbol X_{\mathrm{data}}-\boldsymbol Z$ is observable, whereas its conditional expectation $u$ is not.

\begin{proposition}[Flow matching]\label{proposition:fm:identity}
Suppose $p_0$ has finite second moment. The velocity field $u$ minimizes $\mathcal{L}_{\mathrm{FM}}$, uniquely up to equality almost everywhere with respect to $\mathrm dt\,p_t(\mathrm d\boldsymbol x)$. For every measurable $v\colon(0,1)\times\mathbb{R}^d\to\mathbb{R}^d$ with finite flow-matching risk,
\begin{equation*}
\mathcal{L}_{\mathrm{FM}}(v)-\mathcal{L}_{\mathrm{FM}}(u)=\int_0^1\|v(t,\cdot)-u(t,\cdot)\|_{L^2(p_t)}^2\,\mathrm{d}t.
\end{equation*}
\end{proposition}

\par Given i.i.d.\ samples $\{(t^i,\boldsymbol Z^i,\boldsymbol X_{\mathrm{data}}^i)\}_{i=1}^n$ from $\mathrm{unif}(0,1)\otimes\gamma\otimes p_0$, and writing $\boldsymbol X_{t^i}^i\coloneqq t^i\boldsymbol X_{\mathrm{data}}^i+(1-t^i)\boldsymbol Z^i$, the empirical counterpart of the population risk is
\begin{equation*}
\widehat{\mathcal{L}}_{\mathrm{FM}}(v)\coloneqq\frac{1}{n}\sum_{i=1}^n\|v(t^i,\boldsymbol X_{t^i}^i)-(\boldsymbol X_{\mathrm{data}}^i-\boldsymbol Z^i)\|_2^2.
\end{equation*}
Let $\mathcal V$ be a class of velocity fields and let $\widehat v\in\mathcal V$ minimize this loss up to a tolerance $\epsilon_{\mathrm{opt}}$. For finite risks and finite expected uniform deviation, Proposition~\ref{proposition:fm:identity} and the empirical-risk comparison give
\begin{equation}\label{eq:structural:velocity:oracle}
\begin{aligned}
\mathcal E_{\mathrm{FM}}^2
&\coloneqq \mathbb E\int_0^1\|\widehat v(t,\cdot)-u(t,\cdot)\|_{L^2(p_t)}^2\,\mathrm dt\\
&\leq
\underbrace{\inf_{v\in\mathcal V}\int_0^1
\|v(t,\cdot)-u(t,\cdot)\|_{L^2(p_t)}^2\,\mathrm dt}_{\mathcal E_{\mathrm{app}}^2}+2\underbrace{\mathbb E\sup_{v\in\mathcal V}
\bigl|\mathcal L_{\mathrm{FM}}(v)-\widehat{\mathcal L}_{\mathrm{FM}}(v)\bigr|}_{\mathcal E_{\mathrm{sta}}^2}
+\epsilon_{\mathrm{opt}}.
\end{aligned}
\end{equation}
Here the expectation is over the sampled endpoint pairs and times, and $\epsilon_{\mathrm{opt}}$ is on the same squared-risk scale. Approximation and statistical estimates must control the tails of the interpolant and regression response.

Let $\widehat{\psi}$ denote the flow map generated by the estimated velocity field, that is,
\begin{equation*}
\begin{aligned}
\frac{\mathrm{d}}{\mathrm{d}t}\widehat{\psi}(t,\boldsymbol z)&=\widehat{v}(t,\widehat{\psi}(t,\boldsymbol z)), \quad t\in(0,1), \\
\widehat{\psi}(0,\boldsymbol z)&=\boldsymbol z,
\end{aligned}
\end{equation*}
and let $\widehat{\gamma}$ be the distribution actually used to initialize the ODE. Since the two trajectories solve ODEs with different vector fields, their discrepancy is controlled by Gr\"onwall's inequality. In addition to the ideal transport assumptions above, assume that the learned ODE generates a unique flow on every interval $[0,t]$, $t<1$, and that both flows have finite second moments under their respective initial laws. Suppose that the estimated field satisfies
\[
\|\widehat v(s,\boldsymbol x)-\widehat v(s,\boldsymbol y)\|_2
\leq L_{\mathrm{Lip},v}(s)\|\boldsymbol x-\boldsymbol y\|_2,
\qquad
\int_0^t L_{\mathrm{Lip},v}(s)\,\mathrm ds<\infty\quad(t<1),
\]
for a nonnegative measurable function $L_{\mathrm{Lip},v}$ and all $\boldsymbol x,\boldsymbol y\in\mathbb R^d$, and that the velocity-error integrals below are finite.

\begin{proposition}[Stability of the probability flow ODE]\label{proposition:gronwall}
For every $\boldsymbol z\in\mathbb{R}^d$ and every $t\in(0,1)$,
\begin{equation*}
\|\widehat{\psi}(t,\boldsymbol z)-\psi(t,\boldsymbol z)\|_2\leq\int_0^t\exp\Bigl(\int_r^t L_{\mathrm{Lip},v}(\tau)\,\mathrm{d}\tau\Bigr)\|\widehat{v}(r,\psi(r,\boldsymbol z))-u(r,\psi(r,\boldsymbol z))\|_2\,\mathrm{d}r.
\end{equation*}
Consequently, by coupling the initial variables $(\widehat{\boldsymbol Z}, \boldsymbol Z)$ from $\widehat{\gamma}$ and $\gamma$ and integrating the path-wise estimate, writing $\Lambda_t\coloneqq\exp\bigl(\int_0^t L_{\mathrm{Lip},v}(\tau)\,\mathrm{d}\tau\bigr)$, we obtain
\begin{equation*}
W_2^2\bigl(\widehat{\psi}(t,\cdot)_{\#}\widehat{\gamma},\,p_t\bigr)\leq2\Lambda_t^2\Bigl(W_2^2\bigl(\widehat{\gamma},\gamma\bigr)+t\int_0^t\|\widehat{v}(r,\cdot)-u(r,\cdot)\|_{L^2(p_r)}^2\,\mathrm{d}r\Bigr).
\end{equation*}
\end{proposition}

\par The amplification factor $\Lambda_t$ is the price of integrating the ODE and dictates early stopping. To make its endpoint behavior explicit, suppose in addition that $L_{\mathrm{Lip},v}(t)\leq C_{\mathrm{Lip}}/(1-t)$ near $t=1$ for some $C_{\mathrm{Lip}}>0$. Terminating the integration at $t=1-\delta$ then gives $\Lambda_{1-\delta}\lesssim\delta^{-C_{\mathrm{Lip}}}$, yielding a polynomial amplification factor in $1/\delta$. The remaining bias is bounded via $\boldsymbol X_{1-\delta}-\boldsymbol X_{\mathrm{data}}=\delta(\boldsymbol Z-\boldsymbol X_{\mathrm{data}})$, giving
\begin{equation}\label{eq:early:stopping:fm}
W_2(p_{1-\delta},p_0)\leq\delta\,\mathbb{E}\bigl[\|\boldsymbol X_{\mathrm{data}}-\boldsymbol Z\|_2^2\bigr]^{1/2}\lesssim\delta(\sqrt{d}+m_2),
\end{equation}
where $m_2^2=\mathbb{E}[\|\boldsymbol X_{\mathrm{data}}\|_2^2]$. 

\par Let $\widehat p^{\mathrm{num}}$ be the numerical output at time $1-\delta$, initialized from $\widehat\gamma$, and define
\[
\mathcal E_{\mathrm{disc}}^{\mathrm{FM}}
:=W_2\bigl(\widehat\psi(1-\delta,\cdot)_\sharp\widehat\gamma,\widehat p^{\mathrm{num}}\bigr).
\]
Thus this term compares the exact and numerical learned flows with the same initialization. The triangle inequality and Proposition~\ref{proposition:gronwall}, with $1-\delta\leq1$, give the following decomposition. The early-stopping term is controlled by~\eqref{eq:early:stopping:fm}, while \eqref{eq:structural:velocity:oracle} controls the expected integrated velocity error:
\begin{equation}\label{eq:fm_w2_decomp}
\begin{aligned}
W_2(p_0, \widehat p^{\mathrm{num}})
\;\leq\;& \underbrace{W_2(p_0, p_{1-\delta})}_{\text{early stopping}} 
+ \underbrace{\sqrt{2}\Lambda_{1-\delta} W_2(\widehat{\gamma}, \gamma)}_{\text{initialization}} \\
&+ \underbrace{\sqrt{2}\Lambda_{1-\delta} \left(\int_0^{1-\delta}\|\widehat{v}(r,\cdot)-u(r,\cdot)\|_{L^2(p_r)}^2\,\mathrm{d}r\right)^{1/2}}_{\text{velocity estimation}} 
+ \underbrace{\mathcal E_{\mathrm{disc}}^{\mathrm{FM}}}_{\text{discretization}}.
\end{aligned}
\end{equation}

The stability bound is conditional on the trained field. If $L_{\mathrm{Lip},v}$ is a deterministic bound valid uniformly over the velocity class, taking expectations in \eqref{eq:fm_w2_decomp} and applying Cauchy--Schwarz bounds the velocity-estimation term by $\sqrt{2}\Lambda_{1-\delta}(\mathcal E_{\mathrm{app}}^2+2\mathcal E_{\mathrm{sta}}^2+\epsilon_{\mathrm{opt}})^{1/2}$.

Spatial Lipschitz control limits the amplification of velocity errors, while Euler discretization also requires control of temporal regularity and the magnitude of the learned field. Gao et al.~\cite[Section~4.1]{gao2024convergence} make these roles explicit by separating velocity estimation, Euler discretization, and early stopping. The relevant regularity bounds may deteriorate as the stopping time approaches the data endpoint, so reducing the early-stopping bias may require a finer time grid to maintain numerical accuracy.

\par The sharper rate of \cite{fukumizu2025flow} concerns more general Gaussian interpolants, written in the present time convention as $\boldsymbol X_t=a_t\boldsymbol X_{\mathrm{data}}+b_t\boldsymbol Z$, with noise scale $b_t\asymp(1-t)^\kappa$ near the data endpoint. The analysis uses a time partition and interval-specific velocity networks. Its empirical objectives average over the same $n$ data samples on every interval, but integrate over time and the conditional Gaussian noise, rather than drawing only one time--noise pair per data sample as in the empirical objective above. These choices allow approximation and statistical errors to be controlled separately across time scales.

\begin{theorem}[Flow matching convergence; \citealp{fukumizu2025flow}, Theorem 9]\label{thm:fm_main}
Let $d\geq2$ and $1\le r\le2$. Suppose Assumptions (A1)--(A5) of \cite{fukumizu2025flow} hold with Besov density smoothness $\beta>0$ and interpolation parameter $\kappa\geq1/2$. These assumptions require a target density supported on $[-1,1]^d$, bounded above and away from zero there, with higher Besov smoothness of order greater than $\max\{6\beta-1,1\}$ in a boundary layer. They also require power-law behavior of $1-a_t$ and $b_t$ near the data endpoint, upper and positive lower bounds on $a_t^2+b_t^2$, and derivative bounds on the schedules, including a polylogarithmic bound on their integrated squared derivatives when $\kappa=1/2$. The spatial derivative of the conditional mean $\mathbb E[\boldsymbol X_{\mathrm{data}}\mid\boldsymbol X_t=\boldsymbol x]$ must be uniformly bounded along the truncated path. Apply the time partition, interval-specific network classes, and empirical minimization procedure specified in its Theorem~9, using the same $n$ data samples for each interval. In the present time convention, let $\widehat p_{1-\delta_n}$ denote the exact-flow output, where $\delta_n$ corresponds to the early-stopping parameter $T_0$ chosen in the reference so that the truncation bias is of smaller order than the rate below. For every sufficiently small fixed $\zeta>0$,
\begin{equation*}
\mathbb E\,W_r(\widehat p_{1-\delta_n},p_0) \lesssim n^{-\frac{\beta+(2\kappa)^{-1}-\zeta}{2\beta+d}}.
\end{equation*}
\end{theorem}
Here $\kappa$ describes the decay of the interpolation noise scale near the data endpoint and is distinct from the Lipschitz coefficient $C_{\mathrm{Lip}}$ above. The almost minimax rate corresponds to $\kappa=1/2$, rather than standard linear interpolation where $\kappa=1$.

\par The two analyses both estimate a field by regression and then control its effect on the sampling dynamics. In each case, a conditional expectation minimizes the population quadratic loss, whose excess risk measures an $L^2$ field error. Approximation and statistical bounds, together with an optimization residual, control this error; stability of the sampling dynamics then converts it into distributional error. The tools and error metrics differ: the reverse-SDE analysis uses Girsanov's theorem to obtain KL and TV bounds, with the factor $\delta^{-1}$ arising when the preceding weighted denoising loss is used to control the unweighted integrated squared score error; the ODE analysis uses Gr\"onwall's inequality to obtain Wasserstein bounds, with Lipschitz amplification $\Lambda_{1-\delta}$.

\begin{remark}
In addition to diffusion models and flow-based models, several other dynamics-based generative models have been proposed and analyzed. For example,~\cite{wang2021Deep,jiao2024latent} constructed generative models based on the Schr{\"o}dinger bridge, while~\cite{Gao2019deep,Gao2022Deep} developed methods derived from gradient flows in the space of probability measures. The latter methods estimate density ratios or density differences to construct the transport field. Furthermore, since both the reverse SDE and the probability flow ODE require many evaluations of the learned field during sampling, distillation methods have been developed that compress the trajectory into one or a few steps; an error analysis of this approach is established in~\cite{ding2024characteristic}.
\end{remark}

\begin{remark}[Conditional generative models]
The framework outlined above extends naturally to conditional sampling, where one aims to generate from the conditional distribution of a response given a covariate rather than from a single fixed target. The conditional F{\"o}llmer flow of~\cite{chang2025deep} realizes this idea through an ODE whose velocity field is, once more, a conditional expectation, and its analysis follows the same two steps: an oracle inequality for the estimated velocity field, followed by a stability estimate for the induced flow map. Beyond sampling itself, draws from the estimated conditional law can approximate conditional means and quantiles and support the construction of prediction regions. The statistical guarantees for such inference require controlling the effect of estimating the conditional law; see the prediction and bootstrap analyses in~\cite{Jiao2025Model,chang2026deep}.
\end{remark}

\begin{remark}[Generative models for inverse problems]
A pre-trained generative model can also serve as a learned prior for inverse problems, where the goal is to sample from the posterior distribution of the signal given indirect and noisy measurements~\cite{chang2025provable,chang2026inference}. A closely related application is data assimilation, in which the posterior must be updated sequentially as observations arrive along the trajectory of a dynamical system; see~\cite{Ding2026Nonlinear}.
\end{remark}

\subsection{Reinforcement learning}\label{section:app:rl}

Reinforcement learning (RL) studies how an agent chooses actions through interaction with an environment to maximize its expected cumulative reward. Fitted $Q$-iteration (FQI) holds the previous action-value estimate fixed to construct regression targets from the observed reward and the highest estimated value at the next state. Its final policy selects an action of highest fitted value~\cite{Feng2023Over}. In the approximate policy iteration (API) method studied here, each round first selects a greedy policy from the current action-value estimate, then holds that policy fixed while estimating the return from following it after the initial action~\cite{Jiao2025Deepapproximate}. Actor--critic (AC) methods also combine policy evaluation and improvement. The critic estimates the value of the actor's policy, and the actor uses that estimate to update its policy parameters, typically by a policy-gradient step~\cite{sutton2018reinforcement}. The following analysis focuses on FQI and API.

Consider a discounted Markov decision process $(\mathcal S,\mathcal A,P,\mathcal R,\gamma)$, where $P$ is a state-transition kernel, $\mathcal R(\cdot\mid\boldsymbol x,a)$ is a reward distribution on $[0,R_{\max}]$, and $0<\gamma<1$. Throughout this subsection, $\gamma$ denotes the discount factor. Write $\bar R(\boldsymbol x,a)$ for its mean and assume that measurable greedy policies exist. For a policy $\pi$, let $Q^\pi(\boldsymbol x,a)=\mathbb E_\pi[\sum_{t\geq0}\gamma^tR_t\mid\boldsymbol X_0=\boldsymbol x,A_0=a]$. Then $0\leq Q^\pi\leq V_{\max}:=R_{\max}/(1-\gamma)$, and the optimal action-value function $Q^*$ satisfies
\begin{equation}\label{eq:rl:bellman}
Q^*=\mathcal T^*Q^*,
\qquad
(\mathcal T^*Q)(\boldsymbol x,a)
=\bar R(\boldsymbol x,a)+\gamma\mathbb E_{\boldsymbol X'\sim P(\cdot\mid\boldsymbol x,a)}
\max_{a'\in\mathcal A}Q(\boldsymbol X',a').
\end{equation}
For a fixed policy $\pi$, $Q^\pi$ describes the return from taking action $a$ at state $\boldsymbol x$ and following $\pi$ thereafter. Greedy improvement selects actions using this long-term value estimate.
The Bellman optimality operator $\mathcal T^*$ is a $\gamma$-contraction in the supremum norm. Throughout this subsection, $\varrho$ denotes the training state--action distribution, $\eta_{\mathrm{ev}}$ the evaluation distribution, and $d$ the dimension of the Euclidean state--action input when approximation on $[0,1]^d$ is invoked.

One-step transition data provide the responses for fitting a Bellman update in FQI and for estimating a policy-specific Bellman residual in API. In each case the fit is measured under $\varrho$. To bound the value gap of the final policy, the analysis propagates these errors across discounted iterations and accounts for how well $\varrho$ covers the state--action pairs reached from $\eta_{\mathrm{ev}}$.

\subsubsection{Fitted \texorpdfstring{$Q$}{Q}-iteration with dependent data}\label{section:app:rl:fqi}

FQI fits Bellman updates by supervised regression. In an early batch formulation, Ernst et al.~\cite{ernst2005tree} fitted these updates with regression trees, reusing the same transition data across iterations. The sampling protocol of Feng et al.~\cite[Algorithm~1]{Feng2023Over} uses $J$ mutually independent batches $\mathcal D_j$, each containing $n$ stationary, $\beta$-mixing transitions $\boldsymbol Z_{j,i}=(\boldsymbol X_{j,i},A_{j,i},R_{j,i},\boldsymbol X'_{j,i})$. Each transition records a state--action pair, its realized reward, and the next state under the data-collection process. These one-step observations provide regression inputs and targets, even though the target policy's value is not directly observed. The total sample count is $nJ$. Initialize $\widehat Q_0=0$ and use a network class $\mathcal F_Q$ with outputs in $[0,V_{\max}]$. Set
\begin{equation}\label{eq:fqi_step}
\begin{aligned}
Y_{j,i}&=R_{j,i}+\gamma\max_{a'}\widehat Q_{j-1}(\boldsymbol X'_{j,i},a'),\\
\widehat{\mathcal L}_j(Q)&=\frac1n\sum_{i=1}^n
\bigl(Q(\boldsymbol X_{j,i},A_{j,i})-Y_{j,i}\bigr)^2.
\end{aligned}
\end{equation}
An approximate fit $\widehat Q_j\in\mathcal F_Q$ satisfies $\widehat{\mathcal L}_j(\widehat Q_j)\leq\inf_{Q\in\mathcal F_Q}\widehat{\mathcal L}_j(Q)+\epsilon_{\mathrm{opt},j}$. The output policy is $\pi_J(\boldsymbol x)\in\argmax_a\widehat Q_J(\boldsymbol x,a)$; exact empirical minimization corresponds to $\epsilon_{\mathrm{opt},j}=0$.

Conditional on the earlier batches, the regression target is $\mathcal T^*\widehat Q_{j-1}$. Indeed, under the stationary law with state--action marginal $\varrho$,
\[
\mathbb E\!\left[Y_j\mid \boldsymbol X=\boldsymbol x,A=a\right]
=\bar R(\boldsymbol x,a)
+\gamma\mathbb E_{\boldsymbol X'\sim P(\cdot\mid \boldsymbol x,a)}
\max_{a'}\widehat Q_{j-1}(\boldsymbol X',a')
=\mathcal T^*\widehat Q_{j-1}(\boldsymbol x,a).
\]
Let $\mathcal L_j$ be the conditional population counterpart of $\widehat{\mathcal L}_j$. Conditional-mean orthogonality gives
\begin{equation}\label{eq:rl:regression:identity}
\mathcal L_j(Q)-\mathcal L_j(\mathcal T^*\widehat Q_{j-1})
=\|Q-\mathcal T^*\widehat Q_{j-1}\|_{L^2(\varrho)}^2.
\end{equation}
Conditional on the earlier batches, the empirical-risk comparison of Section~\ref{sec:decomp} applies to the fixed Bellman target as in~\cite[Lemma~20]{Feng2023Over}. Taking expectations gives
\begin{equation}\label{eq:rl:one:step}
\begin{aligned}
\mathbb E\|\widehat Q_j & -\mathcal T^*\widehat Q_{j-1}\|_{L^2(\varrho)}^2 \\
&\leq \mathbb E
\underbrace{\inf_{Q\in\mathcal F_Q}\|Q-\mathcal T^*\widehat Q_{j-1}\|_{L^2(\varrho)}^2}
_{\mathcal E_{\mathrm{app},j}^{\mathrm{FQI}}} +2\mathbb E
\underbrace{\sup_{Q\in\mathcal F_Q}|\mathcal L_j(Q)-\widehat{\mathcal L}_j(Q)|}
_{\mathcal E_{\mathrm{sta},j}^{\mathrm{FQI}}}
+\epsilon_{\mathrm{opt},j}.
\end{aligned}
\end{equation}
The approximation term depends on the earlier batches through $\widehat Q_{j-1}$. A uniform smoothness assumption on the Bellman targets allows the same approximation bound to hold at every iteration.

To convert regression error into policy error, coverage of the evaluation distributions by the training distribution is essential. Let $P^\pi$ denote the state--action transition kernel under policy $\pi$ and define
\begin{equation}\label{eq:rl:concentrability}
\begin{aligned}
c_m&=\sup_{\pi_1,\ldots,\pi_m}
\left\|\frac{\mathrm d(\eta_{\mathrm{ev}}P^{\pi_1}\cdots P^{\pi_m})}{\mathrm d\varrho}\right\|_{L^2(\varrho)},\\
C_{\eta_{\mathrm{ev}},\varrho}&=(1-\gamma)^2\sum_{m\geq1}m\gamma^{m-1}c_m<\infty,
\end{aligned}
\end{equation}
where failure of absolute continuity gives $c_m=\infty$. Let $e_j=\widehat Q_j-\mathcal T^*\widehat Q_{j-1}$ be the stage-$j$ fitting error. For any distribution $\nu_m=\eta_{\mathrm{ev}}P^{\pi_1}\cdots P^{\pi_m}$ reached after $m$ steps, Cauchy--Schwarz gives
\[
\int |e_j|\,\mathrm d\nu_m
=\int |e_j|\,\frac{\mathrm d\nu_m}{\mathrm d\varrho}\,\mathrm d\varrho
\leq c_m\|e_j\|_{L^2(\varrho)}.
\]
Thus each residual in the discounted expansion of the policy value gap can be bounded in the training norm. Applying~\eqref{eq:rl:one:step} to these residuals yields
\begin{align}
\mathbb E\|Q^*-Q^{\pi_J}\|_{L^1(\eta_{\mathrm{ev}})}
\leq &\frac{2\gamma C_{\eta_{\mathrm{ev}},\varrho}}{(1-\gamma)^2}
\max_{1\leq j\leq J}
\Bigl(\mathbb E\mathcal E_{\mathrm{app},j}^{\mathrm{FQI}}
+2\mathbb E\mathcal E_{\mathrm{sta},j}^{\mathrm{FQI}}
+\epsilon_{\mathrm{opt},j}\Bigr)^{1/2} +\frac{4\gamma^{J+1}R_{\max}}{(1-\gamma)^2}.
\label{eq:rl:total:structural}
\end{align}
The policy-error bound follows from the discounted residual expansion underlying Proposition~19 of Feng et al.~\cite{Feng2023Over} and Theorem~6.1 of the full version of Fan et al.~\cite{Fan2020DQN}. Taking expectations term by term before bounding the weighted sum gives the maximum of the expected single-stage errors in~\eqref{eq:rl:total:structural}. The square root converts squared regression risk into an $L^2$ residual, while the geometric term accounts for the finite number of Bellman iterations.

Different regression classes can be used in this propagation argument. For independent transitions, the sparse-network analysis of Fan et al.~\cite{Fan2020DQN} gives the benchmark $\widetilde O(n^{-\beta/(2\beta+d_{\mathcal S})})$ for expected $L^1$ policy-value error, apart from the iteration remainder, in its ordinary H\"older specialization for a fixed finite action space and state dimension $d_{\mathcal S}$. Its assumptions include fresh data at each iteration, smooth Bellman images, and finite concentrability. Kang et al.~\cite{kang2023fqi} analyze ReLU FQI under exponential $\alpha$-mixing and control the statistical error through network size. We consider the weight-norm analysis of Feng et al.~\cite{Feng2023Over}, which allows large widths and accounts for $\beta$-mixing within each batch.

Write $\bar\beta_k$ for the mixing coefficients, distinguishing them from the H\"older exponent $\beta$. Partition each batch into $2B_n$ blocks of length $a_n$, with $n=2a_nB_n$ after discarding any remainder. Coupling alternating blocks to independent copies incurs an accumulated cost of order $B_n\bar\beta_{a_n}$~\cite{yu1994rates,Mohri2009RademacherMixing}. For the value class $\mathcal F_Q$, take $Q=V_{\max}\operatorname{clip}_{[0,1]}(f)$ with $f\in\mathcal{NN}_{\mathrm{ReLU}}(W,L,\mathcal K)$, where clipping projects a scalar onto $[0,1]$. For normalized Bellman targets in a fixed H\"older ball of order $\beta$ and the width and depth specified below, the single-stage terms satisfy
\begin{align}
\mathcal E_{\mathrm{app},j}^{\mathrm{FQI}}
&\lesssim V_{\max}^2\mathcal K^{-2\beta/(d+1)},\label{eq:rl:approx}\\
\mathbb E\mathcal E_{\mathrm{sta},j}^{\mathrm{FQI}}
&\lesssim V_{\max}^2\left[
\frac{\mathcal K\sqrt{L+2+\log(d+1)}}{\sqrt{B_n}}
+B_n\bar\beta_{a_n}\right].
\label{eq:fqi_stat}
\end{align}
For~\eqref{eq:rl:approx}, Theorem~\ref{thm:A3} gives a uniform approximation error of order $\mathcal K^{-\beta/(d+1)}$ for the normalized Bellman target. Since the target takes values in $[0,1]$, clipping can only reduce this error. Squaring and restoring the factor $V_{\max}^2$ gives the stated bound. For~\eqref{eq:fqi_stat}, the blocking argument of Feng et al.~\cite[Theorem~21]{Feng2023Over} uses the Rademacher bound of Lemma~\ref{lem:S6}, which also applies after clipping by contraction. The statistical bound is independent of width. Balancing the two terms that depend on $\mathcal K$ up to logarithmic factors and substituting into~\eqref{eq:rl:total:structural} gives the following result.

\begin{theorem}[Norm-controlled FQI]\label{thm:rl_fqi}
Suppose that the state--action inputs lie in $[0,1]^d$ and that the normalized Bellman targets $(\mathcal T^*\widehat Q_{j-1})/V_{\max}$ belong uniformly to a fixed H\"older ball of order $\beta>0$ on $[0,1]^d$. Assume independent batches of stationary $\beta$-mixing transitions and the concentrability condition~\eqref{eq:rl:concentrability}. Use exact empirical minimization over the value class $\mathcal F_Q$ above, with
\[
\mathcal K\asymp B_n^{(d+1)/(4\beta+2d+2)},\qquad
W\gtrsim\mathcal K^{(2d+\beta)/(2d+2)},\qquad
L\asymp\log(2+\mathcal K).
\]
For sufficiently large $B_n$, this depth meets the lower bound in Theorem~\ref{thm:A3}, and the greedy policy satisfies
\begin{align}
\mathbb E\|Q^*-Q^{\pi_J}\|_{L^1(\eta_{\mathrm{ev}})}
\lesssim{}&
\frac{\gamma C_{\eta_{\mathrm{ev}},\varrho}R_{\max}}{(1-\gamma)^3}
\left[B_n^{-\beta/(4\beta+2d+2)}
\{\log(dB_n)\}^{1/4}+\sqrt{B_n\bar\beta_{a_n}}\right]\notag\\
&+\frac{\gamma^{J+1}R_{\max}}{(1-\gamma)^2}.
\label{eq:rl:feng:rate}
\end{align}
\end{theorem}

If $\bar\beta_k\leq C_{\mathrm{mix}}e^{-b_{\mathrm{mix}}k^\eta}$, choosing $a_n\asymp(\log n)^\tau$ with $\tau\eta>1$ makes the mixing penalty smaller than every polynomial rate; $B_n\asymp n/(\log n)^\tau$ then preserves the sample exponent up to logarithms. Choosing $J\asymp\log n$ with a sufficiently large constant makes the iteration remainder no larger than the learning term. The resulting rate is $\widetilde O(n^{-\beta/(4\beta+2d+2)})$; replacing the per-stage sample size $n$ by the total transition count $N_{\mathrm{tot}}=nJ$ changes only logarithmic factors. Low-dimensional support can improve the rate under the manifold and dynamical assumptions discussed in \citep{Feng2023Over}.

\subsubsection{Approximate policy iteration}\label{section:app:rl:api}

Approximate policy iteration (API) alternates greedy improvement with evaluation of the resulting policy. In least-squares policy iteration, Lagoudakis and Parr~\cite{lagoudakis2003least} used linear features to estimate the action values of a fixed policy from transition data, then chose the next policy greedily. In the neural method studied here, policy evaluation seeks a fixed point of the same policy-specific Bellman operator
\[
(\mathcal T^\pi Q)(\boldsymbol x,a)
=\bar R(\boldsymbol x,a)+\gamma\mathbb E_{\boldsymbol X'\mid\boldsymbol x,a}
\int_{\mathcal A}Q(\boldsymbol X',a')\,\pi(\mathrm da'\mid\boldsymbol X').
\]
Fix $J\geq2$. Starting from $\widehat Q_0=0$, form $\pi_j$ greedily from $\widehat Q_{j-1}$ and approximately evaluate $\pi_j$, for $1\leq j<J$; the output $\pi_J$ is greedy with respect to $\widehat Q_{J-1}$. We assess this policy through the value gap $\|Q^*-Q^{\pi_J}\|_{L^2(\eta_{\mathrm{ev}})}$. For the single-stage population analysis, fix $\pi_j$ and a candidate $Q$. The response $Y_j(Q)=R+\gamma\int Q(\boldsymbol X',a')\,\pi_j(\mathrm da'\mid\boldsymbol X')$ depends on this candidate, not on the frozen previous estimate. Although
\[
\mathbb E\!\left[Y_j(Q)\mid \boldsymbol X=\boldsymbol x,A=a\right]
=\mathcal T^{\pi_j}Q(\boldsymbol x,a),
\]
the squared prediction loss does not reduce to the squared Bellman residual alone. Its conditional decomposition contains
\[
\mathbb E\!\left[(Q(\boldsymbol X,A)-Y_j(Q))^2\right]
=\|Q-\mathcal T^{\pi_j}Q\|_{L^2(\varrho)}^2
+\mathbb E\,\mathrm{Var}\!\left(Y_j(Q)\mid \boldsymbol X,A\right),
\]
and the variance term depends on $Q$. Minimizing the ordinary squared prediction loss can therefore trade off Bellman residual against a $Q$-dependent noise level, rather than directly minimizing the residual.

Antos et al.~\cite{antos2008learning} addressed this double-sampling difficulty in fitted policy iteration from a single trajectory by introducing an auxiliary function. Jiao et al.~\cite{Jiao2025Deepapproximate} use an auxiliary network $O$ evaluated against the same response as $Q$. The two squared losses share the conditional variance of $Y_j(Q)$, so their difference cancels this term at the population level. At evaluation stage $j$, use $n$ transitions $\{(\boldsymbol X_{j,i},A_{j,i},R_{j,i},\boldsymbol X'_{j,i})\}_{i=1}^n$ from a stationary data-collection process with state--action marginal $\varrho$, and let $Y_{j,i}(Q)$ be the corresponding realization of $Y_j(Q)$. Let $\mathcal C_n$ be their ReLU convolutional-network class with outputs bounded in absolute value by $V_{\max}$ and define the empirical minimax loss
\begin{equation}\label{eq:rl:api:loss}
\begin{aligned}
\widehat\ell_j(Q,O)
=\frac1n\sum_{i=1}^n\Big[
\bigl(Q(\boldsymbol X_{j,i},A_{j,i})-Y_{j,i}(Q)\bigr)^2-\bigl(O(\boldsymbol X_{j,i},A_{j,i})-Y_{j,i}(Q)\bigr)^2
\Big].
\end{aligned}
\end{equation}
The estimate approximately minimizes $\sup_{O\in\mathcal C_n}\widehat\ell_j(Q,O)$ over $Q\in\mathcal C_n$, within tolerance $\epsilon_{\mathrm{opt},j}$. For fixed $\pi_j$, let $\ell_j$ be the population counterpart under this data-collection law. The two conditional variances cancel, giving
\begin{equation}\label{eq:rl:api:identity}
\ell_j(Q,O)
=\|Q-\mathcal T^{\pi_j}Q\|_{L^2(\varrho)}^2
-\|O-\mathcal T^{\pi_j}Q\|_{L^2(\varrho)}^2.
\end{equation}
For fixed $Q$, maximizing over $O$ amounts to least-squares fitting of the conditional Bellman mean $\mathcal T^{\pi_j}Q$. The minimax formulation thus provides a criterion based on transition samples for Bellman-residual minimization, with the variance cancellation occurring at the population level. An unrestricted auxiliary class would recover the squared Bellman residual exactly after maximization. A neural auxiliary class contributes an additional approximation error.

Comparing population and empirical minimax objectives twice, and using \eqref{eq:rl:api:identity}, gives the single-stage error decomposition
\begin{equation}\label{eq:rl:api:decomp}
\begin{aligned}
\|\widehat Q_j-\mathcal T^{\pi_j}\widehat Q_j\|_{L^2(\varrho)}^2 &\leq
\underbrace{\inf_{Q\in\mathcal C_n}\|Q-\mathcal T^{\pi_j}Q\|_{L^2(\varrho)}^2+\sup_{Q\in\mathcal C_n}\inf_{O\in\mathcal C_n}
\|O-\mathcal T^{\pi_j}Q\|_{L^2(\varrho)}^2}_{\mathcal E_{\mathrm{app},j}^{\mathrm{API}}}\\
& ~~~~ +2\underbrace{\sup_{Q,O\in\mathcal C_n}
|\ell_j(Q,O)-\widehat\ell_j(Q,O)|}
_{\mathcal E_{\mathrm{sta},j}^{\mathrm{API}}}
+\epsilon_{\mathrm{opt},j}.
\end{aligned}
\end{equation}
This is the two-network counterpart of \eqref{eq:rl:one:step}; see~\cite[Lemma~4.2]{Jiao2025Deepapproximate} for the source decomposition. Uniform approximation of $Q^{\pi_j}$ controls the first term in $\mathcal E_{\mathrm{app},j}^{\mathrm{API}}$, since $\|Q-\mathcal T^{\pi_j}Q\|_\infty\leq(1+\gamma)\|Q-Q^{\pi_j}\|_\infty$. Approximation of the Bellman images $\mathcal T^{\pi_j}\mathcal C_n$ controls the second. The regularity assumptions thus apply to both value functions and Bellman images. The statistical term requires uniform empirical-process control over the pair $(Q,O)$, with its size determined by the network-class complexity and the dependence of the transition data. The analysis in~\cite{Jiao2025Deepapproximate} uses entrywise parameter bounds and geometrically $C$-mixing data, whereas the preceding FQI analysis uses product-norm control and $\beta$-mixing.

Residuals measured under $\varrho$ affect policy performance along future state--action trajectories. Discounting controls their accumulation, while concentrability accounts for the change from the training distribution to the distributions visited by these trajectories.

For this analysis define $c_m^{\mathrm{API}}$ by replacing the $L^2(\varrho)$ density-ratio norm in \eqref{eq:rl:concentrability} with $L^\infty(\varrho)$, and assume $C_{\eta_{\mathrm{ev}},\varrho}^{\mathrm{API}}=(1-\gamma)^2\sum_{m\geq1}m\gamma^{m-1}c_m^{\mathrm{API}}<\infty$. Combining the bound in~\cite[Proposition~4.1]{Jiao2025Deepapproximate} with \eqref{eq:rl:api:decomp} gives the pathwise structural bound
\begin{align}
\|Q^*-Q^{\pi_J}\|_{L^2(\eta_{\mathrm{ev}})}
\leq{}&\frac{2\gamma\sqrt{C_{\eta_{\mathrm{ev}},\varrho}^{\mathrm{API}}}}{(1-\gamma)^2}
\max_{1\leq j<J}
\left(\mathcal E_{\mathrm{app},j}^{\mathrm{API}}
+2\mathcal E_{\mathrm{sta},j}^{\mathrm{API}}
+\epsilon_{\mathrm{opt},j}\right)^{1/2}\notag\\
&+\frac{2\gamma^{(J+2)/2}R_{\max}}{(1-\gamma)^2}.
\label{eq:rl:api:structural}
\end{align}
To state a conditional convergence guarantee, prescribe an iteration count $J=J_n\geq2$. Suppose that a deterministic residual scale $r_n>0$ and a failure probability $0<\tau_n<1$ satisfy, with probability at least $1-\tau_n$,
\begin{equation}\label{eq:rl:api:uniform}
\max_{1\leq j<J_n}
\left(
\mathcal E_{\mathrm{app},j}^{\mathrm{API}}
+2\mathcal E_{\mathrm{sta},j}^{\mathrm{API}}
+\epsilon_{\mathrm{opt},j}
\right)\leq r_n^2.
\end{equation}
This is a simultaneous requirement for the data-dependent policies actually generated by the algorithm. With fresh batches, one way to establish it is to condition each stage's bound on the preceding data, uniformly over the resulting policy, and allocate failure probabilities across all stages. When data are reused, a uniform estimate over the policies generated from those data is one way to account for their dependence. The scale $r_n$ may therefore depend on both $J_n$ and $\tau_n$ through the simultaneous control in \eqref{eq:rl:api:uniform}.

On the same event, substitution into \eqref{eq:rl:api:structural} gives
\begin{equation}\label{eq:rl:api:conditional}
\|Q^*-Q^{\pi_{J_n}}\|_{L^2(\eta_{\mathrm{ev}})}
\leq
\frac{2\gamma\sqrt{C_{\eta_{\mathrm{ev}},\varrho}^{\mathrm{API}}}}{(1-\gamma)^2}\,r_n
+\frac{2\gamma^{(J_n+2)/2}R_{\max}}{(1-\gamma)^2}.
\end{equation}
For fixed $\gamma$ and $R_{\max}$ and uniformly bounded concentrability, $r_n\to0$, $\tau_n\to0$, and $J_n\to\infty$ imply convergence of the policy-value error to zero in probability. The first term collects approximation, statistical, and optimization effects across the evaluation stages, while increasing the number of policy updates reduces the geometric remainder. The conclusion is conditional on \eqref{eq:rl:api:uniform}, leaving its quantitative dependence on the network class and sampling scheme to the corresponding approximation and statistical estimates.

A related approach by Kang et al.~\cite{kang2024dapi} uses ReLU residual networks and an adversarial Bellman \emph{moment}, $\sup_O|\frac1n\sum_i O(\boldsymbol X_{j,i},A_{j,i})[Q(\boldsymbol X_{j,i},A_{j,i})-Y_{j,i}(Q)]|$, under exponentially $\alpha$-mixing sampling. It provides another instance of residual-to-policy error propagation, with a different fitting objective and data-dependence assumption. In both approaches, the outer policy-iteration error is distinct from the stagewise optimization tolerance $\epsilon_{\mathrm{opt},j}$; the latter measures how accurately each neural fitting problem is solved.

\subsection{Representation learning}\label{section:app:repssl}

Representation learning seeks features useful for prediction directly from data, reducing reliance on hand-designed features~\cite{bengio2013representation}. In supervised learning, response labels guide the feature-learning process; in self-supervised learning, training signals are constructed from the observations themselves, for example by learning to distinguish positive and negative pairs of augmented samples~\cite{chen2020simple}. The learned features are then used for downstream tasks such as classification or regression. A central theoretical question is how the representation-learning objective controls performance on these tasks.

We discuss three analyses: supervised deep dimension reduction (DDR), contrastive pretraining followed by a neural downstream head, and adversarial self-supervised learning followed by a class-template classifier~\cite{huang2024ddr,duan2025deep,duan2025advssl}. Each method learns an encoder $\phi:\mathcal X\subseteq\mathbb R^d\to\mathbb R^{d_{\mathrm{enc}}}$. The analyses concern a representation objective for DDR, downstream excess risk for contrastive learning, and the SSL objective and its relation to class-template prediction for Adv-SSL.

In the two self-supervised settings, $n$ denotes the pretraining sample size and $m$ the labeled downstream sample size. Network width, depth, parameter count, and norm budget retain the notation $W,L,\mathfrak D,\mathcal K$ used earlier. Empirical optimization gaps are kept as external residuals where an oracle inequality is displayed; the cited convergence results concern exact empirical optimizers.

\subsubsection{Deep dimension reduction}\label{app:rep_learning}

DDR uses labeled observations $(\boldsymbol X,Y)\sim P_{XY}$, with input marginal $P_X$, to learn a low-dimensional representation that retains information relevant to predicting $Y$. In addition, it regularizes the representation distribution toward the Gaussian reference $\nu_{d_{\mathrm{enc}}}=\mathcal N(\boldsymbol 0,\boldsymbol I_{d_{\mathrm{enc}}})$. The reference provides a simple target distribution with independent coordinates.

Huang et al.\ \cite{huang2024ddr} combine dependence on the response and distributional regularization through
\begin{equation}\label{eq:rep_ddr_loss}
\mathcal L(\phi)
=
-\mathrm{dCov}^{\,2}[\phi(\boldsymbol X),Y]
+\lambda_{\mathrm{lag}}D_f(\phi_\#P_X\,\|\,\nu_{d_{\mathrm{enc}}}),
\qquad \lambda_{\mathrm{lag}}>0,
\end{equation}
where $\mathrm{dCov}^{\,2}$ denotes squared distance covariance \cite{szekely2007measuring} and $D_f$ is an $f$-divergence. The negative sign rewards dependence on $Y$, while the divergence penalizes departure from the reference distribution. We assess the excess objective $\mathcal L(\widehat\phi)-\mathcal L^\star$ of the learned encoder $\widehat\phi$, where $\mathcal L^\star:=\inf_{\phi\in\Phi_{\mathrm{adm}}}\mathcal L(\phi)$ is assumed finite for a prescribed admissible representation family $\Phi_{\mathrm{adm}}$.

For a proper lower-semicontinuous convex divergence generator $f$ with $f(1)=0$ and conjugate $f^*$, the variational identity gives, when $\phi_\#P_X\ll\nu_{d_{\mathrm{enc}}}$,
\[
D_f(\phi_\#P_X\,\|\,\nu_{d_{\mathrm{enc}}})
=
\sup_D J(\phi,D),
\qquad
J(\phi,D)
=
\mathbb E D(\phi(\boldsymbol X))
-\mathbb E_{\boldsymbol Z\sim\nu_{d_{\mathrm{enc}}}}f^*(D(\boldsymbol Z)),
\]
with the supremum over admissible measurable discriminators. This variational form allows the divergence to be approximated using encoded observations and Gaussian reference samples, without explicitly evaluating the representation density. Let $\Phi_n\subseteq\Phi_{\mathrm{adm}}$ and $\mathcal D_n$ be output-bounded ReLU encoder and discriminator classes, respectively. Take $n\geq4$ i.i.d.\ labeled observations $\{(\boldsymbol X_i,Y_i)\}_{i=1}^n$ from $P_{XY}$ and, independently, $n$ i.i.d.\ reference samples $\{\boldsymbol Z_i\}_{i=1}^n$ from $\nu_{d_{\mathrm{enc}}}$. Define
\begin{align*}
\mathcal L_{\mathcal D_n}(\phi)
&=-\mathrm{dCov}^{\,2}[\phi(\boldsymbol X),Y]
+\lambda_{\mathrm{lag}}\sup_{D\in\mathcal D_n}J(\phi,D),\\
\widehat{\mathcal L}_{\mathcal D_n,n}(\phi)
&=-\widehat{\mathrm{dCov}^{\,2}}_n[\phi(\boldsymbol X),Y]
+\lambda_{\mathrm{lag}}\sup_{D\in\mathcal D_n}
\frac1n\sum_{i=1}^n
\big[D(\phi(\boldsymbol X_i))-f^*(D(\boldsymbol Z_i))\big].
\end{align*}
Here $\widehat{\mathrm{dCov}^{\,2}}_n$ is the fourth-order U-statistic estimator. Suppose that $\widehat\phi$ minimizes the displayed empirical objective over $\Phi_n$ up to $\epsilon_{\mathrm{opt}}$. Since restricting the discriminator class gives $\mathcal L_{\mathcal D_n}\leq\mathcal L$, the empirical-risk comparison of Section~\ref{sec:decomp} yields, whenever the risks and the terms below are finite,
\begin{equation}\label{eq:rep_ddr_decomp}
\begin{aligned}
\mathbb E\big[\mathcal L(\widehat\phi)-\mathcal L^\star\big]&\leq
\underbrace{\inf_{\phi\in\Phi_n}\mathcal L(\phi)-\mathcal L^\star}
_{\mathcal E_{\mathrm{app}}^{\mathrm{rep}}}
+\underbrace{\sup_{\phi\in\Phi_n}
\big[\mathcal L(\phi)-\mathcal L_{\mathcal D_n}(\phi)\big]}
_{\mathcal E_{\mathrm{disc}}}\\
& \quad +2\mathbb E
\underbrace{\sup_{\phi\in\Phi_n}
\big|\mathcal L_{\mathcal D_n}(\phi)
-\widehat{\mathcal L}_{\mathcal D_n,n}(\phi)\big|}
_{\Delta_n}
+\underbrace{\mathbb E\epsilon_{\mathrm{opt}}}_{\text{optimization error}}.
\end{aligned}
\end{equation}
The expectations include both samples and any randomness in fitting the encoder. The two approximation contributions require different estimates. For $\mathcal E_{\mathrm{app}}^{\mathrm{rep}}$, network approximation must be transferred to the full objective. The distance-covariance contribution can be analyzed through perturbations of pairwise representation distances, together with suitable moment bounds on the response. The divergence contribution requires stability of the induced distributions for the chosen $f$-divergence. These functional estimates determine which encoder approximation norm and regularity conditions are needed.

To bound $\mathcal E_{\mathrm{disc}}$, one approximates discriminators whose variational values approach the unrestricted supremum, uniformly over $\phi\in\Phi_n$. This requires control of both $D(\phi(\boldsymbol X))$ and $f^*(D(\boldsymbol Z))$. When $f^*$ is Lipschitz on the ranges of the comparison and neural discriminators, approximation of their outputs also controls the conjugate term. The comparison discriminators vary with the encoder, making uniformity over $\Phi_n$ essential.

The statistical deviation separates into a distance-covariance term and two empirical-process terms:
\begin{align*}
\Delta_n
&\leq
\underbrace{\sup_{\phi\in\Phi_n}
\left|\widehat{\mathrm{dCov}^{\,2}}_n[\phi(\boldsymbol X),Y]
-\mathrm{dCov}^{\,2}[\phi(\boldsymbol X),Y]\right|}_{\text{distance-covariance estimation}}\\
&\quad+\lambda_{\mathrm{lag}}
\underbrace{\sup_{\substack{\phi\in\Phi_n\\D\in\mathcal D_n}}
\left|\frac1n\sum_{i=1}^n D(\phi(\boldsymbol X_i))
-\mathbb E D(\phi(\boldsymbol X))\right|}_{\text{data-sample fluctuation}}\\
&\quad+\lambda_{\mathrm{lag}}
\underbrace{\sup_{D\in\mathcal D_n}
\left|\frac1n\sum_{i=1}^n f^*(D(\boldsymbol Z_i))
-\mathbb E_{\boldsymbol Z\sim\nu_{d_{\mathrm{enc}}}}f^*(D(\boldsymbol Z))\right|}_{\text{reference-sample fluctuation}}.
\end{align*}
The first supremum is a U-process: each kernel combines four labeled observations, and U-process symmetrization and entropy bounds account for the overlap between these tuples. The other two terms concern the composition classes $\{D\circ\phi:\phi\in\Phi_n,D\in\mathcal D_n\}$ and $\{f^*\circ D:D\in\mathcal D_n\}$. Under suitable envelope and moment conditions, they can be controlled by the empirical-process tools of Section~\ref{sec:classical_stat}. Thus statistical control involves both networks as well as the labeled and reference samples.

For a fixed encoder class, enlarging the discriminator class can reduce its approximation error while changing the complexity entering the statistical estimate. With both classes fixed, increasing the sample size leaves the approximation terms unchanged, while more accurate empirical optimization reduces the optimization residual. These conclusions concern the representation objective; a downstream prediction guarantee additionally requires a link between this objective and the task risk.

\subsubsection{Contrastive pretraining and downstream prediction}\label{app:ssl}

Contrastive learning encourages agreement between representations of positive pairs and distinguishes them from negative pairs. A theoretical connection to downstream prediction is provided by HaoChen et al.~\cite{haochen2021provable}, who analyze spectral contrastive loss and obtain guarantees for linear classification under assumptions on the data-augmentation structure.

The two-stage procedure of Duan et al.\ \cite{duan2025deep} studies an encoder learned from pairwise similarity data and a neural predictor fitted to the resulting features. Its pretraining observations are i.i.d.\ triples
$\{(\boldsymbol X_i,\boldsymbol X_i',\sigma_i)\}_{i=1}^n$, where $\sigma_i\in\{-1,+1\}$ records whether the pair is positive ($+1$) or negative ($-1$). This label describes the relation between the two inputs, rather than the downstream response of either input. The triples are sampled independently from the same joint law. With $\phi^\star$ denoting the target encoder, the conditional model is
\begin{equation}\label{eq:ssl_classical_likelihood}
P_{\phi^\star}(\sigma\mid\boldsymbol x,\boldsymbol x')
=
\frac{1}{1+\exp\{-\sigma\,\phi^\star(\boldsymbol x)^\top
\phi^\star(\boldsymbol x')\}}.
\end{equation}
The inner product of the encoded inputs is a similarity score: a larger value gives a higher conditional probability to a positive pair. Replacing $\phi^\star$ with a candidate encoder $\phi$ defines $P_\phi$. For an output-bounded ReLU class $\Phi_n$, the pretraining estimator is
\begin{equation}\label{eq:ssl_classical_mle}
\widehat\phi\in\arg\max_{\phi\in\Phi_n}
\frac1n\sum_{i=1}^n
\log P_\phi(\sigma_i\mid\boldsymbol X_i,\boldsymbol X_i').
\end{equation}
Let $P_\phi^{\mathrm{pair}}$ denote the joint law of the pretraining triple, with the common input-pair marginal $P_X\otimes P_X$, where $P_X$ is the input distribution. This law combines two independent inputs from $P_X$ with the conditional label probabilities $P_\phi(\sigma\mid\boldsymbol x,\boldsymbol x')$, and includes both positive and negative pairs. The laws $P_{\phi^\star}^{\mathrm{pair}}$ and $P_\phi^{\mathrm{pair}}$ describe the true and candidate models, respectively. At the population level, the excess negative log-likelihood is the candidate model's expected loss minus that of the true model, both averaged under $P_{\phi^\star}^{\mathrm{pair}}$. The common input marginal cancels in the joint likelihood ratio, giving
\[
\begin{aligned}
\mathbb E_{P_{\phi^\star}^{\mathrm{pair}}}
\log\frac{P_{\phi^\star}(\sigma\mid\boldsymbol X,\boldsymbol X')}
{P_\phi(\sigma\mid\boldsymbol X,\boldsymbol X')}=D_{\mathrm{KL}}(P_{\phi^\star}^{\mathrm{pair}}\,\|\,P_\phi^{\mathrm{pair}})\gtrsim d_{\mathrm{Hel}}^2(P_{\phi^\star}^{\mathrm{pair}},P_\phi^{\mathrm{pair}}),
\end{aligned}
\]
where $d_{\mathrm{Hel}}$ is Hellinger distance. Thus a small population excess loss guarantees a small squared Hellinger distance between the candidate and true triple distributions. Obtaining this population control from the empirical likelihood fit requires finite-sample analysis, using likelihood ratios and bracketing entropy for the induced distribution class together with encoder approximation; see \cite[Supplementary Material, Section~4.3]{duan2025deep}.

A labeled i.i.d.\ sample $\{(\boldsymbol X_j^{\mathrm{lab}},Y_j)\}_{j=1}^m$, independent of the pretraining sample, is then used to fit a norm-constrained ReLU head
$h:\mathbb R^{d_{\mathrm{enc}}}\to\mathbb R$, with $\widehat\phi$ held fixed. Let $h^\star$ denote the target downstream head associated with $\phi^\star$, and write
\[
\mathcal R(h\circ\phi)
=\mathbb E\ell(h(\phi(\boldsymbol X)),Y),
\qquad
\mathcal E(h\circ\phi)
=\mathcal R(h\circ\phi)-\mathcal R(h^\star\circ\phi^\star).
\]
For regression, $\ell$ may be squared loss on a bounded range; for binary classification it may be logistic loss on bounded scores. The theorem below requires boundedness and Lipschitz continuity of the loss on the relevant ranges.

To turn accuracy of the pair distribution into a prediction guarantee, an additional link to the downstream response is needed. The informative property in \cite[Lemma 4.1]{duan2025deep} provides this link. This result uses the paper's latent-variable model linking pretraining to the downstream response, in addition to the similarity likelihood \eqref{eq:ssl_classical_likelihood}: the same target representation enters the downstream regression or classification model through $h^\star$. Let $P_{\phi,h}^{\mathrm{down}}$ denote the induced joint law of input and response in that model. For an encoder $\phi$ in the hypothesis class, an admissible head $h$ can be chosen so that
\[
d_{\mathrm{TV}}\big(P_{\phi,h}^{\mathrm{down}},
P_{\widetilde\phi,h^\star}^{\mathrm{down}}\big)
\leq C_{\mathrm{inf}}\,
d_{\mathrm{Hel}}\big(P_\phi^{\mathrm{pair}},
P_{\widetilde\phi}^{\mathrm{pair}}\big),
\]
where $\widetilde\phi$ is a best $L^2(P_X)$ approximation to $\phi^\star$ in $\Phi_n$. The distributional comparison accommodates the orthogonal invariance of the representation likelihood in \eqref{eq:ssl_classical_likelihood}.

Let $\mathcal H_{\mathcal K}$ be the head class with weight-norm budget $\mathcal K$, and let $\widehat{\mathcal R}_m$ be the sample average of the prediction loss on the labeled data. Let $\widehat h\in\mathcal H_{\mathcal K}$ be the fitted head, with empirical optimization error at most $\epsilon_{\mathrm{opt}}^{\mathrm{head}}$. Conditional on $\widehat\phi$, the empirical-risk comparison of Section~\ref{sec:decomp} gives
\begin{equation}\label{eq:ssl_head_oracle}
\begin{aligned}
\mathcal E(\widehat h\circ\widehat\phi)
\leq{}&
2\underbrace{\sup_{h\in\mathcal H_{\mathcal K}}
\big|\mathcal R(h\circ\widehat\phi)
-\widehat{\mathcal R}_m(h\circ\widehat\phi)\big|}_{\Delta_m(\widehat\phi)\;\text{(downstream statistical error)}}\\
&+\underbrace{\inf_{h\in\mathcal H_{\mathcal K}}\mathcal E(h\circ\widehat\phi)}_{\text{representation and head error}}
+\underbrace{\epsilon_{\mathrm{opt}}^{\mathrm{head}}}_{\text{optimization residual}}.
\end{aligned}
\end{equation}
The first term is controlled by the complexity of the head class on the frozen features. The infimum contains the effect of pretraining as well as head approximation. To bound it, the analysis inserts the downstream law associated with $\widetilde\phi$: boundedness of the loss converts total-variation bounds into differences of expected losses, and the informative property controls the comparison with the learned representation by the pair-distribution Hellinger error. Comparing $\widetilde\phi$ with $\phi^\star$ adds encoder approximation, while restricting the predictor to $\mathcal H_{\mathcal K}$ adds head approximation. The Lipschitz bound supplied by the head's norm constraint also multiplies the encoder approximation error by $\mathcal K$; these contributions are separated in \cite[Supplementary Material, Lemmas~3.1--3.3]{duan2025deep}.

\begin{theorem}[{\cite[Theorem 4.2]{duan2025deep}}]\label{thm:ssl_contrastive}
Consider the model above on $\mathcal X\subseteq[0,1]^d$, with coordinatewise H{\"o}lder-smooth encoder $\phi^\star$ of order $\beta_0$ and a H{\"o}lder-smooth head $h^\star$ of order $0<\beta_1<(d_{\mathrm{enc}}+3)/2$ on a bounded representation domain, rescaled to $[0,1]^{d_{\mathrm{enc}}}$. Assume the bounded-Lipschitz loss conditions of \cite{duan2025deep}, and choose the encoder architecture and the norm-constrained head as in its Theorem 4.2. For exact empirical optimizers, $n\geq m\geq2$, and $0<\delta<1$, with probability at least $1-\delta$,
\begin{equation}\label{eq:ssl_classical_decomp}
\begin{aligned}
\mathcal E(\widehat h\circ\widehat\phi) \lesssim (\log n)^3 n^{-\beta_0/(2\beta_0+d)} +\sqrt{\frac{\log(2/\delta)}{m}} + \left(n^{-\beta_0/(2\beta_0+d)} +\frac{\log m}{\sqrt m}\right)^{2\beta_1/(d_{\mathrm{enc}}+3)}.
\end{aligned}
\end{equation}
Here dimensions, smoothness parameters, loss bounds, and H{\"o}lder radii are held fixed; their dependence is absorbed into the implicit constant.
\end{theorem}

The first term is the pretraining contribution, and the second is the downstream concentration term. The last term results from choosing the head's norm budget jointly with the encoder approximation accuracy. With the fixed quantities absorbed into constants, the $\mathcal K$-dependent contributions in \cite[Supplementary Material, Section~4.5]{duan2025deep} have the form
\[
\underbrace{\mathcal K n^{-\beta_0/(2\beta_0+d)}}_{\text{transferred encoder approximation}}
+\underbrace{\mathcal K\frac{\log m}{\sqrt m}}_{\text{downstream statistical complexity}}
+\underbrace{\mathcal K^{-2\beta_1/(d_{\mathrm{enc}}+3-2\beta_1)}}_{\text{head approximation}}.
\]
Increasing $\mathcal K$ improves head approximation but amplifies encoder approximation and increases downstream statistical complexity. Balancing these contributions gives the power $2\beta_1/(d_{\mathrm{enc}}+3)$ of $n^{-\beta_0/(2\beta_0+d)}+\log m/\sqrt m$ in \eqref{eq:ssl_classical_decomp}. Thus pretraining accuracy affects both the direct pretraining term and the choice of downstream head complexity. If $n\gtrsim m^{(2\beta_0+d)/(2\beta_0)}$, the bound reduces, up to logarithmic factors in $n,m$, to
\[
\mathcal E(\widehat h\circ\widehat\phi)
=
\widetilde O\big(m^{-1/2}+m^{-\beta_1/(d_{\mathrm{enc}}+3)}\big)
\]
for fixed confidence. In the stated smoothness range, $\beta_1/(d_{\mathrm{enc}}+3)<1/2$, so the representation dimension enters the slower-decaying term. The result quantifies how abundant pretraining data can reduce the labeled-data burden in this model.

\subsubsection{Adversarial self-supervised learning}\label{app:ssl_adv}

Adv-SSL \cite{duan2025advssl} learns an encoder from two augmentations $(\boldsymbol X_1,\boldsymbol X_2)$ of each unlabeled source observation. The encoder maps each augmented input to a feature vector that is subsequently used for prediction. We first assess the excess pretraining objective $\mathcal L_{\mathrm{SSL}}(\widehat\phi)-\mathcal L_{\mathrm{SSL}}(\phi^\star)$, where $\phi^\star$ minimizes the population objective over the admissible representation family. The downstream task uses a class-template classifier $Q_{\widehat\phi}$, evaluated by its classification error under the target input--label distribution $P_{\mathrm t}$. Relating these two quantities requires augmentation geometry and source--target transfer conditions.

With fixed $\lambda_{\mathrm{ssl}}>0$, the population objective to be minimized combines agreement between the augmented representations with control of their cross-moment:
\begin{equation}\label{eq:ssl_adv_population}
\begin{aligned}
\mathcal L_{\mathrm{SSL}}(\phi)
&=\mathcal L_{\mathrm{align}}(\phi)+\lambda_{\mathrm{ssl}}\|\boldsymbol C_\phi\|_{\mathrm F}^2,\\
\mathcal L_{\mathrm{align}}(\phi)
&=\mathbb E\|\phi(\boldsymbol X_1)-\phi(\boldsymbol X_2)\|_2^2,\\
\boldsymbol C_\phi&=\mathbb E[\phi(\boldsymbol X_1)\phi(\boldsymbol X_2)^\top]-\boldsymbol I_{d_{\mathrm{enc}}}.
\end{aligned}
\end{equation}
The alignment term encourages the two views of an observation to have similar features, but also vanishes for a constant encoder. The moment penalty adds a requirement on the feature structure: products of the same coordinate across views should have expectation close to one, whereas products of distinct coordinates should have expectation close to zero. For example, the all-zero representation has zero alignment loss but a nonzero moment penalty. The parameter $\lambda_{\mathrm{ssl}}$ balances these two requirements. Let $\widehat{\boldsymbol C}_\phi$ replace the expectation by the average of $n$ independent augmented pairs. For fixed $\phi$, the quadratic plug-in bias is explicit:
\[
\mathbb E\|\widehat{\boldsymbol C}_\phi\|_{\mathrm F}^2
=\|\boldsymbol C_\phi\|_{\mathrm F}^2
+\mathbb E\|\widehat{\boldsymbol C}_\phi-\boldsymbol C_\phi\|_{\mathrm F}^2.
\]
The squared regularizer therefore has a nonnegative plug-in bias.

To rewrite the moment penalty, the paper introduces an auxiliary matrix $\boldsymbol G_{\mathrm{dual}}\in\mathbb R^{d_{\mathrm{enc}}\times d_{\mathrm{enc}}}$. Its entries weight the corresponding entries of $\boldsymbol C_\phi$, so maximizing their Frobenius inner product over a suitable bounded set measures the moment discrepancy for the current encoder. Define
\begin{align*}
\mathcal G_{\mathrm{dual}}(\phi)
&=\{\boldsymbol G_{\mathrm{dual}}:\|\boldsymbol G_{\mathrm{dual}}\|_{\mathrm F}\leq\|\boldsymbol C_\phi\|_{\mathrm F}\},\\
\widehat{\mathcal G}_{\mathrm{dual}}(\phi)
&=\{\boldsymbol G_{\mathrm{dual}}:\|\boldsymbol G_{\mathrm{dual}}\|_{\mathrm F}\leq\|\widehat{\boldsymbol C}_\phi\|_{\mathrm F}\}.
\end{align*}
Cauchy--Schwarz, with equality at $\boldsymbol G_{\mathrm{dual}}=\boldsymbol C_\phi$, gives
\[
\|\boldsymbol C_\phi\|_{\mathrm F}^2
=\sup_{\boldsymbol G_{\mathrm{dual}}\in\mathcal G_{\mathrm{dual}}(\phi)}
\langle \boldsymbol C_\phi,\boldsymbol G_{\mathrm{dual}}\rangle_{\mathrm F}.
\]
This identity preserves the original population objective while expressing the moment penalty as an inner maximization. For a fixed auxiliary matrix, the moment term is a linear functional of $\boldsymbol C_\phi$ and can be estimated by sample averages. Writing
$\mathcal L_{\mathrm{SSL}}(\phi,\boldsymbol G_{\mathrm{dual}})
=\mathcal L_{\mathrm{align}}(\phi)+
\lambda_{\mathrm{ssl}}\langle \boldsymbol C_\phi,\boldsymbol G_{\mathrm{dual}}\rangle_{\mathrm F}$,
the sample objective is
\begin{align*}
\widehat{\mathcal L}_{\mathrm{SSL}}(\phi,\boldsymbol G_{\mathrm{dual}})
={}&\frac1n\sum_{i=1}^n
\|\phi(\boldsymbol X_{i,1})-\phi(\boldsymbol X_{i,2})\|_2^2\\
&+\lambda_{\mathrm{ssl}}
\left\langle
\frac1n\sum_{i=1}^n
\phi(\boldsymbol X_{i,1})\phi(\boldsymbol X_{i,2})^\top
-\boldsymbol I_{d_{\mathrm{enc}}},\boldsymbol G_{\mathrm{dual}}
\right\rangle_{\mathrm F}.
\end{align*}
For an encoder network class $\Phi_n$, the empirical estimator is
\begin{equation}\label{eq:ssl_adv_minimax}
\widehat\phi\in\argmin_{\phi\in\Phi_n}
\sup_{\boldsymbol G_{\mathrm{dual}}\in\widehat{\mathcal G}_{\mathrm{dual}}(\phi)}
\widehat{\mathcal L}_{\mathrm{SSL}}(\phi,\boldsymbol G_{\mathrm{dual}}).
\end{equation}
The encoder is selected by minimizing the sum of the alignment loss and the penalty evaluated by the inner maximization; the auxiliary matrix is needed only during training. For each fixed $(\phi,\boldsymbol G_{\mathrm{dual}})$ the sample objective is unbiased. Maximizing over the sample-dependent dual set recovers the quadratic plug-in objective, so the analysis controls both empirical fluctuations and changes in the dual set.

The following oracle bound controls the population objective of the learned encoder relative to its optimum, with these two statistical contributions entering alongside encoder approximation and optimization error. Suppose $\|\phi(\boldsymbol x)\|_2\leq B_\phi$ for $\phi\in\Phi_n$ and let $B_G=B_\phi^2+\sqrt{d_{\mathrm{enc}}}$. Both dual sets lie in the ball of radius $B_G$, and their radii differ by at most $\|\boldsymbol C_\phi-\widehat{\boldsymbol C}_\phi\|_{\mathrm F}$. Comparing the objectives on this common ball and accounting for the change of radius, the empirical-risk comparison of Section~\ref{sec:decomp} gives, for an approximate solution of \eqref{eq:ssl_adv_minimax},
\begin{equation}\label{eq:ssl_adv_decomp}
\begin{aligned}
\mathcal L_{\mathrm{SSL}}(\widehat\phi)-\mathcal L_{\mathrm{SSL}}(\phi^\star)& \leq \underbrace{\inf_{\phi\in\Phi_n}
[\mathcal L_{\mathrm{SSL}}(\phi)-\mathcal L_{\mathrm{SSL}}(\phi^\star)]}_{\text{encoder approximation}}\\
&\qquad+2\underbrace{\sup_{\substack{\phi\in\Phi_n\\\|\boldsymbol G_{\mathrm{dual}}\|_{\mathrm F}\leq B_G}}
\big|\mathcal L_{\mathrm{SSL}}(\phi,\boldsymbol G_{\mathrm{dual}})
-\widehat{\mathcal L}_{\mathrm{SSL}}(\phi,\boldsymbol G_{\mathrm{dual}})\big|}_{\Delta_{\mathrm{SSL}}\;\text{(empirical fluctuation)}}\\
&\qquad+2\underbrace{\lambda_{\mathrm{ssl}}B_G
\sup_{\phi\in\Phi_n}\|\boldsymbol C_\phi-\widehat{\boldsymbol C}_\phi\|_{\mathrm F}}_{\Delta_{\mathrm{set}}\;\text{(dual-set variation)}}
+\underbrace{\epsilon_{\mathrm{opt}}}_{\text{optimization residual}}.
\end{aligned}
\end{equation}
Let $\beta_0$ denote the coordinatewise H{\"o}lder smoothness of $\phi^\star$. The approximation--statistical balance considered in \cite[Appendix F.4.6]{duan2025advssl} can be stated conditionally as follows. With $B_\phi$, $d_{\mathrm{enc}}$, $\lambda_{\mathrm{ssl}}$, and depth $L$ fixed, suppose that the same norm-constrained encoder class $\Phi_n$, with budget $\mathcal K$, satisfies
\[
\begin{aligned}
\inf_{\phi\in\Phi_n}
[\mathcal L_{\mathrm{SSL}}(\phi)-\mathcal L_{\mathrm{SSL}}(\phi^\star)]
&\lesssim \mathcal K^{-\beta_0/(d+1)},\\
\mathbb E(\Delta_{\mathrm{SSL}}+\Delta_{\mathrm{set}})
&\lesssim \frac{\mathcal K}{\sqrt n},
\end{aligned}
\]
with constants independent of $n$ and $\mathcal K$. The first estimate requires approximation within the constrained encoder class; the second requires uniform empirical-process control under the augmented-data law, including the dual-set term. Taking expectations in \eqref{eq:ssl_adv_decomp} gives
\[
\mathbb E[\mathcal L_{\mathrm{SSL}}(\widehat\phi)-\mathcal L_{\mathrm{SSL}}(\phi^\star)]
\lesssim
\frac{\mathcal K}{\sqrt n}
+\mathcal K^{-\beta_0/(d+1)}
+\mathbb E\epsilon_{\mathrm{opt}}.
\]
For exact empirical minimization, balancing the approximation and statistical terms gives
$\mathcal K\asymp n^{(d+1)/(2(\beta_0+d+1))}$
and an expected excess-objective bound of order $n^{-\beta_0/(2(\beta_0+d+1))}$. If additionally $\mathcal L_{\mathrm{SSL}}(\phi^\star)=0$, the same bound applies to $\mathbb E\mathcal L_{\mathrm{SSL}}(\widehat\phi)$.

Downstream prediction uses the learned encoder with labeled target-domain observations. For each class $k$, average the two augmented representations of its $m_k>0$ labeled observations to obtain a template $\widehat{\boldsymbol\mu}_k$, and classify by
\[
Q_{\widehat\phi}(\boldsymbol x)
=\arg\max_k
\langle\widehat{\boldsymbol\mu}_k,\widehat\phi(\boldsymbol x)\rangle,
\qquad
\mathrm{Err}(Q_{\widehat\phi})
=P_{\mathrm t}\{Q_{\widehat\phi}(\boldsymbol X)\ne Y\}.
\]
On the covered subsets of each class, the augmentation geometry in \cite{duan2025advssl} connects distinct observations through nearby augmented views. Combined with alignment control and encoder Lipschitz regularity, it relates consistency between two views of one observation to concentration of features across observations within a class. Classification then requires the correct class template to have a sufficiently larger inner-product score than competing templates; global cross-moment control alone does not ensure this margin. Finite-sample errors in the templates perturb these scores and must be small relative to the separation.

Transferring the source-domain alignment and moment estimates to the target domain additionally involves class-conditional distribution differences and changes in class proportions. Increasing the source sample size improves pretraining statistical control, whereas increasing the labeled sample size in a class improves estimation of its template. A fixed source--target discrepancy does not disappear merely by collecting more source observations. These geometric, sampling, and transfer considerations are additional to the conditional pretraining-objective bound above.

When the classwise augmentation conditions hold only on subsets of the source and target domains, downstream classification guarantees must account for the uncovered probability mass and use feature separation and template estimation bounds to control errors on the covered subsets.

Taken together, these examples illustrate three different uses of neural approximation theory. DDR approximates an encoder and divergence witnesses; contrastive learning combines a learned encoder with a norm-controlled prediction head; Adv-SSL separates approximation and empirical-process errors in the pretraining objective, with additional geometric and sampling considerations for class-template prediction.

\subsection{Large language models}
\label{section:app:llm}

Many large language models (LLMs) use Transformer architectures~\cite{vaswani2017attention} and are pretrained to predict the next token from its preceding context~\cite{radford2018improving,brown2020language}. At inference time, a prompt can specify a task or provide demonstrations, allowing the model to perform new tasks without updating its parameters~\cite{brown2020language}. In-context learning (ICL) uses these demonstrations, while chain-of-thought (CoT) prompting includes intermediate reasoning steps in the examples~\cite{wei2022chain}.

We first review Transformer approximation as architectural background for LLMs. Sequence regression~\cite{jiao2026approximation} then provides a representative example of how the error-analysis framework connects approximation, statistical complexity, and data dependence to prediction accuracy. We next turn to autoregressive pretraining and prompting within the latent-task model of \cite{jiao2026prompt}. The two settings have different evaluation objects: squared prediction error for numerical responses, and conditional response-probability error for language generation. In the latter setting, the analysis separates learning the response distribution from identifying the task specified by a prompt.

\subsubsection{Transformer approximation and regression error analysis}\label{app:llm_reg}
Universal approximation for Transformers was established in \cite{yun2019transformers} for continuous permutation-equivariant sequence maps with compact support, and, with positional encodings, for general continuous sequence maps on compact domains. The proof uses self-attention to construct a contextual mapping that distinguishes tokens according to their surrounding sequence. Kajitsuka and Sato~\cite{kajitsuka2023transformers} showed that a single self-attention layer with low-rank weight matrices suffices for the relevant contextual mapping; combined with two feed-forward networks, one before and one after attention, it yields universal approximation for continuous permutation-equivariant maps on compact domains. In both results, approximation is measured in $L^p$ with respect to Lebesgue measure on the input domain, with $1\leq p<\infty$.

Quantitative approximation theory further studies the network complexity required for a given accuracy and its dependence on dimension. A different construction was developed in \cite{jiao2025transformers}, where a one-head, one-layer self-attention module combined with feed-forward blocks is used to approximate H\"older-continuous targets with exponent $0<\beta\leq1$ through a Kolmogorov--Arnold representation rather than a contextual-mapping argument. For $L^\infty$ approximation accuracy $\epsilon$ on the input cube, the ReLU--floor construction uses $O(\log(1/\epsilon))$ feed-forward layers with width $O(\epsilon^{-2/\beta}\log(1/\epsilon))$, while richer activation systems permit constant-width constructions. The representation system itself can therefore change the dimension dependence of the approximation problem. More recently, \cite{he2026approximation} derived approximation upper and lower bounds for H\"older classes with standard Transformers using softmax attention, ReLU activations, and residual connections.

To connect approximation to learning from finite data, Jiao et al.~\cite{jiao2026approximation} study nonparametric regression with sequential observations. Let $N_{\mathrm{reg}}$ denote the number of observed time points, $L_{\mathrm{win}}$ a fixed input-window length, and $d_x$ the dimension of each input vector. Consider
\begin{equation}\label{eq:llm_reg_model}
Y_t=f_0(\boldsymbol H_t)+\xi_t,
\qquad
\boldsymbol H_t=(\boldsymbol X_{t-L_{\mathrm{win}}+1},\ldots,\boldsymbol X_t)
\in[0,1]^{d_x\times L_{\mathrm{win}}},
\end{equation}
where $Y_t\in\mathbb R$ and $f_0$ is the conditional mean of $Y_t$ given $\boldsymbol H_t$. We take the noise variables $\xi_t$ to be independent and identically distributed, centered and sub-Gaussian, and independent of the input process. Assume that the input process is stationary and $\beta$-mixing, and let $\Pi$ be the common distribution of its length-$L_{\mathrm{win}}$ windows. The final quantity of interest is the mean squared prediction error $\mathbb E\|\widetilde f-f_0\|_{L^2(\Pi)}^2$, where $\widetilde f$ is the fitted predictor specified below and the outer expectation is over the training observations.

Let $\mathcal F_{\mathrm{Tr}}$ be a class of scalar predictors $f(\boldsymbol H)=\langle T_{\boldsymbol\theta}(\boldsymbol H),\boldsymbol A\rangle_{\mathrm F}$, where $T_{\boldsymbol\theta}$ is a Transformer with softmax attention and ReLU feed-forward blocks, and $\boldsymbol A$ is a trainable matrix of matching dimension. This linear readout converts the Transformer's matrix output into a scalar prediction. From observations $\{(\boldsymbol X_t,Y_t)\}_{t=1}^{N_{\mathrm{reg}}}$, the sliding-window estimator minimizes the empirical squared loss:
\begin{equation}\label{eq:llm_reg_erm}
\widehat f\in\argmin_{f\in\mathcal F_{\mathrm{Tr}}}
\frac{1}{N_{\mathrm{reg}}-L_{\mathrm{win}}+1}
\sum_{t=L_{\mathrm{win}}}^{N_{\mathrm{reg}}}
\bigl(f(\boldsymbol H_t)-Y_t\bigr)^2.
\end{equation}
Overlapping windows are generally dependent, even when the input observations are independent. Stationarity gives them a common population risk,
\[
\mathcal R_{\mathrm{Tr}}(f)
=\mathbb E\bigl[(f(\boldsymbol H)-Y)^2\bigr],
\]
and conditional-mean orthogonality yields
\begin{equation}\label{eq:llm_reg_identity}
\mathcal R_{\mathrm{Tr}}(f)-\mathcal R_{\mathrm{Tr}}(f_0)
=\|f-f_0\|_{L^2(\Pi)}^2.
\end{equation}
Thus, excess squared loss directly measures the prediction error of interest. The reported estimator is the clipped fit $\widetilde f=C_B\widehat f$, where $C_Bf=\max\{-B,\min\{B,f\}\}$ and $B\asymp\log N_{\mathrm{reg}}$. Clipping supplies the bounded output envelope used in the statistical analysis; the empirical minimization in \eqref{eq:llm_reg_erm} is performed before this clipping.

Suppose that $f_0$ belongs to a fixed H\"older ball of order $0<\beta\leq1$ on $[0,1]^{d_x\times L_{\mathrm{win}}}$ and that $\Pi$ has a bounded density there. Use two Transformer blocks with embedding dimension $d_x$, one attention head of dimension one per block, and growing feed-forward width, as in the $L^2$ construction of \cite[Theorem~2.3]{jiao2026approximation}. Let $\mathfrak D$ denote its total parameter count, including the scalar readout. With these architectural dimensions fixed, $\mathfrak D$ is of the same order as the feed-forward width. The clipped least-squares analysis in \cite[Appendix~A.2]{jiao2026approximation}, together with its network-complexity estimate, gives
\begin{equation}\label{eq:llm_reg_decomp}
\begin{aligned}
\mathbb E\|\widetilde f-f_0\|_{L^2(\Pi)}^2
\lesssim{}&
\underbrace{\inf_{f\in\mathcal F_{\mathrm{Tr}}}
\|f-f_0\|_{L^2(\Pi)}^2}_{\text{approximation}}\\
&+\underbrace{\frac{B^2a\,\mathfrak D^2
\log(N_{\mathrm{reg}}^2B)}{N_{\mathrm{reg}}}}_{\text{statistical estimation}}
+\underbrace{\frac{B^2N_{\mathrm{reg}}}{a}\bar\beta_a}_{\text{data dependence}}.
\end{aligned}
\end{equation}
Here $a$ is an auxiliary blocking separation used in the dependence analysis; the empirical estimator still uses all sliding windows in \eqref{eq:llm_reg_erm}. The coefficient $\bar\beta_a$ describes mixing of the window--response process at that separation. With independent noise and fixed $L_{\mathrm{win}}$, this process inherits the decay of the input mixing coefficients, up to a fixed shift in the separation. The bounds hold for sufficiently large $N_{\mathrm{reg}}$, with $d_x$, $L_{\mathrm{win}}$, the H\"older radius, and the noise parameters held fixed. Exact empirical minimization is assumed.

The approximation theorem controls the first term in $L^2$; the density bound transfers Lebesgue approximation to $L^2(\Pi)$. For the stated construction, the approximation and pseudo-dimension estimates are
\[
\inf_{f\in\mathcal F_{\mathrm{Tr}}}\|f-f_0\|_{L^2(\Pi)}^2
\lesssim\mathfrak D^{-2\beta/(d_xL_{\mathrm{win}})},
\qquad
\operatorname{Pdim}(\mathcal F_{\mathrm{Tr}})
\lesssim\mathfrak D^2.
\]
The latter estimate, proved in \cite[Lemma~A.12]{jiao2026approximation} using a computational model that includes the exponential operations in softmax, controls the covering numbers of the clipped class. These covering estimates produce the statistical term in \eqref{eq:llm_reg_decomp}. The dependence term is the cost of replacing separated blocks by independent copies.

For geometrically mixing data, suppose $\bar\beta_a\leq C_{\mathrm{mix}}\exp(-b_{\mathrm{mix}}a^\eta)$ with $\eta>0$. Choose
\[
a\asymp(\log N_{\mathrm{reg}})^{1/\eta},
\qquad
\mathfrak D\asymp
N_{\mathrm{reg}}^{d_xL_{\mathrm{win}}/(2\beta+2d_xL_{\mathrm{win}})},
\]
with a sufficiently large constant in the separation choice. The dependence contribution is then negligible at the displayed scale, and the approximation and statistical bounds give the guarantee of \cite[Theorem~2.7]{jiao2026approximation}:
\begin{equation}\label{eq:llm_reg_rate}
\mathbb E\|\widetilde f-f_0\|_{L^2(\Pi)}^2
\lesssim
N_{\mathrm{reg}}^{-\beta/(\beta+d_xL_{\mathrm{win}})}
(\log N_{\mathrm{reg}})^{3+1/\eta}.
\end{equation}
The decomposition identifies two tradeoffs: increasing the network size improves approximation but enlarges the statistical term, while increasing the blocking separation reduces dependence at a statistical cost. The window length enters the approximation dimension $d_xL_{\mathrm{win}}$, so the guarantee is stated for a fixed window. The cited rate is not minimax optimal; the paper attributes this limitation to the quadratic parameter dependence of its complexity bound. Sharper statistical complexity estimates could therefore improve the guarantee without changing the approximation construction.

This regression example connects Transformer expressivity to finite-sample prediction error for numerical responses. Autoregressive language modelling uses token histories and a logarithmic loss, so its probability error requires a separate analysis. In the conditional-probability model considered next, a finite-history exact-representation result supplies the realizability used to eliminate the memorization term in \eqref{eq:llm_gen_mem}.

\subsubsection{Autoregressive pretraining and prompt comprehension}\label{app:llm_pretraining}
We follow the latent-task framework of \cite{jiao2026prompt} to introduce pretraining error and zero-shot task ambiguity; Section~\ref{app:llm_understand} then examines ICL and CoT in the same model. We first compare the trained model's response probabilities with those of the true marginal response distribution, both conditioned on the same prompt. For prompt comprehension, we then compare the learned probabilities with those conditional on the task with the largest posterior probability, separating pretraining error from task ambiguity.

Let
\[
\mathcal D_{N_{\mathrm{seq}}} = \{\boldsymbol d^{(i)}\}_{i=1}^{N_{\mathrm{seq}}}
\]
be a pretraining corpus consisting of $N_{\mathrm{seq}}$ sequences of common length $L_{\mathrm{seq}}$, generated under the true data-generating distribution $p^*$. We write $N_{\mathrm{tok}}:=N_{\mathrm{seq}}L_{\mathrm{seq}}$ for the corresponding number of token-level prediction examples. This counts the prediction terms in the empirical loss; these terms need not correspond to independent observations. Throughout this and the next subsection, $\mathcal P$ denotes a full prompt, $x$ denotes the target query contained in that prompt, $y$ denotes a single predicted answer or token, and $\boldsymbol y$ denotes a response sequence of length $L_{\mathrm{res}}$. Let $\mathcal V$ denote the finite vocabulary. Writing $t_j^{(i)}\in\mathcal V$ for the $j$th token and $\boldsymbol h_j^{(i)}$ for its preceding context in the $i$th sequence, the ideal autoregressive empirical risk minimizer is
\begin{equation}
\label{eq:llm_erm}
\widehat p_{\mathrm{ERM}} \in \arg\min_{p\in\mathcal F} -\frac{1}{N_{\mathrm{tok}}} \sum_{i=1}^{N_{\mathrm{seq}}}\sum_{j=1}^{L_{\mathrm{seq}}} \log p\left(t_j^{(i)} \mid \boldsymbol h_j^{(i)}\right),
\end{equation}
where $\mathcal F$ (playing the role of the hypothesis class $\mathcal{F}_{NN}$ of Section~\ref{sec:decomp}) is a Transformer class. Assume that both the empirical and population risk minima over $\mathcal F$ are attained. The actual trained model is denoted by $\widehat p$. Autoregressive training thus estimates a conditional probability distribution through next-token prediction.

At the population level, negative log-likelihood has the usual Kullback--Leibler interpretation. If $p^*$ denotes the true data-generating distribution, then
\[
\begin{aligned}
\mathcal R(p)
= \mathbb E_{(t,\boldsymbol h)\sim p^*} \left[-\log p(t\mid\boldsymbol h)\right] = \mathbb E_{\boldsymbol h\sim p^*} D_{\mathrm{KL}}\left(p^*(\cdot\mid\boldsymbol h) \,\|\, p(\cdot\mid\boldsymbol h)\right) + C_{p^*},
\end{aligned}
\]
where $C_{p^*}$ does not depend on $p$. Hence minimizing the population negative log-likelihood amounts to minimizing the conditional KL divergence averaged over the history distribution ($\mathcal R(p)$ plays the role of the population risk $\mathcal{L}$ in the framework of Section~\ref{sec:decomp}). Because the prompting analysis compares responses under a specified prompt, we express the next comparison directly on response probabilities, using a pointwise triangle inequality. Converting the average conditional KL discrepancy into a uniform probability-error bound would additionally require tracking the probabilities of the relevant histories.

We assume that $\epsilon_{\mathrm{opt}}^{\mathrm{pred}}$ bounds $|\widehat p(\boldsymbol y\mid\mathcal P)-\widehat p_{\mathrm{ERM}}(\boldsymbol y\mid\mathcal P)|$ over the prompt--response pairs considered below, including single-token responses. This prediction-level residual is distinct from the empirical-risk optimization gap $\epsilon_{\mathrm{opt}}$ used earlier and is retained as an external quantity. Introducing a population risk minimizer $\bar p\in\arg\min_{p\in\mathcal F}\mathcal R(p)$ gives the triangle inequality
\begin{align}
&\left| \widehat p(\boldsymbol y\mid\mathcal P) - p^*(\boldsymbol y\mid\mathcal P) \right| \label{eq:llm_gen_mem} \\
&\quad\leq \underbrace{\left| \widehat p_{\mathrm{ERM}}(\boldsymbol y\mid\mathcal P) - \bar p(\boldsymbol y\mid\mathcal P) \right|}_{\mathcal E_{\mathrm{gen}}}
+ \underbrace{\left| \bar p(\boldsymbol y\mid\mathcal P) - p^*(\boldsymbol y\mid\mathcal P) \right|}_{\mathcal E_{\mathrm{mem}}} +  \underbrace{\left| \widehat p(\boldsymbol y\mid\mathcal P)-\widehat p_{\mathrm{ERM}}(\boldsymbol y\mid\mathcal P)\right|}_{\leq\epsilon_{\mathrm{opt}}^{\mathrm{pred}}}. \nonumber
\end{align}
For each prompt--response pair, the memorization error $\mathcal E_{\mathrm{mem}}$ measures the discrepancy between the population minimizer and the true response probability, while the generalization error $\mathcal E_{\mathrm{gen}}$ compares the empirical and population minimizers.

All prompts and response prefixes considered below have positive probability under $p^*$, and the true next-token probabilities on the relevant finite token--history space are positive. The number of non-padding tokens in the full prompt, plus $L_{\mathrm{res}}$, must not exceed $L_{\mathrm{seq}}$.

To control the memorization error, one can show that a sufficiently wide Transformer can exactly realize prescribed conditional probability vectors on a finite collection of input sequences satisfying the tokenwise-separability conditions of the memorization theorem. Under Assumption~9 of \cite{jiao2026prompt}, a depth-one Transformer with width of order $|\mathcal V|^{L_{\mathrm{seq}}+1}$ realizes the relevant finite-history conditionals. If $\mathcal F$ contains such a network realizing the true conditionals at every positive-probability history entering the population risk, the population minimum attains zero average KL discrepancy. The population minimizer $\bar p$ therefore agrees with the true conditionals on these histories, so that the memorization component vanishes for the prompt--response pairs specified above:
\[
\mathcal E_{\mathrm{mem}}=0.
\]
Next-token errors enter response probabilities through the autoregressive product. Writing $\boldsymbol h_j=\mathcal P\circ\boldsymbol y_{<j}$ for the prompt followed by the first $j-1$ response tokens, a telescoping expansion gives
\begin{equation}
\label{eq:llm_gen_bound}
\left| \widehat p(\boldsymbol y\mid\mathcal P) - p^*(\boldsymbol y\mid\mathcal P) \right|
\leq \sum_{j=1}^{L_{\mathrm{res}}}
\left|\widehat p(y_j\mid\boldsymbol h_j)-p^*(y_j\mid\boldsymbol h_j)\right|.
\end{equation}
Here $\circ$ denotes sequence concatenation and $y_j$ is the $j$th response token. Thus a uniform next-token probability error contributes at most $L_{\mathrm{res}}$ times that error to a response probability. For the prompting analyses below, collect the prediction error of the trained model into
\begin{equation}
\label{eq:llm_pretraining_error}
\mathcal E_{\mathrm{pre}}
:=
\sup_{\mathcal P,\boldsymbol y}
\left|\widehat p(\boldsymbol y\mid\mathcal P)
-p^*(\boldsymbol y\mid\mathcal P)\right|,
\end{equation}
where the supremum ranges over the prompt--response pairs just described, including the prompting schemes below. Taking the supremum in \eqref{eq:llm_gen_mem} bounds this quantity by the corresponding uniform approximation and estimation discrepancies, together with the prediction-level training residual $\epsilon_{\mathrm{opt}}^{\mathrm{pred}}$. We retain this quantity as the prediction-error input to the prompting analysis.

The latent-task model describes prompt comprehension through posterior inference over tasks. A query may be compatible with several different intents, and the appropriate response depends on which latent task generated or is intended by the prompt. The model introduces a latent task $\theta_{\mathrm{task}}\in\Theta_{\mathrm{task}}$ in a discrete task space, with prior $p^*(\theta_{\mathrm{task}})$, and obtains the marginal document distribution by integrating out this unobserved variable. The task posteriors are quantities of this latent-task model used to analyze how the prompt identifies a task. For a target query $x$, the dominant inferred task is
\[
\theta_{\mathrm{task},x} \in \arg\max_{\theta_{\mathrm{task}}\in\Theta_{\mathrm{task}}}p^*(\theta_{\mathrm{task}}\mid x).
\]
In zero-shot prompting, the prompt consists of the query $x$ without demonstrations. For a response sequence $\boldsymbol y$, we measure the difference between the trained model's response probability and the true response probability conditional on the dominant task by
\[
\mathcal E_{\text{zero-shot}}(x,\boldsymbol y)
\coloneqq
\left|
\widehat p(\boldsymbol y\mid x)
-p^*(\boldsymbol y\mid x,\theta_{\mathrm{task},x})
\right|.
\]
The task ambiguity of the query is defined as
\begin{equation}
\label{eq:llm_A_theta}
\Delta_{\Theta_{\mathrm{task}}}(x) := 1-p^*(\theta_{\mathrm{task},x}\mid x).
\end{equation}
This is a posterior quantity: it is small when the prompt concentrates most of its posterior mass on a single task and large when several latent tasks remain plausible. The marginal response distribution satisfies
\begin{equation}
\label{eq:llm_latent_mixture}
p^*(\boldsymbol y\mid x)
=
\sum_{\theta_{\mathrm{task}}\in\Theta_{\mathrm{task}}}
p^*(\boldsymbol y\mid x,\theta_{\mathrm{task}})
p^*(\theta_{\mathrm{task}}\mid x).
\end{equation}
Subtracting the response distribution associated with the dominant task and using that conditional probabilities lie in $[0,1]$ gives
\begin{align}
&\left| p^*(\boldsymbol y\mid x) - p^*(\boldsymbol y\mid x,\theta_{\mathrm{task},x}) \right| \nonumber\\
&\quad=
\left|
\sum_{\theta_{\mathrm{task}}\neq\theta_{\mathrm{task},x}}
p^*(\theta_{\mathrm{task}}\mid x)
\left[
p^*(\boldsymbol y\mid x,\theta_{\mathrm{task}})
-p^*(\boldsymbol y\mid x,\theta_{\mathrm{task},x})
\right]
\right| \nonumber\\
&\quad\leq
\sum_{\theta_{\mathrm{task}}\neq\theta_{\mathrm{task},x}}
p^*(\theta_{\mathrm{task}}\mid x)
=\Delta_{\Theta_{\mathrm{task}}}(x).
\label{eq:llm_q_diff}
\end{align}
Combining \eqref{eq:llm_q_diff} with \eqref{eq:llm_pretraining_error} gives the zero-shot decomposition used in the proof of Theorem~12 of \cite{jiao2026prompt}:
\begin{equation}
\label{eq:llm_zero_shot_bound}
\mathcal E_{\text{zero-shot}}(x,\boldsymbol y)
\leq
\mathcal E_{\mathrm{pre}}
+\Delta_{\Theta_{\mathrm{task}}}(x).
\end{equation}
This is the total zero-shot error bound: the first term measures prediction of the marginal response distribution, and the second measures ambiguity about the task. Together with \eqref{eq:llm_gen_mem}, it separates the approximation, finite-corpus estimation, and prediction-level training contributions from the additional ambiguity introduced by the query.

As an illustrative refinement, impose the additional interpretation that the same dominant task is represented by a pair $(\theta_x^\star,\theta_y^\star)$: a carrier state described by the query and the state produced by a transformation. Define
\begin{equation}
\label{eq:llm_refined_ambiguities}
\Delta_{\mathrm C}(x)
:=1-p^*(\theta_x^\star\mid x),
\qquad
\Delta_{\mathrm{trans}}(x)
:=1-p^*(\theta_y^\star\mid\theta_x^\star,x).
\end{equation}
These quantities measure uncertainty about the carrier state and the resulting state conditional on it, respectively. The chain rule for conditional probabilities gives
\begin{align}
1-p^*(\theta_x^\star,\theta_y^\star\mid x)
&=1-p^*(\theta_x^\star\mid x)\,
p^*(\theta_y^\star\mid\theta_x^\star,x)\nonumber\\
&=\Delta_{\mathrm C}(x)
+\bigl[1-\Delta_{\mathrm C}(x)\bigr]
\Delta_{\mathrm{trans}}(x)\leq \Delta_{\mathrm C}(x)+\Delta_{\mathrm{trans}}(x).
\label{eq:llm_refined_ambiguity_bound}
\end{align}
Under this paired-task interpretation, the zero-shot error therefore satisfies
\[
\mathcal E_{\text{zero-shot}}(x,\boldsymbol y)
\leq
\underbrace{\mathcal E_{\mathrm{pre}}}_{\text{pretraining error}}
+\underbrace{\Delta_{\mathrm C}(x)}_{\text{carrier-state ambiguity}}
+\underbrace{\Delta_{\mathrm{trans}}(x)}_{\text{transformation-outcome ambiguity}}.
\]
This separates uncertainty about the state described by the query from uncertainty about the transformation outcome conditional on that state. The prompting results below use the overall task ambiguity $\Delta_{\Theta_{\mathrm{task}}}(x)$.

\subsubsection{In-context learning and chain-of-thought}\label{app:llm_understand}
We retain the model and notation of Section~\ref{app:llm_pretraining} to analyze how demonstrations and intermediate reasoning steps affect the error bounds. To reduce the task ambiguity in the zero-shot bound, in-context learning keeps the trained model fixed and adds demonstrations to the prompt. We evaluate its response probabilities against those for the query conditioned on its dominant latent task.

Let
\[
\mathcal P_{\mathrm{ICL}} = (x^{(1)},y^{(1)},\ldots,x^{(m)},y^{(m)},x)
\]
denote a prompt containing $m\geq1$ demonstrations followed by the target query, with delimiter tokens between examples suppressed in this notation. Here $y^{(i)}$ denotes a demonstration response, which may contain multiple tokens. The demonstrations and query share a latent task, and each demonstration has the same dominant task $\theta_{\mathrm{task},x}$.

For a response sequence $\boldsymbol y$, define the ICL response-probability error by
\[
\mathcal E_{\mathrm{ICL}}(\mathcal P_{\mathrm{ICL}},\boldsymbol y)
\coloneqq
\left|\widehat p(\boldsymbol y\mid\mathcal P_{\mathrm{ICL}})
-p^*(\boldsymbol y\mid x,\theta_{\mathrm{task},x})\right|.
\]

The delimiter and nearly-Markov assumptions control how much earlier examples affect token probabilities conditional on this task: removing earlier demonstrations changes each conditional log token probability by at most $\phi_{\mathrm M}$. Under the task-consistency, delimiter, nearly-Markov, and prior-balance assumptions (Assumptions~13--16 of \cite{jiao2026prompt}), the posterior calculation in the proof of Theorem~17 gives
\begin{equation}
\label{eq:llm_icl_bound}
\begin{aligned}
&\mathcal E_{\mathrm{ICL}}(\mathcal P_{\mathrm{ICL}},\boldsymbol y) \\
&\quad  \leq \underbrace{\mathcal E_{\mathrm{pre}}}_{\text{pretraining error}}
+\underbrace{L_{\mathrm{res}}\phi_{\mathrm M}}_{\text{residual dependence}} +\underbrace{\left(e^{2L_{\mathrm{seq}}\phi_{\mathrm M}}
c\epsilon_{\mathrm{ICL}}\right)^m
\Delta_{\Theta_{\mathrm{task}}}(x)}_{\text{task ambiguity}}.
\end{aligned}
\end{equation}
Under the nearly-Markov condition, the residual effect of earlier demonstrations is controlled token by token; accumulating these bounds across the response gives the factor $L_{\mathrm{res}}$. Here $c$ controls imbalance in the task prior, and the ambiguity coefficient depends on both the demonstrations and the query:
\[
\epsilon_{\mathrm{ICL}}
=\frac{1}{1-\Delta_{\Theta_{\mathrm{task}}}(x)}
\max_{1\leq i\leq m}
\frac{\Delta_{\Theta_{\mathrm{task}}}(x^{(i)}\circ y^{(i)})}
{1-\Delta_{\Theta_{\mathrm{task}}}(x^{(i)}\circ y^{(i)})}.
\]
Under these assumptions, the posterior calculation bounds task-likelihood ratios across the demonstrations by a product of factors controlled by their ambiguity ratios. A common upper bound, incorporating the query factor in $\epsilon_{\mathrm{ICL}}$, the prior-imbalance correction $c$, and the residual-dependence correction $e^{2L_{\mathrm{seq}}\phi_{\mathrm M}}$, yields the $m$th power in \eqref{eq:llm_icl_bound}. The task-ambiguity term admits a geometrically decaying upper bound in $m$ if the base
\[
e^{2L_{\mathrm{seq}}\phi_{\mathrm M}}c\epsilon_{\mathrm{ICL}}
\]
stays uniformly below one as demonstrations are added within the context-length limit. The pretraining contribution and nearly-Markov remainder remain additive.

Chain-of-thought prompts include observable intermediate reasoning steps. In the model, these texts provide evidence about a compositional task represented by a sequence of simpler latent tasks:
\[
t_{\mathrm{task}}^\star = \theta_{\mathrm{task},1}^\star \circ\cdots\circ \theta_{\mathrm{task},L_{\mathrm{CoT}}}^\star.
\]
Here $L_{\mathrm{CoT}}$ is the chain length. Let $\mathcal P_{\mathrm{CoT}}$ denote a prompt with demonstrations and intermediate reasoning steps, and let $\widetilde p^*$ be the compositional-task law in \cite{jiao2026prompt}. We evaluate the model's response probabilities against those under $\widetilde p^*$ conditioned on the query and the specified target sequence $t_{\mathrm{task}}^\star$.

Define the CoT response-probability error for this target sequence by
\[
\begin{aligned}
\mathcal E_{\mathrm{CoT}}(\mathcal P_{\mathrm{CoT}},\boldsymbol y;t_{\mathrm{task}}^\star)\coloneqq
\left|\widehat p(\boldsymbol y\mid\mathcal P_{\mathrm{CoT}})
-\widetilde p^*(\boldsymbol y\mid x,t_{\mathrm{task}}^\star)\right|.
\end{aligned}
\]

In this model, the pretraining and compositional laws share the conditional generation law for prompts and responses given a task sequence, but may have different task priors. The pretraining tasks are identified with constant sequences in $\Theta_{\mathrm{task}}^{L_{\mathrm{CoT}}}$, each repeating a single task. This places both task priors on the same space for the distribution-shift comparison. Assume that the prompt has positive probability under both laws and that the target sequence $t_{\mathrm{task}}^\star$ has positive compositional posterior probability.

Following the distribution-shift comparison in Appendix~E of \cite{jiao2026prompt}, define
\[
\begin{aligned}
M_{\mathrm{CoT}}:=\max\left\{
p^*(\mathcal P_{\mathrm{CoT}})^{-1},
\widetilde p^*(\mathcal P_{\mathrm{CoT}})^{-1}
\right\}, \quad  \Delta_{\mathrm{PCoT}}:= \sum_{\substack{\vartheta\in\Theta_{\mathrm{task}}^{L_{\mathrm{CoT}}}\\
p^*(\mathcal P_{\mathrm{CoT}}\mid\vartheta)>0}}
\left|\widetilde p^*(\vartheta)-p^*(\vartheta)\right|.
\end{aligned}
\]
Here $\Delta_{\mathrm{PCoT}}$ measures prior differences among task sequences compatible with the prompt. The factor $M_{\mathrm{CoT}}$ arises from normalizing task posteriors by the prompt probabilities and can be large when the prompt is rare under either law. Comparing the normalized task posteriors bounds the difference between the two marginal response probabilities by $M_{\mathrm{CoT}}\Delta_{\mathrm{PCoT}}$. Expanding the compositional response law over task sequences, as in \eqref{eq:llm_latent_mixture}, bounds its difference from the target-task response by the posterior mass outside $t_{\mathrm{task}}^\star$. If the nearly-Markov log-probability condition also holds for the compositional law at this target task, the triangle inequality gives
\begin{equation}
\label{eq:llm_cot_bound}
\begin{aligned}
& \mathcal E_{\mathrm{CoT}}(\mathcal P_{\mathrm{CoT}},\boldsymbol y;t_{\mathrm{task}}^\star) \\
& \quad \leq
\underbrace{\mathcal E_{\mathrm{pre}}}_{\text{pretraining error}}
+\underbrace{L_{\mathrm{res}}\phi_{\mathrm M}}_{\text{residual dependence}} +\underbrace{M_{\mathrm{CoT}}\Delta_{\mathrm{PCoT}}}_{\text{task-prior shift}}
+\underbrace{\left[1-\widetilde p^*(t_{\mathrm{task}}^\star
\mid\mathcal P_{\mathrm{CoT}})\right]}_{\text{target-sequence ambiguity}}.
\end{aligned}
\end{equation}
The bound separates learning how to predict responses from inferring which task composition the prompt specifies. Pretraining accuracy controls the former, while reasoning steps that concentrate the posterior on the target sequence reduce the ambiguity in the latter. Connecting these two aspects also requires accounting for how earlier context affects the response even after conditioning on the task, and for the change in task prior between pretraining and compositional prompting.

The regression example shows how Transformer approximation and complexity estimates combine with data-dependence control to bound prediction error. In the language-model setting, the pretraining decomposition retains approximation, finite-corpus estimation, and the prediction-level training residual. The same latent-task model then connects this prediction error to zero-shot, ICL, and CoT performance through task ambiguity, residual dependence on earlier demonstrations, and, for CoT, a change in the task distribution.

%% file: sections/horizon.tex
\section{Open Horizons}\label{sec:open_horizons}

This review has examined approaches to error analysis in deep learning by distinguishing approximation, statistical, and optimization errors. The size-controlled and norm-controlled analyses in Chapters~\ref{sec:classical_setting} and~\ref{sec:overparam_setting} balance approximation and statistical errors under exact empirical risk minimization. The former control complexity through network width, depth, and parameter count, while the latter use weight-norm constraints to obtain guarantees for networks with large numbers of parameters. Complete error analysis additionally incorporates a specified training algorithm. Using nonparametric regression and the Deep Ritz Method as running examples, we have explained how approximation, statistical complexity, and optimization estimates combine to bound the error of the algorithm's output.

These analyses show that a complete guarantee requires the component estimates to apply to the trained network under compatible conditions. The networks providing the approximation guarantee must satisfy the constraints underlying the statistical estimates, and the training output must remain within the scope of those estimates. Network architecture, parameter constraints, and training procedures are therefore interdependent. Increasing expressive power may reduce approximation error while enlarging statistical complexity; changing the training constraints may facilitate optimization but invalidate existing statistical bounds. The complete analyses reviewed here illustrate how these requirements can be reconciled, providing a basis for extending their scope, sharpening error bounds, and studying training dynamics in applications.

One important direction is to extend complete error analysis to a broader range of training algorithms and network models. For training algorithms, this entails reducing reliance on projection and accommodating stochastic gradients and adaptive updates. Complete error guarantees for unprojected gradient descent and for stochastic gradient descent have already been obtained in specific neural network regression models~\cite{kohler2025rate,kohler2026sgd}. Without projection, a key question is how to control parameter norms and trajectories through the training dynamics themselves, so that approximation and statistical estimates can still be combined for the trained output. Stochastic gradient noise affects both optimization accuracy and parameter movement, requiring step sizes, batch sizes, and stopping times to account for both effects. For algorithms such as Adam, the analysis must also account for momentum and adaptive rescaling based on past gradients, and determine how they affect optimization error and the complexity of the trained network.

This direction also includes extending the analysis to other architectures and activation functions. The complete analyses in Chapter~\ref{sec:complete_setting} rely substantially on a particular parallel architecture, a linear outer layer, and prescribed initialization constructions. Extending their guarantees to convolutional networks, residual networks, and attention-based architectures calls for arguments that retain the features of these structures. In particular, parameter sharing, residual connections, and attention mechanisms may affect approximation efficiency, statistical complexity, and training behavior. Representing these architectures simply as larger generic networks may conceal their structural advantages in the resulting bounds. Activation functions must also match the regularity required by the task: regression primarily concerns function values, DRM involves first-order spatial derivatives, and strong-form PINN losses may involve higher-order derivatives. Changing the activation therefore requires examining approximation of both functions and their derivatives, the statistical complexity of the corresponding classes, and the smoothness of the training objective.

A second direction is to sharpen approximation and statistical estimates to obtain tighter complete error bounds, with the aim of attaining minimax-optimal convergence rates. Some of the complete guarantees reviewed here remain slower than the minimax benchmark for the corresponding model. On the approximation side, more effective use of target smoothness and structure may yield sharper bounds under prescribed network-size or norm constraints, while reducing losses in parameter and derivative estimates. On the statistical side, the structure of the loss and estimation problem can be used to refine global deviation bounds over the entire candidate class. Localized empirical-process arguments and fast-rate analyses that exploit loss curvature or suitable noise conditions offer routes to more precise statistical estimates. Rebalancing approximation and statistical errors using these estimates, while keeping the optimization contribution below the leading terms, may narrow the gap between existing complete guarantees and the optimal statistical benchmark under the same target class, observation model, and error criterion; see Section~\ref{sec:minimax_lb}.

A third direction is to develop end-to-end error analyses that include the training process in specific applications. Many results in Chapter~\ref{sec:applications} start from empirical risk minimization or retain an optimization residual without estimating its decay. They establish what can be guaranteed once the empirical problem is solved to a given accuracy, leaving the question of how a training algorithm attains that accuracy. This step often depends on the application. For example, the GAN analysis in Section~\ref{section:app:generative:gan} requires control of the empirical minimax objective, whereas a stationarity guarantee does not automatically control its optimality gap. Errors from solving the generator and discriminator problems also need to be accounted for separately. Incorporating optimization analysis into an existing decomposition therefore requires checking that the quantity controlled by the algorithm is the one needed by the error bound.

An application-level analysis must also translate training accuracy into the error of interest. In inverse problems, this translation depends jointly on observation noise, regularization, and inverse stability. In diffusion and flow matching, field-estimation errors propagate through the sampling dynamics and combine with time truncation and numerical discretization errors to determine distributional accuracy. The required optimization accuracy should therefore be related to the target recovery or generation error, accounting for the other error terms and the relevant stability constants. Reinforcement learning similarly requires control of fitting errors accumulated over multiple iterations, while representation learning requires conditions linking the pretraining objective to downstream risk. Connecting these relations to the training process would explain how improvements in the training objective translate into gains on the final task, and distinguish errors that decrease with training from those limited by observation noise, distributional coverage, or transfer conditions.

Error analysis connects network expressivity, finite data, and training dynamics to final accuracy, allowing the performance of deep learning methods to be understood in terms of explicit assumptions and error sources. Identifying the dominant contributions in different settings helps assess the gains that may come from additional data, changes in architecture, or improved training procedures. The methods and open questions discussed in this review provide a basis for developing these connections further, with the continuing goal of making theoretical guarantees more informative about the performance of different methods and more useful in guiding model and training choices.

%% file: sections/appendix.tex
\section{Notation tables}\label{app:notation}

\begingroup\renewcommand{\arraystretch}{1.15}
\begin{longtable}{@{}p{0.18\linewidth} >{\raggedright\arraybackslash}p{0.75\linewidth}@{}}
\caption{List of abbreviations}\\
\toprule
\textbf{Abbreviation} & \textbf{Meaning} \\
\midrule
\endhead
DRM  & Deep Ritz Method \\
ERM  & Empirical Risk Minimizer (or Empirical Risk Minimization) \\
PNN  & Parallel Neural Network \\
PGD  & Projected Gradient Descent \\
PSGD & Projected Stochastic Gradient Descent \\
PDE  & Partial Differential Equation \\
PINN & Physics-Informed Neural Network \\
NTK  & Neural Tangent Kernel \\
GAN  & Generative Adversarial Network \\
IPM  & Integral Probability Metric \\
DSM  & Denoising Score Matching \\
FM   & Flow Matching \\
OU   & Ornstein--Uhlenbeck \\
SDE  & Stochastic Differential Equation \\
ODE  & Ordinary Differential Equation \\
RL   & Reinforcement Learning \\
FQI  & Fitted $Q$-Iteration \\
API  & Approximate Policy Iteration \\
CNN  & Convolutional Neural Network \\
SSL  & Self-Supervised Learning \\
DDR  & Deep Dimension Reduction \\
Adv-SSL & Adversarial Self-Supervised Learning \\
LLM  & Large Language Model \\
MC   & Monte Carlo \\
KL   & Kullback--Leibler \\
TV   & Total Variation \\
ICL  & In-Context Learning \\
CoT  & Chain-of-Thought \\
\bottomrule
\end{longtable}
\endgroup
\medskip
\begingroup\renewcommand{\arraystretch}{1.25}
\begin{longtable}{@{}p{0.27\linewidth} >{\raggedright\arraybackslash}p{0.46\linewidth} >{\raggedright\arraybackslash}p{0.20\linewidth}@{}}
\caption{Neural network architecture}\\
\toprule
\textbf{Symbol} & \textbf{Meaning} & \textbf{Definition / location} \\
\midrule
\endhead
$d$ & Input dimension & \S\ref{nnstructure} \\
$L$ & Network depth (number of affine layers, including the output layer) & \S\ref{nnstructure}, \eqref{eqnn1} \\
$W$ & Maximum hidden-layer width ($L\ge2$) & \S\ref{nnstructure} \\
$N_\ell$ & Number of neurons at layer $\ell$ & \S\ref{nnstructure} \\
$\rho$ & Activation function & \S\ref{nnstructure} \\
$\boldsymbol{A}_\ell, \boldsymbol{b}_\ell$ & Weight matrix and bias at layer $\ell$ & \S\ref{nnstructure} \\
$\boldsymbol{\theta}$ & Parameter vector (all weights and biases) & \S\ref{nnstructure} \\
$\phi_{\boldsymbol{\theta}}$ & Generic neural network function & \S\ref{nnstructure}, \eqref{eqnn1} \\
$\|\cdot\|_{\max}$ & Max absolute entry & \S\ref{nnstructure} \\
$B_{\boldsymbol{\theta}}$ & Bound on $\|\boldsymbol{\theta}\|_\infty$ & \S\ref{nnstructure} \\
$\mathfrak{D}(W, L, d), \mathfrak{D}$ & Padded width-$W$ parameter count ($L\ge2$); parameter count or nonzero-parameter budget, as specified for the network class. In Transformer regression, the total count includes the linear readout & \S\ref{nnstructure}, \eqref{eq:align_dim}, \S\ref{app:llm_reg} \\
$\Theta$ & Parameter space & \S\ref{nnstructure} \\
$\mathcal{NN}(W, L)$ & Generic NN function class & \S\ref{nnstructure} \\
$\mathcal{NN}_\rho(W, L, B_{\boldsymbol{\theta}})$ & Bounded-parameter NN class (size and magnitude) & \S\ref{nnstructure} \\
$\mathcal{NN}_\rho(L, \mathfrak{D}, B_{\boldsymbol{\theta}})$ & NN class by depth, \#~nonzero parameters, magnitude (DRM/tanh) & \S\ref{nnstructure}, \S\ref{sec:A2} \\
$\mathcal{F}_{NN}$ & Hypothesis class (problem-dependent) & \S\ref{nnstructure} \\
$\tilde{\boldsymbol{A}}_\ell$ & Augmented weight matrix (with bias) & Definition~\ref{sec:overparam_setup} \\
$\kappa(\boldsymbol{\theta})$ & Weight-norm product $\prod_\ell \|\tilde{\boldsymbol{A}}_\ell\|$ & Definition~\ref{sec:overparam_setup} \\
$\mathcal{K}$ & Bound on $\kappa(\boldsymbol{\theta})$ & Definition~\ref{sec:overparam_setup}, \S\ref{section:app:generative:gan} \\
$\mathcal{NN}_\rho(W, L, \mathcal{K})$ & Norm-constrained NN class & Definition~\ref{sec:overparam_setup} \\
$\mathcal{PNN}$ & Parallel NN class & \S\ref{sec:new_decomp} \\
$\Phi$ & PNN function & \S\ref{sec:new_decomp} \\
$\mathfrak{m}$ & Number of sub-networks in PNN & \S\ref{sec:new_decomp} \\
$M$ & $\ell^1$-bound on outer coefficients & \S\ref{sec:new_decomp} \\
$c_k$ & Outer coefficients of PNN & \S\ref{sec:new_decomp} \\
$\phi^k$ & $k$-th sub-network of PNN & \S\ref{sec:new_decomp} \\
$\boldsymbol{\theta}_{\mathrm{in}}, \boldsymbol{\theta}_{\mathrm{out}}$ & Inner and outer PNN parameter blocks & \S\ref{sec:new_decomp} \\
$\boldsymbol{\theta}_{\mathrm{total}}$ & Concatenated PNN parameter vector & \S\ref{sec:new_decomp} \\
\bottomrule
\end{longtable}
\endgroup

\begingroup\renewcommand{\arraystretch}{1.25}
\begin{longtable}{@{}p{0.27\linewidth} >{\raggedright\arraybackslash}p{0.46\linewidth} >{\raggedright\arraybackslash}p{0.20\linewidth}@{}}
\caption{Risk and error decomposition}\\
\toprule
\textbf{Symbol} & \textbf{Meaning} & \textbf{Definition / location} \\
\midrule
\endhead
$\ell(y', y)$ & Point-wise loss function & \S\ref{sec:decomp} \\
$\mathcal{L}(\phi)$ & Population (true) risk & \S\ref{sec:decomp} \\
$\widehat{\mathcal{L}}_n(\phi)$ & Empirical risk on $n$ samples & \S\ref{sec:decomp} \\
$\mathcal{L}^*$ & Bayes risk $\inf_\phi \mathcal{L}(\phi)$ & \S\ref{sec:decomp} \\
$\phi^*$ & Bayes-optimal function & \S\ref{sec:decomp} \\
$\phi_{\bar{\boldsymbol{\theta}}}$ & Best in-class approximator & \S\ref{sec:decomp} \\
$\phi_{\hat{\boldsymbol{\theta}}_n}, \widehat{\phi}$ & Empirical risk minimizer (ERM) & \S\ref{sec:decomp} \\
$\phi_{\boldsymbol{\theta}_{\mathcal{A}}}$ & Algorithm output & \S\ref{sec:decomp}, \S\ref{sec:complete_opt} \\
$\mathcal{A}$ & Algorithm (gradient method, PGD, etc.) & \S\ref{sec:decomp} \\
$\mathcal{H}$ & Target-quantity Hilbert space ($L^2(P_X)$ or $H^1(\Omega)$) & \S\ref{sec:running_examples} \\
$\mathcal{E}_{app}$ & Best-in-class approximation error & \S\ref{sec:decomp} \\
$\mathcal{E}_{sta}$ & Sample-dependent uniform population--empirical deviation & \S\ref{sec:decomp} \\
$\mathcal{E}_{sta,\,slow}$ & Twice the expected uniform population--empirical deviation & \S\ref{sec:decomp}, \eqref{eq:sym_decomp} \\
$\mathcal{E}_{sta,\,fast}$ & Expected offset supremum for the fast-rate bound & \S\ref{sec:decomp}, \eqref{eq:asym_decomp} \\
$\mathcal{E}^{\mathrm{REG}}_{sta,\,fast}$ & Regression expected offset supremum & \S\ref{sec:running_examples}, \eqref{eq:reg_offset_def} \\
$\mathcal{E}^{\mathrm{DRM}}_{sta,\,fast}$ & DRM expected offset supremum & \S\ref{sec:running_examples}, \eqref{eq:drm_offset_def} \\
$\mathcal{E}_{opt}$ & ERM-anchored optimization error & \S\ref{sec:decomp} \\
$\mathcal{E}_{opt,-}$ & Approximation-based optimization error & \S\ref{sec:new_decomp}, \eqref{eq:D1} \\
$\mathcal E^{\mathrm c}_{app}$ & Risk gap of the constructed approximating network & \eqref{eq:D1} \\
$\mathcal E^{\mathrm c}_{sta}$ & Twice the uniform population--empirical deviation (sample-dependent) & \eqref{eq:D1} \\
$\mathcal E^{\mathrm{DRM,c}}_{app},\mathcal E^{\mathrm{DRM,c}}_{sta}$ & DRM instances of the constructed-network approximation and statistical terms & Lemma~\ref{lem:D1}, \S\ref{sec:drm_complete} \\
$\bar{\phi}, \bar{f}, \bar{u}$ & Constructed approximating network; regression and DRM instances & \S\ref{sec:new_decomp}, \eqref{eq:bar_phi} \\
\bottomrule
\end{longtable}
\endgroup

\begingroup\renewcommand{\arraystretch}{1.25}
\begin{longtable}{@{}p{0.27\linewidth} >{\raggedright\arraybackslash}p{0.46\linewidth} >{\raggedright\arraybackslash}p{0.20\linewidth}@{}}
\caption{Function spaces, smoothness, and support geometry}\\
\toprule
\textbf{Symbol} & \textbf{Meaning} & \textbf{Definition / location} \\
\midrule
\endhead
$\mathcal{H}^\beta(\Omega, B_{\mathrm{REG}})$ & Hölder ball of smoothness $\beta$, radius $B_{\mathrm{REG}}$ & \S\ref{sec:notation_summary}, \eqref{eq0} \\
$W^{s, p}(\Omega)$ & Sobolev space, smoothness $s$, integrability $p$ & \S\ref{sec:notation_summary} \\
$H^s(\Omega) = W^{s,2}(\Omega)$ & Sobolev space (Hilbert case) & \S\ref{sec:notation_summary} \\
$\beta$ & Hölder smoothness parameter & \S\ref{sec:reg_setup} \\
$s$ & Sobolev smoothness order & \S\ref{sec:drm_setup} \\
$B_{\mathrm{REG}}$ & Hölder ball radius (Regression smoothness bound) & \S\ref{sec:reg_setup} \\
$r$ & H\"older remainder ($r=\beta-s\in(0,1]$, $s=\lceil\beta\rceil-1$) & \S\ref{sec:notation_summary} \\
$\alpha$ & Elliptic regularity exponent ($\partial\Omega \in C^{2+\alpha}$, $0<\alpha<1$) & \S\ref{sec:drm_setup} \\
$B_{\mathrm{DRM}}$ & PDE data sup-norm $\max\{\|h\|_{L^\infty}, \|g\|_{L^\infty}, \|\omega\|_{L^\infty}\}$ (DRM) & \S\ref{sec:drm_setup} \\
$L^2(P_X)$ & $L^2$ space with respect to the design measure $P_X$ & \S\ref{sec:reg_setup} \\
$\mathcal X,d^*,d^\sharp$ & Input support, its upper Minkowski dimension, and a relaxed covering exponent $d^\sharp>d^*$ & \S\ref{sec:data_geometric_statistical_dimension} \\
$\mathcal M,d_{\mathcal M},\rho_{\mathrm{tube}}$ & Support manifold, its dimension, and neighborhood radius; $\mathcal M_{\rho_{\mathrm{tube}}}$ denotes the neighborhood within $[0,1]^d$ & \S\ref{sec:data_geometric_statistical_dimension} \\
$\delta,d_\delta,d_0$ (geometry) & Projection distortion and the construction-specific projection dimensions for approximate-manifold and low-Minkowski support & \S\ref{sec:data_geometric_statistical_dimension} \\
\bottomrule
\end{longtable}
\endgroup

\begingroup\renewcommand{\arraystretch}{1.25}
\begin{longtable}{@{}p{0.27\linewidth} >{\raggedright\arraybackslash}p{0.46\linewidth} >{\raggedright\arraybackslash}p{0.20\linewidth}@{}}
\caption{Running examples}\\
\toprule
\textbf{Symbol} & \textbf{Meaning} & \textbf{Definition / location} \\
\midrule
\endhead
\multicolumn{3}{l}{Regression} \\[1pt]
$f_0$ & True regression function (conditional mean) & \S\ref{sec:reg_setup}, \eqref{eq:reg_model}, \eqref{eq:llm_reg_model} \\
$\boldsymbol{X}, Y$ & Predictor and response & \S\ref{sec:reg_setup} \\
$\xi$ & Noise term & \S\ref{sec:reg_setup} \\
$\sigma$ & Noise bound, sub-Gaussian scale, or Gaussian standard deviation & \S\ref{sec:reg_slow} \\
$P_X$ & Marginal distribution of $\boldsymbol{X}$ & \S\ref{sec:reg_setup} \\
$S_n = \{(\boldsymbol{X}_i, Y_i)\}_{i=1}^n$ & i.i.d.\ training sample & \S\ref{sec:reg_setup} \\
$n$ & Sample size & \S\ref{sec:reg_setup} \\
$\mathcal{L}_{\mathrm{REG}}, \widehat{\mathcal{L}}_{\mathrm{REG}}$ & Population and empirical regression risk & \S\ref{sec:reg_setup} \\
$\widehat{f}_{\boldsymbol{\theta}}$ & Regression ERM estimator & \S\ref{sec:reg_setup}, \eqref{eq:reg_erm} \\
$[\phi]_{\mathcal B}$ & Pointwise clipping to $[-\mathcal B,\mathcal B]$ & \S\ref{sec:reg_slow} \\
\midrule
\multicolumn{3}{l}{DRM (Deep Ritz Method)} \\[1pt]
$\Omega \subset [0,1]^d$ & Bounded domain & \S\ref{sec:drm_setup} \\
$\partial \Omega$ & Domain boundary & \S\ref{sec:drm_setup} \\
$\boldsymbol{n}$ & Outward normal vector & \S\ref{sec:drm_setup}, \eqref{eq:drm_pde} \\
$h, g, \omega$ & PDE source, boundary, coercivity data & \S\ref{sec:drm_setup} \\
$c_0$ & Lower bound on $\omega$ & \S\ref{sec:drm_setup} \\
$u_0$ & Weak PDE solution & \S\ref{sec:drm_setup} \\
$\mathcal{L}_{\mathrm{DRM}}, \widehat{\mathcal{L}}_{\mathrm{DRM}}$ & Energy functional and its MC discretization & \S\ref{sec:drm_setup}, \eqref{eq:drm_energy}, \eqref{eq:drm_emp} \\
$N_{\mathrm{in}}, N_b$ & Interior and boundary MC sample sizes & \S\ref{sec:drm_setup} \\
$S$ & Combined MC sample set $\{\boldsymbol{X}_p\} \cup \{\boldsymbol{Y}_p\}$ & \S\ref{sec:drm_setup} \\
$N$ & MC sample size on each discretization under the equal-size assumption $N_{\mathrm{in}} = N_b = N$ & Theorem~\ref{thm:drm_slow} \\
$T$ (operator) & Trace operator $H^1 \to L^2(\partial\Omega)$ & \S\ref{sec:drm_setup} \\
$\widehat{u}_{\boldsymbol{\theta}}$ & DRM ERM estimator & \S\ref{sec:drm_setup}, \eqref{eq:drm_erm} \\
$c_0 \wedge 1$, $B_{\mathrm{DRM}} \vee 1$ & Energy--norm equivalence constants & \S\ref{sec:drm_setup}, Prop~\ref{prop:drm_eq} \\
\bottomrule
\end{longtable}
\endgroup

\begingroup\renewcommand{\arraystretch}{1.25}
\begin{longtable}{@{}p{0.27\linewidth} >{\raggedright\arraybackslash}p{0.46\linewidth} >{\raggedright\arraybackslash}p{0.20\linewidth}@{}}
\caption{Approximation theorems}\\
\toprule
\textbf{Symbol} & \textbf{Meaning} & \textbf{Definition / location} \\
\midrule
\endhead
$J, I$ & Network parameters ($J$ controls width, $I$ controls depth) & \S\ref{sec:classical_approx}, Thm~\ref{thm:A1} \\
$G$ & Partition-grid resolution & \S\ref{sec:A1} \\
$\mu$ & Arbitrarily small parameter in the tanh approximation theorems & \S\ref{sec:A2}, Thm~\ref{thm:A2} \\
$\vartheta$ & Weight-magnitude exponent of the tanh approximation theorems & Theorem~\ref{thm:A2} \\
$a$ & Bump sharpness parameter ($a = G^\mu$) & \S\ref{sec:A2} proof \\
$\delta$ (tanh) & Tanh Taylor scaling parameter & \S\ref{sec:classical_approx} \\
$\mathcal{NN}_\rho^{\mathfrak{m}}(W, L, \mathcal{K})$ & Linear combinations $\sum_j c_j\phi^j$ with aggregate budget $\sum_j|c_j|\,\kappa(\boldsymbol\theta^j)\leq\mathcal K$ & \S\ref{sec:overparam_approx}, Thm~\ref{thm:A4} \\
$S_{\mathcal K},R_{\mathcal K}$ & Outer-coefficient $\ell^1$ budget and per-sub-network layerwise matrix-norm budget for DRM & \eqref{eq:drm_op_candidate}, Thm~\ref{thm:drm_overparam} \\
$q_d,a_{\mathrm{DRM}}$ & Local DRM abbreviations: $q_d=\lceil\log_2(d+1)\rceil$, $a_{\mathrm{DRM}}=(1-\mu)/(C_1d q_d)$ & Thm~\ref{thm:drm_overparam} and proof \\
\bottomrule
\end{longtable}
\endgroup

\begingroup\renewcommand{\arraystretch}{1.25}
\begin{longtable}{@{}p{0.27\linewidth} >{\raggedright\arraybackslash}p{0.46\linewidth} >{\raggedright\arraybackslash}p{0.20\linewidth}@{}}
\caption{Statistical bounds}\\
\toprule
\textbf{Symbol} & \textbf{Meaning} & \textbf{Definition / location} \\
\midrule
\endhead
$\sigma_i$ & Rademacher signs ($\pm 1$ with prob.\ $1/2$) & \S\ref{sec:classical_stat} \\
$\widehat{\mathfrak{R}}_n(\mathcal{F})$ & Empirical Rademacher complexity & \S\ref{sec:classical_stat}, \eqref{eq:rad_def_classical} \\
$\tilde{\mathfrak{R}}_n(\mathcal{F})$ & Expected symmetric Rademacher complexity (with $|\cdot|$) & \S\ref{sec:classical_stat}, \eqref{eq:rad_def_sym} \\
$\mathfrak{R}_n(\mathcal{F})$ & Expected Rademacher complexity & \S\ref{sec:classical_stat} \\
$\mathcal{C}(\epsilon, \mathcal{F}, \|\cdot\|)$ & $\epsilon$-covering number & \S\ref{sec:classical_stat}, Prop~\ref{prop:dudley} \\
$\mathrm{Pdim}(\mathcal{F})$ & Pseudo-dimension & \S\ref{sec:classical_stat} \\
$\mathcal{B}$ & Uniform function bound; also response and clipping bound in bounded regression & \S\ref{sec:classical_stat} \\
$\mathfrak{a}$ & Statistical prefactor exponent of the tanh-DRM slow bound & \S\ref{sec:drm_slow} \\
$M_g, \Lambda_g, \alpha_g$ & Uniform envelope, Lipschitz, and variance--mean constants of the excess loss & \S\ref{sec:four_conditions} \\
\bottomrule
\end{longtable}
\endgroup

\begingroup\renewcommand{\arraystretch}{1.25}
\begin{longtable}{@{}p{0.27\linewidth} >{\raggedright\arraybackslash}p{0.46\linewidth} >{\raggedright\arraybackslash}p{0.20\linewidth}@{}}
\caption{Optimization (PGD on PNN)}\\
\toprule
\textbf{Symbol} & \textbf{Meaning} & \textbf{Definition / location} \\
\midrule
\endhead
$\eta$ & Projection radius for inner weights $\boldsymbol{\theta}_{\mathrm{in}}$ & \S\ref{sec:pgd_algorithm} \\
$\zeta$ & Projection radius for outer coefficients $\boldsymbol{\theta}_{\mathrm{out}}$ & \S\ref{sec:pgd_algorithm} \\
$A_\eta$ & $\ell^2$-ball constraint set on $\boldsymbol{\theta}_{\mathrm{in}}$ & \S\ref{sec:pgd_algorithm} \\
$B_\zeta$ & $\ell^1$-ball constraint set on $\boldsymbol{\theta}_{\mathrm{out}}$ & \S\ref{sec:pgd_algorithm} \\
$T$ (PGD) & Total number of PGD iterations & \S\ref{sec:pgd_algorithm} \\
$\lambda$ & Step size & \S\ref{sec:pgd_algorithm}\\
$R, Q_{\mathrm{trial}}$ & Replica count and initialization trials per replica, respectively & \S\ref{sec:factorisation}, Lem~\ref{lem:O2} \\
$\delta$ (PGD) & Cover radius in $\ell^\infty$ & \S\ref{sec:complete_opt}, Lem~\ref{lem:O2} \\
$\bar{\mathfrak{m}}, \bar{M}, \bar{W}, \bar{L}, B_{\bar{\boldsymbol{\theta}}}$ & PNN parameters of the comparison network $\bar\phi$ ($\bar f$ or $\bar u$) & \S\ref{sec:new_decomp}, \eqref{eq:bar_phi} \\
$F(\boldsymbol{x},\boldsymbol{y})$ & Generic loss of the convex PGD lemma (Lem~\ref{lem:O2b}) & \S\ref{sec:complete_opt} \\
$V$ & Gradient-norm bound $\|\nabla_{\boldsymbol{y}}F\|_2 \le V$ & \S\ref{sec:complete_opt}, Lem~\ref{lem:O2b} \\
$\iota$ & Free parameter controlling the inner radius $\eta = \epsilon^{-\iota}$ & \S\ref{sec:drm_complete}, Thm~\ref{thm:drm_complete} \\
$\|\cdot\|_{\mathrm{eval}}$ & Uniform value norm for regression; uniform value-and-gradient norm for DRM & Assumption~\ref{assum:L_lip} \\
$\mathcal B_{\mathrm{sub}}(B)$ & Sub-network bound in the evaluation norm at parameter radius $B$ & \eqref{eq:subnet_props} \\
$\mathfrak L(B),\mathfrak L$ & Parameter-Lipschitz bound in the evaluation norm; its value at initialization radius & \eqref{eq:subnet_props}--\eqref{eq:frakturL} \\
$L_{\widehat{\mathcal L}},B_{\ell'}$ & Empirical-loss Lipschitz and sampled-loss derivative bounds & Assumption~\ref{assum:L_lip}, \eqref{eq:grad_bound} \\
\bottomrule
\end{longtable}
\endgroup

\begingroup\renewcommand{\arraystretch}{1.25}
\begin{longtable}{@{}p{0.27\linewidth} >{\raggedright\arraybackslash}p{0.46\linewidth} >{\raggedright\arraybackslash}p{0.20\linewidth}@{}}
\caption{Minimax framework}\\
\toprule
\textbf{Symbol} & \textbf{Meaning} & \textbf{Definition / location} \\
\midrule
\endhead
$\mathfrak{M}_n(\Theta_{\mathrm{stat}}, \mathsf d; \mathcal E_n)$ & Minimax risk over class $\Theta_{\mathrm{stat}}$, metric $\mathsf d$, experiment $\mathcal E_n$ & \S\ref{sec:minimax_lb}, \eqref{eq:minimax-risk-abstract} \\
$\mathsf d$ & Metric in which estimation accuracy is evaluated & \S\ref{sec:minimax_lb} \\
$\theta$ & Statistical experiment parameter ($\theta \in \Theta_{\mathrm{stat}}$) & \S\ref{sec:minimax_lb} \\
$t$ & Target-norm order in the Fano template ($\nu_{\mathrm{loss}} = s-t$) & \S\ref{sec:minimax_lb} \\
$q$ & Observation-operator order in the Fano template ($\nu_{\mathrm{KL}} = s-q$) & \S\ref{sec:minimax_lb} \\
$\nu_{\mathrm{loss}}$, $\nu_{\mathrm{KL}}$ & Separation and distinguishability exponents in the Fano template & \S\ref{sec:minimax_lb}, eqs.~\eqref{eq:fano-separation}--\eqref{eq:fano-kl} \\
$m$ (Fano) & Resolution parameter (bumps per dimension) & \S\ref{sec:minimax_lb} \\
$K_m$ (Fano) & Number of packing elements at resolution $m$ & \eqref{eq:fano-cardinality} \\
$\overline D_{\mathrm{KL},m}$ & Average pairwise KL divergence between one-observation laws & \eqref{eq:fano-kl} \\
\bottomrule
\end{longtable}
\endgroup

\begingroup\renewcommand{\arraystretch}{1.25}
\begin{longtable}{@{}p{0.36\linewidth} >{\raggedright\arraybackslash}p{0.41\linewidth} >{\raggedright\arraybackslash}p{0.16\linewidth}@{}}
\caption{Application-specific notation for Section~\ref{sec:applications}}\\
\toprule
\textbf{Symbol} & \textbf{Meaning} & \textbf{Used in} \\
\midrule
\endhead
\multicolumn{3}{@{}l}{\textit{Scientific machine learning}}\\
$u_0,h_0$ & True state and source in the inverse problem & \S\ref{section:app:sml} \\
$e_y,\boldsymbol e_z,\delta$ & Scalar and gradient observation noises, and their componentwise amplitude bound & \S\ref{section:app:sml} \\
$N_{\mathrm{in}},N_b$ & Interior and boundary sample sizes & \S\ref{section:app:sml} \\
$\mathcal L_{\mathrm{PINN}},\widehat{\mathcal L}_{\mathrm{PINN}}$ & Population and empirical PINN objectives & \S\ref{section:app:sml} \\
$\mathcal F_{NN},\mathcal G_{NN}$ (inverse) & State and source candidate network classes for the inverse problem & \S\ref{section:app:sml}, \eqref{eq:pinn_oracle} \\
$\Delta_{\mathrm{PINN}}(\mathcal F)$,\newline $\mathcal E_{\mathrm{sta}}^{\mathrm{PINN}}$ & Sample-dependent uniform deviation over $\mathcal F$; its expectation over the ERM class $\mathcal F_{NN}$ & \eqref{eq:pinn_decomp} \\
$\mathcal E_{\mathrm{app}}^{\mathrm{PINN}}$,\newline $\mathcal E_{\mathrm{app}}^{\mathrm{PINN,c}}$,\newline $\mathcal E_{\mathrm{app}}^{\mathrm{PINN},\mathrm{inv}}$ & Best-in-class forward approximation; specified comparison-network approximation; joint state--source approximation & \S\ref{section:app:sml}, \eqref{eq:pinn_oracle} \\
$\mathcal F_{NN}^{\mathrm{PSGD}}$,\newline $\mathcal E_{\mathrm{opt},-}^{\mathrm{PINN}}$ & Admissible parallel tanh class and comparison-network optimization gap in the PSGD analysis & \S\ref{section:app:sml} \\
$\lambda_{\mathrm{reg}}$ & Tikhonov regularization parameter in the inverse objective & \S\ref{section:app:sml} \\
$\mathcal L_{\mathrm{PINN}}^{\lambda_{\mathrm{reg}}}$,\newline $\widehat{\mathcal L}_{\mathrm{PINN}}^{\lambda_{\mathrm{reg}}}$ & Regularized population and empirical objectives for inverse source recovery & \eqref{eq:pinn_inv_risk}, \eqref{eq:pinn_inverse_empirical} \\
$\alpha,\mu$ & Sobolev regularity index and approximation slack in inverse recovery & \S\ref{section:app:sml} \\
$B_{\delta,\tau},c_\lambda,\tau$ & Full deterministic risk scale, fixed factor in $\lambda_{\mathrm{reg}}=c_\lambda B_{\delta,\tau}$, and confidence parameter & \eqref{eq:pinn_inverse_risk_scale} \\
\midrule
\multicolumn{3}{@{}l}{\textit{Generative learning}}\\
$p_0,p_z,d_Z$ & Target law, prescribed GAN source law, and latent dimension & \S\ref{section:app:generative:gan} \\
$g,f;\ \mathcal G,\mathcal F,\mathcal H$ & Generator, discriminator; their classes and the benchmark test-function class & \S\ref{section:app:generative:gan} \\
$g_\sharp p_z,\ d_{\mathcal F},d_{\mathcal H}$ & Generated pushforward law and integral probability metrics & \S\ref{section:app:generative:gan} \\
$n,m;\ \widehat p_{0,n},\widehat p_{z,m}$ & Data and source sample sizes; their empirical laws & \S\ref{section:app:generative:gan} \\
$W_{\mathcal G},L_{\mathcal G},W_{\mathcal F},L_{\mathcal F}$ & Generator and discriminator widths and depths & \S\ref{section:app:generative:gan} \\
$\mathcal K_{\mathcal F},\mathcal K_{\mathcal F}^{\mathrm{adm}},B_{\mathcal F}$ & Discriminator Lipschitz radius, sufficient radius from the approximation construction, and uniform output bound & \S\ref{section:app:generative:gan} \\
$\mathcal E(\mathcal H,\mathcal F,\Omega)$ & Uniform approximation error of the benchmark class by discriminators on $\Omega$ & \S\ref{section:app:generative:gan} \\
$\widetilde{\boldsymbol X},\boldsymbol\xi,V_{\mathrm{noise}}$ & Clean low-dimensional component, independent centered additive noise, and its second moment $\mathbb E\|\boldsymbol\xi\|_2^2$ & \S\ref{section:app:generative:gan}, \eqref{eq:gan_intrinsic_noise} \\
$\gamma,\widehat\gamma$ & Standard Gaussian law and the actual FM initialization law & \S\ref{section:app:generative} \\
$\gamma_a$ & Density of $\mathcal N(\boldsymbol 0,a\boldsymbol I_d)$, $a>0$ & \S\ref{section:app:generative:flow_matching} \\
$m_2$ & Root second moment $(\int\|\boldsymbol x\|_2^2\,p_0(\mathrm d\boldsymbol x))^{1/2}$ of the target law & \S\ref{section:app:generative} \\
$t,\tau,\boldsymbol X_t,\boldsymbol Y_\tau,\widehat{\boldsymbol Y}_\tau$ & Forward and reverse times ($\tau=T_{\mathrm{end}}-t$), forward process, exact and learned reverse processes & \eqref{eq:forward:base}, \eqref{eq:base:reversal}, \eqref{eq:approx:reversal} \\
$\boldsymbol Y_\tau^{\mathrm{pf}}$ & Exact-score probability-flow ODE process in reverse time & \eqref{eq:diffusion:probability:flow} \\
$p_t$ & Diffusion forward marginal (data to noise); separately, FM interpolation marginal (noise to data) & \S\ref{section:app:generative:diffusion}, \S\ref{section:app:generative:flow_matching} \\
$\mu_t,\sigma_t,\boldsymbol B_t$ & OU mean coefficient, Gaussian noise scale, and Brownian motion & \S\ref{section:app:generative:diffusion} \\
$s,s^*,\widehat s,\mathcal S$ & Candidate, true, and learned scores, and the score class; $s^*(t,\cdot)=\nabla\log p_t$ & \S\ref{section:app:generative:diffusion} \\
$\mathcal L_{\mathrm{DSM}},\widehat{\mathcal L}_{\mathrm{DSM}},\mathcal E_{\mathrm{score}}^2$ & Population and empirical denoising objectives, and expected weighted integrated score error & \S\ref{section:app:generative:diffusion} \\
$T_{\mathrm{end}},\delta,\widehat p_\delta,\widehat p^{\mathrm{tr}}$ & Forward terminal time and cutoff (reverse sampling stops at $T_{\mathrm{end}}-\delta$), learned output law, and truncated output law & \S\ref{section:app:generative:diffusion} \\
$B^\beta_{p,q};\ \beta,p,q$ & Besov density class; smoothness, integrability, and summability indices & \S\ref{section:app:generative:diffusion} \\
$\bar p_\delta$ & Auxiliary learned reverse-process output initialized from $p_{T_{\mathrm{end}}}$ & \S\ref{section:app:generative:diffusion} \\
$\mathcal E_{\mathrm{early}}^{\mathrm{TV}},\mathcal E_{\mathrm{init}}^{\mathrm{TV}}$ & TV distances $d_{\mathrm{TV}}(p_0,p_\delta)$ and $d_{\mathrm{TV}}(p_{T_{\mathrm{end}}},\gamma)$ & \S\ref{section:app:generative:diffusion} \\
$\boldsymbol X_{\mathrm{data}}$ & Target-distributed sample used in the FM interpolant & \S\ref{section:app:generative:flow_matching} \\
$a_t,b_t$ & Data and Gaussian-noise coefficients in the general FM interpolant & \S\ref{section:app:generative:flow_matching} \\
$u,v,\widehat v,\mathcal V$ & True, candidate, and learned FM velocities, and the velocity class & \S\ref{section:app:generative:flow_matching} \\
$\psi,\widehat\psi$ & Ideal and learned ODE flow maps & \S\ref{section:app:generative:flow_matching} \\
$\mathcal L_{\mathrm{FM}},\widehat{\mathcal L}_{\mathrm{FM}},\mathcal E_{\mathrm{FM}}^2$ & Population and empirical flow-matching objectives, and expected integrated squared velocity error & \S\ref{section:app:generative:flow_matching} \\
$\mathcal E_{\mathrm{app}}^2,\mathcal E_{\mathrm{sta}}^2$ & Best-in-class excess risk and expected uniform empirical deviation, on the DSM/FM squared-risk scale & \S\ref{section:app:generative:diffusion}, \S\ref{section:app:generative:flow_matching} \\
$L_{\mathrm{Lip},v}(t),\Lambda_t$ & Spatial Lipschitz bound and its integrated stability amplification & \S\ref{section:app:generative:flow_matching} \\
$C_{\mathrm{Lip}},\kappa,\delta_n$ & Endpoint Lipschitz coefficient, interpolation noise-decay exponent, and FM stopping-time offset & \S\ref{section:app:generative:flow_matching} \\
$\widehat p^{\mathrm{num}},\mathcal E_{\mathrm{disc}}^{\mathrm{TV}},\mathcal E_{\mathrm{disc}}^{\mathrm{FM}}$ & Numerical output law and diffusion/FM numerical-discretization errors in the stated metrics & \S\ref{section:app:generative} \\
$\epsilon_{\mathrm{opt}}$ & Empirical-objective suboptimality retained as an external residual; not squared again in squared-risk bounds & \S\ref{section:app:generative} \\
\midrule
\multicolumn{3}{@{}l}{\textit{Reinforcement learning}}\\
$\mathcal S,\mathcal A,P,\pi,\gamma$ & State and action spaces, transition kernel, policy, and discount factor & \S\ref{section:app:rl} \\
$Q^\pi,Q^*,\mathcal T^\pi,\mathcal T^*$ & Policy and optimal action-value functions, and their Bellman operators & \S\ref{section:app:rl} \\
$\varrho,\eta_{\mathrm{ev}},V_{\max}$ & Sampling and evaluation distributions; value bound $R_{\max}/(1-\gamma)$ & \S\ref{section:app:rl} \\
$n,J,N_{\mathrm{tot}}$ & Transitions per fitting stage, outer-iteration count, and total FQI sample size $N_{\mathrm{tot}}=nJ$ & \S\ref{section:app:rl:fqi} \\
$a_n,B_n,\bar\beta_k$ & Block length, number of blocks in each alternating group with $n=2a_nB_n$, and $\beta$-mixing coefficient at lag $k$ & \S\ref{section:app:rl:fqi} \\
$C_{\eta_{\mathrm{ev}},\varrho},C^{\mathrm{API}}_{\eta_{\mathrm{ev}},\varrho}$ & FQI/API concentrability, using $L^2/L^\infty$ density-ratio norms, respectively & \S\ref{section:app:rl} \\
$\mathcal F_Q,\mathcal C_n,O$ & FQI value class, API CNN class, and auxiliary Bellman-image network & \S\ref{section:app:rl} \\
$\epsilon_{\mathrm{opt},j}$ & Empirical optimization tolerance at fitting stage $j$, distinct from the outer-iteration remainder & \S\ref{section:app:rl} \\
$r_n,\tau_n,J_n$ & API simultaneous residual scale, failure probability, and iteration count at per-stage sample size $n$ & \eqref{eq:rl:api:uniform} \\
\midrule
\multicolumn{3}{@{}l}{\textit{Representation learning}}\\
$\phi,d_{\mathrm{enc}},\Phi_n,\Phi_{\mathrm{adm}}$ & Encoder, representation dimension, neural encoder class, and admissible benchmark family & \S\ref{section:app:repssl} \\
$\nu_{d_{\mathrm{enc}}},D_f,\mathrm{dCov}^{\,2}$ & Gaussian reference, $f$-divergence, and squared distance covariance & \S\ref{app:rep_learning} \\
$\lambda_{\mathrm{lag}},D,\mathcal D_n$ & DDR regularization weight, discriminator, and discriminator class & \S\ref{app:rep_learning} \\
$\mathcal L_{\mathcal D_n},\widehat{\mathcal L}_{\mathcal D_n,n}$ & Population and empirical DDR objectives with restricted discriminators & \S\ref{app:rep_learning} \\
$\mathcal L^\star,\mathcal E_{\mathrm{app}}^{\mathrm{rep}}$ & Infimum over $\Phi_{\mathrm{adm}}$ and encoder approximation gap relative to that benchmark & \eqref{eq:rep_ddr_decomp} \\
$\mathcal E_{\mathrm{disc}},\Delta_n$ & DDR discriminator approximation error and uniform empirical deviation & \S\ref{app:rep_learning} \\
$n,m,\sigma_i$ & Pretraining and labeled downstream sample sizes; contrastive positive/negative pair indicator & \S\ref{section:app:repssl}, \S\ref{app:ssl} \\
$P_\phi^{\mathrm{pair}},P_{\phi,h}^{\mathrm{down}}$ & Joint laws of the pretraining triple and downstream input--response pair & \S\ref{app:ssl} \\
$\phi^\star,\widetilde\phi$ (contrastive) & Target encoder and its best $L^2(P_X)$ approximation in $\Phi_n$ & \S\ref{app:ssl} \\
$h,\widehat h,h^\star,\mathcal H_{\mathcal K}$ & Candidate, fitted, and target downstream heads, and the head class with norm budget $\mathcal K$ & \S\ref{app:ssl} \\
$\beta_0,\beta_1$ & Encoder and head smoothness; $\beta_0$ also denotes encoder smoothness in Adv-SSL & \S\ref{app:ssl}, \S\ref{app:ssl_adv} \\
$\Delta_m(\widehat\phi),\epsilon_{\mathrm{opt}}^{\mathrm{head}}$ & Conditional downstream empirical deviation and head-fitting optimization tolerance & \S\ref{app:ssl} \\
$\lambda_{\mathrm{ssl}},\boldsymbol C_\phi,\widehat{\boldsymbol C}_\phi$ & Adv-SSL regularization weight, population cross-moment discrepancy, and its sample version & \S\ref{app:ssl_adv} \\
$\boldsymbol G_{\mathrm{dual}},\mathcal G_{\mathrm{dual}}(\phi)$,\newline $\widehat{\mathcal G}_{\mathrm{dual}}(\phi)$ & Dual matrix and population/sample-dependent constraint sets & \S\ref{app:ssl_adv} \\
$\Delta_{\mathrm{SSL}},\Delta_{\mathrm{set}}$ & Uniform empirical fluctuation and dual-set mismatch & \S\ref{app:ssl_adv} \\
$\widehat{\boldsymbol\mu}_k,m_k,Q_{\widehat\phi}$ & Estimated class template, labeled class count, and template classifier & \S\ref{app:ssl_adv} \\
$B_\phi,B_G$ & Encoder output-norm bound and enclosing dual-ball radius $B_G=B_\phi^2+\sqrt{d_{\mathrm{enc}}}$ & \S\ref{app:ssl_adv} \\
\midrule
\multicolumn{3}{@{}l}{\textit{Transformers and large language models}}\\
$N_{\mathrm{reg}},L_{\mathrm{win}},d_x$ & Observed time points, fixed input-window length, and input dimension per time point; there are $N_{\mathrm{reg}}-L_{\mathrm{win}}+1$ training windows & \eqref{eq:llm_reg_model}, \eqref{eq:llm_reg_erm} \\
$\boldsymbol H_t,\Pi$ & Input-window matrix and its common distribution & \eqref{eq:llm_reg_model}, \S\ref{app:llm_reg} \\
$T_{\boldsymbol\theta},\boldsymbol A,\mathcal F_{\mathrm{Tr}}$ & Transformer map, trainable linear readout matrix, and resulting scalar-valued predictor class & \S\ref{app:llm_reg} \\
$\mathcal R_{\mathrm{Tr}}$ & Population squared-loss risk for sequence regression & \eqref{eq:llm_reg_identity} \\
$\widehat f,\widetilde f$ (Transformer) & Sliding-window empirical risk minimizer before clipping, and its clipped version used to evaluate prediction error & \eqref{eq:llm_reg_erm}, \S\ref{app:llm_reg} \\
$B,C_B$ (Transformer) & Output-clipping threshold and pointwise clipping operator: $C_Bf=[f]_B$, with $B\asymp\log N_{\mathrm{reg}}$ & \S\ref{app:llm_reg}; clipping convention in \S\ref{sec:reg_slow} \\
$a,\bar\beta_a$ (Transformer) & Auxiliary blocking separation and mixing coefficient of the window--response process; all sliding windows remain in the empirical loss & \eqref{eq:llm_reg_decomp} \\
$C_{\mathrm{mix}},b_{\mathrm{mix}},\eta$ & Constants and exponent in the geometric mixing-decay bound, specified separately for each data process & \S\ref{section:app:rl:fqi}, \S\ref{app:llm_reg} \\
$\mathcal F,\mathcal R$ (pretraining) & Transformer conditional-probability class and population negative-log-likelihood risk & \eqref{eq:llm_erm}, \S\ref{app:llm_pretraining} \\
$\mathcal V,N_{\mathrm{seq}},N_{\mathrm{tok}}$ & Finite vocabulary, corpus sequence count, and token-prediction count $N_{\mathrm{tok}}=N_{\mathrm{seq}}L_{\mathrm{seq}}$; the count does not assert independence & \S\ref{app:llm_pretraining} \\
$L_{\mathrm{seq}},L_{\mathrm{res}}$ & Common corpus sequence length (also the context limit) and predicted response length & \S\ref{app:llm_pretraining} \\
$\mathcal P,x,y,\boldsymbol y$ & Full prompt, target query, single answer/token, and response sequence & \S\ref{app:llm_pretraining} \\
$y^{(i)}$ (ICL) & Demonstration response, possibly containing multiple tokens & \S\ref{app:llm_understand} \\
$\boldsymbol d^{(i)},t_j^{(i)},\boldsymbol h_j^{(i)}$ & Corpus sequence, its $j$th token, and preceding context & \eqref{eq:llm_erm} \\
$p^*,\bar p,\widehat p_{\mathrm{ERM}},\widehat p$ & True law, population minimizer, ideal empirical minimizer, and actual trained model & \S\ref{app:llm_pretraining} \\
$\mathcal E_{\mathrm{mem}},\mathcal E_{\mathrm{gen}}$ & Memorization and generalization discrepancies on the prediction-probability scale & \eqref{eq:llm_gen_mem} \\
$\epsilon_{\mathrm{opt}}^{\mathrm{pred}}$ & Upper bound on the prediction-probability discrepancy from the ideal empirical minimizer, not an empirical-risk gap & \eqref{eq:llm_gen_mem} \\
$\mathcal E_{\mathrm{pre}}$ & Uniform response-probability discrepancy of the trained model over the specified prompt--response pairs, including the prediction-level training residual & \eqref{eq:llm_pretraining_error} \\
$\theta_{\mathrm{task}},\Theta_{\mathrm{task}}$,\newline $\theta_{\mathrm{task},x},\Delta_{\Theta_{\mathrm{task}}}(x)$ & Latent task, discrete task space, task of highest posterior probability given $x$, and posterior task ambiguity & \S\ref{app:llm_pretraining} \\
$\mathcal E_{\text{zero-shot}}(x,\boldsymbol y)$ & Zero-shot response-probability error relative to the dominant task given the query & \eqref{eq:llm_zero_shot_bound} \\
$\theta_x^\star,\theta_y^\star,\Delta_{\mathrm C},\Delta_{\mathrm{trans}}$ & Carrier and transformed states, and their ambiguities under the optional paired-task interpretation & \eqref{eq:llm_refined_ambiguities} \\
$m,\phi_{\mathrm M},\epsilon_{\mathrm{ICL}}$ & Demonstration count, nearly-Markov departure, and ICL ambiguity coefficient & \S\ref{app:llm_understand} \\
$\mathcal E_{\mathrm{ICL}}(\mathcal P_{\mathrm{ICL}},\boldsymbol y)$ & ICL response-probability error relative to the query's dominant task, with demonstrations in the prompt & \eqref{eq:llm_icl_bound} \\
$t_{\mathrm{task}}^\star,L_{\mathrm{CoT}}$ & Target compositional task and reasoning-chain length & \S\ref{app:llm_understand} \\
$\widetilde p^*,M_{\mathrm{CoT}},\Delta_{\mathrm{PCoT}}$ & Compositional-task law; maximum reciprocal prompt probability under the two laws; prior discrepancy over task sequences compatible with the prompt & \S\ref{app:llm_understand} \\
$\mathcal E_{\mathrm{CoT}}(\mathcal P_{\mathrm{CoT}},\boldsymbol y;t_{\mathrm{task}}^\star)$ & CoT response-probability error relative to the specified target task sequence under the compositional-task law & \eqref{eq:llm_cot_bound} \\
\bottomrule
\end{longtable}
\endgroup

\section{Tool dictionary and convergence rate index}\label{app:tool_index}

\begingroup\renewcommand{\arraystretch}{1.25}
\begin{longtable}{@{}l|>{\raggedright\arraybackslash}p{3.0in}|>{\raggedright\arraybackslash}p{1.4in}@{}}
\caption{Tool dictionary}\\
\toprule
\textbf{Label} & \textbf{Named result} & \textbf{Location} \\
\midrule
\endhead
A1 & ReLU approximation of Hölder functions & Theorem~\ref{thm:A1} \\
A2 & Gühring--Raslan tanh Sobolev approximation & Theorem~\ref{thm:A2} \\
A3 & Norm-controlled ReLU approximation (Jiao et al.) & Theorem~\ref{thm:A3} \\
A4 & Norm-controlled tanh approximation (Ding et al.) & Theorem~\ref{thm:A4} \\
A5 & PNN tanh approximation (Jiao et al.) & Theorem~\ref{thm:A5} \\
A5r & PNN ReLU approximation (Kohler--Cho--Krzy\.{z}ak) & Theorem~\ref{thm:A5r} \\
S1 & Symmetrization + Massart bound for a finite loss net & Prop~\ref{prop:slow_engine} \\
S2 & ReLU pseudo-dimension bound & Prop~\ref{prop:relu_pdim} \\
S3 & Regression slow-rate Rademacher bound & Prop~\ref{prop:S3} \\
S4 & Fast-rate Rademacher via peeling (regression) & Prop~\ref{prop:S4} \\
S5 & DRM slow-rate Rademacher bound & Prop~\ref{prop:S5} \\
S6 & Golowich norm-controlled Rademacher bound & Lemma~\ref{lem:S6} \\
S7 & PNN-to-sub-network Rademacher reduction & Lemma~\ref{lem:S7} \\
S8 & Slow/fast statistical engine under Assumptions~\ref{assum:bounded}--\ref{assum:var} & Props~\ref{prop:slow_engine}, \ref{prop:fast_engine} \\
O1 & NTK and Polyak--Łojasiewicz convergence (background) & \S\ref{sec:opt_landscape} \\
O2 & Combined initialization and iteration bound & Lemma~\ref{lem:O2} \\
\bottomrule
\end{longtable}
\endgroup

\begingroup\renewcommand{\arraystretch}{1.25}
\begin{longtable}{@{}>{\raggedright\arraybackslash}p{0.95in}|>{\raggedright\arraybackslash}p{0.75in}|>{\raggedright\arraybackslash}p{1.35in}|>{\raggedright\arraybackslash}p{1.5in}@{}}
\caption{Convergence rate index. Unless indicated otherwise, only the leading algebraic orders are displayed, with logarithmic factors suppressed. See the referenced theorems for the assumptions and expectation or high-probability formulations.}\\
\toprule
\textbf{Theorem} & \textbf{Reference} & \textbf{Rate} & \textbf{Building blocks} \\
\midrule
\endhead
\multicolumn{4}{@{}l}{\emph{Classical size-controlled (\S\ref{sec:classical_setting}, $\mathcal{E}_{opt} = 0$)}} \\[1pt]
Regression slow rate & Thm \ref{thm:reg_slow} & $n^{-2\beta/(4\beta+d)}$ & A1 + S3, symmetric decomp \eqref{eq:reg_sym_decomp} \\
DRM slow rate & Thm \ref{thm:drm_slow} & $N^{-(1-\mu)/(C d \log(d+1))}$ & A2 ($s=2$) + S5 + Prop \ref{prop:drm_eq} \\
Regression fast rate & Thm \ref{thm:reg_fast} & $n^{-2\beta/(2\beta+d)}$ & A1 + S4 (peeling), asymmetric decomp \eqref{eq:reg_asym_decomp} \\
DRM fast rate (external) & Thm \ref{thm:drm_fast} & $N^{-2(s-1)/(2(s-1)+d)}$ & $\mathrm{ReLU}^3$ + localized Rademacher + Dirichlet energy coercivity \cite{lu2021machine} \\
\midrule
\multicolumn{4}{@{}l}{\emph{Over-parameterized (\S\ref{sec:overparam_setting}, $\mathcal{E}_{opt} = 0$)}} \\[1pt]
Regression over-param & Thm \ref{thm:reg_overparam} & $n^{-\beta/(d+1+2\beta)}$ & A3 + S6 \\
DRM over-param & Thm \ref{thm:drm_overparam} & $N^{-c_\mu/d^2}$ (mean-square) & Tanh construction + layerwise derivative bound \eqref{eq:drm_op_uniform_stat} on class \eqref{eq:drm_op_candidate} \\
\midrule
\multicolumn{4}{@{}l}{\emph{Complete three-error (\S\ref{sec:complete_setting})}} \\[1pt]
Regression complete & Thm \ref{thm:reg_complete} & $n^{-\beta/(2\beta+d)}$ (up to $\log^5 n$) & A5r + S7 + O2 \\
DRM complete & Thm \ref{thm:drm_complete} & $\tilde{\mathcal{O}}(N^{-\tilde C_3^{-1}})$ & A5 + S7 + O2 + Prop \ref{prop:drm_eq} \\
\midrule
\multicolumn{4}{@{}l}{\emph{Risk lower bounds under the stated observation models (\S\ref{sec:minimax_lb})}} \\[1pt]
Hölder regression & Thm \ref{thm:holder-regression-lower-bound} & $n^{-2\beta/(2\beta+d)}$ & Fano's method \\
Elliptic solution recovery (noisy RHS) & Thm \ref{thm:noisy-rhs-lower-bound} & $n^{-2(s-1)/(d+2s-4)}$ & Fano's method \\
\bottomrule
\end{longtable}
\endgroup